\documentclass[twoside]{article}

\usepackage{tocbibind}
\usepackage{fancyhdr,bm,mathtools,graphicx,leftidx}
\usepackage{amssymb,amsmath}
\usepackage{marvosym}
\usepackage{mathstyle}
\usepackage{hieroglf}
\usepackage[T1]{fontenc}
\usepackage{eufrak}
\usepackage{romanbar}
\usepackage[dvipsnames]{xcolor}
\usepackage{ifthen}
\usepackage{tikz-cd}
\usepackage[title]{appendix}
\tikzcdset{
arrow style=tikz,
diagrams={>={Stealth[scale=1]}}
}
\usepackage[T1]{fontenc}
\usepackage{textcomp}

\usepackage{stmaryrd}
\usepackage[mathscr]{euscript}
\usepackage{chngcntr}
\usepackage[style]{fncychap}
\usepackage{latexsym}

\def\Ref#1{[{equation~\ref{#1}}]}

\newcounter{chapterr}
\counterwithin{section}{chapterr}

\newcounter{definition}
\counterwithin{definition}{chapterr}
\def\definition
{\refstepcounter{definition}\vskip0.4\baselineskip\noindent{{\textbf{Definition~}}}{$\bf\arabic{chapterr}.$}{$\bf\arabic{definition}$}:\hskip 0.4\parindent}
\newcounter{theorem}
\counterwithin{theorem}{chapterr}
\def\theorem
{\refstepcounter{theorem}\vskip0.4\baselineskip\noindent{{\textbf{Theorem~}}}{$\bf\arabic{chapterr}.$}{$\bf\arabic{theorem}$}:\hskip 0.4\parindent}
\newcounter{symbol}
\counterwithin{symbol}{chapterr}

\newcounter{QMaxiom}
\counterwithin{QMaxiom}{chapterr}

\newcounter{corollary}
\counterwithin{corollary}{chapterr}
\def\corollary
{\refstepcounter{corollary}\vskip0.4\baselineskip\noindent{{\textbf{Corollary~}}}{$\bf\arabic{chapterr}.$}{$\bf\arabic{corollary}$}:\hskip 0.4\parindent}
\newcounter{lemma}
\counterwithin{lemma}{chapterr}
\def\lemma
{\refstepcounter{lemma}\vskip0.4\baselineskip\noindent{{\textbf{Lemma~}}}{$\bf\arabic{chapterr}.$}{$\bf\arabic{lemma}$}:\hskip 0.4\parindent}
\newcounter{proposition}
\counterwithin{proposition}{chapterr}
\def\proposition
{\refstepcounter{proposition}\vskip0.4\baselineskip\noindent{{\textbf{Proposition~}}}{$\bf\arabic{chapterr}.$}{$\bf\arabic{proposition}$}:\hskip 0.4\parindent}
\newcounter{remark}
\counterwithin{remark}{chapterr}
\def\remark
{\refstepcounter{remark}\vskip0.4\baselineskip\noindent{{\textbf{Remark~}}}{$\bf\arabic{chapterr}.$}{$\bf\arabic{remark}$}:\hskip 0.4\parindent}
\def\proof
{\vskip0.3\baselineskip\noindent{\textbf{proof}}\hskip.08\baselineskip:\hskip 0.4\parindent}

\newcounter{fixed}
\counterwithin{fixed}{chapterr}
\def\fixed
{\refstepcounter{fixed}\vskip0.4\baselineskip\noindent{{\textbf{Fixed Objects~}}}{$\bf\arabic{chapterr}.$}{$\bf\arabic{fixed}$}:\hskip 0.4\parindent}
\newcounter{exercise}
\counterwithin{exercise}{chapterr}
\def\exercise
{\refstepcounter{exercise}\vskip0.4\baselineskip\noindent{{\textbf{Exercise~}}}{$\bf\arabic{chapterr}.$}{$\bf\arabic{exercise}$}:\hskip 0.4\parindent}
\def\caution{\textasteriskcentered\hskip0.25\baselineskip}
\def\V{V}

\def\R{\mathbb{R}}

\newcommand{\p}[1]{{#1^{\prime}}}

\def\R{\mathbb{R}}

\def\f{{\rm F}}

\def\card{{\textsf{card}}}
\def\endef{{\flushright{\noindent$\blacksquare$\\}}\noindent}
\def\endthm{{\flushright{\noindent$\square$\\}}\noindent}
\def\endcor{{\flushright{\noindent$\square$\\}}\noindent}

\def\endlem{{\flushright{\noindent$\square$\\}}\noindent}
\def\endpro{{\flushright{\noindent$\square$\\}}\noindent}

\def\endlem{{\flushright{\noindent$\square$\\}}\noindent}
\def\endfixed{{\flushright{\noindent$\Diamond$\\}}\noindent}
\def\endremark{\vskip0.5\baselineskip}
\def\endexercise{\vskip0.5\baselineskip}

\def\then{\Rightarrow}
\def\thenn{\Leftrightarrow}

\def\({\left(}
\def\){\right)}
\def\[{\left[}
\def\]{\right]}

\newcommand{\CSs}[1]{{\mathcal{P}}\(#1\)}

\newcommand{\Card}[1]{\left| #1\right|}
\newcommand{\CarD}[1]{\card\(#1\)}
\def\cardeq{\overset{\underset{\mathrm{card}}{}}{=}}

\newcommand{\union}[1]{\bigcup#1}
\newcommand{\intersection}[1]{\bigcap#1}

\newcommand{\Union}[3]{\bigcup_{#1\in#2}#3}

\newcommand{\Dproduct}[3]{\prod_{{#1}\in{#2}}{#3}}
\newcommand{\cmp}[2]{#1\circ#2}
\newcommand{\Func}[2]{{\textsf{F}}\bpair{#1}{#2}}
\newcommand{\surFunc}[2]{{{\mathfrak{s}}\textsf{F}}\opair{#1}{#2}}
\newcommand{\IF}[2]{{\mathsf{B\hphantom{}F}}\opair{#1}{#2}}

\newcommand{\domain}[1]{{\mathsf{dom}}\bbsingle{#1}}
\newcommand{\codomain}[1]{{\mathsf{codom}}\bbsingle{#1}}
\newcommand{\funcimage}[1]{{\mathsf{img}}\bbsingle{#1}}
\newcommand{\EqR}[1]{{\mathsf{EqR}}\bbsingle{#1}}
\newcommand{\EqClass}[2]{{#1}\left/{#2}\right.}
\newcommand{\pEqclass}[2]{\[{#1}\]_{#2}}
\newcommand{\PEqclass}[2]{{\mathsf{EqC}}\bbpair{#1}{#2}}
\newcommand{\Cprod}[2]{{#1}\times{#2}}

\newcommand{\OR}[2]{#1\thinspace\lor\thinspace#2}
\newcommand{\AND}[2]{#1\thinspace\land\thinspace#2}
\newcommand*{\suchthat}{\;\ifnum\currentgrouptype=16 \middle\fi|\;}
\newcommand{\Foreach}[2]{\forall\thinspace#1\in#2\negthinspace:\thickspace}
\newcommand{\Exists}[2]{\exists\thinspace#1\in#2\negthinspace:\thickspace}

\newcommand{\defset}[3]{\left\{#1\in#2\thickspace:\thickspace#3\right\}}

\newcommand{\defSet}[2]{\left\{#1\thickspace|\thickspace#2\right\}}
\newcommand{\negation}[1]{\neg{#1}}
\def\eqdef{\overset{\underset{\mathrm{def}}{}}{=}}

\newcommand{\resd}[1]{{\mathfrak{res}}{\mathsf{D}}_{#1}}
\newcommand{\rescd}[1]{{\mathfrak{res}}{\mathsf{C\negthinspace D}}_{#1}}
\newcommand{\res}[1]{{\mathfrak{res}}_{#1}}
\newcommand{\finv}[1]{{#1}^{-1}}
\def\index{{\mathscr{I}}}

\def\empty{\varnothing}

\newcommand{\opair}[2]{\(#1,\thinspace #2\)}

\newcommand{\bpair}[2]{\negthinspace\left(#1,\thinspace #2\right)}
\newcommand{\bbpair}[2]{\negthinspace\left(#1;\thinspace #2\right)}
\newcommand{\bbsingle}[1]{\negthinspace\left({#1}\right)}
\newcommand{\cpair}[2]{\[#1,\thinspace #2\]}

\newcommand{\compl}[2]{#1\setminus#2}

\newcommand{\seta}[1]{\left\{#1\right\}}
\newcommand{\func}[2]{#1\(#2\)}

\def\U{U}

\def\asubset{A}

\def\point{p}
\def\cf{f}

\newcommand{\image}[1]{#1^{\rightarrow}}
\newcommand{\pimage}[1]{#1^{\leftarrow}}

\newcommand{\function}[3]{#1:\thinspace#2\to#3}

\newcommand{\binary}[2]{#1,\thinspace#2}

\newcommand{\identity}[1]{{\rm{Id}}_{#1}}

\newcommand{\VVS}[1]{\VS_{#1}}

\newcommand{\Vspan}[1]{{\mathsf{span}}_{#1}}
\def\VS{{\mathbb{V}}}
\newcommand{\NVS}[1]{{\mathscr{V}}_{#1}}

\newcommand{\Lin}[2]{{\mathsf{L}}\bpair{#1}{#2}}
\newcommand{\Linisom}[2]{{\mathsf{LIsom}}\bpair{#1}{#2}}

\newcommand{\VLin}[2]{{\mathbb{L}}\bpair{#1}{#2}}

\newcommand{\Det}[1]{{\mathrm{det}}_{#1}}
\newcommand{\triple}[3]{\opair{#1}{\binary{#2}{#3}}}

\newcommand{\mtuple}[2]{\(\suc{#1}{#2}\)}

\newcommand{\abs}[1]{\left|{#1}\right|}

\newcommand{\dualV}[1]{{#1}^{\star}}

\newcommand{\suc}[2]{{#1},\ldots,\thinspace{#2}}

\newcommand{\funcprod}[2]{{#1}{\overset{\underset{\mathrm{\mathfrak{f}}}{}}{\times}}{#2}}
\newcommand{\Injection}[2]{{\mathrm{Inj}}_{{#1}\to{#2}}}

\newcommand{\Man}[1]{\mathscr{M}_{#1}}

\newcommand{\manprod}[2]{{#1}{\overset{\underset{\mathrm{\mathfrak{m}}}{}}{\times}}{#2}}

\newcommand{\maxatlases}[3]{{\mathsf{maxAtl}}^{\(#1\)}\bpair{#2}{#3}}
\newcommand{\atlases}[3]{{\mathsf{Atl}}^{\(#1\)}\bpair{#2}{#3}}

\newcommand{\maxatlasgen}[3]{{\mathfrak{maxAtl}}^{\(#1\)}_{\opair{#2}{#3}}}
\newcommand{\M}[1]{\mathrm{M}_{#1}}
\newcommand{\maxatlas}[1]{{\boldsymbol{\mathscr{A}}}_{#1}}
\newcommand{\atlas}[1]{\mathscr{A}_{#1}}
\newcommand{\difclass}[1]{\mathsf{C}^{#1}}

\newcommand{\avecf}[1]{{\mathscr{X}}_{#1}}
\newcommand{\avecff}[1]{{\mathscr{Y}}_{#1}}
\newcommand{\tanspace}[2]{{\mathsf{T}}_{#1}{#2}}
\newcommand{\Tanspace}[2]{{\mathbb{T}}_{#1}{#2}}

\newcommand{\tanbun}[1]{{\mathsf{T}}#1}
\newcommand{\Tanbun}[1]{{\mathbf{T}}#1}
\newcommand{\basep}[1]{\pi_{#1}}
\newcommand{\mapdifclass}[3]{{\mathsf{C}}^{#1}\bpair{#2}{#3}}

\newcommand{\tanspaceiso}[3]{{\boldsymbol{\theta}}_{#1}^{\opair{#2}{#3}}}
\newcommand{\derr}[1]{{\mathsf{D}}{#1}}

\newcommand{\Eucbase}[2]{{\boldsymbol{e}}^{\(#1\)}_{#2}}
\newcommand{\deltaf}[2]{\delta_{{#1},{#2}}}

\newcommand{\zerovec}[1]{{\boldsymbol{0}}_{#1}}
\newcommand{\liebracket}[3]{\left\llbracket#1,~#2\right\rrbracket_{#3}}

\newcommand{\Vdual}[1]{{#1}^{\star}}

\newcommand{\Tensors}[3]{{\mathbf{T}}_{#1}^{#2}\bbsingle{#3}}

\newcommand{\Dsum}[1]{\bigoplus#1}
\newcommand{\directsum}[2]{{#1}\oplus{#2}}

\newcommand{\tensor}[1]{\overset{\underset{\mathrm{#1}}{}}{\otimes}}
\newcommand{\vsbase}[1]{{\mathbf{e}}_{#1}}

\newcommand{\quadruple}[4]{\opair{#1}{\binary{#2}{\binary{#3}{#4}}}}

\newcommand{\subman}[2]{\left. #1\right|_{#2}}

\newcommand{\reS}[2]{\left. #1\right|_{#2}}

\newcommand{\fbmorb}[1]{\widehat{#1}}

\newcommand{\quintuple}[5]{\opair{#1}{\binary{#2}{\binary{#3}{\binary{#4}{#5}}}}}
\newcommand{\vbundle}[1]{{\mathbb{V}}_{#1}}
\newcommand{\vbtotal}[1]{{\mathscr{E}}_{#1}}
\newcommand{\vbbase}[1]{{\mathscr{B}}_{#1}}
\newcommand{\vbfiber}[1]{{\mathbb{X}}_{#1}}
\newcommand{\vbprojection}[1]{{\boldsymbol{\pi}}_{#1}}
\newcommand{\vTot}[1]{{\mathrm{E}}_{#1}}
\newcommand{\vB}[1]{{\mathrm{B}}_{#1}}

\newcommand{\vbatlas}[1]{{\mathscr{A}}_{#1}}

\newcommand{\GL}[2]{{\mathsf{GL}}_{#2}\bbsingle{#1}}
\newcommand{\fibervecs}[2]{{\mathbf{Fib}}_{#1}\bbsingle{#2}}
\newcommand{\vbsec}[1]{\sigma_{#1}}
\newcommand{\vbsections}[1]{{\boldsymbol{\Gamma}}\bbsingle{#1}}

\newcommand{\vbtensorbundle}[3]{{\mathbb{TB}}^{\opair{#1}{#2}}\bbsingle{#3}}
\newcommand{\TF}[3]{{\mathsf{TF}}^{\infty}_{\opair{#1}{#2}}\bbsingle{#3}}

\newcommand{\VTF}[3]{{\mathbb{TF}}^{\infty}_{\opair{#1}{#2}}\bbsingle{#3}}
\newcommand{\DVTF}[1]{\overset{\underset{\mathrm{VTF}}{}}{\bigotimes}#1}

\newcommand{\vbmorphisms}[2]{{\mathsf{VBMor}}\bpair{#1}{#2}}
\newcommand{\vbisomorphisms}[2]{{\mathsf{VBIsom}}\bpair{#1}{#2}}
\newcommand{\VBpullback}[3]{{#1}^{\star}_{\opair{#2}{#3}}}

\newcommand{\vbDualbundle}[1]{{#1}^{\star}}

\newcommand{\myitem}[1]{${\fontsize{6.65}{7}\selectfont{\textbf{#1}}}$}
\newcommand{\connection}[1]{\nabla_{#1}}
\newcommand{\connections}[1]{{\mathsf{CN}}\bbsingle{#1}}
\newcommand{\vectorfields}[1]{{\boldsymbol{\mathfrak{X}}}\bbsingle{#1}}
\newcommand{\oneforms}[1]{{\boldsymbol{\mathfrak{X}}}^{\star}\bbsingle{#1}}
\newcommand{\smoothmaps}[1]{\difclass{\infty}\bbsingle{#1}}
\newcommand{\lieder}[2]{\mathcal{L}_{#1}{#2}}
\newcommand{\con}[2]{\nabla_{#1}{#2}}
\newcommand{\rescon}[2]{{#1}^{\(#2\)}}
\newcommand{\vbsecc}[1]{\zeta_{#1}}
\newcommand{\TBcon}[1]{\widetilde{#1}}
\newcommand{\tbcon}[2]{\widetilde{\nabla}_{#1}{#2}}
\newcommand{\connectionform}[3]{\leftidx{^{#1}}\Lambda^{#2}_{#3}}
\newcommand{\connectionformch}[4]{\leftidx{^{#1}}\Lambda^{#2}_{{#3}{#4}}}

\newcommand{\TFF}[3]{{\mathfrak{TF}}^{\infty}_{\opair{#1}{#2}}\bbsingle{#3}}
\newcommand{\seccof}[1]{s^{#1}}
\newcommand{\veccof}[1]{X^{#1}}
\newcommand{\veccoff}[1]{Y^{#1}}
\newcommand{\localframevecf}[1]{\partial_{#1}}
\newcommand{\localframeoneform}[1]{{\mathrm{dx}}_{#1}}
\newcommand{\Localframevecf}[2]{\partial_{#1}^{#2}}
\newcommand{\ltchart}[1]{\widetilde{#1}}
\newcommand{\Christoffel}[4]{\leftidx{^{#1}}\Gamma^{#2}_{{#3}{#4}}}
\newcommand{\interval}[1]{\mathscr{I}_{#1}}
\newcommand{\Curve}[1]{\gamma_{#1}}
\newcommand{\valongc}[1]{\xi_{#1}}
\newcommand{\Valongc}[2]{\xi_{#1}^{#2}}
\newcommand{\valongcurves}[2]{{\mathsf{V}\negthinspace\mathfrak{c}}^{\infty}\bbpair{#1}{#2}}
\newcommand{\parallelvalongcurves}[2]{{\mathscr{P}\mathsf{V}\negthinspace\mathfrak{c}}^{\infty}\bbpair{#1}{#2}}
\newcommand{\covder}[2]{{\mathcal{D}}_{#1}^{#2}}
\newcommand{\indexset}[1]{\mathfrak{I}_{#1}}
\newcommand{\pvfalongc}[4]{{\mathscr{PV}}_{#4}\bbpair{#1}{\binary{#2}{#3}}}
\newcommand{\ptransport}[4]{{\mathfrak{P}}_{#4}\bbpair{#1}{\binary{#2}{#3}}}
\newcommand{\framealongc}[1]{E_{#1}}
\newcommand{\geodesics}[2]{{\mathsf{GD}}\bbpair{#1}{#2}}
\newcommand{\dualEucbase}[2]{{\overline{\boldsymbol{e}}}^{#2}_{\(#1\)}}
\newcommand{\maxgeodesic}[3]{{\mathfrak{g}}\bbpair{\binary{#1}{#2}}{#3}}
\newcommand{\vectorspace}[1]{\mathbb{E}_{#1}}
\newcommand{\scalarprod}[1]{{\mathcal{S}}_{#1}}
\newcommand{\scalarprodf}[3]{{\mathcal{S}}_{#1}\opair{#2}{#3}}
\newcommand{\scalarprodflat}[1]{{#1}^{\flat}}
\newcommand{\scalarprodsharp}[1]{{#1}^{\sharp}}
\newcommand{\scalarprodmatrix}[2]{{\[#1\]}_{#2}}
\newcommand{\dualvsbase}[1]{{\overline{\mathbf{e}}}^{#1}}
\newcommand{\scalarprodindex}[2]{{\mathrm{index}}_{#1}\bbsingle{#2}}
\newcommand{\scalarprods}[2]{{\mathsf{SC}}_{#2}\bbsingle{#1}}
\newcommand{\idmatrix}[1]{{\mathbf{I}}_{#1\times #1}}
\newcommand{\zeromatrix}[2]{{\boldsymbol{0}}_{#1\times #2}}
\newcommand{\metrictensor}[1]{{\mathbf{g}}_{#1}}
\newcommand{\metrictensors}[2]{{\mathsf{SRM}}_{#2}\bbsingle{#1}}
\newcommand{\mtflat}[1]{{#1}^{\flat}}
\newcommand{\mtsharp}[1]{{#1}^{\sharp}}
\newcommand{\vfmetricprodmap}[2]{\langle,\rangle_{\opair{#1}{#2}}}
\newcommand{\vfmetricprod}[4]{\left\langle{#3},\thinspace{#4}\right\rangle_{\opair{#1}{#2}}}
\newcommand{\vfmetricproduct}[2]{\left\langle{#1},\thinspace{#2}\right\rangle}
\newcommand{\Rderivative}[2]{{\mathnormal{d}}_{#1}{#2}}
\newcommand{\grad}[1]{{\mathrm{grad}}_{#1}}
\newcommand{\aoneform}[1]{\eta_{#1}}
\newcommand{\oneform}[1]{\omega{#1}}
\newcommand{\TTF}[3]{{\mathscr{TF}}^{\infty}_{\opair{#1}{#2}}\bbsingle{#3}}
\newcommand{\TFequiv}[1]{\overline{#1}}
\newcommand{\TTFequiv}[1]{\overline{\overline{#1}}}
\newcommand{\TFequivv}[1]{\widetilde{#1}}
\newcommand{\TTFequivv}[1]{\widetilde{#1}}
\newcommand{\Torsion}[1]{{\mathrm{Tor}}_{#1}}
\newcommand{\LCconnection}[2]{\overset{\underset{\mathrm{\mathfrak{\opair{#1}{#2}}}}{}}{\nabla}}
\newcommand{\LCcon}[4]{{\overset{\underset{\mathrm{\mathfrak{\opair{#1}{#2}}}}{}}{\nabla}}_{#3}#4}
\newcommand{\metrictensorchart}[3]{{\mathrm{g}}^{#1}_{#2#3}}
\newcommand{\metrictensorchartinv}[3]{{\widetilde{\mathrm{g}}}^{#1}_{#2#3}}
\newcommand{\connectionpb}[1]{{#1}^{\hookleftarrow}}
\newcommand{\sectionpushforward}[1]{{#1}^{\hookrightarrow}}
\newcommand{\tfconnection}[3]{\leftidx{_{#2}^{#3}}{#1}}
\newcommand{\tfcon}[4]{\leftidx{_{#1}^{#2}}\nabla_{#3}#4}
\newcommand{\vstrace}[1]{{\mathrm{tr}}_{#1}}
\newcommand{\Tr}[1]{{\mathrm{Tr}}_{#1}}
\newcommand{\contraction}[5]{{\mathbf{C}}{#5}_{\opair{#1}{#3}}^{\opair{#2}{#4}}}
\newcommand{\mendomorphism}[2]{{\mathsf{End}}_{#1}{#2}}
\newcommand{\curvature}[1]{{\mathbf{R}}_{#1}}
\newcommand{\mcurvature}[1]{{\mathrm{R}}_{#1}}
\newcommand{\triplet}[3]{\binary{#1}{\binary{#2}{#3}}}
\newcommand{\Quadruple}[4]{\binary{#1}{\triplet{#2}{#3}{#4}}}
\newcommand{\bbinary}[2]{#1;\thinspace#2}

\newcommand{\pCurve}[1]{\Gamma_{#1}}
\newcommand{\pcurve}[3]{{#1}^{\bbsingle{#2}}_{#3}}
\newcommand{\pcovder}[3]{\leftidx{^{#3}}{{\boldsymbol{\mathcal{D}}}}_{#1}^{#2}}
\newcommand{\Riemcurvature}[2]{{\mathrm{Riem}}_{\opair{#1}{#2}}}
\newcommand{\Riemcurvaturecoef}[5]{{\mathrm{Rm}}^{#1}_{{#2}{#3}{#4}{#5}}}
\newcommand{\Ricci}[2]{{\mathrm{Ric}}_{\opair{#1}{#2}}}
\newcommand{\Riccicoef}[3]{{\mathrm{Rc}}^{#1}_{{#2}{#3}}}
\newcommand{\scalarcurvature}[2]{{\mathrm{Sc}}_{\opair{#1}{#2}}}
\newcommand{\semiEucspace}[2]{\R^{#1}_{#2}}
\newcommand{\semiEucmetric}[2]{{\mathfrak{E}}^{#1}_{#2}}
\newcommand{\Eucstandardframe}[2]{{\mathbf{E}}^{\(#1\)}_{#2}}
\newcommand{\isometries}[1]{{\mathbf{Isom}}\bbsingle{#1}}
\newcommand{\speed}[3]{{\mathrm{Speed}}_{\opair{#1}{#2}}\bbsingle{#3}}
\newcommand{\length}[5]{{\mathrm{Length}}_{\opair{#1}{#2}}\bbpair{#3}{\binary{#4}{#5}}}
\def\TT{\mathcal{T}}

\newcommand{\maprule}[3]{#1:\thinspace#2\mapsto#3}
\def\timeder{\frac{\mathrm{d}}{\mathrm{dt}}}
\def\varfill{\dotfill}

\def\toclevel@section{1}\def\toclevel@subection{2}
\addtocontents{toc}{\string\let\string\l@section\string\l@subsection}
\def\toclevel@subsection{2}\def\toclevel@subsubection{3}
\addtocontents{toc}{\string\let\string\l@subsection\string\l@subsubsection}
\usepackage{tocloft}
\newcommand{\chapteR}[1]{\cleardoublepage
{\refstepcounter{chapterr}\vskip\baselineskip\centering{\fontsize{21}{21}\selectfont${\bf\Roman{chapterr}}$
\vskip0.6\baselineskip{\fontsize{21}{21}\selectfont{\textbf{#1}}}}
\vskip 5.9\baselineskip}
\addcontentsline{toc}{0}
{\protect\vskip0.5\baselineskip\noindent\bf{\Roman{chapterr}\hskip0.5\baselineskip#1\hspace{\fill}}}\par
\fancyhead[LO]{\ifthenelse{\value{chapterr}=0}{#1}{$\bf\Roman{chapterr}$\hskip0.7\baselineskip#1}}
}
\newcommand{\Bibliography}[1]{\vskip0.5\baselineskip\centering{\huge\bf{References}}\vskip \baselineskip
\addcontentsline{toc}{0}
{\protect\vskip0.5\baselineskip\noindent\bf{References}\hspace{\fill}}\par
\fancyhead[LO,RE]{\bf{References}#1}
}

\def\refthm#1{[{theorem~\ref{#1}}]}
\def\reflem#1{[{lemma~\ref{#1}}]}
\def\refdef#1{[{definition~\ref{#1}}]}
\def\refcor#1{[{corollary~\ref{#1}}]}
\def\refpro#1{[{proposition~\ref{#1}}]}

\newcommand{\mathleft}{\@fleqntrue\@mathmargin0pt}
\renewcommand{\sectionmark}[1]{\ifthenelse{\value{section}=0}{\markright{#1}{}}
{\markright{${\arabic{section}}$ #1}{}}}
\makeatletter
\def\cleardoublepage{\clearpage\if@twoside \ifodd\c@page\else
\hbox{}
\vspace*{\fill}
\vspace{\fill}
\thispagestyle{empty}
\newpage
\if@twocolumn\hbox{}\newpage\fi\fi\fi}
\makeatother
\fancypagestyle{plain}{
\fancyhf{}
\fancyhead{}
\fancyhead[R]{\thepage}
\fancyhead[L]{\leftmark}
}
\newcommand{\newsymp}[1]{{#1}\equiv}
\newcommand{\newsymb}[1]{{#1}\equiv}

\newcommand{\SET}[1]{{\mathscr{S}}_{#1}}
\newcommand{\FUNCTION}[1]{{{f}}_{#1}}

\newcommand{\prop}[1]{{\mathfrak{p}}\llparenthesis{#1}\rrparenthesis}
\def\Prop{\mathfrak{p}}
\newcommand{\propos}[1]{{\mathbf{\mathfrak{p}}}_{#1}}
\def\dummy{\centerdot}
\usepackage[hang,flushmargin]{footmisc}
\renewcommand{\footnoterule}{%
  \kern 20pt
  \hrule width \textwidth height 0.5pt
  \kern 5pt
}
\def\quotl{``}
\def\quotr{"}
\usepackage[inner=3.8cm,outer=3cm,bottom=3.9cm,twoside]{geometry}
\begin{document}
\thispagestyle{empty}
\noindent
{\\ \\\textbf{\fontsize{30}{30}\selectfont
{\textsf{Semi-Riemannian}}}}
\\[0.8\baselineskip]{\textbf{\fontsize{30}{30}\selectfont
{\textsf{Structures on Manifolds}}}}
\\[8\baselineskip]
\noindent
{\fontsize{21}{21}\selectfont
{\textsf{Farzad Shahi}}
}
\\[11\baselineskip]
\noindent
{\fontsize{11}{11}\selectfont
{\textrm{Version:} 1.00}
}
%%%%%%%%%%%%%%%%%%%%%%%%%%%%%%%%%%%%%%%%%
\newpage
\thispagestyle{empty}
\noindent
{\fontsize{9.4}{9.4}\selectfont
{\underline{{\bf\textsf{Title:}} Semi-Riemannian Structures on Manifolds}}}\\
{\fontsize{9.4}{9.4}\selectfont
{\underline{{\bf\textsf{Author:}} Farzad Shahi}}}\\
{{\bf\textsf{email}}}:~{\texttt{shahi.farzad@gmail.com}}
\vskip 4\baselineskip
\noindent
{\fontsize{9.4}{9.4}\selectfont
{\underline{{\bf\textsf{Version:}} 1.00}}
}\\
\noindent
{\fontsize{9.4}{9.4}\selectfont
{\bf\textsf{2022}}}
\vskip 5\baselineskip
\noindent
{\fontsize{9.4}{9.4}\selectfont
{\bf\textsf{Abstract:}}
This text serves as an introduction to the basic concepts of the theory of semi-Riemannian geometry
on real finite-dimensional manifolds without boundary.
}
%%%%%%%%%%%%%%%%%%%%%%%%%%%%%%%%%%%%%%%%%%
\newpage
\thispagestyle{empty}
\noindent
{\fontsize{20}{20}\selectfont
{\bf Preface}}
\vskip1.5\baselineskip
\noindent
Euclidean spaces may be regarded as the natural settings that give rise to fundamental and intuitive geometric concepts.
The simplest case is when the whole of an Euclidean space is regarded as the ambient space.
But, taken to consider a Euclidean space of some dimension, like the familiar three-dimensional space,
there is a natural quest for generalizing the geometric concepts of interest from that Euclidean space
to its submanifolds, like hypersurfaces of the three-dimensional Euclidean space.\\
The geometric concepts of distance and angle between stright lines in a Euclidean space can be
expressed in terms of the standard inner-product operation on a Euclidean space.
These geometric ideas can be generalized to arbitrary vector-spaces by abstracting the standard inner-product operation
on Euclidean spaces to introduce the general notion of an inner-product on a vector-space.
One further level of abstraction generalizes the notion of inner product on finite-dimensional vector-spaces to
smooth manifolds by bringing it from local to global. This generalized notion of an inner-product on a smooth manifold
is famously known as the concept of a $\quotl$Riemannian metric$\quotr$ on a smooth manifold. Precisely, a
Riemannian metric on a manifold is an assignment of an inner-product, in the sense of vector-spaces,
to every tangent-space of the manifold, which is additionally assumed to be smooth. So, since
an inner-product on a vector-space is a bilinear form on it, a Riemannian metric on a manifold is actually
a specific type of a $2$-covariant smooth tensor field on its tengent-bundle, which yields a positive-definite
bilinear form at the tangent-space of each point. A manifold equipped with a Riemannian metric on it is called a
$\quotl$Riemannian manifold$\quotr$.\\
Two other well-known geometric notions of Euclidean spaces are the velocity and acceleration of a smooth
parametrized curve in a Euclidean spaces, which are of crucial importance in studying the motion of particles
in the three-dimensional space. The velocity of a smooth curve in a Euclidean space
is a smooth map from the domain of the curve (assumed to be an open interval of $\R$, regarded as a time interval) to the Euclidean
space under consideration, and its acceleration is defined to be the derived function of velocity, which again
is a smooth map from the domain of the curve to the Euclidean space.
Fixing a smooth curve, the velocity function of the curve is a specific example of the general class of
smooth maps from the domain of the curve to the underlying Euclidean space, whose instances are called
the smooth vector fields along the curve on the underlying Euclidean space.
There is an equivalent way of looking at
these notions in Euclidean spaces. When the ambient Euclidean space is regarded to be equipped with its canonical
differentiable structure, the velocity function of a smooth curve or any smooth vector field along it on a Euclidean space
can be thought of as a smooth map from the domain of the curve to the tangent-bundle of the Euclidean space,
sending each instant of the time interval to the tangent-space of the point that the curve passes through at that instant.
Although this way of thinking about the velocity function or any other smooth vector field along the curve
is absolutely unnecessary in the scope of Euclidean spaces,
the natural generalization of the notion of a smooth vector field along a curve from the category of Euclidean spaces
to the category of manifolds stems from it in an straightforward manner. Given a smooth manifold and a smooth curve on it,
a smooth vector field along the curve is defined to be a smooth assignmet of vectors from the tangent-bundle of the manifold
to every instant of the time interval, in such a way that the vector assigned to any instant lies in the tangent-space of the
point that the curve passes through at that instant.\\
Now let us consider the notion of acceleration. As mentioned previously, the acceleration of a smooth curve in a Euclidean space,
by definition, is the derived function of the velocity of that curve. This can be extended to the general problem of differentiating
any smooth vector field along a smooth curve on a Euclidean space, which yields again a smooth vector field along the curve.
When the Euclidean space is viewed as a manifold
and the smooth vector field along the curve is viewed as a smooth map from the domain of the curve to the tangent-bundle
of the Euclidean space, again there is no obstacle to define the differential of the vector field along the curve in this framework,
basically being equivalent to the ordinary differential. The reason is that the tangent-spaces of all points are
naturally identified with the Euclidean space itself, and thus all can be viewed as copies of the underlying Euclidean space,
and therefore the difference of any pair of vectors lying in different tangent-spaces makes sense. But on general manifolds
there is no natural identification of tangent-spaces of different points, and hence the idea of differentiating
smooth vector fields along smooth curves on general manifolds cannot be defined exactly as the same way in the category
of Euclidean manifolds.\\
To deal with the problem of differentiating smooth vector fields along curves on general manifolds,
it is inspiring to first consider the case of submanifolds of Euclidean spaces. Actually based on the
\textit{strong Whitney embedding theorem}, every smooth manifold can be embedded into a Euclidean space, and thus
can be viewed as a submanifold of a Eiclidean space; so, precisely, we mean to consider first the case that
the differentiable structure of the manifold is the inherited one from the canonical differentiable structure of a
particular Euclidean space.
Given a submanifold of a Euclidean space, the tangent-space of the manifold at each point can be identified
with a unique subspace of the tangent-space of the ambient Euclidean space at that point. Hence, the tangent-space of
the manifold at each point can naturally be viewd as a subspace of the ambient Euclidean space itself.
Given a smooth curve on the submanifold, trivially it is also a smooth curve in the ambient Euclidean space.
Additionally, considering the tangent-space at any point of the submanifold as a subspace of the tangent-space of the Euclidean space
at that point, any smooth vector field along the curve on the submanifold is also canonically a
smooth vector field along the curve in the ambient Euclidean space. So, given a smooth vector field along the curve
on the submanifold, differentiating it in the ambient Euclidean space yields a new smooth vector field
along the curve in the Euclidean space. Then by projecting the value of this newly achieved vector field along the curve
at any instant, to the tangent-space of the submanifold at the point through which the curve passes at that instant,
a new smooth vector field along the curve on the submanifold is obtained. Note that the projection operation on the tangent-spaces
stated as above is absolutely determined by the standard inner-product of the ambient Euclidean space. This can be regarded as the natural
way of differentiating smooth vector fields along curves on submanifolds of Euclidean spaces, and the resulting
smooth vector field along the curve on the submanifold thus achieved is called the
$\quotl$covariant derivative of the initially given vector field along the curve on the submanifold$\quotr$.
So, for every smooth curve on the submanifold, an operation is obtained on the set of all vector fields along the curve,
sending each one to its covariant derivative. Such an operation is called the $\quotl$covariant differentiation along the
given curve$\quotr$.
Taking into account all covariant differentiations along all smooth curves on the submanifold, a binary operation
on the set of all smooth vector fields on the submanifold is determined, which is famous for the $\quotl$tangential connection of the submanifold
(with respect to the ambient Euclidean space)$\quotr$.
\cite[Chap. 7]{Boothby} provides a rigorous explanation of this procedure.\\
The notion of connection on submanifolds of an Euclidean space, as mentioned above, satisfies a collection of properties
that are conventionally chosen to define an abstract notion of a connection on a general smooth manifold in an axiomatizing way.
Additionally, an abstract connection on a manifold determines uniquely an operation on the set of all smooth vector fields along a curve on the
manifold via a collection of axioms that are all satisfied in the case of submanifolds of Euclidean spaces explained above.
This operation is also referred to as the $\quotl$covariant differentiation along a given smooth curve on the manifold
relative to the given connection$\quotr$.
So, in general, the covariant differentiation along a curve on a manifold becomes dependent on the choice of
a particular connection on that manifold.\\
Given a submanifold of an Euclidean space, the inner-product of the Euclidean space induces
an inner-product on the tangent-space of the submanifold at every point. Collecting these inner-products on tangent spaces
gives rise to a Riemannian metric on the submanifold which is called the $\quotl$induced Riemannian metric on the submanifold
from that of the ambient Euclidean space$\quotr$.
There are two further properties that the tangential connection of
a submanifold of an Euclidean space satisfies which are famous for the $\quotl$metric compatibility$\quotr$ and
$\quotl$torsion-lessness$\quotr$ properties. The metric compatibility property can be asserted in terms of the
induced Riemannian metric on the submanifold. The interesting problem is that given an arbitrary smooth manifold
and an arbitrary Riemannian metric on it, there exists a unique connection on the manifold that satisfies the
torsion-lessness property together with the metric compatibility with respect to the Riemannian metric on it.
This unique connection is called the $\quotl$Levi-Civita connection of the considered Riemannian manifold$\quotr$.
Another interesting problem is that according to \textit{Nash embedding theorem}, given a Riemannian manifold,
there can be found an Euclidean space such that the considered Riemannian manifold can be embedded into that Euclidean
space and the induced Riemannian metric on it from that of the Euclidean space coincides with the manifold's own Riemannian metric.
Then considering this manifold as a submanifold of that particular Euclidean space, the tangential connection
of this submanifold with respect to this Euclidean space must agree with the Levi-Civita connection of the given
Riemannian manifold. Informally speaking, the Levi-Civita connection of a Riemannian manifold is nothing more than a
tangential connection basically.\\
The notion of Riemannian metrics on manifolds can be generalized to the notion of semi-Riemannian
metrics on manifolds, in the same way that inner-products on vector-spaces are generalized to scalar-products
on vector-spaces. Any semi-Riemannian manifold also accepts a unique connection that satisfies the
torsion-lessness condition and the metric compatibility with respect to the semi-Riemannian metric given to the manifold.
This unique connection is also called the $\quotl$Levi-Civita connection of the given semi-Riemannian manifold$\quotr$.\\
The concept of the acceleration of a smooth curve in a Euclidean space can be generalized to any Riemannian or
semi-Riemannian manifold, or even a manifold equipped with a connection. The generalization takes place
on the basis of the notion of covariant differentiation along the given curve relative to the connection,
and is defined to be the covariant derivative of the velocity of the curve along the curve itself
(relative to the connection). Any smooth curve with vanishing acceleration is defined to be a $\quotl$geodesic$\quotr$.
So, the concept of a geodesic is in some sense the generalization of the idea of the motion on a straight line with constant speed
in a Euclidean space. The theory of classical mechanics actually makes this generalization much more sensible;
the motion of a particle constrained to move on a (frictionless) hypersurface must be a geodesic on the hypersurface in the absence of
any external forces.\\
An important geometric notion in the category of semi-Riemannian manifolds is the idea of curvature.
The importance of this notion relates to its vastly diverse pool of applications and interpretations,
connecting it to other geometrical and also topological concepts, which are absolutely beyond the
reach of this text. A geometric interpretation of curvature which is studied in this text, relates it
to the idea of parallel transportation of vectors of the tangent-bundle of a semi-Riemannian manifold
along closed smooth curves. In particular, vanishing of curvature is equivalent to the fact that
any vector of the tangent-bundle remains intact under the parallel transportation along any smooth closed curve.
Another important fact about curvature that is studied here, is that it is an invariant of local-isometry,
which precisely means that the curvature of two locally-isometric semi-Riemannian manifolds are related naturally.
There are many derivations of the main curvature operation defined initially, each carrying a part of
information inherent in the main curvature, and having its own interpretations, applications, and connections
to other concepts. Two important such derivations are the Ricci and scalar curvatures
that play a direct role in the theory of general relativity.\\
Also it is worthwhile to mention that
the notion of a connection on a manifold can further be generalized to the category of smooth vector bundles.
This generalization introduces a notion of differentiating the sections of a smooth vector bundle with respect
to the vector fields of its base space.\\
%%%%%%%%%%%%%
\vskip0.5\baselineskip
\noindent
{\bf{A brief explanation about chapters:}}
The first chapter begins with introducing the abstract notion of connection in the most general setting
of a smooth vector bundle. Then the abstract notions of covariant differentiation of vector fields along curves,
parallel transportation, and geodesics are introduced and studied on manifolds equipped with a connection.
In this chapter, manifolds are not assumed to be equipped with any semi-Riemannian metric and all geometrical ideas
are introduced and studied relative to an arbitrary connection on the manifold.\\
In the second chapter, the manifolds are assumed to have a semi-Riemannian metric. Apart from studying
some basic aspects of semi-Riemannian metrics on manifolds, the existance of a unique connection on a manifold in alignment
with its semi-Riemannian metric via the metric compatibility property is studied. This connection is usually assumed implicitly
as the considered connection on a semi-Riemannian manifold, and all concepts introduced in the first chapter can
be automatically carried over the semi-Riemannian manifolds.\\
The third chapter is devoted to the introduction of the abstract notion of curvature on a manifold with connection,
a basic geometric interpretation of it, and the introduction of some of its derivations like
Riemann, Ricci, and scalar curvatures on semi-Riemannian manifolds.
\vskip0.5\baselineskip
\noindent
{\bf{Prerequisites:}}
A sufficient knowledge of linear algebra, general topology, basic differential geometry,
and the theory of smooth vector bundles, along some elementary-level knowledge of some algebraic
concepts like rings and modules is assumed.
\vskip\baselineskip
\hfill
{\textsf{Farzad Shahi}}
%%%%%%%%%%%%%%%%%%%%%%%%%%%%
%%%%%%%%%%%%%%%%%%%%%%%%%%%%
%%%%%%%%%%%%%%%%%%%%%%%%%%%%%%%%%%%%%%%%%%
%\tableofcontents
\newpage
\thispagestyle{empty}
%\Thecontents{}
\section*{\fontsize{21}{21}\selectfont\bf{Contents}}
\addtocontents{toc}{\protect\setcounter{tocdepth}{-1}}
\tableofcontents
\addtocontents{toc}{\protect\setcounter{tocdepth}{3}}
%\tableofcontents
\newpage
%%%%%%%%%%%%%%%%%%%%%%%%%%%%%%%%%%%%%%%%%%%%%%%%%%%%%%%%%%
%%%%%%%%%%%%%%%%%%%%%%%%%%%%%%%%%%%%%%%%%%%%%%%%%%%%%%%%%%
\newpage
\chapteR{
Mathematical Symbols
}
\thispagestyle{fancy}
\section*{
Set-theory
}
%%%%%%%%%%%%%%%%%%%%%%%%%%%%%%%%%%%%%%%%%%%%%%%%%%%%%%%%%%%%%%%%%%%%%%%%%%%%%%%%%%%%%%%%
$\newsymb{\empty}$
empty-set
\varfill
$\empty$\\
%%%%%%%%%%%%%%%%%%%%%%%%%%%%%%%%%%%%%%%%%%%%%%%%%%%%%%%%%%%%%%%%%%%%%%%%%%%%%%%%%%%%%%%%
$\newsymp{\SET{1}=\SET{2}}$
$\SET{1}$
equals
$\SET{2}$.
\varfill
$=$\\
%%%%%%%%%%%%%%%%%%%%%%%%%%%%%%%%%%%%%%%%%%%%%%%%%%%%%%%%%%%%%%%%%%%%%%%%%%%%%%%%%%%%%%%%
%%%%%%%%%%%%%%%%%%%%%%%%%%%%%%%%%%%%%%%%%%%%%%%%%%%%%%%%%%%%%%%%%%%%%%%%%%%%%%%%%%%%%%%%
$\newsymp{\SET{1}\in\SET{2}}$
$\SET{1}$
is an element of
$\SET{2}$.
\varfill
$\in$\\
%%%%%%%%%%%%%%%%%%%%%%%%%%%%%%%%%%%%%%%%%%%%%%%%%%%%%%%%%%%%%%%%%%%%%%%%%%%%%%%%%%%%%%%%
$\newsymp{\SET{1}\ni\SET{2}}$
$\SET{1}$
contains
$\SET{2}$.
\varfill
$\ni$\\
%%%%%%%%%%%%%%%%%%%%%%%%%%%%%%%%%%%%%%%%%%%%%%%%%%%%%%%%%%%%%%%%%%%%%%%%%%%%%%%%%%%%%%%%
$\newsymb{\seta{\binary{\SET{1}}{\SET{2}}}}$
the set composed of
$\SET{1}$
and
$\SET{2}$
\varfill
$\seta{\binary{\dummy}{\dummy}}$\\
%%%%%%%%%%%%%%%%%%%%%%%%%%%%%%%%%%%%%%%%%%%%%%%%%%%%%%%%%%%%%%%%%%%%%%%%%%%%%%%%%%%%%%%%
$\newsymb{\defset{\SET{1}}{\SET{2}}{\prop{\SET{1}}}}$
all elements of
$\SET{2}$
having the property
$\Prop$
\varfill
$\defset{\dummy}{\dummy}{\dummy}$\\
%%%%%%%%%%%%%%%%%%%%%%%%%%%%%%%%%%%%%%%%%%%%%%%%%%%%%%%%%%%%%%%%%%%%%%%%%%%%%%%%%%%%%%%%
$\newsymb{\union{\SET{}}}$
union of all elements of
$\SET{}$
\varfill
$\bigcup$\\
%%%%%%%%%%%%%%%%%%%%%%%%%%%%%%%%%%%%%%%%%%%%%%%%%%%%%%%%%%%%%%%%%%%%%%%%%%%%%%%%%%%%%%%%
$\newsymb{\intersection{\SET{}}}$
intersection of all elements of
$\SET{}$
\varfill
$\bigcap$\\
%%%%%%%%%%%%%%%%%%%%%%%%%%%%%%%%%%%%%%%%%%%%%%%%%%%%%%%%%%%%%%%%%%%%%%%%%%%%%%%%%%%%%%%%
$\newsymb{\CSs{\SET{}}}$
power-set of
$\SET{}$
\varfill
$\CSs{\dummy}$\\
%%%%%%%%%%%%%%%%%%%%%%%%%%%%%%%%%%%%%%%%%%%%%%%%%%%%%%%%%%%%%%%%%%%%%%%%%%%%%%%%%%%%%%%%
$\newsymb{\Dproduct{\alpha}{\index}{\SET{\alpha}}}$
Cartesian-product of the collection of indexed sets
${\seta{\SET{\alpha}}}_{\alpha\in\index}$
\varfill
$\prod$\\
%%%%%%%%%%%%%%%%%%%%%%%%%%%%%%%%%%%%%%%%%%%%%%%%%%%%%%%%%%%%%%%%%%%%%%%%%%%%%%%%%%%%%%%%
$\newsymp{\SET{1}\subseteq\SET{2}}$
$\SET{1}$
is a subset of
$\SET{2}$.
\varfill
$\subseteq$\\
%%%%%%%%%%%%%%%%%%%%%%%%%%%%%%%%%%%%%%%%%%%%%%%%%%%%%%%%%%%%%%%%%%%%%%%%%%%%%%%%%%%%%%%%
$\newsymp{\SET{1}\supseteq\SET{2}}$
$\SET{1}$
includes
$\SET{2}$.
\varfill
$\supseteq$\\
%%%%%%%%%%%%%%%%%%%%%%%%%%%%%%%%%%%%%%%%%%%%%%%%%%%%%%%%%%%%%%%%%%%%%%%%%%%%%%%%%%%%%%%%
$\newsymp{\SET{1}\subset\SET{2}}$
$\SET{1}$
is a proper subset of
$\SET{2}$.
\varfill
$\subset$\\
%%%%%%%%%%%%%%%%%%%%%%%%%%%%%%%%%%%%%%%%%%%%%%%%%%%%%%%%%%%%%%%%%%%%%%%%%%%%%%%%%%%%%%%%
$\newsymp{\SET{1}\supset\SET{2}}$
$\SET{1}$
properly includes
$\SET{2}$.
\varfill
$\supset$\\
%%%%%%%%%%%%%%%%%%%%%%%%%%%%%%%%%%%%%%%%%%%%%%%%%%%%%%%%%%%%%%%%%%%%%%%%%%%%%%%%%%%%%%%%
$\newsymb{\SET{1}\cup\SET{2}}$
union of
$\SET{1}$
and
$\SET{2}$
\varfill
$\cup$\\
%%%%%%%%%%%%%%%%%%%%%%%%%%%%%%%%%%%%%%%%%%%%%%%%%%%%%%%%%%%%%%%%%%%%%%%%%%%%%%%%%%%%%%%%
$\newsymb{\SET{1}\cap\SET{2}}$
intersection of
$\SET{1}$
and
$\SET{2}$
\varfill
$\cap$\\
%%%%%%%%%%%%%%%%%%%%%%%%%%%%%%%%%%%%%%%%%%%%%%%%%%%%%%%%%%%%%%%%%%%%%%%%%%%%%%%%%%%%%%%%
$\newsymb{\SET{1}\times\SET{2}}$
Cartesian-product of
$\SET{1}$
and
$\SET{2}$
\varfill
$\times$\\
%%%%%%%%%%%%%%%%%%%%%%%%%%%%%%%%%%%%%%%%%%%%%%%%%%%%%%%%%%%%%%%%%%%%%%%%%%%%%%%%%%%%%%%%
$\newsymb{\compl{\SET{1}}{\SET{2}}}$
the relative complement of
$\SET{2}$
with respect to
$\SET{1}$
\varfill
$\setminus$\\
%%%%%%%%%%%%%%%%%%%%%%%%%%%%%%%%%%%%%%%%%%%%%%%%%%%%%%%%%%%%%%%%%%%%%%%%%%%%%%%%%%%%%%%%
$\newsymb{\func{\FUNCTION{}}{\SET{}}}$
value of the function
$\FUNCTION{}$
at
$\SET{}$
\varfill
$\func{\dummy}{\dummy}$\\
%%%%%%%%%%%%%%%%%%%%%%%%%%%%%%%%%%%%%%%%%%%%%%%%%%%%%%%%%%%%%%%%%%%%%%%%%%%%%%%%%%%%%%%%
$\newsymb{\domain{\FUNCTION{}}}$
domain of the function
$\FUNCTION{}$
\varfill
$\domain{\dummy}$\\
%%%%%%%%%%%%%%%%%%%%%%%%%%%%%%%%%%%%%%%%%%%%%%%%%%%%%%%%%%%%%%%%%%%%%%%%%%%%%%%%%%%%%%%%
$\newsymb{\codomain{\FUNCTION{}}}$
codomain of the function
$\FUNCTION{}$
\varfill
$\codomain{\dummy}$\\
%%%%%%%%%%%%%%%%%%%%%%%%%%%%%%%%%%%%%%%%%%%%%%%%%%%%%%%%%%%%%%%%%%%%%%%%%%%%%%%%%%%%%%%%
$\newsymb{\funcimage{\FUNCTION{}}}$
image of the function
$\FUNCTION{}$
\varfill
$\funcimage{\dummy}$\\
%%%%%%%%%%%%%%%%%%%%%%%%%%%%%%%%%%%%%%%%%%%%%%%%%%%%%%%%%%%%%%%%%%%%%%%%%%%%%%%%%%%%%%%%
$\newsymb{\image{\FUNCTION{}}}$
image-map of the function
$\FUNCTION{}$
\varfill
$\image{\dummy}$\\
%%%%%%%%%%%%%%%%%%%%%%%%%%%%%%%%%%%%%%%%%%%%%%%%%%%%%%%%%%%%%%%%%%%%%%%%%%%%%%%%%%%%%%%%
$\newsymb{\pimage{\FUNCTION{}}}$
inverse-image-map of the function
$\FUNCTION{}$
\varfill
$\pimage{\dummy}$\\
%%%%%%%%%%%%%%%%%%%%%%%%%%%%%%%%%%%%%%%%%%%%%%%%%%%%%%%%%%%%%%%%%%%%%%%%%%%%%%%%%%%%%%%%
$\newsymb{\resd{\FUNCTION{}}}$
domain-restriction-map of the function
$\FUNCTION{}$
\varfill
$\resd{\dummy}$\\
%%%%%%%%%%%%%%%%%%%%%%%%%%%%%%%%%%%%%%%%%%%%%%%%%%%%%%%%%%%%%%%%%%%%%%%%%%%%%%%%%%%%%%%%
$\newsymb{\rescd{\FUNCTION{}}}$
codomain-restriction-map of the function
$\FUNCTION{}$
\varfill
$\rescd{\dummy}$\\
%%%%%%%%%%%%%%%%%%%%%%%%%%%%%%%%%%%%%%%%%%%%%%%%%%%%%%%%%%%%%%%%%%%%%%%%%%%%%%%%%%%%%%%%
$\newsymb{\func{\res{\FUNCTION{}}}{\SET{}}}$
domain-restriction and codomain-restriction of\\ the function
$\FUNCTION{}$ to $\SET{}$ and $\func{\image{\FUNCTION{}}}{\SET{}}$, respectively
\varfill
$\res{\dummy}$\\
%%%%%%%%%%%%%%%%%%%%%%%%%%%%%%%%%%%%%%%%%%%%%%%%%%%%%%%%%%%%%%%%%%%%%%%%%%%%%%%%%%%%%%%%
$\newsymb{\Func{\SET{1}}{\SET{2}}}$
the set of all maps from
$\SET{1}$
to
$\SET{2}$
\varfill
$\Func{\dummy}{\dummy}$\\
%%%%%%%%%%%%%%%%%%%%%%%%%%%%%%%%%%%%%%%%%%%%%%%%%%%%%%%%%%%%%%%%%%%%%%%%%%%%%%%%%%%%%%%%
$\newsymb{\IF{\SET{1}}{\SET{2}}}$
the set of all bijective functions from
$\SET{1}$
to
$\SET{2}$
\varfill
$\IF{\dummy}{\dummy}$\\
%%%%%%%%%%%%%%%%%%%%%%%%%%%%%%%%%%%%%%%%%%%%%%%%%%%%%%%%%%%%%%%%%%%%%%%%%%%%%%%%%%%%%%%%
$\newsymb{\finv{\FUNCTION{}}}$
the inverse mapping of the bijective function $\FUNCTION{}$
\varfill
$\finv{\dummy}$\\
%%%%%%%%%%%%%%%%%%%%%%%%%%%%%%%%%%%%%%%%%%%%%%%%%%%%%%%%%%%%%%%%%%%%%%%%%%%%%%%%%%%%%%%%
$\newsymb{\surFunc{\SET{1}}{\SET{2}}}$
the set of all surjective functions from
$\SET{1}$
to
$\SET{2}$
\varfill
$\surFunc{\dummy}{\dummy}$\\
%%%%%%%%%%%%%%%%%%%%%%%%%%%%%%%%%%%%%%%%%%%%%%%%%%%%%%%%%%%%%%%%%%%%%%%%%%%%%%%%%%%%%%%%
$\newsymb{\cmp{\FUNCTION{1}}{\FUNCTION{2}}}$
composition of the function
$\FUNCTION{1}$
with the function
$\FUNCTION{2}$
\varfill
$\cmp{}{}$\\
%%%%%%%%%%%%%%%%%%%%%%%%%%%%%%%%%%%%%%%%%%%%%%%%%%%%%%%%%%%%%%%%%%%%%%%%%%%%%%%%%%%%%%%%
$\newsymb{\Injection{\SET{1}}{\SET{2}}}$
the injection-mapping of the set $\SET{1}$ into the set $\SET{2}$
\varfill
$\Injection{\dummy}{\dummy}$\\
%%%%%%%%%%%%%%%%%%%%%%%%%%%%%%%%%%%%%%%%%%%%%%%%%%%%%%%%%%%%%%%%%%%%%%%%%%%%%%%%%%%%%%%%
$\newsymb{\funcprod{\cf_1}{\cf_2}}$
the function-product of the function $\cf_1$ and $\cf_2$
\varfill
$\Cprod{\dummy}{\dummy}$\\
%%%%%%%%%%%%%%%%%%%%%%%%%%%%%%%%%%%%%%%%%%%%%%%%%%%%%%%%%%%%%%%%%%%%%%%%%%%%%%%%%%%%%%%%
$\newsymb{\EqR{\SET{}}}$
the set of all equivalence relations on the set
$\SET{}$
\varfill
$\EqR{\dummy}$\\
%%%%%%%%%%%%%%%%%%%%%%%%%%%%%%%%%%%%%%%%%%%%%%%%%%%%%%%%%%%%%%%%%%%%%%%%%%%%%%%%%%%%%%%%
$\newsymb{\EqClass{\SET{1}}{\SET{2}}}$
quotient-set of
$\SET{1}$
by the equivalence-relation
$\SET{1}$
\varfill
$\EqClass{\dummy}{\dummy}$\\
%%%%%%%%%%%%%%%%%%%%%%%%%%%%%%%%%%%%%%%%%%%%%%%%%%%%%%%%%%%%%%%%%%%%%%%%%%%%%%%%%%%%%%%%
$\newsymb{\pEqclass{\SET{1}}{\SET{2}}}$
equivalence-class of
$\SET{1}$
by the equivalence-relation
$\SET{2}$
\varfill
$\pEqclass{\dummy}{\dummy}$\\
%%%%%%%%%%%%%%%%%%%%%%%%%%%%%%%%%%%%%%%%%%%%%%%%%%%%%%%%%%%%%%%%%%%%%%%%%%%%%%%%%%%%%%%%
$\newsymb{\PEqclass{\SET{1}}{\SET{2}}}$
equivalence-class of
$\SET{2}$
by the equivalence-relation
$\SET{1}$
\varfill
$\PEqclass{\dummy}{\dummy}$\\
%%%%%%%%%%%%%%%%%%%%%%%%%%%%%%%%%%%%%%%%%%%%%%%%%%%%%%%%%%%%%%%%%%%%%%%%%%%%%%%%%%%%%%%%
$\newsymp{\Card{\SET{1}}\cardeq\Card{\SET{2}}}$
$\IF{\SET{1}}{\SET{2}}$
is non-empty.
\varfill
$\Card{\dummy}\cardeq\Card{\dummy}$\\
%%%%%%%%%%%%%%%%%%%%%%%%%%%%%%%%%%%%%%%%%%%%%%%%%%%%%%%%%%%%%%%%%%%%%%%%%%%%%%%%%%%%%%%%
$\newsymb{\CarD{\SET{}}}$
cardinality of
$\SET{}$
\varfill
$\CarD{\dummy}$
%%%%%%%%%%%%%%%%%%%%%%%%%%%%%%%%%%%%%%%%%%%%%%%%%%%%%%%%%%%%%%%%%%%%%%%%%%%%%%%%%%%%%%%%
%%%%%%%%%%%%%%%%%%%%%%%%%%%%%%%%%%%%%%%%%%%%%%%%%%%%%%%%%%%%%%%%%%%%%%%%%%%%%%%%%%%%%%%%
%%%%%%%%%%%%%%%%%%%%%%%%%%%%%%%%%%%%%%%%%%%%%%%%%%%%%%%%%%%%%%%%%%%%%%%%%%%%%%%%%%%%%%%%
%%%%%%%%%%%%%%%%%%%%%%%%%%%%%%%%%%%%%%%%%%%%%%%%%%%%%%%%%%%%%%%%%%%%%%%%%%%%%%%%%%%%%%%%
%%%%%%%%%%%%%%%%%%%%%%%%%%%%%%%%%%%%%%%%%%%%%%%%%%%%%%%%%%%%%%%%%%%%%%%%%%%%%%%%%%%%%%%%
\section*{
Logic
}
%%%%%%%%%%%%%%%%%%%%%%%%%%%%%%%%%%%%%%%%%%%%%%%%%%%%%%%%%%%%%%%%%%%%%%%%%%%%%%%%%%%%%%%%
$\newsymp{\AND{\propos{1}}{\propos{2}}}$
$\propos{1}$
and
$\propos{2}$.
\varfill
$\AND{}{}$\\
%%%%%%%%%%%%%%%%%%%%%%%%%%%%%%%%%%%%%%%%%%%%%%%%%%%%%%%%%%%%%%%%%%%%%%%%%%%%%%%%%%%%%%%%
$\newsymp{\OR{\propos{1}}{\propos{2}}}$
$\propos{1}$
or
$\propos{2}$.
\varfill
$\OR{}{}$\\
%%%%%%%%%%%%%%%%%%%%%%%%%%%%%%%%%%%%%%%%%%%%%%%%%%%%%%%%%%%%%%%%%%%%%%%%%%%%%%%%%%%%%%%%
$\newsymp{{\propos{1}}\then{\propos{2}}}$
if
$\propos{1}$,
then
$\propos{2}$.
\varfill
$\then$\\
%%%%%%%%%%%%%%%%%%%%%%%%%%%%%%%%%%%%%%%%%%%%%%%%%%%%%%%%%%%%%%%%%%%%%%%%%%%%%%%%%%%%%%%%
$\newsymp{{\propos{1}}\thenn{\propos{2}}}$
$\propos{1}$,
if-and-only-if
$\propos{2}$.
\varfill
$\thenn$\\
%%%%%%%%%%%%%%%%%%%%%%%%%%%%%%%%%%%%%%%%%%%%%%%%%%%%%%%%%%%%%%%%%%%%%%%%%%%%%%%%%%%%%%%%
$\newsymp{\negation{\propos{}}}$
$\propos{1}$,
negation of
$\propos{}$.
\varfill
$\negation{}$\\
%%%%%%%%%%%%%%%%%%%%%%%%%%%%%%%%%%%%%%%%%%%%%%%%%%%%%%%%%%%%%%%%%%%%%%%%%%%%%%%%%%%%%%%%
$\newsymp{\Foreach{\SET{1}}{\SET{2}}{\prop{\SET{1}}}}$
for every
$\SET{2}$
in
$\SET{1}$,
$\prop{\SET{1}}$.
\varfill
$\Foreach{\dummy}{\dummy}\dummy$\\
%%%%%%%%%%%%%%%%%%%%%%%%%%%%%%%%%%%%%%%%%%%%%%%%%%%%%%%%%%%%%%%%%%%%%%%%%%%%%%%%%%%%%%%%
$\newsymp{\Exists{\SET{1}}{\SET{2}}{\prop{\SET{1}}}}$
exists
$\SET{2}$
in
$\SET{1}$ such that
$\prop{\SET{1}}$.
\varfill
$\Foreach{\dummy}{\dummy}\dummy$\\
%%%%%%%%%%%%%%%%%%%%%%%%%%%%%%%%%%%%%%%%%%%%%%%%%%%%%%%%%%%%%%%%%%%%%%%%%%%%%%%%%%%%%%%%
%%%%%%%%%%%%%%%%%%%%%%%%%%%%%%%%%%%%%%%%%%%%%%%%%%%%%%%%%%%%%%%%%%%%%%%%%%%%%%%%%%%%%%%%
%%%%%%%%%%%%%%%%%%%%%%%%%%%%%%%%%%%%%%%%%%%%%%%%%%%%%%%%%%%%%%%%%%%%%%%%%%%%%%%%%%%%%%%%
%%%%%%%%%%%%%%%%%%%%%%%%%%%%%%%%%%%%%%%%%%%%%%%%%%%%%%%%%%%%%%%%%%%%%%%%%%%%%%%%%%%%%%%%
%%%%%%%%%%%%%%%%%%%%%%%%%%%%%%%%%%%%%%%%%%%%%%%%%%%%%%%%%%%%%%%%%%%%%%%%%%%%%%%%%%%%%%%%
%%%%%%%%%%%%%%%%%%%%%%%%%%%%%%%%%%%%%%%%%%%%%%%%%%%%%%%%%%%%%%%%%%%%%%%%%%%%%%%%%%%%%%%%
%%%%%%%%%%%%%%%%%%%%%%%%%%%%%%%%%%%%%%%%%%%%%%%%%%%%%%%%%%%%%%%%%%%%%%%%%%%%%%%%%%%%%%%%
%%%%%%%%%%%%%%%%%%%%%%%%%%%%%%%%%%%%%%%%%%%%%%%%%%%%%%%%%%%%%%%%%%%%%%%%%%%%%%%%%%%%%%%%
\section*{
Linear Algebra
}
%%%%%%%%%%%%%%%%%%%%%%%%%%%%%%%%%%%%%%%%%%%%%%%%%%%%%%%%%%%%%%%%%%%%%%%%%%%%%%%%%%%%%%%%
%$\newsymp{\subvec{\vectorspace{}}{m}}$
%the set of all sets of vectors of all $m$-dimensional vector-subspaces of
%the vector-space $\vectorspace{}$
%\varfill
%$\subvec{\dummy}{\dummy}$\\
%%%%%%%%%%%%%%%%%%%%%%%%%%%%%%%%%%%%%%%%%%%%%%%%%%%%%%%%%%%%%%%%%%%%%%%%%%%%%%%%%%%%%%%%
$\newsymp{\func{\Vspan{\vectorspace{}}}{\asubset}}$
the vector-subspace of the vector-space $\VVS{}$ spanned by the subset $\asubset$
of the set of all vectors of $\vectorspace{}$
\varfill
$\func{\Vspan{\dummy}}{\dummy}$\\
%%%%%%%%%%%%%%%%%%%%%%%%%%%%%%%%%%%%%%%%%%%%%%%%%%%%%%%%%%%%%%%%%%%%%%%%%%%%%%%%%%%%%%%%
%$\newsymp{\ovecbasis{\VVS{}}}$
%the set of all ordered-bases of the vector-space $\VVS{}$
%\varfill
%$\ovecbasis{\dummy}$\\
%%%%%%%%%%%%%%%%%%%%%%%%%%%%%%%%%%%%%%%%%%%%%%%%%%%%%%%%%%%%%%%%%%%%%%%%%%%%%%%%%%%%%%%%
$\newsymp{\Lin{\vectorspace{1}}{\vectorspace{2}}}$
the set of all linear maps from the vector-space $\vectorspace{1}$ to the vector-space
$\vectorspace{2}$
\varfill
$\Lin{\dummy}{\dummy}$\\
%%%%%%%%%%%%%%%%%%%%%%%%%%%%%%%%%%%%%%%%%%%%%%%%%%%%%%%%%%%%%%%%%%%%%%%%%%%%%%%%%%%%%%%%
$\newsymp{\Linisom{\vectorspace{1}}{\vectorspace{2}}}$
the set of all linear isomorphisms from the vector-space $\vectorspace{1}$ to the
vector-space $\vectorspace{2}$
\varfill
$\Linisom{\dummy}{\dummy}$\\
%%%%%%%%%%%%%%%%%%%%%%%%%%%%%%%%%%%%%%%%%%%%%%%%%%%%%%%%%%%%%%%%%%%%%%%%%%%%%%%%%%%%%%%%
$\newsymp{\GL{\vectorspace{}}{}}$
the set of all linear isomorphisms from the vector-space $\vectorspace{}$ to itself
\varfill
$\GL{\dummy}{}$\\
%%%%%%%%%%%%%%%%%%%%%%%%%%%%%%%%%%%%%%%%%%%%%%%%%%%%%%%%%%%%%%%%%%%%%%%%%%%%%%%%%%%%%%%%
$\newsymp{\VLin{\vectorspace{1}}{\vectorspace{2}}}$
the canonical vector-space of all linear maps from the vector-space
$\vectorspace{1}$ to the vector-space $\vectorspace{2}$
\varfill
$\VLin{\dummy}{\dummy}$\\
%%%%%%%%%%%%%%%%%%%%%%%%%%%%%%%%%%%%%%%%%%%%%%%%%%%%%%%%%%%%%%%%%%%%%%%%%%%%%%%%%%%%%%%%
$\newsymp{\Det{n}}$
the determinant function on the set of all square $\R$-matrices of degree $n$
\varfill
$\Det{\dummy}$\\
%%%%%%%%%%%%%%%%%%%%%%%%%%%%%%%%%%%%%%%%%%%%%%%%%%%%%%%%%%%%%%%%%%%%%%%%%%%%%%%%%%%%%%%%
$\newsymp{\directsum{\vectorspace{1}}{\vectorspace{2}}}$
the direct sum of the vector spaces $\vectorspace{1}$ and $\vectorspace{2}$ over
the same field
\varfill
$\directsum{\dummy}{\dummy}$\\
%%%%%%%%%%%%%%%%%%%%%%%%%%%%%%%%%%%%%%%%%%%%%%%%%%%%%%%%%%%%%%%%%%%%%%%%%%%%%%%%%%%%%%%%
$\newsymp{\Tensors{r}{s}{\vectorspace{}}}$
the set of all $r$-covariant and $s$-contravariant tensors on the vector-space
$\vectorspace{}$
\varfill
$\Tensors{\dummy}{\dummy}{\dummy}$\\
%%%%%%%%%%%%%%%%%%%%%%%%%%%%%%%%%%%%%%%%%%%%%%%%%%%%%%%%%%%%%%%%%%%%%%%%%%%%%%%%%%%%%%%%
$\newsymp{T_1\tensor{\vectorspace{}}T_2}$
the tensor-product of the tensors $T_1$ and $T_2$ on the vector-space
$\vectorspace{}$
\varfill
$\dummy\tensor{\dummy}\dummy$
%%%%%%%%%%%%%%%%%%%%%%%%%%%%%%%%%%%%%%%%%%%%%%%%%%%%%%%%%%%%%%%%%%%%%%%%%%%%%%%%%%%%%%%%
%%%%%%%%%%%%%%%%%%%%%%%%%%%%%%%%%%%%%%%%%%%%%%%%%%%%%%%%%%%%%%%%%%%%%%%%%%%%%%%%%%%%%%%%
%%%%%%%%%%%%%%%%%%%%%%%%%%%%%%%%%%%%%%%%%%%%%%%%%%%%%%%%%%%%%%%%%%%%%%%%%%%%%%%%%%%%%%%%
%%%%%%%%%%%%%%%%%%%%%%%%%%%%%%%%%%%%%%%%%%%%%%%%%%%%%%%%%%%%%%%%%%%%%%%%%%%%%%%%%%%%%%%%
%%%%%%%%%%%%%%%%%%%%%%%%%%%%%%%%%%%%%%%%%%%%%%%%%%%%%%%%%%%%%%%%%%%%%%%%%%%%%%%%%%%%%%%%
%%%%%%%%%%%%%%%%%%%%%%%%%%%%%%%%%%%%%%%%%%%%%%%%%%%%%%%%%%%%%%%%%%%%%%%%%%%%%%%%%%%%%%%%
%%%%%%%%%%%%%%%%%%%%%%%%%%%%%%%%%%%%%%%%%%%%%%%%%%%%%%%%%%%%%%%%%%%%%%%%%%%%%%%%%%%%%%%%
%%%%%%%%%%%%%%%%%%%%%%%%%%%%%%%%%%%%%%%%%%%%%%%%%%%%%%%%%%%%%%%%%%%%%%%%%%%%%%%%%%%%%%%%
\section*{Differential Geometry}
%%%%%%%%%%%%%%%%%%%%%%%%%%%%%%%%%%%%%%%%%%%%%%%%%%%%%%%%%%%%%%%%%%%%%%%%%%%%%%%%%%%%%%%%
$\newsymp{\atlases{r}{\M{}}{\NVS{}}}$
the set of all atlases of differentiablity class $\difclass{r}$
on the set $\M{}$
modeled on the Banach-space $\NVS{}$
\varfill
$\atlases{\dummy}{\dummy}{\dummy}$\\
%%%%%%%%%%%%%%%%%%%%%%%%%%%%%%%%%%%%%%%%%%%%%%%%%%%%%%%%%%%%%%%%%%%%%%%%%%%%%%%%%%%%%%%%
$\newsymp{\maxatlases{r}{\M{}}{\NVS{}}}$
the set of all maximal-atlases of differentiablity class $\difclass{r}$
on the set $\M{}$
modeled on the Banach-space $\NVS{}$
\varfill
$\maxatlases{\dummy}{\dummy}{\dummy}$\\
%%%%%%%%%%%%%%%%%%%%%%%%%%%%%%%%%%%%%%%%%%%%%%%%%%%%%%%%%%%%%%%%%%%%%%%%%%%%%%%%%%%%%%%%
$\newsymp{\func{\maxatlasgen{r}{\M{}}{\NVS{}}}{\atlas{}}}$
the maximal-atlas of differentiablity class $\difclass{r}$
on the set $\M{}$ modeled on
the Banach-space $\NVS{}$, generated by the atlas $\atlas{}$ in
$\atlases{r}{\M{}}{\NVS{}}$
\varfill
$\func{\maxatlasgen{\dummy}{\dummy}{\dummy}}{\dummy}$\\
%%%%%%%%%%%%%%%%%%%%%%%%%%%%%%%%%%%%%%%%%%%%%%%%%%%%%%%%%%%%%%%%%%%%%%%%%%%%%%%%%%%%%%%%
$\newsymp{\manprod{\Man{1}}{\Man{2}}}$
the manifold-product of the manifolds $\Man{1}$ and $\Man{2}$
\varfill
$\manprod{\dummy}{\dummy}$\\
%%%%%%%%%%%%%%%%%%%%%%%%%%%%%%%%%%%%%%%%%%%%%%%%%%%%%%%%%%%%%%%%%%%%%%%%%%%%%%%%%%%%%%%%
$\newsymp{\subman{\Man{}}{S}}$
the (regular) submanifold $S$ of the manifold $\Man{}$ endowed with its canonically
inherited differentiable structure
\varfill
$\subman{\dummy}{\dummy}$\\
%%%%%%%%%%%%%%%%%%%%%%%%%%%%%%%%%%%%%%%%%%%%%%%%%%%%%%%%%%%%%%%%%%%%%%%%%%%%%%%%%%%%%%%%
$\newsymp{\tanspace{\point}{\Man{}}}$
tangent space at point $\point$ of manifold $\Man{}$
\varfill
$\tanspace{\dummy}{\dummy}$\\
%%%%%%%%%%%%%%%%%%%%%%%%%%%%%%%%%%%%%%%%%%%%%%%%%%%%%%%%%%%%%%%%%%%%%%%%%%%%%%%%%%%%%%%%
$\newsymp{\tanspace{\point}{\Man{}}}$
tangent space at point $\point$ of manifold $\Man{}$ endowed with its linear structure
\varfill
$\tanspace{\dummy}{\dummy}$\\
%%%%%%%%%%%%%%%%%%%%%%%%%%%%%%%%%%%%%%%%%%%%%%%%%%%%%%%%%%%%%%%%%%%%%%%%%%%%%%%%%%%%%%%%
$\newsymp{\Tanbun{\Man{}}}$
tangent-bundle of the manifold $\Man{}$
\varfill
$\Tanbun{\dummy}$\\
%%%%%%%%%%%%%%%%%%%%%%%%%%%%%%%%%%%%%%%%%%%%%%%%%%%%%%%%%%%%%%%%%%%%%%%%%%%%%%%%%%%%%%%%
$\newsymp{\basep{\Man{}}}$
the projection map of the tangent-bundle of the manifold $\Man{}$
\varfill
$\basep{\dummy}$\\
%%%%%%%%%%%%%%%%%%%%%%%%%%%%%%%%%%%%%%%%%%%%%%%%%%%%%%%%%%%%%%%%%%%%%%%%%%%%%%%%%%%%%%%%
$\newsymp{\tanspaceiso{\point}{\Man{}}{\phi}}$
the canonical correspondence from the tangent-space at point $\point$ of
the manifold $\Man{}$ to the Banach-space it is modeled on, via the chart $\phi$
\varfill
$\oneforms{\dummy}$\\
%%%%%%%%%%%%%%%%%%%%%%%%%%%%%%%%%%%%%%%%%%%%%%%%%%%%%%%%%%%%%%%%%%%%%%%%%%%%%%%%%%%%%%%%
$\newsymp{\mapdifclass{\infty}{\Man{1}}{\Man{2}}}$
the set of all smooth maps from the manifold $\Man{1}$ to the manifold $\Man{2}$
\varfill
$\mapdifclass{\dummy}{\dummy}{\dummy}$\\
%%%%%%%%%%%%%%%%%%%%%%%%%%%%%%%%%%%%%%%%%%%%%%%%%%%%%%%%%%%%%%%%%%%%%%%%%%%%%%%%%%%%%%%%
$\newsymp{\smoothmaps{\Man{}}}$
the set of all smooth maps from the manifold $\Man{1}$ to $\R$
\varfill
$\smoothmaps{\dummy}$\\
%%%%%%%%%%%%%%%%%%%%%%%%%%%%%%%%%%%%%%%%%%%%%%%%%%%%%%%%%%%%%%%%%%%%%%%%%%%%%%%%%%%%%%%%
$\newsymp{\vectorfields{\Man{}}}$
the set of all smooth vector fields on the manifold $\Man{}$
\varfill
$\vectorfields{\dummy}$\\
%%%%%%%%%%%%%%%%%%%%%%%%%%%%%%%%%%%%%%%%%%%%%%%%%%%%%%%%%%%%%%%%%%%%%%%%%%%%%%%%%%%%%%%%
$\newsymp{\oneforms{\Man{}}}$
the set of all smooth co-vector fields (one-forms) on the manifold $\Man{}$
\varfill
$\oneforms{\dummy}$\\
%%%%%%%%%%%%%%%%%%%%%%%%%%%%%%%%%%%%%%%%%%%%%%%%%%%%%%%%%%%%%%%%%%%%%%%%%%%%%%%%%%%%%%%%
$\newsymp{\derr{f}}$
the tangent-map of the smooth map $f$ between two manifolds
\varfill
$\derr{\dummy}$\\
%%%%%%%%%%%%%%%%%%%%%%%%%%%%%%%%%%%%%%%%%%%%%%%%%%%%%%%%%%%%%%%%%%%%%%%%%%%%%%%%%%%%%%%%
$\newsymp{\lieder{\avecf{}}{f}}$
Lie-derivative of the real-valued smooth map $f$ on a manifold with respect
to the smooth vector field $\avecf{}$ on it
\varfill
$\lieder{\dummy}{\dummy}$\\
%%%%%%%%%%%%%%%%%%%%%%%%%%%%%%%%%%%%%%%%%%%%%%%%%%%%%%%%%%%%%%%%%%%%%%%%%%%%%%%%%%%%%%%%
$\newsymp{\liebracket{\avecf{1}}{\avecf{2}}{\Man{}}}$
Lie-bracket of the smooth vector fields $\avecf{1}$ and $\avecf{2}$ on the
manifold $\Man{}$
\varfill
$\liebracket{\dummy}{\dummy}{\dummy}$\\
%%%%%%%%%%%%%%%%%%%%%%%%%%%%%%%%%%%%%%%%%%%%%%%%%%%%%%%%%%%%%%%%%%%%%%%%%%%%%%%%%%%%%%%%
$\newsymp{\fibervecs{\vbundle{}}{\point}}$
the fiber at the point $\point$ of the smooth vector bundle $\vbundle{}$
\varfill
$\fibervecs{\dummy}{\dummy}$\\
%%%%%%%%%%%%%%%%%%%%%%%%%%%%%%%%%%%%%%%%%%%%%%%%%%%%%%%%%%%%%%%%%%%%%%%%%%%%%%%%%%%%%%%%
$\newsymp{\vbsections{\vbundle{}}}$
the set of all smooth sections of the smooth vector bundle $\vbundle{}$
\varfill
$\vbsections{\dummy}$\\
%%%%%%%%%%%%%%%%%%%%%%%%%%%%%%%%%%%%%%%%%%%%%%%%%%%%%%%%%%%%%%%%%%%%%%%%%%%%%%%%%%%%%%%%
$\newsymp{\vbtensorbundle{r}{s}{\vbundle{}}}$
the $\opair{r}{s}$-tensor-bundle of the smooth vector bundle $\vbundle{}$
\varfill
$\vbtensorbundle{\dummy}{\dummy}{\dummy}$\\
%%%%%%%%%%%%%%%%%%%%%%%%%%%%%%%%%%%%%%%%%%%%%%%%%%%%%%%%%%%%%%%%%%%%%%%%%%%%%%%%%%%%%%%%
$\newsymp{\TF{r}{s}{\vbundle{}}}$
the set of all $\opair{r}{s}$ tensor fields of the smooth vector bundle $\vbundle{}$
\varfill
$\TF{\dummy}{\dummy}{\dummy}$
%%%%%%%%%%%%%%%%%%%%%%%%%%%%%%%%%%%%%%%%%%%%%%%%%%%%%%%%%%%%%%%%%%%%%%%%%%%%%%%%%%%%%%%%
$\newsymp{\VTF{r}{s}{\vbundle{}}}$
the vector-space of $\opair{r}{s}$ tensor fields of the smooth vector bundle
	$\vbundle{}$
\varfill
$\TF{\dummy}{\dummy}{\dummy}$
%%%%%%%%%%%%%%%%%%%%%%%%%%%%%%%%%%%%%%%%%%%%%%%%%%%%%%%%%%%%%%%%%%%%%%%%%%%%%%%%%%%%%%%%
$\newsymp{T_1\tensor{\vbundle{}}T_2}$
the tensor-product of the tensor fields $T_1$ and $T_2$ on the smooth vector bundle
$\vbundle{}$
\varfill
$\dummy\tensor{\dummy}\dummy$\\
%%%%%%%%%%%%%%%%%%%%%%%%%%%%%%%%%%%%%%%%%%%%%%%%%%%%%%%%%%%%%%%%%%%%%%%%%%%%%%%%%%%%%%%%
$\newsymp{\DVTF{\vbundle{}}}$
the universal tensor algebra of the smooth vector bundle $\vbundle{}$
$\vbundle{}$
\varfill
$\vbundle{\dummy}$\\
%%%%%%%%%%%%%%%%%%%%%%%%%%%%%%%%%%%%%%%%%%%%%%%%%%%%%%%%%%%%%%%%%%%%%%%%%%%%%%%%%%%%%%%%
$\newsymp{\vbmorphisms{\vbundle{1}}{\vbundle{2}}}$
the set of all smooth vector bundle morphisms from the vector bundle $\vbundle{1}$
to the vector bundle $\vbundle{2}$
$\vbundle{}$
\varfill
$\vbmorphisms{\dummy}{\dummy}$\\
%%%%%%%%%%%%%%%%%%%%%%%%%%%%%%%%%%%%%%%%%%%%%%%%%%%%%%%%%%%%%%%%%%%%%%%%%%%%%%%%%%%%%%%%
$\newsymp{\vbisomorphisms{\vbundle{1}}{\vbundle{2}}}$
the set of all smooth vector bundle isomorphisms from the vector bundle $\vbundle{1}$
to the vector bundle $\vbundle{2}$
$\vbundle{}$
\varfill
$\vbisomorphisms{\dummy}{\dummy}$\\
%%%%%%%%%%%%%%%%%%%%%%%%%%%%%%%%%%%%%%%%%%%%%%%%%%%%%%%%%%%%%%%%%%%%%%%%%%%%%%%%%%%%%%%%
$\newsymp{\func{\VBpullback{f}{r}{s}}{T}}$
the action of the $\opair{r}{s}$-pullback of the smooth vector bundle morphism $f$
between two vector bundles on the $\opair{r}{s}$ tensor field $T$
\varfill
$\func{\VBpullback{\dummy}{\dummy}{\dummy}}{\dummy}$
%%%%%%%%%%%%%%%%%%%%%%%%%%%%%%%%%%%%%%%%%%%%%%%%%%%%%%%%%%%%%%%%%%%%%%%%%%%%%%%%%%%%%%%%
%%%%%%%%%%%%%%%%%%%%%%%%%%%%%%%%%%%%%%%%%%%%%%%%%%%%%%%%%%%%%%%%%%%%%%%%%%%%%%%%%%%%%%%%
%%%%%%%%%%%%%%%%%%%%%%%%%%%%%%%%%%%%%%%%%%%%%%%%%%%%%%%%%%%%%%%%%%%%%%%%%%%%%%%%%%%%%%%%
%%%%%%%%%%%%%%%%%%%%%%%%%%%%%%%%%%%%%%%%%%%%%%%%%%%%%%%%%%%%%%%%%%%%%%%%%%%%%%%%%%%%%%%%
\section*{Mathematical Environments}
$\newsymp{\blacksquare}$
end of definition
\varfill
$\blacksquare$\\
%%%%%%%%%%%%%%%%%%%%%%%%%%%%%%%%%%%%%%%%%%%%%%%%%%%%%%%%%%%%%%%%%%%%%%%%%%%%%%%%%%%%%%%%
$\newsymp{\square}$
end of theorem, lemma, proposition, or corollary
\varfill
$\square$\\
%%%%%%%%%%%%%%%%%%%%%%%%%%%%%%%%%%%%%%%%%%%%%%%%%%%%%%%%%%%%%%%%%%%%%%%%%%%%%%%%%%%%%%%%
$\newsymp{\Diamond}$
end of the introduction of new fixed objects
\varfill
$\Diamond$
\chapteR{Connections}
\thispagestyle{fancy}
\section{Connections on Smooth Vector Bundles}
%%%%%%%%%%%%%%%%%%%%%%%%%%%%%%%%%%%%%%%%%%%%%%%%%%%%%%%%%%%%%%%%%%%%%%%%%%%%%%%%%%%%%%%%%%%%%%%%%%%%%%%%%%%%%%%%%%%%%%%%%%%%%%%%
\fixed
$\vbundle{}=\quintuple{\vbtotal{}}{\vbprojection{}}{\vbbase{}}{\vbfiber{}}{\vbatlas{}}$ is fixed as a real smooth vector bundle
of rank $d$,
where $\vbtotal{}=\opair{\vTot{}}{\maxatlas{\vTot{}}}$ and
$\vbbase{}=\opair{\vB{}}{\maxatlas{\vB{}}}$ are $\difclass{\infty}$ manifolds
modeled on the Banach-spaces $\R^{n_{\vTot{}}}$ and $\R^{n_{\vB{}}}$, respectively.
\endfixed
%%%%%%%%%%%%%%%%%%%%%%%%%%%%%%%%%%%%%%%%%%%%%%%%%%%%%%%%%%%%%%%%%%%%%%%%%%%%%%%%%%%%%%%%%%%%%%%%%%%%%%%%%%%%%%%%%%%%%%%%%%%%%%%%
\definition\label{defconnection}
A map $\function{\connection{}}{\Cprod{\vectorfields{\vbbase{}}}{\vbsections{\vbundle{}}}}{\vbsections{\vbundle{}}}$
is referred to as a $\quotl$(Koszul) connection on the smooth vector bundle $\vbundle{}$$\quotr$ if it satisfies the following conditions.
$\vectorfields{\vbbase{}}$ and $\vbsections{\vbundle{}}$ denote the set of all
smooth vector fields on the manifold $\vbbase{}$, and the set of all smooth sections of the vector bundle $\vbundle{}$,
respectively, considered to be endowed with their canonical module structures.
\begin{itemize}
\item[\myitem{CON~1.}]
\begin{align}
&\Foreach{\opair{f_1}{f_2}}{\Cprod{\smoothmaps{\vbbase{}}}{\smoothmaps{\vbbase{}}}}
\Foreach{\opair{\avecf{1}}{\avecf{2}}}{\Cprod{\vectorfields{\vbbase{}}}{\vectorfields{\vbbase{}}}}
\Foreach{\vbsec{}}{\vbsections{\vbundle{}}}\cr
&\func{\connection{}}{\binary{f_1\avecf{1}+f_2\avecf{2}}{\vbsec{}}}=
f_1\func{\connection{}}{\binary{\avecf{1}}{\vbsec{}}}+
f_2\func{\connection{}}{\binary{\avecf{2}}{\vbsec{}}}.
\end{align}
\item[\myitem{CON~2.}]
\begin{align}
\Foreach{\avecf{}}{\vectorfields{\vbbase{}}}
\Foreach{\opair{\vbsec{1}}{\vbsec{2}}}{\Cprod{\vbsections{\vbundle{}}}{\vbsections{\vbundle{}}}}
\func{\connection{}}{\binary{\avecf{}}{\vbsec{1}+\vbsec{2}}}=
\func{\connection{}}{\binary{\avecf{}}{\vbsec{1}}}+\func{\connection{}}{\binary{\avecf{}}{\vbsec{2}}}.
\end{align}
\item[\myitem{CON~3.}]
\begin{align}
\Foreach{\avecf{}}{\vectorfields{\vbbase{}}}
\Foreach{f}{\smoothmaps{\vbbase{}}}
\Foreach{\vbsec{}}{\vbsections{\vbundle{}}}
\func{\connection{}}{\binary{\avecf{}}{f\vbsec{}}}=
\(\lieder{\avecf{}}f\)\vbsec{}+
f\func{\connection{}}{\binary{\avecf{}}{\vbsec{}}},
\end{align}
where, $\lieder{\avecf{}}f$ denotes the Lie-derivative of $f$ with respect to $\avecf{}$.
\end{itemize}
The set of all connections on $\vbundle{}$ will be denoted by $\connections{\vbundle{}}$.
$\func{\connection{}}{\binary{\avecf{}}{\vbsec{}}}$ will also be alternatively denoted by
$\con{\avecf{}}{\vbsec{}}$.
\endef
%%%%%%%%%%%%%%%%%%%%%%%%%%%%%%%%%%%%%%%%%%%%%%%%%%%%%%%%%%%%%%%%%%%%%%%%%%%%%%%%%%%%%%%%%%%%%%%%%%%%%%%%%%%%%%%%%%%%%%%%%%%%%%%%
\fixed
$\connection{}$ is fixed as an element of $\connections{\vbundle{}}$.
\endfixed
%%%%%%%%%%%%%%%%%%%%%%%%%%%%%%%%%%%%%%%%%%%%%%%%%%%%%%%%%%%%%%%%%%%%%%%%%%%%%%%%%%%%%%%%%%%%%%%%%%%%%%%%%%%%%%%%%%%%%%%%%%%%%%%%
\lemma
Let $\avecf{}\in\vectorfields{\vbbase{}}$ be the everywhere-vanishing vector field. For every
$\vbsec{}\in\vbsections{\vbundle{}}$ the section $\con{\avecf{}}{\vbsec{}}$ vanishes everywhere,
that is $\func{\[\con{\avecf{}}{\vbsec{}}\]}{\point}=\zerovec{}$ for every $\point\in\vB{}$.
\proof
Let $\function{f}{\vB{}}{\R}$ be the trivial smooth map, that is $\func{f}{\point}=0$ for every $\point\in\vB{}$.
Then, trivially $f\avecf{}=\avecf{}$, and thus,
\begin{equation}
\con{\avecf{}}{\vbsec{}}=\con{f\avecf{}}{\vbsec{}}=f\con{\avecf{}}{\vbsec{}}=0.
\end{equation}
\endlem
%%%%%%%%%%%%%%%%%%%%%%%%%%%%%%%%%%%%%%%%%%%%%%%%%%%%%%%%%%%%%%%%%%%%%%%%%%%%%%%%%%%%%%%%%%%%%%%%%%%%%%%%%%%%%%%%%%%%%%%%%%%%%%%%
\lemma\label{locallityofconnectionwrtvectorfields}
Let $\avecf{}\in\vectorfields{\vbbase{}}$ and $\vbsec{}\in\vbsections{\vbundle{}}$.
Let $\point$ be a point of $\vbbase{}$.
If there exists an open neighborhood $\U$ of $\point$ in $\vbbase{}$ such that
$\Foreach{x}{\U}\func{\avecf{}}{x}=\zerovec{}$, then $\func{\[\con{\avecf{}}{\vbsec{}}\]}{\point}=\zerovec{}$.
\proof
Assume that there exists an open neighborhood $\U$ of $\point$ in $\vbbase{}$ such that
$\Foreach{x}{\U}\func{\avecf{}}{x}=\zerovec{}$. According to \cite[page 28, Chapter 1, Lemma 1.69]{JLee},
there exists a smooth $\R$-valued map $f\in\smoothmaps{\vbbase{}}$ such that $\func{f}{\point}=1$ and
$\Foreach{x}{\(\compl{\vB{}}{\U}\)}\func{f}{x}=0$. Clearly, $f\avecf{}=0$, and thus $\displaystyle\con{f\avecf{}}{\vbsec{}}=0$,
and in particular $\displaystyle\func{\[\con{f\avecf{}}{\vbsec{}}\]}{\point}=\zerovec{}$. On the other hand,
$\displaystyle\func{\[\con{f\avecf{}}{\vbsec{}}\]}{\point}=\func{f}{\point}\func{\[\con{\avecf{}}{\vbsec{}}\]}{\point}=
\func{\[\con{\avecf{}}{\vbsec{}}\]}{\point}$. Therefore, $\func{\[\con{\avecf{}}{\vbsec{}}\]}{\point}=\zerovec{}$.
\endlem
%%%%%%%%%%%%%%%%%%%%%%%%%%%%%%%%%%%%%%%%%%%%%%%%%%%%%%%%%%%%%%%%%%%%%%%%%%%%%%%%%%%%%%%%%%%%%%%%%%%%%%%%%%%%%%%%%%%%%%%%%%%%%%%%
\lemma\label{locallityofconnectionwrtvectorsections}
Let $\avecf{}\in\vectorfields{\vbbase{}}$ and $\vbsec{}\in\vbsections{\vbundle{}}$.
Let $\point$ be a point of $\vbbase{}$.
If there exists an open neighborhood $\U$ of $\point$ in $\vbbase{}$ such that
$\Foreach{x}{\U}\func{\vbsec{}}{x}=\zerovec{}$, then $\func{\[\con{\avecf{}}{\vbsec{}}\]}{\point}=\zerovec{}$.
\proof
there exists an open neighborhood $\U$ of $\point$ in $\vbbase{}$ such that
$\Foreach{x}{\U}\func{\vbsec{}}{x}=\zerovec{}$.
There exists a smooth $\R$-valued map $f\in\smoothmaps{\vbbase{}}$ such that $\func{f}{\point}=1$ and
$\Foreach{x}{\(\compl{\vB{}}{\U}\)}\func{f}{x}=0$. So, $f\avecf{}=0$ and hence
$\displaystyle\func{\[\con{f\avecf{}}{\vbsec{}}\]}{\point}=\zerovec{}$. Furthermore, considering that
$\func{\vbsec{}}{\point}=\zerovec{}$ and $\func{f}{\point}=1$,
\begin{align}
\func{\[\con{f\avecf{}}{\vbsec{}}\]}{\point}=
\func{\[\lieder{\avecf{}}{f}\]}{\point}\func{\vbsec{}}{\point}+
\func{f}{\point}\func{\[\con{\avecf{}}{\vbsec{}}\]}{\point}=
\func{\[\con{\avecf{}}{\vbsec{}}\]}{\point}.
\end{align}
Thus, $\func{\[\con{\avecf{}}{\vbsec{}}\]}{\point}=\zerovec{}$.
\endlem
%%%%%%%%%%%%%%%%%%%%%%%%%%%%%%%%%%%%%%%%%%%%%%%%%%%%%%%%%%%%%%%%%%%%%%%%%%%%%%%%%%%%%%%%%%%%%%%%%%%%%%%%%%%%%%%%%%%%%%%%%%%%%%%%
\corollary
Let $\avecf{1}$ and $\avecf{2}$ be a pair of smooth vector fields on $\vbbase{}$, and let
$\vbsec{1}$ and $\vbsec{2}$ be a pair of smooth sections of the vector bundle $\vbundle{}$, such that
$\reS{\avecf{1}}{\U}=\reS{\avecf{2}}{\U}$ and $\reS{\vbsec{1}}{\U}=\reS{\vbsec{2}}{\U}$ for some non-empty open set $\U$
of $\vbbase{}$. Then $\Foreach{\point}{\U}\func{\[\con{\avecf{1}}{\vbsec{1}}\]}{\point}=
\func{\[\con{\avecf{2}}{\vbsec{2}}\]}{\point}$.
\endcor
%%%%%%%%%%%%%%%%%%%%%%%%%%%%%%%%%%%%%%%%%%%%%%%%%%%%%%%%%%%%%%%%%%%%%%%%%%%%%%%%%%%%%%%%%%%%%%%%%%%%%%%%%%%%%%%%%%%%%%%%%%%%%%%%
\remark
We remind that given a regular submanifold of $\M{}$ of $\vbbase{}$, we denote by $\subman{\vbbase{}}{\Man{}}$ the submanifold of
$\vbbase{}$ with the set of points $\M{}$ endowed with the
differentiable structure canonically inherited from that of $\vbbase{}$. We also denote by $\subman{\vbundle{}}{\M{}}$
the restriction of the smooth vector bundle $\vbundle{}$ to $\M{}$, as defined in \cite{ShahiVB}.
Moreover, given any section $\vbsec{}\in\vbsections{\vbundle{}}$, as usual, we denote by $\reS{\vbsec{}}{\M{}}$
the restriction of $\vbsec{}$ to $\M{}$, which is a smooth section of the vector bundle $\subman{\vbundle{}}{\M{}}$.
\endremark
%%%%%%%%%%%%%%%%%%%%%%%%%%%%%%%%%%%%%%%%%%%%%%%%%%%%%%%%%%%%%%%%%%%%%%%%%%%%%%%%%%%%%%%%%%%%%%%%%%%%%%%%%%%%%%%%%%%%%%%%%%%%%%%%
\definition\label{defrestrictedconnection}
Let $\U$ be a non-empty open subset of $\vbbase{}$. We define the map
$\function{\rescon{\connection{}}{\U}}{\Cprod{\vectorfields{\subman{\vbbase{}}{\U}}}{\vbsections{\subman{\vbundle{}}{\U}}}}
{\vbsections{\subman{\vbundle{}}{\U}}}$ in the following way.\\
Let $\avecf{}\in\vectorfields{\subman{\vbbase{}}{\U}}$ and $\vbsec{}\in\vbsections{\subman{\vbundle{}}{\U}}$.
Given any $\point\in\U$, choose an open $\V\subset\U$, a vector field $\avecff{}\in\vectorfields{\vbbase{}}$,
and a smooth section $\vbsecc{}\in\vbsections{\vbundle{}}$,
such that $\point\in\V$ and $\reS{\avecff{}}{\V}=\reS{\avecf{}}{\V}$, $\reS{\vbsecc{}}{\V}=\reS{\vbsec{}}{\V}$.
The existence of such triple $\triple{\V}{\avecff{}}{\vbsecc{}}$ is guaranteed, and multiplying $\avecf{}$ and $\vbsec{}$ by
appropriate cut-off functions provides a class of straightforward examples. Now define,
\begin{equation}
\func{\[\func{\rescon{\connection{}}{\U}}{\binary{\avecf{}}{\vbsec{}}}\]}{\point}\eqdef
\func{\[\con{\avecff{}}{\vbsecc{}}\]}{\point}.
\end{equation}
This definition is permissible, because according to \reflem{locallityofconnectionwrtvectorfields} and
\reflem{locallityofconnectionwrtvectorsections}, $\func{\[\con{\avecff{}}{\vbsecc{}}\]}{\point}$ is
independent of the choice of the triple $\triple{\V}{\avecff{}}{\vbsecc{}}$ with the desired properties stated above.
\endef
%%%%%%%%%%%%%%%%%%%%%%%%%%%%%%%%%%%%%%%%%%%%%%%%%%%%%%%%%%%%%%%%%%%%%%%%%%%%%%%%%%%%%%%%%%%%%%%%%%%%%%%%%%%%%%%%%%%%%%%%%%%%%%%%
\remark
Note that given a vector field $\avecf{}\in\vectorfields{\subman{\vbbase{}}{\U}}$ and a vector field
$\avecff{}\in\vectorfields{\vbbase{}}$, when $\U$ is an open submanifold of $\vbbase{}$, considering that the tangent-space
of $\subman{\vbbase{}}{\U}$ at any $\point\in\U$ is canonically identified with the tangent-space of $\vbbase{}$
at $\point$, $\func{\avecf{}}{\point}$ and $\func{\avecff{}}{\point}$ can be thought of as vectors of a same vector-space.
Thus the question of coincidence of $\avecf{}$ and $\avecff{}$ in a region of $\U$ is meaningful.
\endremark
%%%%%%%%%%%%%%%%%%%%%%%%%%%%%%%%%%%%%%%%%%%%%%%%%%%%%%%%%%%%%%%%%%%%%%%%%%%%%%%%%%%%%%%%%%%%%%%%%%%%%%%%%%%%%%%%%%%%%%%%%%%%%%%%
\theorem
Let $\U$ be a non-empty open subset of $\vbbase{}$. $\rescon{\connection{}}{\U}$ is a connection on
the restricted vector bundle $\subman{\vbundle{}}{\U}$, that is
$\rescon{\connection{}}{\U}\in\connections{\subman{\vbundle{}}{\U}}$.
\proof
Let $\avecf{}$ and $\avecf{1}$ be a pair of elements of $\vectorfields{\subman{\vbbase{}}{\U}}$, and
$\vbsec{}$ and $\vbsec{1}$ a pair of elements of $\vbsections{\subman{\vbundle{}}{\U}}$. Also, let
$f$ and $f_1$ be a pair of elements of $\smoothmaps{\subman{\vbbase{}}{\U}}$.\\
Let $\point\in\U$. It is clearly known that there exists a quadruple $\quadruple{\V}{\avecff{}}{\vbsecc{}}{g}$
with $\V$ an open subset of $\vbbase{}$, $\avecff{}\in\vectorfields{\vbbase{}}$,
$\vbsecc{}\in\vbsections{\vbundle{}}$, and $g\in\smoothmaps{\vbbase{}}$ such that
$\point\in\V$ and $\reS{\avecff{}}{\V}=\reS{\avecf{}}{\V}$,
$\reS{\vbsecc{}}{\V}=\reS{\vbsec{}}{\V}$, and $\reS{g}{\V}=\reS{f}{\V}$.
Similarly, there exists a quadruple $\quadruple{\V_1}{\avecff{1}}{\vbsecc{1}}{g_1}$
with $\V_1$ an open subset of $\vbbase{1}$, $\avecff{1}\in\vectorfields{\vbbase{}}$,
$\vbsecc{1}\in\vbsections{\vbundle{}}$, and $g_1\in\smoothmaps{\vbbase{}}$ such that
$\point\in\V_1$ and $\reS{\avecff{1}}{\V_1}=\reS{\avecf{1}}{\V_1}$,
$\reS{\vbsecc{1}}{\V_1}=\reS{\vbsec{1}}{\V_1}$, and $\reS{g_1}{\V_1}=\reS{f_1}{\V_1}$.
Now, by defining $W:=\V\cap\V_1$, trivially $\reS{\avecff{}}{W}=\reS{\avecf{}}{W}$,
$\reS{\vbsecc{}}{W}=\reS{\vbsec{}}{W}$, $\reS{g}{W}=\reS{f}{W}$, $\reS{\avecff{1}}{W}=\reS{\avecf{1}}{W}$,
$\reS{\vbsecc{1}}{W}=\reS{\vbsec{1}}{W}$, and $\reS{g_1}{W}=\reS{f_1}{W}$, and hence
$\reS{\(g\avecff{}+g_1\avecff{1}\)}{W}=\reS{\(f\avecf{}+f_1\avecf{1}\)}{W}$. Therefore,
\begin{align}
\func{\[\func{\rescon{\connection{}}{\U}}{\binary{f\avecf{}+f_1\avecf{1}}{\vbsec{}}}\]}{\point}&=
\func{\[\con{\(g\avecff{}+g_1\avecff{1}\)}{\vbsecc{}}\]}{\point}\cr
&=\func{g}{\point}\func{\[\con{\avecff{}}{\vbsecc{}}\]}{\point}+
\func{g_1}{\point}\func{\[\con{\avecff{1}}{\vbsecc{1}}\]}{\point}\cr
&=\func{f}{\point}\func{\[\func{\rescon{\connection{}}{\U}}{\binary{\avecf{}}{\vbsec{}}}\]}{\point}+
\func{f_1}{\point}\func{\[\func{\rescon{\connection{}}{\U}}{\binary{\avecf{1}}{\vbsec{1}}}\]}{\point}\cr
&=\func{\[f\func{\rescon{\connection{}}{\U}}{\binary{\avecf{}}{\vbsec{}}}+
f_1\func{\rescon{\connection{}}{\U}}{\binary{\avecf{1}}{\vbsec{1}}}\]}{\point}.
\end{align}
Moreover, trivially $\reS{\(\vbsecc{}+\vbsecc{1}\)}{W}=\reS{\vbsec{}+\vbsec{1}}{W}$, and thus,
\begin{align}
\func{\[\func{\rescon{\connection{}}{\U}}{\binary{\avecf{}}{\vbsec{}+\vbsec{1}}}\]}{\point}&=
\func{\[\con{\avecff{}}{\(\vbsecc{}+\vbsecc{1}\)}\]}{\point}\cr
&=\func{\[\con{\avecff{}}{\vbsecc{}}\]}{\point}+
\func{\[\con{\avecff{1}}{\vbsecc{1}}\]}{\point}\cr
&=\func{\[\func{\rescon{\connection{}}{\U}}{\binary{\avecf{}}{\vbsec{}}}\]}{\point}+
\func{\[\func{\rescon{\connection{}}{\U}}{\binary{\avecf{}}{\vbsec{1}}}\]}{\point}\cr
&=\func{\[\func{\rescon{\connection{}}{\U}}{\binary{\avecf{}}{\vbsec{}}}+
\func{\rescon{\connection{}}{\U}}{\binary{\avecf{}}{\vbsec{1}}}\]}{\point}.
\end{align}
Now consider the Lie derivatives $\lieder{\avecf{}}{f}$ and $\lieder{\avecff{}}{g}$. There is a notational
subtlety, and the point is that $\lieder{\avecf{}}{f}$ and $\lieder{\avecff{}}{g}$ are Lie derivatives
on the manifolds $\subman{\vbbase{}}{\U}$ and $\vbbase{}$. Since $\subman{\vbbase{}}{\U}$ is an open
submanifol of $\vbbase{}$, $f$ coincides with $g$ in a neighborhood of $\point$,
and $\func{\avecf{}}{\point}=\func{\avecff{}}{\point}$,
it is easy to verify that,
\begin{equation}
\func{\[\lieder{\avecf{}}{f}\]}{\point}=\func{\[\lieder{\avecff{}}{g}\]}{\point}.
\end{equation}
Therefore, considering that $\reS{f\vbsec{}}{W}=\reS{g\vbsecc{}}{W}$,
\begin{align}
\func{\[\func{\rescon{\connection{}}{\U}}{\binary{\avecf{}}{f\vbsec{}}}\]}{\point}&=
\func{\[\con{\avecff{}}{\(g\vbsecc{}\)}\]}{\point}\cr
&=\func{\[\lieder{\avecff{}}{g}\]}{\point}\func{\vbsecc{}}{\point}+
\func{g}{\point}\func{\[\con{\avecff{}}{\vbsecc{}}\]}{\point}\cr
&=\func{\[\lieder{\avecf{}}{f}\]}{\point}\func{\vbsec{}}{\point}+
\func{f}{\point}\func{\[\func{\rescon{\connection{}}{\U}}{\binary{\avecf{}}{\vbsec{}}}\]}{\point}\cr
&=\func{\[\(\lieder{\avecf{}}{f}\)\vbsec{}+f\func{\rescon{\connection{}}{\U}}{\binary{\avecf{}}{\vbsec{}}}\]}{\point}.
\end{align}
Ultimately, according to \refdef{defconnection}, it is seen that $\rescon{\connection{}}{\U}$
is a connection on $\subman{\vbundle{}}{\U}$.
\endthm
%%%%%%%%%%%%%%%%%%%%%%%%%%%%%%%%%%%%%%%%%%%%%%%%%%%%%%%%%%%%%%%%%%%%%%%%%%%%%%%%%%%%%%%%%%%%%%%%%%%%%%%%%%%%%%%%%%%%%%%%%%%%%%%%
\remark
Given a non-empty open subset $\U$ of $\Man{}$, we will refer to $\rescon{\connection{}}{\U}$ as the
$\quotl$restriction of the connection $\connection{}$ to the open subset $\U$ of $\Man{}$$\quotr$.
\endremark
%%%%%%%%%%%%%%%%%%%%%%%%%%%%%%%%%%%%%%%%%%%%%%%%%%%%%%%%%%%%%%%%%%%%%%%%%%%%%%%%%%%%%%%%%%%%%%%%%%%%%%%%%%%%%%%%%%%%%%%%%%%%%%%%
\theorem\label{thmnaturalityofconnection}
Let $\V$ and $\U$ be open subsets of $\vbbase{}$ such that $\V\subseteq\U$. For every
$\avecf{}\in\vectorfields{\subman{\vbbase{}}{\U}}$ and $\vbsec{}\in\vbsections{\subman{\vbundle{}}{\U}}$,
\begin{equation}
\func{\rescon{\connection{}}{\V}}{\binary{\reS{\avecf{}}{\V}}{\reS{\vbsec{}}{\V}}}=
\reS{\func{\rescon{\connection{}}{\U}}{\binary{\avecf{}}{\vbsec{}}}}{\V}.
\end{equation}
\proof
It is trivial and left as an exercise.
\endthm
%%%%%%%%%%%%%%%%%%%%%%%%%%%%%%%%%%%%%%%%%%%%%%%%%%%%%%%%%%%%%%%%%%%%%%%%%%%%%%%%%%%%%%%%%%%%%%%%%%%%%%%%%%%%%%%%%%%%%%%%%%%%%%%%
\theorem
Let $\avecf{}\in\vectorfields{\vbbase{}}$, and $\func{\avecf{}}{\point}=\zerovec{}$ for some $\point\in\vB{}$.
\begin{equation}
\Foreach{\vbsec{}}{\vbsections{\vbundle{}}}
\func{\[\con{\avecf{}}{\vbsec{}}\]}{\point}=\zerovec{}.
\end{equation}
\proof
Choose a chart $\opair{\U}{\phi}$ of $\vbbase{}$ around the point $\point$, and let $\mtuple{\avecff{1}}{\avecff{n}}$
be the local frame-field corresponded to this chart. So, for each $i$, $\avecff{i}\in\vectorfields{\subman{\vbbase{}}{\U}}$,
and for every $x\in\U$, $\mtuple{\func{\avecff{1}}{x}}{\func{\avecff{n}}{x}}$ is an ordered basis of $\tanspace{x}{\vbbase{}}$.
So, there exists smooth maps $f_i\in\smoothmaps{\subman{\vbbase{}}{\U}}$ such that,
\begin{equation}
\reS{\avecf{}}{\U}=\sum_{i=1}^{n}f_{i}\avecff{i},
\end{equation}
and for every $i\in\seta{\suc{1}{n}}$, $\func{f_i}{\point}=0$.
Therefore, according to \refdef{defconnection} and \refdef{defrestrictedconnection},
\begin{align}
\func{\[\con{\avecf{}}{\vbsec{}}\]}{\point}&=
\func{\[\func{\rescon{\connection{}}{\U}}{\binary{\reS{\avecf{}}{\U}}{\reS{\vbsec{}}{\U}}}\]}{\point}\cr
&=\func{\[\func{\rescon{\connection{}}{\U}}{\binary{\sum_{i=1}^{n}f_{i}\avecff{i}}{\reS{\vbsec{}}{\U}}}\]}{\point}\cr
&=\sum_{i=1}^{n}\func{f_i}{\point}
\func{\[\func{\rescon{\connection{}}{\U}}{\binary{\avecff{i}}{\reS{\vbsec{}}{\U}}}\]}{\point}\cr
&=\zerovec{}.
\end{align}
\endthm
%%%%%%%%%%%%%%%%%%%%%%%%%%%%%%%%%%%%%%%%%%%%%%%%%%%%%%%%%%%%%%%%%%%%%%%%%%%%%%%%%%%%%%%%%%%%%%%%%%%%%%%%%%%%%%%%%%%%%%%%%%%%%%%%
\corollary\label{corconnectionpointwisewrtvectorfields}
Let $\avecf{1}$ and $\avecf{2}$ be a pair of vector fields on $\vbbase{}$ such that
$\func{\avecf{1}}{\point}=\func{\avecf{2}}{\point}$ for some $\point\in\vB{}$.
\begin{equation}
\Foreach{\vbsec{}}{\vbsections{\vbundle{}}}
\func{\[\con{\avecf{1}}{\vbsec{}}\]}{\point}=
\func{\[\con{\avecf{2}}{\vbsec{}}\]}{\point}.
\end{equation}
\endcor
%%%%%%%%%%%%%%%%%%%%%%%%%%%%%%%%%%%%%%%%%%%%%%%%%%%%%%%%%%%%%%%%%%%%%%%%%%%%%%%%%%%%%%%%%%%%%%%%%%%%%%%%%%%%%%%%%%%%%%%%%%%%%%%%
\definition\label{defcovariantderivativeinducedbyconnection}
The map $\function{\TBcon{\connection{}}}{\Cprod{\tanbun{\vbbase{}}}{\vbsections{\vbundle{}}}}{\vbtotal{}}$
is defined as the following.\\
Given a $v\in\tanspace{\point}{\vbbase{}}$ for some $\point\in\vB{}$, and a $\vbsec{}\in\vbsections{\vbundle{}}$,
choose a vector field $\avecf{}\in\vectorfields{\vbbase{}}$ such that $\func{\avecf{}}{\point}=v$ and,
\begin{equation}
\func{\TBcon{\connection{}}}{\binary{v}{\vbsec{}}}\eqdef
\func{\[\con{\avecf{}}{\vbsec{}}\]}{\point}.
\end{equation}
This definition is permissible according to \refcor{corconnectionpointwisewrtvectorfields}.
$\TBcon{\connection{}}$ is called the $\quotl$covariant derivative on $\vbundle{}$ induced by the
connection $\connection{}$$\quotr$.
$\func{\TBcon{\connection{}}}{\binary{v}{\vbsec{}}}$ can alternatively be denoted by
$\tbcon{v}{\vbsec{}}$.
\endef
%%%%%%%%%%%%%%%%%%%%%%%%%%%%%%%%%%%%%%%%%%%%%%%%%%%%%%%%%%%%%%%%%%%%%%%%%%%%%%%%%%%%%%%%%%%%%%%%%%%%%%%%%%%%%%%%%%%%%%%%%%%%%%%%
\remark
Given a non-empty open subset $\U$ of $\vbbase{}$, the covariant derivative on the restricted vector bundle
$\subman{\vbundle{}}{\U}$ induced by the connection $\rescon{\connection{}}{\U}$, that is
$\TBcon{\rescon{\connection{}}{\U}}$, will simply be denoted by $\rescon{\TBcon{\connection{}}}{\U}$, which is called
the $\quotl$restriction of the covariant derivative $\TBcon{\connection{}}$ to $\U$$\quotr$. It can be easily verified that,
the restriction of covariant derivatives to the open sets possess the same
naturality properties as those of the connections. Thus it is plausible to denote $\rescon{\TBcon{\connection{}}}{\U}$
simply by $\TBcon{\connection{}}$, by an abuse of notation, when there is no ambiguity about the open set $\U$
to which the covariant derivative is restricted.
\endremark
%%%%%%%%%%%%%%%%%%%%%%%%%%%%%%%%%%%%%%%%%%%%%%%%%%%%%%%%%%%%%%%%%%%%%%%%%%%%%%%%%%%%%%%%%%%%%%%%%%%%%%%%%%%%%%%%%%%%%%%%%%%%%%%%
\remark
Given a manifold $\Man{}$ of dimension $m$, it is a trivial fact that the tangent-bundle $\Tanbun{\Man{}}$
possesses a canonical smooth vector bundle structure, having $\tanbun{\Man{}}$, $\Man{}$, $\basep{\Man{}}$,
and $\R^{m}$ as its total space, base space, projection map, and fiber, respectively. We denote this vector bundle
by the same notation $\Tanbun{\Man{}}$, for convenience.
For every $\point\in\Man{}$, the fiber space of $\Tanbun{\Man{}}$ over $\point$ is simply the tangent-space of $\Man{}$
at $\point$, that is $\tanspace{\point}{\Man{}}$. Since the set of sections of $\Tanbun{}$ as a vector bundle
coincides with the set of smooth vector fields on $\Man{}$, any connection $\connection{}$ on $\Tanbun{\Man{}}$
must be a map $\function{\connection{}}{\Cprod{\vectorfields{\Man{}}}{\vectorfields{\Man{}}}}{\vectorfields{\Man{}}}$
satisfying the conditions of \refdef{defconnection}.\\
Any connection on the tangent-bundle of $\Man{}$ will be called a $\quotl$(affine) connection on the manifold $\Man{}$$\quotr$.
\endremark
%%%%%%%%%%%%%%%%%%%%%%%%%%%%%%%%%%%%%%%%%%%%%%%%%%%%%%%%%%%%%%%%%%%%%%%%%%%%%%%%%%%%%%%%%%%%%%%%%%%%%%%%%%%%%%%%%%%%%%%%%%%%%%%%
%%%%%%%%%%%%%%%%%%%%%%%%%%%%%%%%%%%%%%%%%%%%%%%%%%%%%%%%%%%%%%%%%%%%%%%%%%%%%%%%%%%%%%%%%%%%%%%%%%%%%%%%%%%%%%%%%%%%%%%%%%%%%%%%
%%%%%%%%%%%%%%%%%%%%%%%%%%%%%%%%%%%%%%%%%%%%%%%%%%%%%%%%%%%%%%%%%%%%%%%%%%%%%%%%%%%%%%%%%%%%%%%%%%%%%%%%%%%%%%%%%%%%%%%%%%%%%%%%
%%%%%%%%%%%%%%%%%%%%%%%%%%%%%%%%%%%%%%%%%%%%%%%%%%%%%%%%%%%%%%%%%%%%%%%%%%%%%%%%%%%%%%%%%%%%%%%%%%%%%%%%%%%%%%%%%%%%%%%%%%%%%%%%
%%%%%%%%%%%%%%%%%%%%%%%%%%%%%%%%%%%%%%%%%%%%%%%%%%%%%%%%%%%%%%%%%%%%%%%%%%%%%%%%%%%%%%%%%%%%%%%%%%%%%%%%%%%%%%%%%%%%%%%%%%%%%%%%
%%%%%%%%%%%%%%%%%%%%%%%%%%%%%%%%%%%%%%%%%%%%%%%%%%%%%%%%%%%%%%%%%%%%%%%%%%%%%%%%%%%%%%%%%%%%%%%%%%%%%%%%%%%%%%%%%%%%%%%%%%%%%%%%
%%%%%%%%%%%%%%%%%%%%%%%%%%%%%%%%%%%%%%%%%%%%%%%%%%%%%%%%%%%%%%%%%%%%%%%%%%%%%%%%%%%%%%%%%%%%%%%%%%%%%%%%%%%%%%%%%%%%%%%%%%%%%%%%
%%%%%%%%%%%%%%%%%%%%%%%%%%%%%%%%%%%%%%%%%%%%%%%%%%%%%%%%%%%%%%%%%%%%%%%%%%%%%%%%%%%%%%%%%%%%%%%%%%%%%%%%%%%%%%%%%%%%%%%%%%%%%%%%
%%%%%%%%%%%%%%%%%%%%%%%%%%%%%%%%%%%%%%%%%%%%%%%%%%%%%%%%%%%%%%%%%%%%%%%%%%%%%%%%%%%%%%%%%%%%%%%%%%%%%%%%%%%%%%%%%%%%%%%%%%%%%%%%
%%%%%%%%%%%%%%%%%%%%%%%%%%%%%%%%%%%%%%%%%%%%%%%%%%%%%%%%%%%%%%%%%%%%%%%%%%%%%%%%%%%%%%%%%%%%%%%%%%%%%%%%%%%%%%%%%%%%%%%%%%%%%%%%
%%%%%%%%%%%%%%%%%%%%%%%%%%%%%%%%%%%%%%%%%%%%%%%%%%%%%%%%%%%%%%%%%%%%%%%%%%%%%%%%%%%%%%%%%%%%%%%%%%%%%%%%%%%%%%%%%%%%%%%%%%%%%%%%
\section{Connection Forms}
\textit{
Let $\opair{\U}{\phi}$ be a local trivialization of the vector bundle $\vbundle{}$. So,
$\function{\phi}{\func{\pimage{\phi}}{\U}}{\Cprod{\U}{\R^d}}$ is a diffeomorphism from
$\subman{\vbtotal{}}{\func{\pimage{\phi}}{\U}}$ to $\manprod{\subman{\vbbase{}}{\U}}{\R^d}$.
We denote by $\mtuple{\vbsecc{1}}{\vbsecc{d}}$ the system of local frame fields of the vector bundle
$\vbundle{}$ corresponded to this local trivialization. Hence, for each $i\in\seta{\suc{1}{d}}$,
$\function{\vbsecc{i}}{\U}{\func{\pimage{\phi}}{\U}}$ is the section of the restricted vector bundle
$\subman{\vbundle{}}{\U}$ defined by
\begin{equation}
\Foreach{\point}{\U}\func{\vbsecc{i}}{\point}\eqdef\func{\finv{\phi}}{\binary{\point}{\Eucbase{n}{i}}},
\end{equation}
where $\mtuple{\Eucbase{n}{1}}{\Eucbase{n}{n}}$ denotes the standard orthonormal ordered basis of $\R^d$.
Moreover, it is obvious that for every $\point\in\U$, the set $\seta{\suc{\func{\vbsecc{1}}{\point}}{\func{\vbsecc{1}}{\point}}}$
forms a basis of the vector-space $\fibervecs{\vbundle{}}{\point}=\fibervecs{\subman{\vbundle{}}{\U}}{\point}$,
that is the fiber space of the vector bundle $\vbundle{}$ at $\point$.\\
Given a vector field $\avecff{}\in\vectorfields{\subman{\vbbase{}}{\U}}$, for each $j\in\seta{\suc{1}{d}}$,
$\func{\rescon{\connection{}}{\U}}{\binary{\avecff{}}{\vbsecc{j}}}$ has a unique decomposition as
$\displaystyle\func{\rescon{\connection{}}{\U}}{\binary{\avecff{}}{\vbsecc{j}}}=
\sum_{i=1}^{d}\func{\connectionform{\phi}{i}{j}}{\avecff{}}\vbsecc{i}$, where for each $i$ and $j$
in $\seta{\suc{1}{d}}$, $\function{\func{\connectionform{\phi}{i}{j}}{\avecff{}}}{\U}{\R}$ is smooth,
that is $\func{\connectionform{\phi}{i}{j}}{\avecff{}}\in\smoothmaps{\subman{\vbbase{}}{\U}}$.\\
Thus, given a local trivialization $\opair{\U}{\phi}$, we are provided with mappings
$\function{\connectionform{\phi}{i}{j}}{\vectorfields{\subman{\vbbase{}}{\U}}}{\smoothmaps{\subman{\vbbase{}}{\U}}}$
for each $i$ and $j$ in $\seta{\suc{1}{d}}$, as defined above. It is easy to show that every
$\connectionform{\phi}{i}{j}$ is a module-homomorphism ($\smoothmaps{\subman{\vbbase{}}{\U}}$-linear)
map, when $\vectorfields{\subman{\vbbase{}}{\U}}$ and $\smoothmaps{\subman{\vbbase{}}{\U}}$ are considered to be endowed
with their canonical $\smoothmaps{\subman{\vbbase{}}{\U}}$-module structures; it is actually a
trivial consequence of the $\smoothmaps{\subman{\vbbase{}}{\U}}$-linearity of $\rescon{\connection{}}{\U}$ with respect
to the first variable. Therefore, each $\connectionform{\phi}{i}{j}$ is actually a $1$-covariant tensor field on
$\subman{\vbbase{}}{\U}$, or equivalently a $1$-form on $\subman{\vbbase{}}{\U}$, that is an element of
$\TFF{1}{0}{\subman{\vbbase{}}{\U}}$.
}
%%%%%%%%%%%%%%%%%%%%%%%%%%%%%%%%%%%%%%%%%%%%%%%%%%%%%%%%%%%%%%%%%%%%%%%%%%%%%%%%%%%%%%%%%%%%%%%%%%%%%%%%%%%%%%%%%%%%%%%%%%%%%%%%
\definition
Let $\opair{\U}{\phi}$ be a local trivialization of the vector bundle $\vbundle{}$, and let
$\mtuple{\vbsecc{1}}{\vbsecc{d}}$ denote the system of local frame fields of the vector bundle
$\vbundle{}$ corresponded to this local trivialization. For each $i$ and $j$ in $\seta{\suc{1}{d}}$,
we define the map
$\function{\connectionform{\phi}{i}{j}}{\vectorfields{\subman{\vbbase{}}{\U}}}{\smoothmaps{\subman{\vbbase{}}{\U}}}$ in the
following way.
For every $\avecff{}\in\vectorfields{\subman{\vbbase{}}{\U}}$, $\func{\connectionform{\phi}{i}{j}}{\avecff{}}$
is defined to be the unique element of the ring $\smoothmaps{\subman{\vbbase{}}{\U}}$ that is the coefficient of
$\vbsecc{i}$ in the expansion
$\displaystyle\func{\rescon{\connection{}}{\U}}{\binary{\avecff{}}{\vbsecc{j}}}=
\sum_{i=1}^{d}\func{\connectionform{\phi}{i}{j}}{\avecff{}}\vbsecc{i}$.\\
When there is no ambiguity about the local trivialization $\phi$, $\connectionform{\phi}{i}{j}$
can simply be denoted by $\connectionform{}{i}{j}$.
\endef
%%%%%%%%%%%%%%%%%%%%%%%%%%%%%%%%%%%%%%%%%%%%%%%%%%%%%%%%%%%%%%%%%%%%%%%%%%%%%%%%%%%%%%%%%%%%%%%%%%%%%%%%%%%%%%%%%%%%%%%%%%%%%%%%
\theorem
Let $\opair{\U}{\phi}$ be a local trivialization of the vector bundle $\vbundle{}$.
For each $i$ and $j$ in $\seta{\suc{1}{d}}$, $\connectionform{\phi}{i}{j}$ is a module-homomorphism
from the $\smoothmaps{\subman{\vbbase{}}{\U}}$-module 
$\vectorfields{\subman{\vbbase{}}{\U}}$ to $\smoothmaps{\subman{\vbbase{}}{\U}}$,
and hence a $\opair{1}{0}$ tensor field on the restricted vector bundle
$\subman{\vbundle{}}{\U}$.
\proof
Let $\mtuple{\vbsecc{1}}{\vbsecc{d}}$ denote the system of local frame fields of the vector bundle
$\vbundle{}$ corresponded to the local trivialization $\opair{\U}{\phi}$.
Let $\avecff{}$ and $\avecff{1}$ be a pair of vector fields of $\subman{\vbbase{}}{\U}$,
and let $f\in\smoothmaps{\subman{\vbbase{}}{\U}}$. Clearly, for every $j\in\seta{\suc{1}{d}}$,
\begin{align}
\sum_{i=1}^{d}\func{\connectionform{}{i}{j}}{f\avecff{}+\avecff{1}}\vbsecc{i}&=
\func{\rescon{\connection{}}{\U}}{\binary{f\avecff{}+\avecff{1}}{\vbsecc{j}}}\cr
&=f\func{\rescon{\connection{}}{\U}}{\binary{\avecff{}}{\vbsecc{j}}}+
\func{\rescon{\connection{}}{\U}}{\binary{\avecff{1}}{\vbsecc{j}}}\cr
&=f\[\sum_{i=1}^{d}\func{\connectionform{}{i}{j}}{\avecff{}}\vbsecc{i}\]+
\sum_{i=1}^{d}\func{\connectionform{}{i}{j}}{\avecff{1}}\vbsecc{i}\cr
&=\sum_{i=1}^{d}\(f\func{\connectionform{}{i}{j}}{\avecff{}}+\func{\connectionform{}{i}{j}}{\avecff{1}}\)\vbsecc{i},
\end{align}
and therefore, for every $i$ and $j$ in $\seta{\suc{1}{d}}$,
\begin{equation}
\func{\connectionform{}{i}{j}}{f\avecff{}+\avecff{1}}=
f\func{\connectionform{}{i}{j}}{\avecff{}}+\func{\connectionform{}{i}{j}}{\avecff{1}}.
\end{equation}
\endthm
%%%%%%%%%%%%%%%%%%%%%%%%%%%%%%%%%%%%%%%%%%%%%%%%%%%%%%%%%%%%%%%%%%%%%%%%%%%%%%%%%%%%%%%%%%%%%%%%%%%%%%%%%%%%%%%%%%%%%%%%%%%%%%%%
\theorem\label{thmconnectionformsexpansion0}
Let $\opair{\U}{\phi}$ be a local trivialization of the vector bundle $\vbundle{}$, and let
$\mtuple{\vbsecc{1}}{\vbsecc{d}}$ denote the system of local frame fields of the vector bundle
$\vbundle{}$ corresponded to this local trivialization.
Let $\vbsec{}\in\vbsections{\vbundle{}}$ and $\avecf{}\in\vectorfields{\vbbase{}}$.
Let $\displaystyle\reS{\vbsec{}}{\U}=\sum_{i=1}^{d}\seccof{i}\vbsecc{i}$.
%$\displaystyle\avecf{}=\sum_{\mu=1}^{n}\veccof{\mu}{\avecff{\mu}}$.
\begin{equation}
\reS{\(\con{\avecf{}}{\vbsec{}}\)}{\U}=
\sum_{i=1}^{d}\[\lieder{\reS{\avecf{}}{\U}}{\seccof{i}}+
\sum_{j=1}^{d}\func{\connectionform{}{i}{j}}{\reS{\avecf{}}{\U}}\seccof{j}\]\vbsecc{i}.
\end{equation}
\proof
According to \refthm{thmnaturalityofconnection} and \refdef{defconnection},
\begin{align}
\reS{\(\con{\avecf{}}{\vbsec{}}\)}{\U}&=
\func{\rescon{\connection{}}{\U}}{\binary{\reS{\avecf{}}{\U}}{\reS{\vbsec{}}{\U}}}\cr
&=\func{\rescon{\connection{}}{\U}}{\binary{\reS{\avecf{}}{\U}}{\sum_{i=1}^{d}\seccof{i}\vbsecc{i}}}\cr
&=\sum_{i=1}^{d}\(\lieder{\reS{\avecf{}}{\U}}{\seccof{i}}\)\vbsecc{i}+
\sum_{j=1}^{d}\seccof{j}\[\func{\rescon{\connection{}}{\U}}{\binary{\reS{\avecf{}}{\U}}{\vbsecc{j}}}\]\cr
&=\sum_{i=1}^{d}\(\lieder{\reS{\avecf{}}{\U}}{\seccof{i}}\)\vbsecc{i}+
\sum_{j=1}^{d}\seccof{j}\[\sum_{i=1}^{d}\func{\connectionform{}{i}{j}}{\reS{\avecf{}}{\U}}\vbsecc{i}\]\cr
&=\sum_{i=1}^{d}\[\lieder{\reS{\avecf{}}{\U}}{\seccof{i}}+
\sum_{j=1}^{d}\func{\connectionform{}{i}{j}}{\reS{\avecf{}}{\U}}\seccof{j}\]\vbsecc{i}.
\end{align}
\endthm
%%%%%%%%%%%%%%%%%%%%%%%%%%%%%%%%%%%%%%%%%%%%%%%%%%%%%%%%%%%%%%%%%%%%%%%%%%%%%%%%%%%%%%%%%%%%%%%%%%%%%%%%%%%%%%%%%%%%%%%%%%%%%%%%
\definition\label{defconnectionformsch}
Choose a small enough open set $\U$ of $\vbbase{}$, and let $\opair{\U}{\phi}$ be a local trivialization of
$\vbundle{}$ and $\opair{\U}{\psi}$ be a chart of the manifold $\vbbase{}$. Let
$\mtuple{\localframevecf{1}}{\localframevecf{n}}$
denote the local frame of vector fields corresponded to the chart $\psi$. For each $i$ and $j$ in $\seta{\suc{1}{d}}$
and each $\mu$ in $\seta{\suc{1}{n}}$, we define,
\begin{equation}
\connectionformch{\opair{\phi}{\psi}}{i}{j}{\mu}:=
\func{\connectionform{\phi}{i}{j}}{\localframevecf{\mu}}.
\end{equation}
When the choices of the local trivialization $\phi$ and the chart $\psi$ are clear enough,
$\connectionformch{\opair{\phi}{\psi}}{i}{j}{\mu}$ can be simply denoted by
$\connectionformch{}{i}{j}{\mu}$.
\endef
%%%%%%%%%%%%%%%%%%%%%%%%%%%%%%%%%%%%%%%%%%%%%%%%%%%%%%%%%%%%%%%%%%%%%%%%%%%%%%%%%%%%%%%%%%%%%%%%%%%%%%%%%%%%%%%%%%%%%%%%%%%%%%%%
\theorem
Choose a small enough open set $\U$ of $\vbbase{}$, and let $\opair{\U}{\phi}$ be a local trivialization of
$\vbundle{}$ and $\opair{\U}{\psi}$ be a chart of the manifold $\vbbase{}$. Let
$\mtuple{\localframevecf{1}}{\localframevecf{n}}$
denote the local frame of vector fields corresponded to the chart $\psi$, and let
$\mtuple{\vbsecc{1}}{\vbsecc{d}}$ denote the system of local frame fields of the vector bundle
$\vbundle{}$ corresponded to this local trivialization.
Let $\vbsec{}\in\vbsections{\vbundle{}}$ and $\avecf{}\in\vectorfields{\vbbase{}}$.
Let $\displaystyle\reS{\avecf{}}{\U}=\sum_{\mu=1}^{n}\veccof{\mu}\localframevecf{\mu}$, and
let $\displaystyle\reS{\vbsec{}}{\U}=\sum_{i=1}^{d}\seccof{i}\vbsecc{i}$.
\begin{equation}
\reS{\(\con{\avecf{}}{\vbsec{}}\)}{\U}=
\sum_{i=1}^{d}\[\sum_{j=1}^{d}\sum_{\mu=1}^{n}\(\veccof{\mu}\(\lieder{\localframevecf{\mu}}{\seccof{i}}\)+
\veccof{\mu}\connectionformch{}{i}{j}{\mu}\seccof{j}\)\]\vbsecc{i}.
\end{equation}
\proof
According to \refdef{defconnection}, \refthm{thmconnectionformsexpansion0}, and \refdef{defconnectionformsch},
it is straightforward.
\endthm
%%%%%%%%%%%%%%%%%%%%%%%%%%%%%%%%%%%%%%%%%%%%%%%%%%%%%%%%%%%%%%%%%%%%%%%%%%%%%%%%%%%%%%%%%%%%%%%%%%%%%%%%%%%%%%%%%%%%%%%%%%%%%%%%
\textit{
Let $\Man{}$ be a manifold with dimension $m$, and let $\opair{\U}{\phi}$ be a chart of $\Man{}$.
It is known that corresponded to $\phi$ is a canonical local trivialization of the tangent-bundle
$\Tanbun{\Man{}}$ regarded as a vector bundle, which we will simply denote by $\ltchart{\phi}$.
Precisely, $\function{\ltchart{\phi}}{\func{\pimage{\basep{}}}{\U}}{\Cprod{\U}{\R^{m}}}$ is such that
$\func{\ltchart{\phi}}{v}=\opair{\point}{\func{\tanspaceiso{\point}{\Man{}}{\phi}}{v}}$, where $\point=\func{\basep{}}{v}$.
It is obvious that system of local frame fields of the local trivialization $\ltchart{\phi}$ of the
vector bundle $\Tanbun{\Man{}}$ coincides with the local frame of vector fields associated with the chart
$\phi$ of the manifold $\Man{}$.
}
%%%%%%%%%%%%%%%%%%%%%%%%%%%%%%%%%%%%%%%%%%%%%%%%%%%%%%%%%%%%%%%%%%%%%%%%%%%%%%%%%%%%%%%%%%%%%%%%%%%%%%%%%%%%%%%%%%%%%%%%%%%%%%%%
\definition\label{defChristoffelSymbols}
Let $\Man{}$ be a manifold with dimension $m$, and let $\opair{\U}{\phi}$ be a chart of $\Man{}$.
Let $\mtuple{\localframevecf{1}}{\localframevecf{m}}$
denote the local frame of vector fields on $\Man{}$ corresponded to the chart $\phi$. For every $i$, $j$, and $k$ in
$\seta{\suc{1}{m}}$, we define,
\begin{equation}
\Christoffel{\phi}{i}{j}{k}:=\connectionform{\opair{\ltchart{\phi}}{\phi}}{i}{j}{k}=
\func{\connectionform{\ltchart{\phi}}{i}{j}}{\localframevecf{k}}.
\end{equation}
When the chart $\phi$ is clear enough, we will denote $\Christoffel{\phi}{i}{j}{k}$ simply by
$\Christoffel{}{i}{j}{k}$. $\Christoffel{}{i}{j}{k}$-s are called the $\quotl$Christoffel symbols of
the connection $\connection{}$ on $\Man{}$ with respect to the chart $\opair{\U}{\phi}$.\\
Obviously, each $\Christoffel{}{i}{j}{k}$ is a smooth map, that is an element of $\smoothmaps{\subman{\Man{}}{\U}}$, and,
\begin{equation}
\func{\rescon{\connection{}}{\U}}{\binary{\localframevecf{k}}{\localframevecf{j}}}=
\sum_{i=1}^{m}\Christoffel{}{i}{j}{k}\localframevecf{i}.
\end{equation}
Here, we can simply denote $\func{\rescon{\connection{}}{\U}}{\binary{\localframevecf{k}}{\localframevecf{j}}}$
by $\con{\localframevecf{k}}{\localframevecf{j}}$ for convenience, since the restriction of $\connection{}$
to $\U$ works in almost the same way as the global connection $\connection{}$ itself when the attention is focused on
the region $\U$.
\endef
%%%%%%%%%%%%%%%%%%%%%%%%%%%%%%%%%%%%%%%%%%%%%%%%%%%%%%%%%%%%%%%%%%%%%%%%%%%%%%%%%%%%%%%%%%%%%%%%%%%%%%%%%%%%%%%%%%%%%%%%%%%%%%%%
\corollary\label{corlocalexpressionofconnection}
Let $\Man{}$ be a manifold with dimension $m$, and let $\opair{\U}{\phi}$ be a chart of $\Man{}$.
Let $\mtuple{\localframevecf{1}}{\localframevecf{m}}$
denote the local frame of vector fields on $\Man{}$ corresponded to the chart $\phi$.
Let $\avecf{}$ and $\avecff{}$ be elements of $\vectorfields{\Man{}}$, and let
$\displaystyle\reS{\avecf{}}{\U}=\sum_{i=1}^{m}\veccof{i}{\localframevecf{i}}$ and
$\displaystyle\reS{\avecff{}}{\U}=\sum_{i=1}^{m}\veccoff{i}{\localframevecf{i}}$.
\begin{equation}
\reS{\(\con{\avecf{}}{\avecff{}}\)}{\U}=
\sum_{i=1}^{m}\[\sum_{j=1}^{m}\sum_{k=1}^{m}\(\veccof{k}\(\lieder{\localframevecf{k}}{\veccoff{i}}\)+
\veccof{k}\Christoffel{}{i}{j}{k}\veccoff{j}\)\]\localframevecf{i}.
\end{equation}
\endthm
%%%%%%%%%%%%%%%%%%%%%%%%%%%%%%%%%%%%%%%%%%%%%%%%%%%%%%%%%%%%%%%%%%%%%%%%%%%%%%%%%%%%%%%%%%%%%%%%%%%%%%%%%%%%%%%%%%%%%%%%%%%%%%%%
%%%%%%%%%%%%%%%%%%%%%%%%%%%%%%%%%%%%%%%%%%%%%%%%%%%%%%%%%%%%%%%%%%%%%%%%%%%%%%%%%%%%%%%%%%%%%%%%%%%%%%%%%%%%%%%%%%%%%%%%%%%%%%%%
%%%%%%%%%%%%%%%%%%%%%%%%%%%%%%%%%%%%%%%%%%%%%%%%%%%%%%%%%%%%%%%%%%%%%%%%%%%%%%%%%%%%%%%%%%%%%%%%%%%%%%%%%%%%%%%%%%%%%%%%%%%%%%%%
%%%%%%%%%%%%%%%%%%%%%%%%%%%%%%%%%%%%%%%%%%%%%%%%%%%%%%%%%%%%%%%%%%%%%%%%%%%%%%%%%%%%%%%%%%%%%%%%%%%%%%%%%%%%%%%%%%%%%%%%%%%%%%%%
%%%%%%%%%%%%%%%%%%%%%%%%%%%%%%%%%%%%%%%%%%%%%%%%%%%%%%%%%%%%%%%%%%%%%%%%%%%%%%%%%%%%%%%%%%%%%%%%%%%%%%%%%%%%%%%%%%%%%%%%%%%%%%%%
%%%%%%%%%%%%%%%%%%%%%%%%%%%%%%%%%%%%%%%%%%%%%%%%%%%%%%%%%%%%%%%%%%%%%%%%%%%%%%%%%%%%%%%%%%%%%%%%%%%%%%%%%%%%%%%%%%%%%%%%%%%%%%%%
%%%%%%%%%%%%%%%%%%%%%%%%%%%%%%%%%%%%%%%%%%%%%%%%%%%%%%%%%%%%%%%%%%%%%%%%%%%%%%%%%%%%%%%%%%%%%%%%%%%%%%%%%%%%%%%%%%%%%%%%%%%%%%%%
%%%%%%%%%%%%%%%%%%%%%%%%%%%%%%%%%%%%%%%%%%%%%%%%%%%%%%%%%%%%%%%%%%%%%%%%%%%%%%%%%%%%%%%%%%%%%%%%%%%%%%%%%%%%%%%%%%%%%%%%%%%%%%%%
%%%%%%%%%%%%%%%%%%%%%%%%%%%%%%%%%%%%%%%%%%%%%%%%%%%%%%%%%%%%%%%%%%%%%%%%%%%%%%%%%%%%%%%%%%%%%%%%%%%%%%%%%%%%%%%%%%%%%%%%%%%%%%%%
%%%%%%%%%%%%%%%%%%%%%%%%%%%%%%%%%%%%%%%%%%%%%%%%%%%%%%%%%%%%%%%%%%%%%%%%%%%%%%%%%%%%%%%%%%%%%%%%%%%%%%%%%%%%%%%%%%%%%%%%%%%%%%%%
%%%%%%%%%%%%%%%%%%%%%%%%%%%%%%%%%%%%%%%%%%%%%%%%%%%%%%%%%%%%%%%%%%%%%%%%%%%%%%%%%%%%%%%%%%%%%%%%%%%%%%%%%%%%%%%%%%%%%%%%%%%%%%%%
\section{Pullbacks of Connections}
%%%%%%%%%%%%%%%%%%%%%%%%%%%%%%%%%%%%%%%%%%%%%%%%%%%%%%%%%%%%%%%%%%%%%%%%%%%%%%%%%%%%%%%%%%%%%%%%%%%%%%%%%%%%%%%%%%%%%%%%%%%%%%%%
\fixed
$\vbundle{}=\quintuple{\vbtotal{}}{\vbprojection{}}{\vbbase{}}{\vbfiber{}}{\vbatlas{}}$ is fixed as a real smooth vector bundle
of rank $d$,
where $\vbtotal{}=\opair{\vTot{}}{\maxatlas{\vTot{}}}$ and
$\vbbase{}=\opair{\vB{}}{\maxatlas{\vB{}}}$ are $\difclass{\infty}$ manifolds
modeled on the Banach-spaces $\R^{n_{\vTot{}}}$ and $\R^{n_{\vB{}}}$, respectively.
Also, $\vbundle{1}=\quintuple{\vbtotal{1}}{\vbprojection{1}}{\vbbase{1}}{\vbfiber{1}}{\vbatlas{1}}$
is fixed as a smooth vector bundle of rank $d_1$.
\endfixed
%%%%%%%%%%%%%%%%%%%%%%%%%%%%%%%%%%%%%%%%%%%%%%%%%%%%%%%%%%%%%%%%%%%%%%%%%%%%%%%%%%%%%%%%%%%%%%%%%%%%%%%%%%%%%%%%%%%%%%%%%%%%%%%%
\textit{
The definition of the pullback of a connection from a smooth vector bundle to another isomorphic one, relies on the notion of
pushforward of a section of a vector bundle to that of another isomorphic vector bundle. Suppose that the smooth vector bundles
$\vbundle{}$ and $\vbundle{1}$ are isomorphic with $\function{f}{\vbtotal{1}}{\vbtotal{}}$ as an isomorphism between them
along the diffeomorphism $\function{\fbmorb{f}}{\vbbase{1}}{\vbbase{}}$. It is known from the theory of vector bundles that
$\fbmorb{f}$ is the unique diffeomorphism from $\vbbase{1}$ to $\vbbase{}$ such that
$\cmp{\vbprojection{}}{f}=\cmp{\fbmorb{f}}{\vbprojection{1}}$.
There exists a map $\function{\sectionpushforward{f}}{\vbsections{\vbundle{1}}}{\vbsections{\vbundle{}}}$, defined as,
\begin{align}
\Foreach{\vbsec{}}{\vbsections{\vbundle{1}}}
\Foreach{\point}{\vbbase{}}
\func{\[\func{\sectionpushforward{f}}{\vbsec{}}\]}{\point}\eqdef
\func{f}{\func{\vbsec{}}{\func{\finv{\fbmorb{f}}}{\point}}}.
\end{align}
It is easy to show that $\sectionpushforward{f}$ is a linear-isomorphism from $\vbsections{\vbundle{1}}$ to
$\vbsections{\vbundle{}}$ endowed with their canonical linear structures.
Moreover, given $h\in\smoothmaps{\vbbase{1}}$  and $\vbsec{}\in\vbsections{\vbundle{1}}$, it is obvious that
$\func{\sectionpushforward{f}}{h\vbsec{}}=\(\cmp{h}{\finv{\fbmorb{f}}}\)\func{\sectionpushforward{f}}{\vbsec{}}$.\\
In the particular case of the tangent-bundles of manifolds, we slightly alter the notation.
Given a pair of diffeomorphic manifolds $\Man{}$ and $\Man{1}$, with $\function{F}{\Man{1}}{\Man{}}$
a diffeomorphism from $\Man{1}$ to $\Man{}$, we define the map
$\function{\sectionpushforward{F}}{\vectorfields{\Man{1}}}{\vectorfields{\Man{}}}$ as
$\sectionpushforward{F}:=\sectionpushforward{\(\derr{F}\)}$, where $\derr{F}$ denotes the tangent-map of $F$.
So,
\begin{align}
\Foreach{\avecf{}}{\vectorfields{\Man{1}}}
\Foreach{\point}{\Man{}}
\func{\[\func{\sectionpushforward{F}}{\avecf{}}\]}{\point}\eqdef
\func{\derr{F}}{\func{\avecf{}}{\func{\finv{F}}{\point}}}.
\end{align}
It is easy to check that given $\avecf{}\in\vectorfields{\Man{1}}$, $h\in\smoothmaps{\Man{}}$,
$\lieder{\func{\sectionpushforward{F}}{\avecf{}}}{\(\cmp{h}{\finv{F}}\)}=
\cmp{\(\lieder{\avecf{}}{h}\)}{\finv{F}}$.
}
%%%%%%%%%%%%%%%%%%%%%%%%%%%%%%%%%%%%%%%%%%%%%%%%%%%%%%%%%%%%%%%%%%%%%%%%%%%%%%%%%%%%%%%%%%%%%%%%%%%%%%%%%%%%%%%%%%%%%%%%%%%%%%%%
\definition
Suppose that the smooth vector bundles $\vbundle{}$ and $\vbundle{1}$ are isomorphic, and let
$\function{f}{\vbtotal{1}}{\vbtotal{}}$ be a vector bundle isomorphism from $\vbundle{1}$ to $\vbundle{}$
along the diffeomorphism $\function{\fbmorb{f}}{\vbbase{1}}{\vbbase{}}$,
that is an element of $\vbisomorphisms{\vbundle{1}}{\vbundle{2}}$.
We define the map $\function{\connectionpb{f}}{\connections{\vbundle{}}}
{\Func{\Cprod{\vectorfields{\vbbase{1}}}{\vbsections{\vbundle{1}}}}{\vbsections{\vbundle{1}}}}$ as,
\begin{align}
&\Foreach{\connection{}}{\connections{\vbundle{}}}
\Foreach{\avecf{}}{\vectorfields{\vbbase{1}}}
\Foreach{\vbsec{}}{\vbsections{\vbundle{1}}}\cr
&\func{\[\func{\connectionpb{f}}{\connection{}}\]}{\binary{\avecf{}}{\vbsec{}}}\eqdef
\func{\sectionpushforward{\(\finv{f}\)}}{\con{\func{\sectionpushforward{\fbmorb{f}}}{\avecf{}}}
{\func{\sectionpushforward{f}}{\vbsec{}}}}.
\end{align}
\endef
%%%%%%%%%%%%%%%%%%%%%%%%%%%%%%%%%%%%%%%%%%%%%%%%%%%%%%%%%%%%%%%%%%%%%%%%%%%%%%%%%%%%%%%%%%%%%%%%%%%%%%%%%%%%%%%%%%%%%%%%%%%%%%%%
\lemma
Suppose that the smooth vector bundles $\vbundle{}$ and $\vbundle{1}$ are isomorphic, and let
$\function{f}{\vbtotal{1}}{\vbtotal{}}$ be a vector bundle isomorphism from $\vbundle{1}$ to $\vbundle{}$
along the diffeomorphism $\function{\fbmorb{f}}{\vbbase{1}}{\vbbase{}}$,
that is an element of $\vbisomorphisms{\vbundle{1}}{\vbundle{2}}$. For every connection $\connection{}$
on the smooth vector bundle $\vbundle{}$, $\func{\connectionpb{f}}{\connection{}}$ is a connection on
the smooth vector bundle $\vbundle{1}$.
\proof
Let $\binary{\avecf{}}{\avecf{1}}\in\vectorfields{\vbbase{1}}$, $\binary{\vbsec{}}{\vbsec{1}}\in\vbsections{\vbundle{1}}$,
and $\binary{h}{h_1}\in\smoothmaps{\vbbase{1}}$.
\begin{itemize}
\item[\myitem{pr-1.}]
\begin{align}
&~~~~\func{\[\func{\connectionpb{f}}{\connection{}}\]}{\binary{h\avecf{}+h_1\avecf{1}}{\vbsec{}}}\cr
&=\func{\sectionpushforward{\(\finv{f}\)}}{\con{\func{\sectionpushforward{\fbmorb{f}}}{h\avecf{}+h_1\avecf{1}}}{\func{\sectionpushforward{f}}{\vbsec{}}}}\cr
&=\func{\sectionpushforward{\(\finv{f}\)}}{\con{\(\cmp{h}{\finv{\fbmorb{f}}}\)\func{\sectionpushforward{\fbmorb{f}}}{\avecf{}}+
\(\cmp{h_1}{\finv{\fbmorb{f}}}\)\func{\sectionpushforward{\fbmorb{f}}}{\avecf{1}}}
{\func{\sectionpushforward{f}}{\vbsec{}}}}\cr
&=\func{\sectionpushforward{\(\finv{f}\)}}{\con{\(\cmp{h}{\finv{\fbmorb{f}}}\)\func{\sectionpushforward{\fbmorb{f}}}{\avecf{}}+
\(\cmp{h_1}{\finv{\fbmorb{f}}}\)\func{\sectionpushforward{\fbmorb{f}}}{\avecf{1}}}{\func{\sectionpushforward{f}}{\vbsec{}}}}\cr
&=\func{\sectionpushforward{\(\finv{f}\)}}{\(\cmp{h}{\finv{\fbmorb{f}}}\)\[\con{\func{\sectionpushforward{\fbmorb{f}}}{\avecf{}}}
{\func{\sectionpushforward{f}}{\vbsec{}}}\]+\(\cmp{h_1}{\finv{\fbmorb{f}}}\)\[\con{\func{\sectionpushforward{\fbmorb{f}}}{\avecf{1}}}
{\func{\sectionpushforward{f}}{\vbsec{}}}\]}\cr
&=h
%\(\cmp{\cmp{h}{\finv{\fbmorb{f}}}}{\finv{\fbmorb{\finv{f}}}}\)
\[\func{\sectionpushforward{\(\finv{f}\)}}{\con{\func{\sectionpushforward{\fbmorb{f}}}{\avecf{}}}
{\func{\sectionpushforward{f}}{\vbsec{}}}}\]+
h_1\[\func{\sectionpushforward{\(\finv{f}\)}}{\con{\func{\sectionpushforward{\fbmorb{f}}}{\avecf{1}}}
{\func{\sectionpushforward{f}}{\vbsec{}}}}\]\cr
&=h\func{\[\func{\connectionpb{f}}{\connection{}}\]}{\binary{\avecf{}}{\vbsec{}}}+
h_1\func{\[\func{\connectionpb{f}}{\connection{}}\]}{\binary{\avecf{1}}{\vbsec{}}}.
\end{align}
where, we have used the fact that $\finv{\fbmorb{\finv{f}}}=\fbmorb{f}$.
\item[\myitem{pr-2.}]
\begin{align}
\func{\[\func{\connectionpb{f}}{\connection{}}\]}{\binary{\avecf{}}{\vbsec{}+\vbsec{1}}}&=
\func{\sectionpushforward{\(\finv{f}\)}}{\con{\func{\sectionpushforward{\fbmorb{f}}}{\avecf{}}}
{\func{\sectionpushforward{f}}{\vbsec{}+\vbsec{1}}}}\cr
&=\func{\sectionpushforward{\(\finv{f}\)}}{\con{\func{\sectionpushforward{\fbmorb{f}}}{\avecf{}}}
{\[\func{\sectionpushforward{f}}{\vbsec{}}+\func{\sectionpushforward{f}}{\vbsec{1}}\]}}\cr
&=\func{\sectionpushforward{\(\finv{f}\)}}{\con{\func{\sectionpushforward{\fbmorb{f}}}{\avecf{}}}
{\func{\sectionpushforward{f}}{\vbsec{}}}+
\con{\func{\sectionpushforward{\fbmorb{f}}}{\avecf{}}}
{\func{\sectionpushforward{f}}{\vbsec{1}}}}\cr
&=\func{\sectionpushforward{\(\finv{f}\)}}{\con{\func{\sectionpushforward{\fbmorb{f}}}{\avecf{}}}
{\func{\sectionpushforward{f}}{\vbsec{}}}}+
\func{\sectionpushforward{\(\finv{f}\)}}{\con{\func{\sectionpushforward{\fbmorb{f}}}{\avecf{}}}
{\func{\sectionpushforward{f}}{\vbsec{1}}}}\cr
&=\func{\[\func{\connectionpb{f}}{\connection{}}\]}{\binary{\avecf{}}{\vbsec{}}}+
\func{\[\func{\connectionpb{f}}{\connection{}}\]}{\binary{\avecf{}}{\vbsec{1}}}.
\end{align}
\item[\myitem{pr-3.}]
\begin{align}
\func{\[\func{\connectionpb{f}}{\connection{}}\]}{\binary{\avecf{}}{h\vbsec{}}}&=
\func{\sectionpushforward{\(\finv{f}\)}}{\con{\func{\sectionpushforward{\fbmorb{f}}}{\avecf{}}}
{\func{\sectionpushforward{f}}{h\vbsec{}}}}\cr
&=\func{\sectionpushforward{\(\finv{f}\)}}{\con{\func{\sectionpushforward{\fbmorb{f}}}{\avecf{}}}
{\[\(\cmp{h}{\finv{\fbmorb{f}}}\)\func{\sectionpushforward{f}}{\vbsec{}}\]}}\cr
&=\func{\sectionpushforward{\(\finv{f}\)}}{\[\lieder{\func{\sectionpushforward{\fbmorb{f}}}{\avecf{}}}{\(\cmp{h}{\finv{\fbmorb{f}}}\)}\]
\func{\sectionpushforward{f}}{\vbsec{}}+
\(\cmp{h}{\finv{\fbmorb{f}}}\)\[\con{\func{\sectionpushforward{\fbmorb{f}}}{\avecf{}}}
{\func{\sectionpushforward{f}}{\vbsec{}}}\]}\cr
&=\(\cmp{\[\lieder{\func{\sectionpushforward{\fbmorb{f}}}{\avecf{}}}{\(\cmp{h}{\finv{\fbmorb{f}}}\)}\]}{\fbmorb{f}}\)
\[\func{\sectionpushforward{\(\finv{f}\)}}{\func{\sectionpushforward{f}}{\vbsec{}}}\]+
h\func{\sectionpushforward{\(\finv{f}\)}}{\con{\func{\sectionpushforward{\fbmorb{f}}}{\avecf{}}}
{\func{\sectionpushforward{f}}{\vbsec{}}}}\cr
&=\(\lieder{\avecf{}}{h}\)\vbsec{}+
h\func{\[\func{\connectionpb{f}}{\connection{}}\]}{\binary{\avecf{}}{\vbsec{}}}.
\end{align}
\end{itemize}
\endlem
%%%%%%%%%%%%%%%%%%%%%%%%%%%%%%%%%%%%%%%%%%%%%%%%%%%%%%%%%%%%%%%%%%%%%%%%%%%%%%%%%%%%%%%%%%%%%%%%%%%%%%%%%%%%%%%%%%%%%%%%%%%%%%%%
\definition
Suppose that the smooth vector bundles $\vbundle{}$ and $\vbundle{1}$ are isomorphic, and let
$\function{f}{\vbtotal{1}}{\vbtotal{}}$ be a vector bundle isomorphism from $\vbundle{1}$ to $\vbundle{}$
along the diffeomorphism $\function{\fbmorb{f}}{\vbbase{1}}{\vbbase{}}$,
that is an element of $\vbisomorphisms{\vbundle{1}}{\vbundle{2}}$.
The map $\function{\connectionpb{f}}{\connections{\vbundle{}}}{\connections{\vbundle{1}}}$ defined as,
\begin{align}
&\Foreach{\connection{}}{\connections{\vbundle{}}}
\Foreach{\avecf{}}{\vectorfields{\vbbase{1}}}
\Foreach{\vbsec{}}{\vbsections{\vbundle{1}}}\cr
&\func{\[\func{\connectionpb{f}}{\connection{}}\]}{\binary{\avecf{}}{\vbsec{}}}\eqdef
\func{\sectionpushforward{\(\finv{f}\)}}{\con{\func{\sectionpushforward{\fbmorb{f}}}{\avecf{}}}
{\func{\sectionpushforward{f}}{\vbsec{}}}},
\end{align}
is referred to as the $\quotl$pullback of the connection $\connection{}$ relative to $f$$\quotr$.
\endef
%%%%%%%%%%%%%%%%%%%%%%%%%%%%%%%%%%%%%%%%%%%%%%%%%%%%%%%%%%%%%%%%%%%%%%%%%%%%%%%%%%%%%%%%%%%%%%%%%%%%%%%%%%%%%%%%%%%%%%%%%%%%%%%%
%%%%%%%%%%%%%%%%%%%%%%%%%%%%%%%%%%%%%%%%%%%%%%%%%%%%%%%%%%%%%%%%%%%%%%%%%%%%%%%%%%%%%%%%%%%%%%%%%%%%%%%%%%%%%%%%%%%%%%%%%%%%%%%%
%%%%%%%%%%%%%%%%%%%%%%%%%%%%%%%%%%%%%%%%%%%%%%%%%%%%%%%%%%%%%%%%%%%%%%%%%%%%%%%%%%%%%%%%%%%%%%%%%%%%%%%%%%%%%%%%%%%%%%%%%%%%%%%%
%%%%%%%%%%%%%%%%%%%%%%%%%%%%%%%%%%%%%%%%%%%%%%%%%%%%%%%%%%%%%%%%%%%%%%%%%%%%%%%%%%%%%%%%%%%%%%%%%%%%%%%%%%%%%%%%%%%%%%%%%%%%%%%%
%%%%%%%%%%%%%%%%%%%%%%%%%%%%%%%%%%%%%%%%%%%%%%%%%%%%%%%%%%%%%%%%%%%%%%%%%%%%%%%%%%%%%%%%%%%%%%%%%%%%%%%%%%%%%%%%%%%%%%%%%%%%%%%%
%%%%%%%%%%%%%%%%%%%%%%%%%%%%%%%%%%%%%%%%%%%%%%%%%%%%%%%%%%%%%%%%%%%%%%%%%%%%%%%%%%%%%%%%%%%%%%%%%%%%%%%%%%%%%%%%%%%%%%%%%%%%%%%%
%%%%%%%%%%%%%%%%%%%%%%%%%%%%%%%%%%%%%%%%%%%%%%%%%%%%%%%%%%%%%%%%%%%%%%%%%%%%%%%%%%%%%%%%%%%%%%%%%%%%%%%%%%%%%%%%%%%%%%%%%%%%%%%%
%%%%%%%%%%%%%%%%%%%%%%%%%%%%%%%%%%%%%%%%%%%%%%%%%%%%%%%%%%%%%%%%%%%%%%%%%%%%%%%%%%%%%%%%%%%%%%%%%%%%%%%%%%%%%%%%%%%%%%%%%%%%%%%%
%%%%%%%%%%%%%%%%%%%%%%%%%%%%%%%%%%%%%%%%%%%%%%%%%%%%%%%%%%%%%%%%%%%%%%%%%%%%%%%%%%%%%%%%%%%%%%%%%%%%%%%%%%%%%%%%%%%%%%%%%%%%%%%%
%%%%%%%%%%%%%%%%%%%%%%%%%%%%%%%%%%%%%%%%%%%%%%%%%%%%%%%%%%%%%%%%%%%%%%%%%%%%%%%%%%%%%%%%%%%%%%%%%%%%%%%%%%%%%%%%%%%%%%%%%%%%%%%%
%%%%%%%%%%%%%%%%%%%%%%%%%%%%%%%%%%%%%%%%%%%%%%%%%%%%%%%%%%%%%%%%%%%%%%%%%%%%%%%%%%%%%%%%%%%%%%%%%%%%%%%%%%%%%%%%%%%%%%%%%%%%%%%%
%%%%%%%%%%%%%%%%%%%%%%%%%%%%%%%%%%%%%%%%%%%%%%%%%%%%%%%%%%%%%%%%%%%%%%%%%%%%%%%%%%%%%%%%%%%%%%%%%%%%%%%%%%%%%%%%%%%%%%%%%%%%%%%%
\section{Covariant Derivative of Vector Fields Along Curves}
\textit{
The notion of covariant derivative of vector fields along curves is derived from the notion of an affine connection
on a manifold. This can be generalized naturally to the case of Koszul connections on smooth vector bundles in a straightforward
manner, in order to acquire a similar notion of differentiating sections of a vector bundle along the curves of its base manifold.
Since this generalization is absolutely natural and straightforward, and since the case of differentiating the
vector fields of a manifold with an affine connection along curves lies in the center of our attention, it is preferred
to introduce and study the former case in detail here.
}
%%%%%%%%%%%%%%%%%%%%%%%%%%%%%%%%%%%%%%%%%%%%%%%%%%%%%%%%%%%%%%%%%%%%%%%%%%%%%%%%%%%%%%%%%%%%%%%%%%%%%%%%%%%%%%%%%%%%%%%%%%%%%%%%
\fixed
$\Man{}$ is fixed as a manifold with dimension $m$, and $\connection{}$ as an affine connection on $\Man{}$.
\endfixed
%%%%%%%%%%%%%%%%%%%%%%%%%%%%%%%%%%%%%%%%%%%%%%%%%%%%%%%%%%%%%%%%%%%%%%%%%%%%%%%%%%%%%%%%%%%%%%%%%%%%%%%%%%%%%%%%%%%%%%%%%%%%%%%%
\definition\label{defvectorfieldsalongcurves}
Let $\function{\Curve{}}{\interval{}}{\Man{}}$ be a smooth curve on $\Man{}$
for some connected open interval $\interval{}$ of $\R$.
The map $\function{\valongc{}}{\interval{}}{\tanbun{\Man{}}}$ is called a $\quotl$smooth vector field on $\Man{}$
along the curve $\Curve{}$$\quotr$, if $\valongc{}$ is smooth, that is an element of
$\mapdifclass{\infty}{\interval{}}{\Tanbun{\Man{}}}$, and,
\begin{equation}
\Foreach{t}{\interval{}}
\func{\valongc{}}{t}\in\tanspace{\func{\Curve{}}{t}}{\Man{}}.
\end{equation}
The set of all smooth vector fields on $\Man{}$ along $\Curve{}$ will be denoted by $\valongcurves{\Man{}}{\Curve{}}$.\\
\begin{itemize}
\item
The addition operation on $\valongcurves{\Man{}}{\Curve{}}$ is defined as,
\begin{equation}
\Foreach{\opair{\valongc{1}}{\valongc{2}}}{\Cprod{\valongcurves{\Man{}}{\Curve{}}}{\valongcurves{\Man{}}{\Curve{}}}}
\Foreach{t}{\interval}
\func{\[\valongc{1}+\valongc{2}\]}{t}\eqdef\func{\valongc{1}}{t}+\func{\valongc{2}}{t}.
\end{equation}
\item
The scalar-multiplication operation on $\valongcurves{\Man{}}{\Curve{}}$ is defined as,
\begin{equation}
\Foreach{\valongc{}}{\valongcurves{\Man{}}{\Curve{}}}
\Foreach{c}{\R}
\Foreach{t}{\interval{}}
\func{\[c\valongc{}\]}{t}\eqdef c\func{\valongc{}}{t}.
\end{equation}
\item
The $\smoothmaps{\interval{}}$-multiplication operation on $\valongcurves{\Man{}}{\Curve{}}$ is defined as,
\begin{equation}
\Foreach{\valongc{}}{\valongcurves{\Man{}}{\Curve{}}}
\Foreach{f}{\smoothmaps{\interval{}}}
\Foreach{t}{\interval{}}
\func{\[f\valongc{}\]}{t}\eqdef\func{f}{t}\func{\valongc{}}{t}.
\end{equation}
\end{itemize}
\endef
%%%%%%%%%%%%%%%%%%%%%%%%%%%%%%%%%%%%%%%%%%%%%%%%%%%%%%%%%%%%%%%%%%%%%%%%%%%%%%%%%%%%%%%%%%%%%%%%%%%%%%%%%%%%%%%%%%%%%%%%%%%%%%%%
\theorem
Let $\function{\Curve{}}{\interval{}}{\Man{}}$ be a smooth curve on $\Man{}$
for some connected open interval $\interval{}$ of $\R$.
\begin{itemize}
\item
$\valongcurves{\Man{}}{\Curve{}}$ together with its addition and scalar-multiplication defined in
\refdef{defvectorfieldsalongcurves} forms an $\R$-vector-space, which is referred to as the $\quotl$canonical
linear structure of $\valongcurves{\Man{}}{\Curve{}}$$\quotr$.
\item
$\valongcurves{\Man{}}{\Curve{}}$ together with its addition and $\smoothmaps{\interval{}}$-multiplication
defined in \refdef{defvectorfieldsalongcurves}, forms a module over the commutative ring $\smoothmaps{\interval{}}$,
which is referred to as the $\quotl$canonical module structure of $\valongcurves{\Man{}}{\Curve{}}$$\quotr$.
\end{itemize}
By an abuse of notation, both the canonical linear and module structures of the set of vector fields on $\Man{}$ along $\Curve{}$
will be denoted simply by the same notation $\valongcurves{\Man{}}{\Curve{}}$, for convenience.
\proof
It is trivial.
\endthm
%%%%%%%%%%%%%%%%%%%%%%%%%%%%%%%%%%%%%%%%%%%%%%%%%%%%%%%%%%%%%%%%%%%%%%%%%%%%%%%%%%%%%%%%%%%%%%%%%%%%%%%%%%%%%%%%%%%%%%%%%%%%%%%%
\definition
Let $\function{\Curve{}}{\interval{}}{\Man{}}$ be a smooth curve on $\Man{}$
for some connected open interval $\interval{}$ of $\R$.
An element $\valongc{}$ is called an $\quotl$extendible vector field along $\Curve{}$$\quotr$ if
there exists $\avecf{}\in\vectorfields{\Man{}}$ such that,
\begin{equation}
\valongc{}=\cmp{\avecf{}}{\Curve{}}.
%\Foreach{t}{\interval{}}
%\func{\avecf{}}{\func{\Curve{}}{t}}=\func{\valongc{}}{t}.
\end{equation}
Such a vector field $\avecf{}$ is called an $\quotl$extension of $\valongc{}$$\quotr$.
\endef
%%%%%%%%%%%%%%%%%%%%%%%%%%%%%%%%%%%%%%%%%%%%%%%%%%%%%%%%%%%%%%%%%%%%%%%%%%%%%%%%%%%%%%%%%%%%%%%%%%%%%%%%%%%%%%%%%%%%%%%%%%%%%%%%
\theorem\label{thmcovarianderivativealongcurves}
Let $\function{\Curve{}}{\interval{}}{\Man{}}$ be a smooth curve on $\Man{}$
for some connected open interval $\interval{}$ of $\R$.
There exists a unique map
$\function{\covder{\connection{}}{\Curve{}}}{\valongcurves{\Man{}}{\Curve{}}}{\valongcurves{\Man{}}{\Curve{}}}$
such that the following conditions are satisfied.
\begin{itemize}
\item[\myitem{CD~1.}]
$\covder{\connection{}}{\Curve{}}$ is $\R$-linear, that is,
\begin{align}
&\Foreach{\opair{x}{y}}{\R^2}
\Foreach{\opair{\valongc{1}}{\valongc{2}}}
{\Cprod{\valongcurves{\Man{}}{\Curve{}}}
{\valongcurves{\Man{}}{\Curve{}}}}\cr
&~\func{\covder{\connection{}}{\Curve{}}}{x\valongc{1}+y\valongc{2}}=
x\func{\covder{\connection{}}{\Curve{}}}{\valongc{1}}+
y\func{\covder{\connection{}}{\Curve{}}}{\valongc{2}}.
\end{align}
\item[\myitem{CD~2.}]
\begin{equation}
\Foreach{f}{\smoothmaps{\interval{}}}
\Foreach{\valongc{}}{\valongcurves{\Man{}}{\Curve{}}}
\func{\covder{\connection{}}{\Curve{}}}{f\valongc{}}=\dot{f}\valongc{}+
f\func{\covder{\connection{}}{\Curve{}}}{\valongc{}},
\end{equation}
where $\dot{f}$ denotes the derived function of $f$.
\item[\myitem{CD~3.}]
If $\valongc{}$ is extendible with the extension $\avecf{}\in\vectorfields{\Man{}}$, then,
\begin{equation}
\Foreach{t}{\interval{}}
\func{\[\func{\covder{\connection{}}{\Curve{}}}{\valongc{}}\]}{t}=
\tbcon{\func{\dot{\Curve{}}}{t}}{\avecf{}},
\end{equation}
where $\func{\dot{\Curve{}}}{t}$ denotes the value of the derivative of $\Curve{}$ at $t$,
which is an element of $\tanspace{\func{\Curve{}}{t}}{\Man{}}$.
\end{itemize}
\proof
We break down the proof into parts.
\begin{itemize}
\item[\myitem{pr-1.}]
We first show that, given a smooth curve $\function{\Curve{}}{\interval{}}{\Man{}}$,
if there exists such an operator $\covder{\connection{}}{\Curve{}}$ satisfying the conditions above,
then for any sub-interval $O$ of $\interval{}$,
there exists an operator
$\function{\covder{}{O}}{\valongcurves{\Man{}}{\reS{\Curve{}}{O}}}
{\valongcurves{\Man{}}{\reS{\Curve{}}{O}}}$ that coincides with $\covder{\connection{}}{\Curve{}}$
in the region $O$, and for which the same properties are hold.\\
Suppose that such an operator $\covder{\connection{}}{\Curve{}}$ exists.
Given a sub-interval $\opair{t_0-\delta}{t_0+\delta}$ of $\interval{}$,
and bearing in mind the arguments related to cut-off functions, \myitem{CD~2} implies in
a straightforward manner that $\func{\[\func{\covder{\connection{}}{\Curve{}}}{\valongc{1}}\]}{t_0}=
\func{\[\func{\covder{\connection{}}{\Curve{}}}{\valongc{2}}\]}{t_0}$ for any $\valongc{1}$ and $\valongc{2}$
in $\valongcurves{\Man{}}{\Curve{}}$ that coincide in $\opair{t_0-\delta}{t_0+\delta}$.
So given any connected open subset $O$ of $\interval{}$, we can define an operator
$\function{\covder{}{O}}{\valongcurves{\Man{}}{\reS{\Curve{}}{O}}}
{\valongcurves{\Man{}}{\reS{\Curve{}}{O}}}$ in the way that for any
$\valongc{O}\in\valongcurves{\Man{}}{\reS{\Curve{}}{O}}$ and for any $t\in O$,
$\func{\[\func{\covder{}{O}}{\valongc{O}}\]}{t}\eqdef\func{\[\func{\covder{\connection{}}{\Curve{}}}{\valongc{}}\]}{t}$
for any $\valongc{}\in\valongcurves{\Man{}}{\Curve{}}$ that coincides with $\valongc{O}$ in a neighborhood of $t$.
Furthermore, it can be easily checked that $\covder{}{O}$ satisfies the conditions
\myitem{CD~1}, \myitem{CD~2}, \myitem{CD~3}. Also, clearly for every $\valongc{}\in\valongcurves{\Man{}}{\Curve{}}$,
$\func{\covder{}{O}}{\reS{\valongc{}}{O}}=\reS{\[\func{\covder{\connection{}}{\Curve{}}}{\valongc{}}\]}{O}$.
%%%%%%%%%%%%%%%%%%%%%%%%
\item[\myitem{pr-2.}]
Second, we show that given a smooth curve $\function{\Curve{}}{\interval{}}{\Man{}}$,
if there exists such an operator $\covder{\connection{}}{\Curve{}}$ satisfying the conditions
\myitem{CD~1}, \myitem{CD~2}, \myitem{CD~3}, it is unique.\\
Suppose that such an operator $\covder{\connection{}}{\Curve{}}$ exists.
Let $t$ be an arbitrary element of $\interval{}$, and find a chart $\opair{\U}{\phi}$ of $\Man{}$ such that
$\func{\Curve{}}{t}\in\U$.
%$\defSet{\func{\Curve{}}{t}}{t\in\cpair{-\varepsilon}{\varepsilon}}\subseteq\U$.
Let $\mtuple{\avecff{1}}{\avecff{m}}$ denote the local frame of vector fields
associated with the chart $\opair{\U}{\phi}$. Furthermore, find an open set $\V\subseteq\U$ containing $\func{\Curve{}}{t}$
and an $m$-tuple $\mtuple{\avecf{1}}{\avecf{m}}$ of global vector fields such that $\reS{\avecf{i}}{\V}=\reS{\avecff{i}}{\V}$
for every $i\in\seta{\suc{1}{m}}$. Clearly, there exists a positive integer $\varepsilon$ such that
$\defSet{\func{\Curve{}}{t}}{t\in\cpair{t-\varepsilon}{t+\varepsilon}}\subseteq\V$. Let $O:=\opair{t-\varepsilon}{t+\varepsilon}$.
It is obvious that $\cmp{\avecf{i}}{\reS{\Curve{}}{O}}\in\valongcurves{\Man{}}{\reS{\Curve{}}{O}}$.
Let $\valongc{}\in\valongcurves{\Man{}}{\Curve{}}$. Clearly,
$\displaystyle\reS{\valongc{}}{O}=\sum_{i=1}^{m}f_{i}\(\cmp{\avecf{i}}{\reS{\Curve{}}{O}}\)=
\sum_{i=1}^{m}f_{i}\(\cmp{\avecff{i}}{\reS{\Curve{}}{O}}\)$ for unique $m$-tuple
$\mtuple{f_1}{f_m}$ of real-valued smooth maps on $O$. Therefore, considering that $\covder{}{O}$,
as defined in \myitem{pr-3.},
satisfies \myitem{CD~1}, \myitem{CD~2}, \myitem{CD~3}, and each element $\cmp{\avecf{i}}{\reS{\Curve{}}{O}}$ of
$\valongcurves{\Man{}}{\reS{\Curve{}}{O}}$ is extendible with the extension $\avecf{i}$,
\begin{align}
\func{\[\func{\covder{\connection{}}{\Curve{}}}{\valongc{}}\]}{t}&=
\func{\[\func{\covder{}{O}}{\reS{\valongc{}}{O}}\]}{t}\cr
&=\func{\[\func{\covder{}{O}}{\sum_{i=1}^{m}f_{i}\(\cmp{\avecf{i}}{\reS{\Curve{}}{O}}\)}\]}{t}\cr
&=\sum_{i=1}^{m}\(\func{\dot{f}_{i}}{t}\func{\avecf{i}}{t}+
\func{f_i}{t}\func{\[\func{\covder{}{O}}{\cmp{\avecf{i}}{\reS{\Curve{}}{O}}}\]}{t}\)\cr
&=\sum_{i=1}^{m}\(\func{\dot{f}_{i}}{t}\func{\avecf{i}}{t}+
\func{f_i}{t}\[\tbcon{\func{\dot{\Curve{}}}{t}}{\avecf{i}}\]\)\cr
&=\sum_{i=1}^{m}\(\func{\dot{f}_{i}}{t}\func{\avecff{i}}{t}+
\func{f_i}{t}\[\tbcon{\func{\dot{\Curve{}}}{t}}{\avecff{i}}\]\).
\end{align}
Since the choice of $t$ in $\interval{}$ was arbitrary, this argument reveals that
given a smooth curve $\function{\Curve{}}{\interval{}}{\Man{}}$,
there exists at most one operator $\covder{\connection{}}{\Curve{}}$ satisfying the conditions
\myitem{CD~1}, \myitem{CD~2}, \myitem{CD~3}.
%%%%%%%%%%%%%%%%%%%%%%%%
\item[\myitem{pr-3.}]
Let $\opair{\U}{\phi}$ be a chart of $\Man{}$ that intersects the image of a given curve
$\function{\Curve{}}{\interval{}}{\Man{}}$.
Let $\mtuple{\avecff{1}}{\avecff{m}}$ denote the local frame of vector fields
associated with the chart $\opair{\U}{\phi}$. Let $t_0$ be an element of $\interval{}$ and
$\varepsilon$ be a positive integer such that
$\defSet{\func{\Curve{}}{t}}{t\in\cpair{t_0-\varepsilon}{t_0+\varepsilon}}\subseteq\U$.
Let $O:=\opair{t_0-\varepsilon}{t_0+\varepsilon}$. Here, we define the operator
$\function{\covder{}{O}}{\valongcurves{\Man{}}{\reS{\Curve{}}{O}}}
{\valongcurves{\Man{}}{\reS{\Curve{}}{O}}}$ in the following way. Let
$\valongc{}$ be an element of $\valongcurves{\Man{}}{\reS{\Curve{}}{O}}$, and let
$\displaystyle\valongc{}=\sum_{i=1}^{m}f_{i}\(\cmp{\avecff{i}}{\reS{\Curve{}}{O}}\)$
for unique $m$-tuple $\mtuple{f_1}{f_m}$ of real-valued smooth maps on $O$.
\begin{equation}
\Foreach{t}{O}
\func{\[\func{\covder{}{O}}{\valongc{}}\]}{t}\eqdef
\sum_{i=1}^{m}\(\func{\dot{f}_{i}}{t}\func{\avecff{i}}{t}+
\func{f_i}{t}\[\tbcon{\func{\dot{\Curve{}}}{t}}{\avecff{i}}\]\).
\end{equation}
It can be easily checked that this operator satisfies the conditions
\myitem{CD~1}, \myitem{CD~2}, \myitem{CD~3} for the curve $\reS{\Curve{}}{O}$.
%%%%%%%%%%%%%%%%%%%%%%%%
\item[\myitem{pr-4.}]
Let $\function{\Curve{}}{\interval{}}{\Man{}}$ be a curve on $\Man{}$.
Now, let $\opair{\U_1}{\phi_1}$ and $\opair{\U_2}{\phi_2}$ be charts of $\Man{}$
that intersect the image of the curve $\Curve{}$, and let $t_1$ and $t_2$ be elements of $\interval{}$ and
$\varepsilon_1$ and $\varepsilon_2$ be positive integers such that
$\defSet{\func{\Curve{}}{t}}{t\in\cpair{t_1-\varepsilon_1}{t_1+\varepsilon_1}}\subseteq\U_1$ and
$\defSet{\func{\Curve{}}{t}}{t\in\cpair{t_2-\varepsilon_2}{t_2+\varepsilon_2}}\subseteq\U_2$. Furthermore,
suppose that $O_1:=\opair{t_1-\varepsilon_1}{t_1+\varepsilon_1}$ and
$O_2:=\opair{t_2-\varepsilon_2}{t_2+\varepsilon_2}$ intersect, and denote this intersection by $O$.
Also, let $\mtuple{\avecf{1}}{\avecf{m}}$ and $\mtuple{\avecff{1}}{\avecff{m}}$ denote the local frame of vector fields
associated with the charts $\opair{\U_1}{\phi_1}$ and $\opair{\U_1}{\phi_1}$, respectively.
Since the operators $\covder{}{O_1}$ and $\covder{}{O_2}$, as defined in \myitem{pr-3}, satisfy
the conditions \myitem{CD~1}, \myitem{CD~2}, \myitem{CD~3} for the
curves $\reS{\Curve{}}{O_1}$ and $\reS{\Curve{}}{O_2}$, respectively, and since $O$ is simultaneously a
sub-interval of $O_1$ and $O_2$, the result of \myitem{pr-1} implies that $\covder{}{O_1}$ induces an operator
$\function{\covder{1}{O}}{\valongcurves{\Man{}}{\reS{\Curve{}}{O}}}{\valongcurves{\Man{}}{\reS{\Curve{}}{O}}}$
and $\covder{}{O_2}$ induces an operator
$\function{\covder{2}{O}}{\valongcurves{\Man{}}{\reS{\Curve{}}{O}}}{\valongcurves{\Man{}}{\reS{\Curve{}}{O}}}$
such that $\covder{i}{O}$ coincides with $\covder{}{O_i}$ in the region $O$, and satisfies the conditions
\myitem{CD~1}, \myitem{CD~2}, \myitem{CD~3} for the curve $\reS{\Curve{}}{O}$,
for $i\in\seta{\binary{1}{2}}$. But the result of \myitem{pr-3} implies that $\covder{1}{O}=\covder{2}{O}$.
Therefore, it becomes clear that, given $\valongc{}\in\valongcurves{\Man{}}{\Curve{}}$ with
$\displaystyle\reS{\valongc{}}{O_1}=\sum_{i=1}^{m}f_{i}\(\cmp{\avecf{i}}{\reS{\Curve{}}{O_1}}\)$ and
$\displaystyle\reS{\valongc{}}{O_2}=\sum_{i=1}^{m}g_{i}\(\cmp{\avecff{i}}{\reS{\Curve{}}{O_2}}\)$,
\begin{align}
\Foreach{t}{O}
\func{\[\func{\covder{}{O_1}}{\reS{\valongc{}}{O_1}}\]}{t}&=
\sum_{i=1}^{m}\(\func{\dot{f}_{i}}{t}\func{\avecf{i}}{t}+
\func{f_i}{t}\[\tbcon{\func{\dot{\Curve{}}}{t}}{\avecf{i}}\]\)\cr
&=\sum_{i=1}^{m}\(\func{\dot{g}_{i}}{t}\func{\avecff{i}}{t}+
\func{g_i}{t}\[\tbcon{\func{\dot{\Curve{}}}{t}}{\avecff{i}}\]\)\cr
&=\func{\[\func{\covder{}{O_2}}{\reS{\valongc{}}{O_2}}\]}{t}.
\end{align}
%%%%%%%%%%%%%%%%%%%%%%%%
\item[\myitem{pr-5.}]
Given a curve $\function{\Curve{}}{\interval{}}{\Man{}}$ on $\Man{}$,
the results of \myitem{pr-3} and \myitem{pr-4} clearly imply the existance of an operator
$\function{\covder{\connection{}}{\Curve{}}}{\valongcurves{\Man{}}{\Curve{}}}{\valongcurves{\Man{}}{\Curve{}}}$
satisfying the conditions \myitem{CD~1}, \myitem{CD~2}, and \myitem{CD~3}. This can be done as the following.\\
Given $\valongc{}\in\valongcurves{\Man{}}{\Curve{}}$, and any $t\in\interval{}$,
choose a chart $\opair{\U}{\phi}$ of $\Man{}$ containing $\func{\Curve{}}{t}$ and choose a positive integer
$\varepsilon$ such that
$\defSet{\func{\Curve{}}{s}}{s\in\cpair{t-\varepsilon}{t+\varepsilon}}\subseteq\U$, and for every $s$ in
$\opair{t-\varepsilon}{t+\varepsilon}$,
$\displaystyle\func{\valongc{}}{s}=\sum_{i=1}^{m}\func{f_{i}}{s}\func{\avecf{i}}{\func{\Curve{}}{s}}$, where
$\mtuple{\avecf{1}}{\avecf{m}}$ denotes the local frame of vector fields
associated with the chart $\opair{\U}{\phi}$. According to \myitem{pr-4}, it is permissible to define,
\begin{equation}
\func{\[\func{\covder{\connection{}}{\Curve{}}}{\valongc{}}\]}{t}\eqdef
\sum_{i=1}^{m}\(\func{\dot{f}_{i}}{t}\func{\avecf{i}}{t}+
\func{f_i}{t}\[\tbcon{\func{\dot{\Curve{}}}{t}}{\avecf{i}}\]\).
\end{equation}
According to \myitem{pr-3}, this definition satisfies the properties \myitem{CD~1}, \myitem{CD~2}, and \myitem{CD~3}.
\end{itemize}
\endthm
%%%%%%%%%%%%%%%%%%%%%%%%%%%%%%%%%%%%%%%%%%%%%%%%%%%%%%%%%%%%%%%%%%%%%%%%%%%%%%%%%%%%%%%%%%%%%%%%%%%%%%%%%%%%%%%%%%%%%%%%%%%%%%%%
\definition\label{defcovarianderivativealongcurves}
Let $\function{\Curve{}}{\interval{}}{\Man{}}$ be a smooth curve on $\Man{}$
for some connected open interval $\interval{}$ of $\R$.
The unique operator from $\valongcurves{\Man{}}{\Curve{}}$ to $\valongcurves{\Man{}}{\Curve{}}$
satisfying the conditions of \refthm{thmcovarianderivativealongcurves}, will be denoted by
$\covder{\connection{}}{\Curve{}}$ and referred to as the $\quotl$covariant derivative operator of the vector fields on $\Man{}$
along the curve $\Curve{}$, relative to the connection $\connection{}$$\quotr$.
Given $\valongc{}\in\valongcurves{\Man{}}{\Curve{}}$,
$\func{\covder{\connection{}}{\Curve{}}}{\valongc{}}$ will be called the $\quotl$covariant derivative of $\valongc{}$
along the curve $\Curve{}$ with respect to the connection $\connection{}$$\quotr$.
When the underlying connection is clear enough, $\covder{\connection{}}{\Curve{}}$ can simply be denoted by
$\covder{}{\Curve{}}$, and if in addition there is no ambiguity about the curve under consideration, it can more simply
be denoted by $\covder{}{}$.
\endef
%%%%%%%%%%%%%%%%%%%%%%%%%%%%%%%%%%%%%%%%%%%%%%%%%%%%%%%%%%%%%%%%%%%%%%%%%%%%%%%%%%%%%%%%%%%%%%%%%%%%%%%%%%%%%%%%%%%%%%%%%%%%%%%%
\corollary\label{corcovarianderivativealongcurves}
Let $\function{\Curve{}}{\interval{}}{\Man{}}$ be a smooth curve on $\Man{}$
for some connected open interval $\interval{}$ of $\R$.
Let $\valongc{}\in\valongcurves{\Man{}}{\Curve{}}$, and
let $\opair{\U}{\phi}$ be a chart of $\Man{}$ that intersects the image of $\Curve{}$. Let $O$ be a sub-interval of
$\interval{}$ such that
$\defSet{\func{\Curve{}}{t}}{t\in O}\subseteq\U$. Let $\mtuple{f_1}{f_m}$ be the unique $m$-tuple of smooth functions on $O$
such that for every $t\in O$,
$\displaystyle\func{\valongc{}}{t}=\sum_{i=1}^{m}\func{f_{i}}{t}\func{\avecf{i}}{\func{\Curve{}}{t}}$, where
$\mtuple{\avecf{1}}{\avecf{m}}$ denotes the local frame of vector fields
associated with the chart $\opair{\U}{\phi}$.
\begin{equation}
\Foreach{t}{O}
\func{\[\func{\covder{\connection{}}{\Curve{}}}{\valongc{}}\]}{t}=
\sum_{i=1}^{m}\(\func{\dot{f}_{i}}{t}\func{\avecf{i}}{t}+
\func{f_i}{t}\[\tbcon{\func{\dot{\Curve{}}}{t}}{\avecf{i}}\]\).
\end{equation}
\endcor
%%%%%%%%%%%%%%%%%%%%%%%%%%%%%%%%%%%%%%%%%%%%%%%%%%%%%%%%%%%%%%%%%%%%%%%%%%%%%%%%%%%%%%%%%%%%%%%%%%%%%%%%%%%%%%%%%%%%%%%%%%%%%%%%
\theorem\label{thmlocalexpressionofcovariantderivativeofvectorfieldsalongcurves}
Let $\function{\Curve{}}{\interval{}}{\Man{}}$ be a smooth curve on $\Man{}$
for some connected open interval $\interval{}$ of $\R$.
Let $\valongc{}\in\valongcurves{\Man{}}{\Curve{}}$, and
let $\opair{\U}{\phi}$ be a chart of $\Man{}$ that intersects the image of $\Curve{}$. Let $O$ be a sub-interval of
$\interval{}$ such that
$\defSet{\func{\Curve{}}{t}}{t\in O}\subseteq\U$. Let $\mtuple{f_1}{f_m}$ be the unique $m$-tuple of smooth functions on $O$
such that for every $t\in O$,
$\displaystyle\func{\valongc{}}{t}=\sum_{i=1}^{m}\func{f_{i}}{t}\func{\localframevecf{i}}{\func{\Curve{}}{t}}$, where
$\mtuple{\localframevecf{1}}{\localframevecf{m}}$ denotes the local frame of vector fields
associated with the chart $\opair{\U}{\phi}$. Let $\Curve{i}:=\cmp{\dualEucbase{m}{i}}{\cmp{\phi}{\Curve{}}}$ for every $i$ in
$\seta{\suc{1}{m}}$, where $\mtuple{\dualEucbase{m}{1}}{\dualEucbase{m}{m}}$ denotes the dual of the
standard orthonormal ordered-basis of $\R^m$.
\begin{equation}
\Foreach{t}{O}
\func{\[\func{\covder{\connection{}}{\Curve{}}}{\valongc{}}\]}{t}=
\sum_{i=1}^{m}\[\func{\dot{f}_{i}}{t}+\sum_{j=1}^{m}\sum_{k=1}^{m}\func{{\dot{\Curve{}}}_{j}}{t}
\func{\Christoffel{\phi}{i}{j}{k}}{\func{\Curve{}}{t}}\func{f_k}{t}\]\func{\localframevecf{i}}{\func{\Curve{}}{t}},
\end{equation}
or equivalently,
\begin{equation}
\reS{\[\func{\covder{\connection{}}{\Curve{}}}{\valongc{}}\]}{O}=
\reS{\(\sum_{i=1}^{m}\[\dot{f}_{i}+\sum_{j=1}^{m}\sum_{k=1}^{m}{\dot{\Curve{}}}_{j}
\(\cmp{\Christoffel{\phi}{i}{j}{k}}{\Curve{}}\)f_k\]\(\cmp{\localframevecf{i}}{\Curve{}}\)\)}{O}.
\end{equation}
\proof
It is an immediate consequence of \refdef{defcovariantderivativeinducedbyconnection}, \refcor{corlocalexpressionofconnection},
and \refcor{corcovarianderivativealongcurves}
\endthm
%%%%%%%%%%%%%%%%%%%%%%%%%%%%%%%%%%%%%%%%%%%%%%%%%%%%%%%%%%%%%%%%%%%%%%%%%%%%%%%%%%%%%%%%%%%%%%%%%%%%%%%%%%%%%%%%%%%%%%%%%%%%%%%%
\theorem
Let $\point$ be a point of the manifold $\Man{}$, and let $v\in\tanspace{\point}{\Man{}}$.
Let $\avecf{}$ and $\avecff{}$ be elements of $\vectorfields{\Man{}}$. If there exists a curve
$\function{\Curve{}}{\interval{}}{\Man{}}$ for some open interval $\interval{}$ of $\R$, such that for some
$s\in\interval{}$, $\func{\Curve{}}{s}=\point$ and $\func{\dot{\Curve{}}}{s}=v$, and further
$\cmp{\avecf{}}{\Curve{}}=\cmp{\avecff{}}{\Curve{}}$, then,
\begin{equation}
\tbcon{v}{\avecf{}}=\tbcon{v}{\avecff{}}.
\end{equation}
\proof
Clearly $\valongc{}:=\cmp{\avecf{}}{\Curve{}}=\cmp{\avecff{}}{\Curve{}}$ is an element of $\valongcurves{\Man{}}{\Curve{}}$,
that is a smooth vector field on $\Man{}$ along $\Curve{}$. It is also trivial that both $\avecf{}$ and $\avecff{}$
are extensions of $\valongc{}$. Thus,
\begin{align}
\func{\[\func{\covder{\connection{}}{\Curve{}}}{\valongc{}}\]}{s}=
\tbcon{\func{\dot{\Curve{}}}{s}}{\avecf{}}=
\tbcon{v}{\avecf{}},
\end{align}
and,
\begin{align}
\func{\[\func{\covder{\connection{}}{\Curve{}}}{\valongc{}}\]}{s}=
\tbcon{\func{\dot{\Curve{}}}{s}}{\avecff{}}=
\tbcon{v}{\avecff{}}.
\end{align}
Therefore, $\tbcon{v}{\avecf{}}=\tbcon{v}{\avecff{}}$.
\endthm
%%%%%%%%%%%%%%%%%%%%%%%%%%%%%%%%%%%%%%%%%%%%%%%%%%%%%%%%%%%%%%%%%%%%%%%%%%%%%%%%%%%%%%%%%%%%%%%%%%%%%%%%%%%%%%%%%%%%%%%%%%%%%%%%
\remark
It has been revealed that given a pair of smooth vector fields $\avecf{}$ and $\avecff{}$ on $\Man{}$,
and a point $\point$ of $\Man{}$, the value of the vector field $\con{\avecf{}}{\avecff{}}$ at point $\point$,
that is $\func{\[\con{\avecf{}}{\avecff{}}\]}{\point}$, depends only on the value of $\avecf{}$ at $\point$
and on the values of $\avecff{}$ at the points in the image of some curve passing through the point $\point$ whose velocity
at $\point$ coincides with $\func{\avecf{}}{\point}$.
\endremark
%%%%%%%%%%%%%%%%%%%%%%%%%%%%%%%%%%%%%%%%%%%%%%%%%%%%%%%%%%%%%%%%%%%%%%%%%%%%%%%%%%%%%%%%%%%%%%%%%%%%%%%%%%%%%%%%%%%%%%%%%%%%%%%%
\definition\label{defparallelvectorfieldsalongcurves}
Let $\function{\Curve{}}{\interval{}}{\Man{}}$ be a curve on $\Man{}$ for some open interval $\interval{}$ of $\R$.
Let $\valongc{}\in\valongcurves{\Man{}}{\Curve{}}$. $\valongc{}$ is called a $\quotl$parallel (smooth) vector field
on $\Man{}$ along the curve $\Curve{}$ relative to the connection $\connection{}$$\quotr$,
if $\func{\covder{\connection{}}{\Curve{}}}{\valongc{}}=0$, that is,
\begin{equation}
\Foreach{t}{\interval{}}
\func{\[\func{\covder{\connection{}}{\Curve{}}}{\valongc{}}\]}{t}=\zerovec{}.
\end{equation}
The set of all parallel vector fields on $\Man{}$ along $\Curve{}$ with respect to the connection $\connection{}$
will be denoted by $\parallelvalongcurves{\connection{}}{\Curve{}}$.
\endef
%%%%%%%%%%%%%%%%%%%%%%%%%%%%%%%%%%%%%%%%%%%%%%%%%%%%%%%%%%%%%%%%%%%%%%%%%%%%%%%%%%%%%%%%%%%%%%%%%%%%%%%%%%%%%%%%%%%%%%%%%%%%%%%%
\corollary\label{corparallelvfalongcurveexpresion}
Let $\function{\Curve{}}{\interval{}}{\Man{}}$ be a curve on $\Man{}$ for some open interval $\interval{}$ of $\R$.
Let $\valongc{}\in\valongcurves{\Man{}}{\Curve{}}$. The following statements are equivalent.
\begin{itemize}
\item[\myitem{1.}]
$\valongc{}$ is a parallel vector field on $\Man{}$ along the curve $\Curve{}$, that is
$\valongc{}\in\parallelvalongcurves{\connection{}}{\Curve{}}$.
\item[\myitem{2.}]
For every chart
$\opair{\U}{\phi}$ of $\Man{}$ that intersects the image of $\Curve{}$, and every sub-interval $O$ of
$\interval{}$ such that $\defSet{\func{\Curve{}}{t}}{t\in O}\subseteq\U$,
\begin{align}
&\Foreach{i}{\seta{\suc{1}{m}}}\Foreach{t}{O}\cr
&\func{\dot{f}_{i}}{t}=
-\sum_{k=1}^{m}\sum_{j=1}^{m}\func{{\dot{\Curve{}}}_{j}}{t}
\func{\Christoffel{\phi}{i}{j}{k}}{\func{\Curve{}}{t}}\func{f_k}{t},
\end{align}
or equivalently,
\begin{equation}
\dot{f_i}=-\sum_{k=1}^{m}\sum_{j=1}^{m}
\reS{{\dot{\Curve{}}}_j}{O}
\(\cmp{\Christoffel{\phi}{i}{j}{k}}{\reS{\Curve{}}{O}}\)f_k,
\end{equation}
where, $\mtuple{f_1}{f_m}$ denotes the unique $m$-tuple of smooth functions on $O$
such that for every $t\in O$,
$\displaystyle\func{\valongc{}}{t}=\sum_{i=1}^{m}\func{f_{i}}{t}\func{\localframevecf{i}}{\func{\Curve{}}{t}}$,
$\mtuple{\localframevecf{1}}{\localframevecf{m}}$ denoting the local frame of vector fields
associated with the chart $\opair{\U}{\phi}$, and $\Curve{i}:=\cmp{\dualEucbase{m}{i}}{\cmp{\phi}{\Curve{}}}$ for every $i$ in
$\seta{\suc{1}{m}}$.
\item[\myitem{3.}]
Given a family $\seta{\opair{\U_\alpha}{\phi_\alpha}}_{\alpha\in\indexset{}}$ of charts of $\Man{}$ such that
$\funcimage{\Curve{}}\subseteq\Union{\alpha}{\indexset{}}{\U_\alpha}$ and
$\U_\alpha\cap\funcimage{\Curve{}}\neq\emptyset$, and given a family of
open sub-intervals $\seta{O_\alpha}_{\alpha\in\indexset{}}$ of the interval $\interval{}$ such that
$\defSet{\func{\Curve{}}{t}}{t\in O_\alpha}\subseteq\U_\alpha$, for every $\alpha\in\indexset{}$,
%%%%%%
\begin{align}
&\Foreach{i}{\seta{\suc{1}{m}}}\Foreach{t}{O_\alpha}\cr
&\func{\dot{f}_{i}^{\alpha}}{t}=
-\sum_{k=1}^{m}\sum_{j=1}^{m}\func{{\dot{\Curve{}}}_{j}}{t}
\func{\Christoffel{\phi_\alpha}{i}{j}{k}}{\func{\Curve{}}{t}}\func{f_k^\alpha}{t},
\end{align}
or equivalently,
\begin{equation}
\dot{f_i^\alpha}=-\sum_{k=1}^{m}\sum_{j=1}^{m}
\reS{{\dot{\Curve{}}}_j}{O_\alpha}
\(\cmp{\Christoffel{\phi_\alpha}{i}{j}{k}}{\reS{\Curve{}}{O}}\)f_k^\alpha,
\end{equation}
where, $\mtuple{f_1^\alpha}{f_m^\alpha}$ denotes the unique $m$-tuple of smooth functions on $O_\alpha$
such that for every $t\in O_\alpha$,
$\displaystyle\func{\valongc{}}{t}=\sum_{i=1}^{m}\func{f_{i}^\alpha}{t}\func{\Localframevecf{i}{\alpha}}{\func{\Curve{}}{t}}$,
$\mtuple{\Localframevecf{1}{\alpha}}{\Localframevecf{m}{\alpha}}$ denoting the local frame of vector fields
associated with the chart $\opair{\U_\alpha}{\phi_\alpha}$,
and $\Curve{i}:=\cmp{\dualEucbase{m}{i}}{\cmp{\phi}{\Curve{}}}$ for every $i$ in
$\seta{\suc{1}{m}}$.
\end{itemize}
\endcor
%%%%%%%%%%%%%%%%%%%%%%%%%%%%%%%%%%%%%%%%%%%%%%%%%%%%%%%%%%%%%%%%%%%%%%%%%%%%%%%%%%%%%%%%%%%%%%%%%%%%%%%%%%%%%%%%%%%%%%%%%%%%%%%%
\definition
A smooth vector field $\avecf{}$ on $\Man{}$ is said to be a $\quotl$parallel vector field on $\Man{}$ relative
to the connection $\connection{}$$\quotr$ if for every smooth curve $\Curve{}$ on $\Man{}$,
$\cmp{\avecf{}}{\Curve{}}$ is a parallel vector field along $\Curve{}$ relative to $\connection{}$.
\endef
%%%%%%%%%%%%%%%%%%%%%%%%%%%%%%%%%%%%%%%%%%%%%%%%%%%%%%%%%%%%%%%%%%%%%%%%%%%%%%%%%%%%%%%%%%%%%%%%%%%%%%%%%%%%%%%%%%%%%%%%%%%%%%%%
\corollary
Let $\avecf{}\in\vectorfields{\Man{}}$. $\avecf{}$ is a parallel vector field along $\Curve{}$ relative to $\connection{}$
if and only if
\begin{equation}
\Foreach{\avecff{}}{\vectorfields{\Man{}}}
\con{\avecff{}}{\avecf{}}=0.
\end{equation}
\endcor
%%%%%%%%%%%%%%%%%%%%%%%%%%%%%%%%%%%%%%%%%%%%%%%%%%%%%%%%%%%%%%%%%%%%%%%%%%%%%%%%%%%%%%%%%%%%%%%%%%%%%%%%%%%%%%%%%%%%%%%%%%%%%%%%
%%%%%%%%%%%%%%%%%%%%%%%%%%%%%%%%%%%%%%%%%%%%%%%%%%%%%%%%%%%%%%%%%%%%%%%%%%%%%%%%%%%%%%%%%%%%%%%%%%%%%%%%%%%%%%%%%%%%%%%%%%%%%%%%
%%%%%%%%%%%%%%%%%%%%%%%%%%%%%%%%%%%%%%%%%%%%%%%%%%%%%%%%%%%%%%%%%%%%%%%%%%%%%%%%%%%%%%%%%%%%%%%%%%%%%%%%%%%%%%%%%%%%%%%%%%%%%%%%
%%%%%%%%%%%%%%%%%%%%%%%%%%%%%%%%%%%%%%%%%%%%%%%%%%%%%%%%%%%%%%%%%%%%%%%%%%%%%%%%%%%%%%%%%%%%%%%%%%%%%%%%%%%%%%%%%%%%%%%%%%%%%%%%
%%%%%%%%%%%%%%%%%%%%%%%%%%%%%%%%%%%%%%%%%%%%%%%%%%%%%%%%%%%%%%%%%%%%%%%%%%%%%%%%%%%%%%%%%%%%%%%%%%%%%%%%%%%%%%%%%%%%%%%%%%%%%%%%
%%%%%%%%%%%%%%%%%%%%%%%%%%%%%%%%%%%%%%%%%%%%%%%%%%%%%%%%%%%%%%%%%%%%%%%%%%%%%%%%%%%%%%%%%%%%%%%%%%%%%%%%%%%%%%%%%%%%%%%%%%%%%%%%
%%%%%%%%%%%%%%%%%%%%%%%%%%%%%%%%%%%%%%%%%%%%%%%%%%%%%%%%%%%%%%%%%%%%%%%%%%%%%%%%%%%%%%%%%%%%%%%%%%%%%%%%%%%%%%%%%%%%%%%%%%%%%%%%
%%%%%%%%%%%%%%%%%%%%%%%%%%%%%%%%%%%%%%%%%%%%%%%%%%%%%%%%%%%%%%%%%%%%%%%%%%%%%%%%%%%%%%%%%%%%%%%%%%%%%%%%%%%%%%%%%%%%%%%%%%%%%%%%
%%%%%%%%%%%%%%%%%%%%%%%%%%%%%%%%%%%%%%%%%%%%%%%%%%%%%%%%%%%%%%%%%%%%%%%%%%%%%%%%%%%%%%%%%%%%%%%%%%%%%%%%%%%%%%%%%%%%%%%%%%%%%%%%
%%%%%%%%%%%%%%%%%%%%%%%%%%%%%%%%%%%%%%%%%%%%%%%%%%%%%%%%%%%%%%%%%%%%%%%%%%%%%%%%%%%%%%%%%%%%%%%%%%%%%%%%%%%%%%%%%%%%%%%%%%%%%%%%
%%%%%%%%%%%%%%%%%%%%%%%%%%%%%%%%%%%%%%%%%%%%%%%%%%%%%%%%%%%%%%%%%%%%%%%%%%%%%%%%%%%%%%%%%%%%%%%%%%%%%%%%%%%%%%%%%%%%%%%%%%%%%%%%
\section{Parallel Transport}
\fixed
$\Man{}$ is fixed as a manifold with dimension $m$, and $\connection{}$ as an affine connection on $\Man{}$.
\endfixed
%%%%%%%%%%%%%%%%%%%%%%%%%%%%%%%%%%%%%%%%%%%%%%%%%%%%%%%%%%%%%%%%%%%%%%%%%%%%%%%%%%%%%%%%%%%%%%%%%%%%%%%%%%%%%%%%%%%%%%%%%%%%%%%%
\textit{
The following lemma is a classical result of the theory of ordinary differential equations, concerning the
existance and uniqueness of the solution of a linear system of first-order ordinary differential equations.
The proof of this lemma
is not of concern here, and we accept its truth on faith. This lemma plays a crucial role in proving the following theorem.
}
%%%%%%%%%%%%%%%%%%%%%%%%%%%%%%%%%%%%%%%%%%%%%%%%%%%%%%%%%%%%%%%%%%%%%%%%%%%%%%%%%%%%%%%%%%%%%%%%%%%%%%%%%%%%%%%%%%%%%%%%%%%%%%%%
\lemma\label{lemlinearODEs}
Let $\interval{}$ be an open interval of $\R$ and $t_0\in\interval{}$, and let $N$ be a positive integer.
Let $\function{\lambda_{ij}}{\interval{}}{\R}$
be a real-valued smooth map on $\interval{}$, for every $i$ and $j$ in $\seta{\suc{1}{N}}$.
Let $\mtuple{c_1}{c_N}\in\R^N$.
There exists a unique
$N$-tuple $\mtuple{f_1}{f_N}$ of real-valued smooth functions on $\interval{}$ satisfying the conditions,
\begin{equation}
\begin{pmatrix}
{\dot{f}}_1\\
\vdots\\
{\dot{f}}_N
\end{pmatrix}=
\begin{pmatrix}
\lambda_{11} & \cdots & \lambda_{1N}\\
\vdots & \ddots & \vdots\\
\lambda_{N1} & \cdots & \lambda_{NN}
\end{pmatrix}
\begin{pmatrix}
f_1\\
\vdots\\
f_N
\end{pmatrix},
\end{equation}
and,
\begin{equation}
\mtuple{\func{f_1}{t_0}}{\func{f_N}{t_0}}=\mtuple{c_1}{c_N},
\end{equation}
where, ${\dot{f}}_{i}$ denotes the derived function of $f_i$.\\
Moreover, if $\bf{f}^{\bf{c}}$ denotes the unique solution $\mtuple{f_1}{f_N}$ of this linear system of ODE-s
with the initial condition ${\bf{c}}=\mtuple{c_1}{c_N}\in\R^{N}$ at $t_0$, for any choice of $\bf{c}$,
then $\bf{f}^{\bf{c}}$ is linear with respect to $\bf{c}$, that is,
\begin{equation}
\Foreach{x}{\R}
\Foreach{\opair{{\bf{c}}_1}{{\bf{c}}_2}}{\Cprod{\R^N}{\R^N}}
{\bf{f}}^{x{\bf{c}}_1+{\bf{c}}_2}=x{\bf{f}}^{{\bf{c}}_1}+{\bf{f}}^{{\bf{c}}_2}.
\end{equation}
\endlem
%%%%%%%%%%%%%%%%%%%%%%%%%%%%%%%%%%%%%%%%%%%%%%%%%%%%%%%%%%%%%%%%%%%%%%%%%%%%%%%%%%%%%%%%%%%%%%%%%%%%%%%%%%%%%%%%%%%%%%%%%%%%%%%%
\theorem\label{thmexistanceanduniquenessofparalleltransport}
Let $\function{\Curve{}}{\interval{}}{\Man{}}$ be a smooth curve on $\Man{}$ for some open interval $\interval{}$ of $\R$.
Let $t_0\in\interval{}$ and $v\in\tanspace{\func{\Curve{}}{t_0}}{\Man{}}$. There exists a unique element $\valongc{}$
of $\parallelvalongcurves{\connection{}}{\Curve{}}$, that is a unique parallel vector field on $\Man{}$ along the curve $\Curve{}$,
such that $\func{\valongc{}}{t_0}=v$.
\proof
\begin{itemize}
\item[\myitem{pr-1.}]
We first show that if the image of the curve $\Curve{}$ is included in a chart of $\Man{}$,
then there exists a unique parallel vector field on $\Man{}$ along the curve $\Curve{}$,
such that $\func{\valongc{}}{t_0}=v$.\\
Suppose that $\func{\image{\Curve{}}}{\interval{}}$ is included in the domain of a chart $\opair{\U}{\phi}$.
Let $\mtuple{\localframevecf{1}}{\localframevecf{m}}$ denote the system of local frame fields associated with
the chart $\opair{\U}{\phi}$, and let $\Curve{i}:=\cmp{\dualEucbase{m}{i}}{\cmp{\phi}{\Curve{}}}$.
Let $\mtuple{f_1}{f_m}$ denote the unique $m$-tuple of real-valued smooth functions on $\interval{}$
satisfying the linear system of first-order ordinary differential equations,
\begin{equation}\label{eqthmexistanceanduniquenessofparalleltransport1}
\Foreach{i}{\seta{\suc{1}{m}}}
\dot{f_i}=-\sum_{k=1}^{m}\sum_{j=1}^{m}
{\dot{\Curve{}}}_j
\(\cmp{\Christoffel{\phi}{i}{j}{k}}{\Curve{}}\)f_k,
\end{equation}
with the initial values,
\begin{equation}
\Foreach{i}{\seta{\suc{1}{m}}}
\func{f_i}{t_0}=v_i,
\end{equation}
where $\displaystyle v=\sum_{i=1}^{m}v_i\func{\localframevecf{i}}{\func{\Curve{}}{t_0}}$
is considered to be the unique expansion of $v$ in terms of the basis
$\mtuple{\func{\localframevecf{1}}{\func{\Curve{}}{t_0}}}{\func{\localframevecf{m}}{\func{\Curve{}}{t_0}}}$
of $\tanspace{\func{\Curve{}}{t_0}}{\Man{}}$. Now define,
\begin{equation}
\valongc{}:=\sum_{i=1}^{m}f_i\(\cmp{\localframevecf{i}}{\Curve{}}\).
\end{equation}
Clearly, $\valongc{}\in\valongcurves{\Man{}}{\Curve{}}$ such that $\func{\valongc{}}{t_0}=v$. Furthermore,
according to \refcor{corparallelvfalongcurveexpresion}, \Ref{eqthmexistanceanduniquenessofparalleltransport1}
implies that $\func{\covder{\connection{}}{\Curve{}}}{\valongc{}}=0$.
So, there exists a parallel vector field $\valongc{}$ on $\Man{}$ along $\Curve{}$ with the initial condition
$\func{\valongc{}}{t_0}=v$.\\
Now let $\valongc{1}$ be a parallel vector field along $\Curve{}$ with the initial condition
$\func{\valongc{1}}{t_0}=v$. Let $\mtuple{g_1}{g_m}$ be the unique $m$-tuple of smooth real-valued
functions on $\R$ such that $\displaystyle\valongc{1}=\sum_{i=1}^{m}g_i\(\cmp{\localframevecf{i}}{\Curve{}}\)$.
Since $\func{\covder{\connection{}}{\Curve{}}}{\valongc{1}}=0$, according to \refcor{corparallelvfalongcurveexpresion},
\begin{equation}
\Foreach{i}{\seta{\suc{1}{m}}}
\dot{g_i}=-\sum_{k=1}^{m}\sum_{j=1}^{m}
{\dot{\Curve{}}}_j
\(\cmp{\Christoffel{\phi}{i}{j}{k}}{\Curve{}}\)g_k,
\end{equation}
and moreover, considering the initial value condition, inevitably,
\begin{equation}
\Foreach{i}{\seta{\suc{1}{m}}}
\func{g_i}{t_0}=v_i.
\end{equation}
Since both $\mtuple{f_1}{f_m}$ and $\mtuple{g_1}{g_m}$ satisfy the same linear system of first order ODE-s
with the same initial value condition, \reflem{lemlinearODEs} immediately implies that $g_i=f_i$
for every $i\in\seta{\suc{1}{m}}$, and thus $\valongc{1}=\valongc{}$.
%%%%%%%%%%%%%%%%%%%%%%%%%
\item[\myitem{pr-2.}]
Now suppose that there cannot be found any chart whose domain includes the image of $\Curve{}$.
Consider first the portion $\defset{t}{\interval{}}{t>t_0}$ of $\interval{}$.
%%%%%%%%
Let $T\in\interval{}$ and $T>t_0$. Choose  a finite sequence of real numbers
$s_0<s_1<\cdots<s_{2N+1}$ within $\interval{}$ such that $s_0<t_0<s_1$ and $s_{2N}<T<s_{2N+1}$,
and further for every $\alpha\in\seta{\suc{0}{2N-1}}$,
$\func{\image{\Curve{}}}{\cpair{s_{\alpha}}{s_{\alpha+2}}}$ is included in the domain of a chart
$\opair{\U_\alpha}{\phi_\alpha}$. This can be achieved according to the fact that the interval
$\cpair{t_0}{T}$ is compact, and hence every open cover of it has definitely a finite sub-cover.
Let $\mtuple{\Localframevecf{1}{\alpha}}{\Localframevecf{m}{\alpha}}$
denote the system of local frame fields corresponded to the chart $\opair{\U_\alpha}{\phi_\alpha}$,
for every $\alpha\in\seta{\suc{0}{2N-1}}$.
Choose some $\tau_\alpha$ in
$\opair{s_\alpha}{s_{\alpha+1}}$ for every $\alpha\in\seta{\suc{1}{2N-1}}$.
For every $\alpha\in\seta{\suc{0}{2N-1}}$
and every $i$ and $k$ in $\seta{\suc{1}{m}}$, we define,
\begin{equation}\label{eqthmexistanceanduniquenessofparalleltransport11}
\lambda_{ik}^{\alpha}:=\sum_{j=1}^{m}\reS{{\dot{\Curve{}}}_j}{\opair{s_\alpha}{s_{\alpha+2}}}
\(\cmp{\Christoffel{\phi_{\alpha}}{i}{j}{k}}{\reS{\Curve{}}{\opair{s_\alpha}{s_{\alpha+2}}}}\),
\end{equation}
Obviously, each $\lambda_{ik}^{\alpha}$ is a real-valued smooth function on the open interval
$\opair{s_\alpha}{s_{\alpha+2}}$.
We now define inductively the finite sequence of $m$-tuples of real-valued smooth functions
$\seta{\mtuple{\function{f_1^\alpha}{\opair{s_0}{s_2}}{\R}}{\function{f_{m}^\alpha}{\opair{s_{2N-1}}{s_{2N+1}}}{\R}}}_{\alpha=1}^{2N-1}$
%for every $\alpha\in\seta{\suc{0}{2N-1}}$
as the following. Let $\mtuple{f_1^0}{f_m^0}$ be the unique solution of the linear system of first-order ODE-s,
\begin{equation}
\begin{pmatrix}
{\dot{f}}_1^0\\
\vdots\\
{\dot{f}}_N^0
\end{pmatrix}=
\begin{pmatrix}
\lambda_{11}^0 & \cdots & \lambda_{1N}^0\\
\vdots & \ddots & \vdots\\
\lambda_{N1}^0 & \cdots & \lambda_{NN}^0
\end{pmatrix}
\begin{pmatrix}
f_1^0\\
\vdots\\
f_N^0
\end{pmatrix},
\end{equation}
with the initial condition,
\begin{equation}
v=\sum_{i=1}^{m}\func{f_i^0}{t_0}\func{\Localframevecf{i}{0}}{\func{\Curve{}}{t_0}}.
\end{equation}
Moreover, for every $\alpha\in\seta{\suc{0}{2N-2}}$,
%assuming that $\mtuple{f_1^\alpha}{f_m^\alpha}$ is known,
let $\mtuple{f_1^{\alpha+1}}{f_m^{\alpha+1}}$ be the unique solution of the linear system of first-order ODE-s,
\begin{equation}
\begin{pmatrix}
{\dot{f}}_1^{\alpha+1}\\
\vdots\\
{\dot{f}}_N^{\alpha+1}
\end{pmatrix}=
\begin{pmatrix}
\lambda_{11}^{\alpha+1} & \cdots & \lambda_{1N}^{\alpha+1}\\
\vdots & \ddots & \vdots\\
\lambda_{N1}^{\alpha+1} & \cdots & \lambda_{NN}^{\alpha+1}
\end{pmatrix}
\begin{pmatrix}
f_1^{\alpha+1}\\
\vdots\\
f_N^{\alpha+1}
\end{pmatrix},
\end{equation}
with the initial condition,
\begin{equation}
v_{\alpha+1}=\sum_{i=1}^{m}\func{f_i^{\alpha+1}}{\tau_{\alpha+1}}
\func{\Localframevecf{i}{\alpha+1}}{\func{\Curve{}}{\tau_{\alpha+1}}},
\end{equation}
where,
\begin{equation}
v_{\alpha+1}:=\sum_{i=1}^{m}\func{f_i^{\alpha}}{\tau_{\alpha+1}}
\func{\Localframevecf{i}{\alpha}}{\func{\Curve{}}{\tau_{\alpha+1}}}.
\end{equation}
Now, for every $\alpha\in\seta{\suc{0}{2N-1}}$, we define the map
$\function{\valongc{\alpha}}{\opair{s_\alpha}{s_{\alpha+2}}}{\tanbun{\Man{}}}$ as,
\begin{equation}
\Foreach{t}{\opair{s_\alpha}{s_{\alpha+2}}}
\func{\valongc{\alpha}}{t}\eqdef
\sum_{i=1}^{m}\func{f_i^\alpha}{t}\func{\Localframevecf{i}{\alpha}}{\func{\Curve{}}{t}}.
\end{equation}
It is trivial that each $\valongc{\alpha}$ is a smooth vector field along the curve
$\reS{\Curve{}}{\opair{s_{\alpha}}{s_{\alpha+2}}}$.
Moreover, based on \refcor{corparallelvfalongcurveexpresion}, and considering the way in which each $\valongc{\alpha}$
is defined, it is completely clear that each $\valongc{\alpha}$ is a parallal vector field along the curve
$\reS{\Curve{}}{\opair{s_{\alpha}}{s_{\alpha+2}}}$, that is an element of
$\parallelvalongcurves{\connection{}}{\reS{\Curve{}}{\opair{s_{\alpha}}{s_{\alpha+2}}}}$.\\
Now, for every $\alpha\in\seta{\suc{0}{2N-2}}$, consider the intersection of the domains of
$\valongc{\alpha}$ and $\valongc{\alpha+1}$, that is $\opair{s_{\alpha+1}}{s_{\alpha+2}}$.
According to the naturality of covariant derivative of vector fields along curves with respect to restricting of curves
to open sub-intervals of their domains, and according to \refdef{defparallelvectorfieldsalongcurves},
it is clear that the restriction of a parallel vector field along a curve to some open sub-interval
of the domain of that curve is again a parallel vector field along the restriction of that curve to that sub-interval.
Therefore, the restriction of each $\valongc{\alpha}$ and $\valongc{\alpha+1}$ to $\opair{s_{\alpha+1}}{s_{\alpha+2}}$
is a parallel vector field along the smooth curve $\reS{\Curve{}}{\opair{s_{\alpha+1}}{s_{\alpha+2}}}$.
Furthermore, knowing that the image of $\reS{\Curve{}}{\opair{s_{\alpha+1}}{s_{\alpha+2}}}$ is included in
the domain of some chart of $\Man{}$, for example any of the charts $\opair{\U_\alpha}{\phi_\alpha}$
and $\opair{\U_{\alpha+1}}{\phi_{\alpha+1}}$, the result of \myitem{pr-1} implies that there exists just
one parallel vector field along $\reS{\Curve{}}{\opair{s_{\alpha+1}}{s_{\alpha+2}}}$. Thus,
$\valongc{\alpha}$ and $\valongc{\alpha+1}$ coincide in the intersection of their domains, that is,
\begin{equation}
\Foreach{\alpha}{\seta{\suc{0}{2N-2}}}
\reS{\valongc{\alpha}}{\opair{s_{\alpha+1}}{s_{\alpha+2}}}=
\reS{\valongc{\alpha+1}}{\opair{s_{\alpha+1}}{s_{\alpha+2}}}.
\end{equation}
Hence, we can define a map $\function{\valongc{}}{\opair{s_0}{s_{2N+1}}}{\tanbun{\Man{}}}$ as,
\begin{equation}
\Foreach{\alpha}{\seta{\suc{0}{2N-1}}}
\Foreach{t}{\opair{s_\alpha}{s_{\alpha+2}}}
\func{\valongc{}}{t}\eqdef\func{\valongc{\alpha}}{t}.
\end{equation}
Trivially, $\valongc{}$ is smooth, and furthermore considering form of the ODE-s governing the local expressions of
$\valongc{}$ for some family of charts whose domains cover the image of $\reS{\Curve{}}{\opair{s_0}{s_{2N+1}}}$,
\refcor{corparallelvfalongcurveexpresion} immediately implies that $\valongc{}$ is a parallel vector field
on $\Man{}$ along the curve $\reS{\Curve{}}{\opair{s_0}{s_{2N+1}}}$, that is
$\valongc{}\in\parallelvalongcurves{\connection{}}{\reS{\Curve{}}{\opair{s_{0}}{s_{2N+1}}}}$.\\
If $T<t_0$, then we can again construct a parallel vector field $\valongc{}$ on $\Man{}$ along $\reS{\Curve{}}{\V}$
for some open sub-interval $\V$ of $\interval{}$ containing $t_0$ and $T$ with $\func{\valongc{}}{t_0}=v$,
in a way absolutely similar to that of the case $T<t_0$ explained above.\\
Ultimately, it is seen that for any $T\in\interval{}$ we can construct and fix a parallel
vector field $\Valongc{}{T}$ along $\reS{\Curve{}}{\V_T}$ for some open sub-interval $\V_T$
of $\interval{}$ containing $t_0$ and $T$ with $\func{\Valongc{}{T}}{t_0}=v$.
%%%%%%%%%%%%%%%%%%%%%%%%%
\item[\myitem{pr-3.}]
Now, we show that given any smooth curve $\function{\Curve{1}}{\interval{1}}{\Man{}}$, if there exists a
parallel smooth vector field $\valongc{}$ along $\Curve{1}$ such that $\func{\valongc{}}{t_0}=v$
for some $t_0\in\interval{1}$ and some $v\in\tanspace{\func{\Curve{1}}{t_0}}{\Man{}}$, it is unique.\\
Let $\valongc{}$ be such a parallel smooth vector field along a smooth curve $\function{\Curve{1}}{\interval{1}}{\Man{}}$.
Choose an arbitrary $T$ in the open interval $\interval{1}$. Similar to \myitem{pr-2},
Choose  a finite sequence of real numbers
$s_0<s_1<\cdots<s_{2N+1}$ within $\interval{}$ such that $s_0<t_0<s_1$ and $s_{2N}<T<s_{2N+1}$,
and further for every $\alpha\in\seta{\suc{0}{2N-1}}$,
$\func{\image{\Curve{}}}{\cpair{s_{\alpha}}{s_{\alpha+2}}}$ is included in the domain of a chart
$\opair{\U_\alpha}{\phi_\alpha}$. Let $\mtuple{\Localframevecf{1}{\alpha}}{\Localframevecf{m}{\alpha}}$
denote the system of local frame fields corresponded to the chart $\opair{\U_\alpha}{\phi_\alpha}$,
for every $\alpha\in\seta{\suc{0}{2N-1}}$.
Choose some $\tau_\alpha$ in
$\opair{s_\alpha}{s_{\alpha+1}}$ for every $\alpha\in\seta{\suc{1}{2N-1}}$. Foe every
$\alpha\in\seta{\suc{0}{2N-1}}$, let $\mtuple{f_1^\alpha}{f_m^\alpha}$ denote the
unique $m$-tuple of real-valued smooth functions on $\opair{s_{\alpha}}{s_{\alpha+2}}$
such that $\displaystyle\reS{\valongc{}}{\opair{s_{\alpha}}{s_{\alpha+2}}}=
\sum_{i=1}^{m}f_i\(\cmp{\Localframevecf{i}{\alpha}}{\reS{\Curve{1}}{\opair{s_{\alpha}}{s_{\alpha+2}}}}\)$.
According to \refcor{corparallelvfalongcurveexpresion}, for each $\alpha\in\seta{\suc{0}{2N-1}}$,
$\mtuple{\Localframevecf{1}{\alpha}}{\Localframevecf{m}{\alpha}}$ must satisfy the linear system of
first-order ODE-s,
\begin{equation}
\begin{pmatrix}
{\dot{f}}_1^{\alpha}\\
\vdots\\
{\dot{f}}_N^{\alpha}
\end{pmatrix}=
\begin{pmatrix}
\lambda_{11}^{\alpha} & \cdots & \lambda_{1N}^{\alpha}\\
\vdots & \ddots & \vdots\\
\lambda_{N1}^{\alpha} & \cdots & \lambda_{NN}^{\alpha}
\end{pmatrix}
\begin{pmatrix}
f_1^{\alpha}\\
\vdots\\
f_N^{\alpha}
\end{pmatrix},
\end{equation}
for the smooth maps $\lambda_{ij}^{\alpha}$ defined similar to what is done in
\Ref{eqthmexistanceanduniquenessofparalleltransport11}.
Furthermore, the liear system of ODE-s stated above must satisfy successively the initial conditions,
\begin{equation}
v=\sum_{i=1}^{m}\func{f_i^0}{t_0}\func{\Localframevecf{i}{0}}{\func{\Curve{}}{t_0}},
\end{equation}
and,
\begin{equation}
\Foreach{\alpha}{\seta{\suc{0}{2N-2}}}
v_{\alpha+1}=\sum_{i=1}^{m}\func{f_i^{\alpha+1}}{\tau_{\alpha+1}}
\func{\Localframevecf{i}{\alpha+1}}{\func{\Curve{}}{\tau_{\alpha+1}}},
\end{equation}
where,
\begin{equation}
\Foreach{\alpha}{\seta{\suc{0}{2N-2}}}
v_{\alpha+1}:=\sum_{i=1}^{m}\func{f_i^{\alpha}}{\tau_{\alpha+1}}
\func{\Localframevecf{i}{\alpha}}{\func{\Curve{}}{\tau_{\alpha+1}}}.
\end{equation}
Now, pay attention that independent of the map $\valongc{}$, the succession of
$m$-tuples of smooth functions $\seta{\mtuple{f_1^\alpha}{f_m^\alpha}}_{\alpha=0}^{2N-1}$
is uniquely determined according to the uniqueness of solutions of the ODE-s with the initial values
as indicated above. So, in particular, since $\displaystyle\func{\valongc{}}{T}=\sum_{i=1}^m\func{f_{i}^{2N-1}}{T}
\func{\Localframevecf{i}{2N-1}}{\func{\Curve{1}}{T}}$, $\func{\valongc{}}{T}$ is uniquely determined.
Moreover, considering that the choice of $T\in\interval{1}$ was arbitrary, we conclude that
$\valongc{}$ must be uniquely determined.
%%%%%%%%%%%%%%%%%%%%%%%%%
\item[\myitem{pr-4.}]
Now, for every $T\in\interval{}$, let $\Valongc{}{T}$ be the parallel
vector field $\Valongc{}{T}$ along $\reS{\Curve{}}{\V_T}$ for some open sub-interval $\V_T$
of $\interval{}$ containing $t_0$ and $T$ with $\func{\Valongc{}{T}}{t_0}=v$, as constructed and fixed in
\myitem{pr-2}. According to the result of \myitem{pr-3}, and considering the naturality of covariant derivative
of vector fields along curves with respect to restricting of curves to open sub-intervals of their domains,
it becomes apparent that for every $T_1$ and $T_2$ in $\interval{}$, $\Valongc{}{T_1}$ and
$\Valongc{}{T_2}$ must coincide in the region that their domains overlap. Hence, it is permissible to define
a map $\function{\valongc{}}{\interval{}}{\tanbun{\Man{}}}$ as,
\begin{equation}
\Foreach{t}{\interval{}}
\func{\valongc{}}{t}\eqdef\func{\Valongc{}{t}}{t}.
\end{equation}
The map $\valongc{}$ is obviously a smooth vector field on $\Man{}$ along the curv $\Curve{}$.
Moreover, considering the arguments of \myitem{pr-2}, according to \refcor{corparallelvfalongcurveexpresion},
it becomes clear that $\valongc{}$ is a parallel smooth vector field on $\Man{}$ along the curve $\Curve{}$,
and thus the existance of an element of $\parallelvalongcurves{\connection{}}{\Curve{}}$ with the desired initial condition
is assured, which additionally must be unique according to the result of \myitem{pr-3}.
\end{itemize}
\endthm
%%%%%%%%%%%%%%%%%%%%%%%%%%%%%%%%%%%%%%%%%%%%%%%%%%%%%%%%%%%%%%%%%%%%%%%%%%%%%%%%%%%%%%%%%%%%%%%%%%%%%%%%%%%%%%%%%%%%%%%%%%%%%%%%
\definition\label{defparalleltransport0}
Let $\function{\Curve{}}{\interval{}}{\Man{}}$ be a smooth curve on $\Man{}$ for some open interval $\interval{}$ of $\R$.
Let $t_0\in\interval{}$ and $v\in\tanspace{\func{\Curve{}}{t_0}}{\Man{}}$. We will denote by
$\pvfalongc{\Curve{}}{t_0}{v}{\connection{}}$
the unique parallel smooth vector field on $\Man{}$ along $\Curve{}$ having the value $v$ at $t_0$.
\endef
%%%%%%%%%%%%%%%%%%%%%%%%%%%%%%%%%%%%%%%%%%%%%%%%%%%%%%%%%%%%%%%%%%%%%%%%%%%%%%%%%%%%%%%%%%%%%%%%%%%%%%%%%%%%%%%%%%%%%%%%%%%%%%%%
\definition\label{defparalleltransport}
Let $\function{\Curve{}}{\interval{}}{\Man{}}$ be a smooth curve on $\Man{}$ for some open interval $\interval{}$ of $\R$.
Let $t_0$ and $t_1$ be a pair of elements of $\interval{}$.
We define the map $\function{\ptransport{\Curve{}}{t_0}{t_1}{\connection{}}}{\tanspace{\func{\Curve{}}{t_0}}{\Man{}}}
{\tanspace{\func{\Curve{}}{t_1}}{\Man{}}}$ as,
\begin{equation}
\Foreach{v}{\tanspace{\func{\Curve{}}{t_0}}{\Man{}}}
\func{\[\ptransport{\Curve{}}{t_0}{t_1}{\connection{}}\]}{v}\eqdef
\func{\[\pvfalongc{\Curve{}}{t_0}{v}{\connection{}}\]}{t_1}.
\end{equation}
$\ptransport{\Curve{}}{t_0}{t_1}{\connection{}}$ is referred to as the $\quotl$parallel transport operator of the
curve $\Curve{}$ (with respect to the connection $\connection{}$ on $\Man{}$) from $t_0$ to $t_1$$\quotr$.
For eny $v\in\tanspace{\func{\Curve{}}{t_0}}{\Man{}}$, $\func{\[\ptransport{\Curve{}}{t_0}{t_1}{\connection{}}\]}{v}$ is
called the $\quotl$parallel transport of $v$ along the curve $\Curve{}$ (with respect to the connection $\connection{}$)
from $t_0$ to $t_1$$\quotr$.
When there is no ambiguity about the connection $\connection{}$, $\ptransport{\Curve{}}{t_0}{t_1}{\connection{}}$
can be simply denoted by $\ptransport{\Curve{}}{t_0}{t_1}{}$.
\endef
%%%%%%%%%%%%%%%%%%%%%%%%%%%%%%%%%%%%%%%%%%%%%%%%%%%%%%%%%%%%%%%%%%%%%%%%%%%%%%%%%%%%%%%%%%%%%%%%%%%%%%%%%%%%%%%%%%%%%%%%%%%%%%%%
\theorem
Let $\function{\Curve{}}{\interval{}}{\Man{}}$ be a smooth curve on $\Man{}$ for some open interval $\interval{}$ of $\R$.
Let $t_0$ and $t_1$ be a pair of elements of $\interval{}$. The operator $\ptransport{\Curve{}}{t_0}{t_1}{\connection{}}$
is a linear-isomorphism from $\tanspace{\func{\Curve{}}{t_0}}{\Man{}}$ to $\tanspace{\func{\Curve{}}{t_1}}{\Man{}}$, and,
\begin{equation}
\finv{\[\ptransport{\Curve{}}{t_0}{t_1}{\connection{}}\]}=\ptransport{\Curve{}}{t_1}{t_0}{\connection{}}.
\end{equation}
\proof
\begin{itemize}
\item[\myitem{pr-1.}]
According to \refdef{defparalleltransport} and the way of globalizing a local vector field along $\Curve{}$ with some
given initial value as explained in the proof of \refthm{thmexistanceanduniquenessofparalleltransport},
the linearity of $\ptransport{\Curve{}}{t_0}{t_1}{\connection{}}$ is clearly seen to be a direct consequence of
the linear behavior of the solution of a linear system of first-order ODE-s with respect to the initial value, as
precisely stated in \reflem{lemlinearODEs}.
\item[\myitem{pr-2.}]
Let $v$ be an arbitrary element of $\tanspace{\func{\Curve{}}{t_0}}{\Man{}}$. According to \refdef{defparalleltransport0},
\begin{align}
\func{\[\ptransport{\Curve{}}{t_1}{t_0}{}\]}
{\func{\[\ptransport{\Curve{}}{t_0}{t_1}{}\]}{v}}&=
\func{\[\pvfalongc{\Curve{}}{t_1}{\func{\[\ptransport{\Curve{}}{t_0}{t_1}{}\]}{v}}{\connection{}}\]}{t_0}\cr
&=\func{\[\pvfalongc{\Curve{}}{t_1}{\func{\[\pvfalongc{\Curve{}}{t_0}{v}{\connection{}}\]}{t_1}}{\connection{}}\]}{t_0}.
\end{align}
Note that $\pvfalongc{\Curve{}}{t_1}{\func{\[\pvfalongc{\Curve{}}{t_0}{v}{\connection{}}\]}{t_1}}{\connection{}}$
is an element of $\parallelvalongcurves{\connection{}}{\Curve{}}$ with the value
$\func{\[\pvfalongc{\Curve{}}{t_0}{v}{\connection{}}\]}{t_1}$ at $t_1$. On the other hand, it is trivial that
$\pvfalongc{\Curve{}}{t_0}{v}{\connection{}}$
is also an element of $\parallelvalongcurves{\connection{}}{\Curve{}}$ with the value
$\func{\[\pvfalongc{\Curve{}}{t_0}{v}{\connection{}}\]}{t_1}$ at $t_1$. So, since such an element of
$\parallelvalongcurves{\connection{}}{\Curve{}}$ is unique, clearly,
\begin{equation}
\pvfalongc{\Curve{}}{t_1}{\func{\[\pvfalongc{\Curve{}}{t_0}{v}{\connection{}}\]}{t_1}}{\connection{}}=
\pvfalongc{\Curve{}}{t_0}{v}{\connection{}},
\end{equation}
and thus,
\begin{align}
\func{\[\ptransport{\Curve{}}{t_1}{t_0}{}\]}
{\func{\[\ptransport{\Curve{}}{t_0}{t_1}{}\]}{v}}&=
\func{\[\pvfalongc{\Curve{}}{t_0}{v}{\connection{}}\]}{t_0}\cr
&=v.
\end{align}
Since the choice of $v$ was arbitrary, we conclude that,
\begin{equation}
\cmp{\[\ptransport{\Curve{}}{t_1}{t_0}{}\]}{\[\ptransport{\Curve{}}{t_0}{t_1}{}\]}=
\identity{\tanspace{\func{\Curve{}}{t_0}}{\Man{}}}.
\end{equation}
Similarly, we must have,
\begin{equation}
\cmp{\[\ptransport{\Curve{}}{t_0}{t_1}{}\]}{\[\ptransport{\Curve{}}{t_1}{t_0}{}\]}=
\identity{\tanspace{\func{\Curve{}}{t_1}}{\Man{}}}.
\end{equation}
Therefore, $\ptransport{\Curve{}}{t_0}{t_1}{}$ is injective and onto, and,
\begin{equation}
\finv{\[\ptransport{\Curve{}}{t_0}{t_1}{}\]}=\ptransport{\Curve{}}{t_1}{t_0}{}.
\end{equation}
\end{itemize}
\endthm
%%%%%%%%%%%%%%%%%%%%%%%%%%%%%%%%%%%%%%%%%%%%%%%%%%%%%%%%%%%%%%%%%%%%%%%%%%%%%%%%%%%%%%%%%%%%%%%%%%%%%%%%%%%%%%%%%%%%%%%%%%%%%%%%
\definition
Let $\function{\Curve{}}{\interval{}}{\Man{}}$ be a smooth curve on $\Man{}$ for some open interval $\interval{}$ of $\R$.
Let $\mtuple{\framealongc{1}}{\framealongc{m}}$ be an $m$-tuple of parallel smooth vector fields
along the curve $\Curve{}$, that is elements of $\parallelvalongcurves{\connection{}}{\Curve{}}$, such that
for every $t\in\interval{}$, $\mtuple{\func{\framealongc{1}}{t}}{\func{\framealongc{m}}{t}}$ is an ordered-basis of
the vector-space $\Tanspace{\func{\Curve{}}{t}}{\Man{}}$. Then,
$\mtuple{\framealongc{1}}{\framealongc{m}}$ is referred to as a $\quotl$parallel (moving) frame along the curve $\Curve{}$
(relative to the connection $\connection{}$ on $\Man{}$)$\quotr$.
\endef
%%%%%%%%%%%%%%%%%%%%%%%%%%%%%%%%%%%%%%%%%%%%%%%%%%%%%%%%%%%%%%%%%%%%%%%%%%%%%%%%%%%%%%%%%%%%%%%%%%%%%%%%%%%%%%%%%%%%%%%%%%%%%%%%
\theorem\label{thmcovariantderivativeintermsofaparallelframe}
Let $\function{\Curve{}}{\interval{}}{\Man{}}$ be a smooth curve on $\Man{}$ for some open interval $\interval{}$ of $\R$.
Let $\mtuple{\framealongc{1}}{\framealongc{m}}$ be a parallel frame along $\Curve{}$, and let $\valongc{}$ be a smooth
smooth vector field along $\Curve{}$, that is an element of $\valongcurves{\Man{}}{\Curve{}}$.
Let $\mtuple{f_1}{f_m}$ denote the unique $m$-tuple of real-valued smooth functions
on $\interval{}$ such that $\displaystyle\valongc{}=\sum_{i=1}^{m}f_i\framealongc{i}$. Then,
\begin{equation}
\func{\covder{\connection{}}{\Curve{}}}{\valongc{}}=\sum_{i=1}^{m}{\dot{f}}_i\framealongc{i}.
\end{equation}
\proof
Since for each $i\in\seta{\suc{1}{m}}$, $\func{\covder{\connection{}}{\Curve{}}}{\framealongc{i}}=0$,
according to the definition of covariant derivative of smooth vector fields along curves,
\begin{align}
\func{\covder{\connection{}}{\Curve{}}}{\valongc{}}&=
\func{\covder{\connection{}}{\Curve{}}}{\sum_{i=1}^{m}f_i\framealongc{i}}\cr
&=\sum_{i=1}^{m}{\dot{f}}_i\framealongc{i}+
\sum_{i=1}^{m}f_i\func{\covder{\connection{}}{\Curve{}}}{\framealongc{i}}\cr
&=\sum_{i=1}^{m}{\dot{f}}_i\framealongc{i}.
\end{align}
\endthm
%%%%%%%%%%%%%%%%%%%%%%%%%%%%%%%%%%%%%%%%%%%%%%%%%%%%%%%%%%%%%%%%%%%%%%%%%%%%%%%%%%%%%%%%%%%%%%%%%%%%%%%%%%%%%%%%%%%%%%%%%%%%%%%%
\theorem
Let $\function{\Curve{}}{\interval{}}{\Man{}}$ be a smooth curve on $\Man{}$ for some open interval $\interval{}$ of $\R$.
Let $t_0\in\interval{}$. Let $\mtuple{v_1}{v_m}$ be an ordered-basis of the vector-space
$\Tanspace{\func{\Curve{}}{t_0}}{\Man{}}$. For every $t\in\interval{}$,
$\mtuple{\func{\[\ptransport{\Curve{}}{t_0}{t}{\connection{}}\]}{v_1}}
{\func{\[\ptransport{\Curve{}}{t_0}{t}{\connection{}}\]}{v_1}}=
\mtuple{\func{\[\pvfalongc{\Curve{}}{t_0}{v_1}{\connection{}}\]}{t}}
{\func{\[\pvfalongc{\Curve{}}{t_0}{v_m}{\connection{}}\]}{t}}$ is an ordered-basis of
$\Tanspace{\func{\Curve{}}{t}}{\Man{}}$, and hence $\mtuple{\pvfalongc{\Curve{}}{t_0}{v_1}{\connection{}}}
{\pvfalongc{\Curve{}}{t_0}{v_m}{\connection{}}}$ is a parallel frame along $\Curve{}$.
\proof
For every $t\in\interval{}$, $\ptransport{\Curve{}}{t_0}{t}{\connection{}}$ is a linear-isomorphism from
$\Tanspace{\func{\Curve{}}{t_0}}{\Man{}}$ to $\Tanspace{\func{\Curve{}}{t}}{\Man{}}$, and hence carries any set of
$m$ linearly-independent vectors to a set of $m$ linearly-independent vectors.
\endthm
%%%%%%%%%%%%%%%%%%%%%%%%%%%%%%%%%%%%%%%%%%%%%%%%%%%%%%%%%%%%%%%%%%%%%%%%%%%%%%%%%%%%%%%%%%%%%%%%%%%%%%%%%%%%%%%%%%%%%%%%%%%%%%%%
\corollary
Let $\function{\Curve{}}{\interval{}}{\Man{}}$ be a smooth curve on $\Man{}$ for some open interval $\interval{}$ of $\R$.
There exists a parallel moving frame along $\Curve{}$ (relative to the connection $\connection{}$ on $\Man{}$).
\endcor
%%%%%%%%%%%%%%%%%%%%%%%%%%%%%%%%%%%%%%%%%%%%%%%%%%%%%%%%%%%%%%%%%%%%%%%%%%%%%%%%%%%%%%%%%%%%%%%%%%%%%%%%%%%%%%%%%%%%%%%%%%%%%%%%
\theorem\label{thmcovariantderivativeintermsofparalleltransport}
Let $\function{\Curve{}}{\interval{}}{\Man{}}$ be a smooth curve on $\Man{}$ for some open interval $\interval{}$ of $\R$.
Let $t_0\in\interval{}$ and let $\valongc{}$ be an element of $\valongcurves{\Man{}}{\Curve{}}$,
that is a smooth vector field along $\Curve{}$.\\
\begin{equation}
\func{\[\func{\covder{\connection{}}{\Curve{}}}{\valongc{}}\]}{t_0}=
\lim_{t\to t_0}\frac{\func{\[\ptransport{\Curve{}}{t}{t_0}{\connection{}}\]}{\func{\valongc{}}{t}}-
\func{\valongc{}}{t_0}}{t-t_0}.
\end{equation}
Note that $\Tanspace{\func{\Curve{}}{t_0}}{\Man{}}$ is assumed to be normed-space,
that is a vector-space endowed with a norm, in order to make the limit well-defined.
But since $\Tanspace{\func{\Curve{}}{t_0}}{\Man{}}$ is a finite-dimensional vector-space,
all of the norms on it are equivalent. Thus, the choice of a specific norm does not make sense,
and $\Tanspace{\func{\Curve{}}{t_0}}{\Man{}}$ is assumed to be endowed with an arbitrary norm.
\proof
We assume that
Choose a parallel moving frame $\mtuple{\framealongc{1}}{\framealongc{m}}$ along $\Curve{}$. Let
$\displaystyle\valongc{}=\sum_{i=1}^{m}f_i\framealongc{i}$. According to
\refthm{thmcovariantderivativeintermsofaparallelframe},
\begin{align}
\displaystyle
\func{\[\func{\covder{\connection{}}{\Curve{}}}{\valongc{}}\]}{t_0}&=
\sum_{i=1}^{m}\func{{\dot{f}}_i}{t_0}\func{\framealongc{i}}{t_0}\cr
&=\sum_{i=1}^{m}\(\lim_{t\to t_0}\frac{\func{f_i}{t}-\func{f_i}{t_0}}{t-t_0}\)\func{\framealongc{i}}{t_0}\cr
&=\lim_{t\to t_0}\frac{\displaystyle\sum_{i=1}^{m}\func{f_i}{t}\func{\framealongc{i}}{t_0}-\sum_{i=1}^{m}\func{f_i}{t_0}
\func{\framealongc{i}}{t_0}}{t-t_0}\cr
&=\lim_{t\to t_0}\frac{\displaystyle\sum_{i=1}^{m}\func{f_i}{t}\func{\framealongc{i}}{t_0}-\func{\valongc{}}{t_0}}{t-t_0}.
\end{align}
On the other hand, considering the linearity of parallel transport operator of $\Curve{}$ from $t$ to $t_0$ for any
$t\in\interval{}$,
\begin{align}
\Foreach{t}{\interval{}}
\func{\[\ptransport{\Curve{}}{t}{t_0}{}\]}{\func{\valongc{}}{t}}&=
\func{\[\ptransport{\Curve{}}{t}{t_0}{}\]}{\sum_{i=1}^{m}\func{f_i}{t}\func{\framealongc{i}}{t}}\cr
&=\sum_{i=1}^{m}\func{f_i}{t}\func{\[\ptransport{\Curve{}}{t}{t_0}{}\]}{\func{\framealongc{i}}{t}}\cr
&=\sum_{i=1}^{m}\func{f_i}{t}\func{\framealongc{i}}{t_0}.
\end{align}
Therefore, the result is established.
\endthm
%%%%%%%%%%%%%%%%%%%%%%%%%%%%%%%%%%%%%%%%%%%%%%%%%%%%%%%%%%%%%%%%%%%%%%%%%%%%%%%%%%%%%%%%%%%%%%%%%%%%%%%%%%%%%%%%%%%%%%%%%%%%%%%%
\corollary
Let $\avecf{}$ and $\avecff{}$ be a pair of elements of $\vectorfields{\Man{}}$, and $\point\in\Man{}$.
For every smooth curve $\function{\Curve{}}{\interval{}}{\Man{}}$ on $\Man{}$ such that
$\func{\Curve{}}{t_0}=\point$ and $\func{\dot{\Curve{}}}{t_0}=\func{\avecf{}}{\point}$ for some $t_0\in\interval{}$,
\begin{equation}
\func{\[\con{\avecf{}}{\avecff{}}\]}{\point}=
\lim_{t\to t_0}\frac{\func{\[\ptransport{\Curve{}}{t}{t_0}{}\]}{\func{\avecff{}}{\func{\Curve{}}{t}}}-
\func{\avecff{}}{\point}}{t-t_0}.
\end{equation}
\endcor
%%%%%%%%%%%%%%%%%%%%%%%%%%%%%%%%%%%%%%%%%%%%%%%%%%%%%%%%%%%%%%%%%%%%%%%%%%%%%%%%%%%%%%%%%%%%%%%%%%%%%%%%%%%%%%%%%%%%%%%%%%%%%%%%
%%%%%%%%%%%%%%%%%%%%%%%%%%%%%%%%%%%%%%%%%%%%%%%%%%%%%%%%%%%%%%%%%%%%%%%%%%%%%%%%%%%%%%%%%%%%%%%%%%%%%%%%%%%%%%%%%%%%%%%%%%%%%%%%
%%%%%%%%%%%%%%%%%%%%%%%%%%%%%%%%%%%%%%%%%%%%%%%%%%%%%%%%%%%%%%%%%%%%%%%%%%%%%%%%%%%%%%%%%%%%%%%%%%%%%%%%%%%%%%%%%%%%%%%%%%%%%%%%
%%%%%%%%%%%%%%%%%%%%%%%%%%%%%%%%%%%%%%%%%%%%%%%%%%%%%%%%%%%%%%%%%%%%%%%%%%%%%%%%%%%%%%%%%%%%%%%%%%%%%%%%%%%%%%%%%%%%%%%%%%%%%%%%
%%%%%%%%%%%%%%%%%%%%%%%%%%%%%%%%%%%%%%%%%%%%%%%%%%%%%%%%%%%%%%%%%%%%%%%%%%%%%%%%%%%%%%%%%%%%%%%%%%%%%%%%%%%%%%%%%%%%%%%%%%%%%%%%
%%%%%%%%%%%%%%%%%%%%%%%%%%%%%%%%%%%%%%%%%%%%%%%%%%%%%%%%%%%%%%%%%%%%%%%%%%%%%%%%%%%%%%%%%%%%%%%%%%%%%%%%%%%%%%%%%%%%%%%%%%%%%%%%
%%%%%%%%%%%%%%%%%%%%%%%%%%%%%%%%%%%%%%%%%%%%%%%%%%%%%%%%%%%%%%%%%%%%%%%%%%%%%%%%%%%%%%%%%%%%%%%%%%%%%%%%%%%%%%%%%%%%%%%%%%%%%%%%
%%%%%%%%%%%%%%%%%%%%%%%%%%%%%%%%%%%%%%%%%%%%%%%%%%%%%%%%%%%%%%%%%%%%%%%%%%%%%%%%%%%%%%%%%%%%%%%%%%%%%%%%%%%%%%%%%%%%%%%%%%%%%%%%
%%%%%%%%%%%%%%%%%%%%%%%%%%%%%%%%%%%%%%%%%%%%%%%%%%%%%%%%%%%%%%%%%%%%%%%%%%%%%%%%%%%%%%%%%%%%%%%%%%%%%%%%%%%%%%%%%%%%%%%%%%%%%%%%
%%%%%%%%%%%%%%%%%%%%%%%%%%%%%%%%%%%%%%%%%%%%%%%%%%%%%%%%%%%%%%%%%%%%%%%%%%%%%%%%%%%%%%%%%%%%%%%%%%%%%%%%%%%%%%%%%%%%%%%%%%%%%%%%
%%%%%%%%%%%%%%%%%%%%%%%%%%%%%%%%%%%%%%%%%%%%%%%%%%%%%%%%%%%%%%%%%%%%%%%%%%%%%%%%%%%%%%%%%%%%%%%%%%%%%%%%%%%%%%%%%%%%%%%%%%%%%%%%
\section{Geodesics}
\fixed
$\Man{}$ is fixed as a manifold with dimension $m$, and $\connection{}$ as an affine connection on $\Man{}$.
\endfixed
%%%%%%%%%%%%%%%%%%%%%%%%%%%%%%%%%%%%%%%%%%%%%%%%%%%%%%%%%%%%%%%%%%%%%%%%%%%%%%%%%%%%%%%%%%%%%%%%%%%%%%%%%%%%%%%%%%%%%%%%%%%%%%%%
\theorem
Let $\function{\Curve{}}{\interval{}}{\Man{}}$ be a smooth curve on $\Man{}$ for some open interval $\interval{}$ of $\R$.
The velocity map of $\Curve{}$ is a smooth vector field on $\Man{}$ along the curve $\Curve{}$, that is
$\dot{\Curve{}}\in\valongcurves{\Man{}}{\Curve{}}$.
\proof
It is trivial.
\endthm
%%%%%%%%%%%%%%%%%%%%%%%%%%%%%%%%%%%%%%%%%%%%%%%%%%%%%%%%%%%%%%%%%%%%%%%%%%%%%%%%%%%%%%%%%%%%%%%%%%%%%%%%%%%%%%%%%%%%%%%%%%%%%%%%
\definition\label{defgeodesiconmanifold}
Let $\function{\Curve{}}{\interval{}}{\Man{}}$ be a smooth curve on $\Man{}$ for some open interval $\interval{}$ of $\R$.
\begin{itemize}
\item
$\func{\covder{\connection{}}{\Curve{}}}{\dot{\Curve{}}}$ is called the $\quotl$acceleration map of the curve $\Curve{}$
with respect to the connection $\connection{}$ on $\Man{}$$\quotr$.
\item
$\Curve{}$ is called a $\quotl$geodesic on $\Man{}$ relative to the connection $\connection{}$$\quotr$ if its
acceleration map is a parallel smooth vector field along $\Curve{}$, that is
$\func{\covder{\connection{}}{\Curve{}}}{\dot{\Curve{}}}\in\parallelvalongcurves{\connection{}}{\Curve{}}$, or in
other words,
\begin{equation}
\Foreach{t}{\interval{}}
\func{\[\func{\covder{\connection{}}{\Curve{}}}{\dot{\Curve{}}}\]}{t}=\zerovec{}.
\end{equation}
\item
The set of all geodesics on $\Man{}$ relative to the connection $\connection{}$ will be denoted by
$\geodesics{\Man{}}{\connection{}}$.
\end{itemize}
\endef
%%%%%%%%%%%%%%%%%%%%%%%%%%%%%%%%%%%%%%%%%%%%%%%%%%%%%%%%%%%%%%%%%%%%%%%%%%%%%%%%%%%%%%%%%%%%%%%%%%%%%%%%%%%%%%%%%%%%%%%%%%%%%%%%
\theorem\label{thmgeodesiclocalexpresion}
Let $\function{\Curve{}}{\interval{}}{\Man{}}$ be a curve on $\Man{}$ for some open interval $\interval{}$ of $\R$.
The following statements are equivalent.
\begin{itemize}
\item[\myitem{1.}]
$\Curve{}$ is a geodesic on $\Man{}$ relative to $\connection{}$, that is
$\func{\covder{\connection{}}{\Curve{}}}{\dot{\Curve{}}}\in\parallelvalongcurves{\connection{}}{\Curve{}}$.
\item[\myitem{2.}]
For every chart
$\opair{\U}{\phi}$ of $\Man{}$ that intersects the image of $\Curve{}$, and every sub-interval $O$ of
$\interval{}$ such that $\defSet{\func{\Curve{}}{t}}{t\in O}\subseteq\U$,
\begin{align}
&\Foreach{i}{\seta{\suc{1}{m}}}\Foreach{t}{O}\cr
&\func{\ddot{\Curve{}}_{i}}{t}+\sum_{k=1}^{m}\sum_{j=1}^{m}
\func{\Christoffel{\phi}{i}{j}{k}}{\func{\Curve{}}{t}}\func{{\dot{\Curve{}}}_{j}}{t}
\func{{\dot{\Curve{}}}_{k}}{t}=0,
\end{align}
or equivalently,
\begin{equation}
\reS{{\ddot{\Curve{}}}_i}{O}+\sum_{k=1}^{m}\sum_{j=1}^{m}
\(\cmp{\Christoffel{\phi}{i}{j}{k}}{\reS{\Curve{}}{O}}\)
\reS{{\dot{\Curve{}}}_j}{O}
\reS{{\dot{\Curve{}}}_k}{O}=0,
\end{equation}
where, $\Curve{i}:=\cmp{\dualEucbase{m}{i}}{\cmp{\phi}{\Curve{}}}$.\\
\caution
These equations are called a $\quotl$system of local equations of the geodesic $\Curve{}$ on $\Man{}$
(relative to $\connection{}$) with respect to the chart $\opair{\U}{\phi}$$\quotr$.
\item[\myitem{3.}]
Given a family $\seta{\opair{\U_\alpha}{\phi_\alpha}}_{\alpha\in\indexset{}}$ of charts of $\Man{}$ such that
$\funcimage{\Curve{}}\subseteq\Union{\alpha}{\indexset{}}{\U_\alpha}$ and
$\U_\alpha\cap\funcimage{\Curve{}}\neq\emptyset$, and given a family of
open sub-intervals $\seta{O_\alpha}_{\alpha\in\indexset{}}$ of the interval $\interval{}$ such that
$\defSet{\func{\Curve{}}{t}}{t\in O_\alpha}\subseteq\U_\alpha$, for every $\alpha\in\indexset{}$,
%%%%%%
\begin{align}
&\Foreach{i}{\seta{\suc{1}{m}}}\Foreach{t}{O_\alpha}\cr
&\func{{\ddot{\Curve{}}}_{i}}{t}+\sum_{k=1}^{m}\sum_{j=1}^{m}
\func{\Christoffel{\phi_\alpha}{i}{j}{k}}{\func{\Curve{}}{t}}\func{{\dot{\Curve{}}}_{j}}{t}
\func{{\dot{\Curve{}}}_{k}}{t}=0,
\end{align}
or equivalently,
\begin{equation}
\reS{{\ddot{\Curve{}}}_i}{O_\alpha}+\sum_{k=1}^{m}\sum_{j=1}^{m}
\(\cmp{\Christoffel{\phi_\alpha}{i}{j}{k}}{\reS{\Curve{}}{O}}\)\reS{{\dot{\Curve{}}}_j}{O_\alpha}
\reS{{\dot{\Curve{}}}_k}{O_\alpha}=0,
\end{equation}
where, $\Curve{i}:=\cmp{\dualEucbase{m}{i}}{\cmp{\phi}{\Curve{}}}$.
\end{itemize}
\proof
It is known that given a chart $\opair{\U}{\phi}$ of $\Man{}$ that intersects the image of
$\Curve{}$, and a sub-interval $O$ of $\interval{}$ such that $\defSet{\func{\Curve{}}{t}}{t\in O}\subseteq\U$,
$\displaystyle\reS{\Curve{}}{O}=\sum_{i=1}^m\reS{\Curve{i}}{O}\(\cmp{\localframevecf{i}}{\reS{\Curve{}}{O}}\)$,
where $\mtuple{\localframevecf{1}}{\localframevecf{m}}$ denotes the system of local frame fields
associated with the chart $\opair{\U}{\phi}$ and
$\Curve{i}:=\cmp{\dualEucbase{m}{i}}{\cmp{\phi}{\Curve{}}}$. Thus, according to
\refcor{corparallelvfalongcurveexpresion} and \refdef{defgeodesiconmanifold}, it is obvious.
\endthm
%%%%%%%%%%%%%%%%%%%%%%%%%%%%%%%%%%%%%%%%%%%%%%%%%%%%%%%%%%%%%%%%%%%%%%%%%%%%%%%%%%%%%%%%%%%%%%%%%%%%%%%%%%%%%%%%%%%%%%%%%%%%%%%%
\textit{
The following lemma concerns the
existance and uniqueness of the solution of a system of second-order ordinary differential equations,
the proof of which is not of concern here, and we accept its truth on faith.
This lemma will be used in the sequel.
}
%%%%%%%%%%%%%%%%%%%%%%%%%%%%%%%%%%%%%%%%%%%%%%%%%%%%%%%%%%%%%%%%%%%%%%%%%%%%%%%%%%%%%%%%%%%%%%%%%%%%%%%%%%%%%%%%%%%%%%%%%%%%%%%%
\lemma\label{lemmaexistanceanduniquenessofgeodesicODEs}
Let $N$ be a positive integer,
and let $\U$ be a non-empty open subset of $\R^{N}$.
Let $\function{\lambda_{jk}^{i}}{\U}{\R}$
be a real-valued smooth map on $\U$, for every $i$, $j$, and $k$ in $\seta{\suc{1}{N}}$.
Let $\mtuple{c_1}{c_N}\in\R^N$ and $\mtuple{a_1}{a_N}\in\R^N$. Let $t_0\in\R$.
There exists an open interval $\interval{}$ of $\R$ containing $t_0$ such that there exists a unique
$N$-tuple $\mtuple{\function{f_1}{\interval{}}{\U}}{\function{f_N}{\interval{}}{\U}}$ of real-valued smooth
functions on $\interval{}$ satisfying the system
of second-order ordinary differential equations,
\begin{equation}
\Foreach{i}{\seta{\suc{1}{N}}}
{\ddot{f}}_i+\sum_{k=1}^{N}\sum_{j=1}^{N}
\(\cmp{\lambda_{jk}^{i}}{\sum_{\mu=1}^{N}f_{\mu}\dualEucbase{N}{\mu}}\){\dot{f}}_j{\dot{f}}_k=0,
\end{equation}
with the initial-value conditions,
\begin{align}
\mtuple{\func{f_1}{t_0}}{\func{f_N}{t_0}}&=\mtuple{c_1}{c_N},\cr
\mtuple{\func{\dot{f}_1}{t_0}}{\func{\dot{f}_N}{t_0}}&=\mtuple{a_1}{a_N}.
\end{align}
where, ${\dot{f}}_{i}$ denotes the derived function of $f_i$, and
$\mtuple{\dualEucbase{N}{1}}{\dualEucbase{N}{N}}$ denotes the dual of the standard ordered-basis of $\R^N$.
\endlem
%%%%%%%%%%%%%%%%%%%%%%%%%%%%%%%%%%%%%%%%%%%%%%%%%%%%%%%%%%%%%%%%%%%%%%%%%%%%%%%%%%%%%%%%%%%%%%%%%%%%%%%%%%%%%%%%%%%%%%%%%%%%%%%%
\lemma\label{lemexsitanceofgeodesics}
Let $\point\in\Man{}$, $v\in\tanspace{\point}{\Man{}}$, and $t_0\in\R$. There exists an open interval $\interval{}$ of $\R$
containing $t_0$ such that there exists a unique geodesic $\function{\Curve{}}{\interval{}}{\Man{}}$
(relative to $\connection{}$) satisfying the condition $\func{\dot{\Curve{}}}{t_0}=v$.
\proof
Let $\opair{\U}{\phi}$ be a chart of $\Man{}$ around $\point$. 
Let $\mtuple{\localframevecf{1}}{\localframevecf{m}}$ be the system of local frame fields associated with
$\opair{\U}{\phi}$. 
According to \reflem{lemmaexistanceanduniquenessofgeodesicODEs},
there exists an open interval $\interval{}$ of $\R$ containing $t_0$ such that
there exists a unique $m$-tuple of real-valued smooth functions
$\mtuple{\function{x_1}{\interval{}}{\R}}{\function{x_m}{\interval{}}{\R}}$
being the solution of the system of second-order ODE-s,
\begin{equation}
\Foreach{i}{\seta{\suc{1}{N}}}
{\ddot{x}}_i+\sum_{k=1}^{N}\sum_{j=1}^{N}
\(\cmp{\Christoffel{\phi}{i}{j}{k}}{\cmp{\finv{\phi}}{\sum_{\mu=1}^{N}x_{\mu}\dualEucbase{N}{\mu}}}\){\dot{x}}_j{\dot{x}}_k=0,
\end{equation}
with the initial conditions,
\begin{align}
\mtuple{\func{x_1}{t_0}}{\func{x_m}{t_0}}&=\mtuple{\func{\phi^1}{\point}}
{\func{\phi^m}{\point}},\cr
\mtuple{\func{\dot{x}_1}{t_0}}{\func{\dot{x}_m}{t_0}}&=\mtuple{v^1}{v^m},
%\mtuple{\func{\dualEucbase{N}{1}}{\func{\tanspaceiso{\point}{\Man{}}{\phi}}{v}}}
%{\func{\dualEucbase{N}{N}}{\func{\tanspaceiso{\point}{\Man{}}{\phi}}{v}}}.
\end{align}
where, $\phi^{i}:=\cmp{\dualEucbase{N}{i}}{\phi}$ for every $i$, and
$\displaystyle v=\sum_{i=1}^{m}v_i\func{\localframevecf{i}}{\point}$.\\
\refthm{thmgeodesiclocalexpresion} and \reflem{lemmaexistanceanduniquenessofgeodesicODEs}
imply that $\displaystyle\Curve{}:=
\cmp{\finv{\phi}}{\sum_{\mu=1}^{N}x_{\mu}\dualEucbase{N}{\mu}}$
is the unique geodesic with domain $\interval{}$ sch that
$\func{\Curve{}}{t_0}=\point$ and $\func{\dot{\Curve{}}}{t_0}=v$.
\endlem
%%%%%%%%%%%%%%%%%%%%%%%%%%%%%%%%%%%%%%%%%%%%%%%%%%%%%%%%%%%%%%%%%%%%%%%%%%%%%%%%%%%%%%%%%%%%%%%%%%%%%%%%%%%%%%%%%%%%%%%%%%%%%%%%
\lemma\label{lemrestrictionofgeodesicsaregeodesics}
Let $\function{\Curve{}}{\interval{}}{\Man{}}$ be a geodesic on $\Man{}$ (relative to $\connection{}$)
for some open interval $\interval{}$ of $\R$. For every open sub-interval $O$ of $\interval{}$,
$\reS{\Curve{}}{O}$ is again a geodesic on $\Man{}$ (relative to $\connection{}$).
\proof
It is trivial, based on \refthm{thmgeodesiclocalexpresion}.
\endlem
%%%%%%%%%%%%%%%%%%%%%%%%%%%%%%%%%%%%%%%%%%%%%%%%%%%%%%%%%%%%%%%%%%%%%%%%%%%%%%%%%%%%%%%%%%%%%%%%%%%%%%%%%%%%%%%%%%%%%%%%%%%%%%%%
\lemma\label{lemcoincidenceofgeodesics}
Let $\function{\Curve{1}}{\interval{}}{\Man{}}$ and $\function{\Curve{2}}{\interval{}}{\Man{}}$ be a pair
of geodesics on $\Man{}$ (relative to $\connection{}$) for some open interval $\interval{}$ of $\R$, and let $t_0\in\interval{}$.
If $\func{\dot{\Curve{1}}}{t_0}=\func{\dot{\Curve{2}}}{t_0}$, then $\Curve{1}=\Curve{2}$.
\proof
Assume that $\Curve{1}\neq\Curve{2}$. Then there exists $\tau\in\interval{}$ such that $\func{\Curve{1}}{\tau}\neq
\func{\Curve{2}}{\tau}$. Trivially, either $\tau>t_0$ or $\tau<t_0$; we consider the former case, since the other one
is not different in essence. Define $S:=\defset{t}{\interval{}}{t\geq t_0,~\func{\Curve{1}}{t}\neq\func{\Curve{2}}{t}}$.
Hence, considering that $\tau\in S$, $S$ is a non-empty subset of $\R$ bounded from below, and hence has an infimum.
Define $s:=\func{\inf}{S}$.\\
If $s=t_0$, then trivially $\func{\dot{\Curve{1}}}{s}=\func{\dot{\Curve{2}}}{s}$; assume otherwise that $s>t_0$.
Then for every $t\in\opair{t_0}{s}$,
$\func{\dot{\Curve{1}}}{t}=\func{\dot{\Curve{2}}}{t}$, and thus according to the continuity of
$\dot{\Curve{1}}$ and $\dot{\Curve{2}}$, $\func{\dot{\Curve{1}}}{s}=\func{\dot{\Curve{2}}}{s}$.
According to \reflem{lemexsitanceofgeodesics}, there exists an open interval $\interval{0}$ of $\R$ containing $s$, such that
there exists a unique geodesic $\function{\Curve{}}{\interval{0}}{\Man{}}$ on $\Man{}$ satisfying
$\func{\dot{\Curve{}}}{s}=\func{\dot{\Curve{1}}}{s}=\func{\dot{\Curve{2}}}{s}$. In addition, according to
\reflem{lemrestrictionofgeodesicsaregeodesics}, each $\reS{\Curve{1}}{\interval{0}}$ and $\reS{\Curve{2}}{\interval{0}}$
must be a geodesics on $\Man{}$ with the same initial condition. Therefore,
$\reS{\Curve{1}}{\interval{0}}=\reS{\Curve{2}}{\interval{0}}$, which implies that $\func{\inf}{S}>s$,
which is a contradiction.\\
Therefore, the earliest assumption cannot be valid, and thus $\Curve{1}=\Curve{2}$.
\endlem
%%%%%%%%%%%%%%%%%%%%%%%%%%%%%%%%%%%%%%%%%%%%%%%%%%%%%%%%%%%%%%%%%%%%%%%%%%%%%%%%%%%%%%%%%%%%%%%%%%%%%%%%%%%%%%%%%%%%%%%%%%%%%%%%
\definition\label{defmaximalgeodesic}
Let $\function{\Curve{}}{\interval{}}{\Man{}}$ be a geodesic on $\Man{}$ (relative to $\connection{}$)
for some open interval $\interval{}$ of $\R$. $\Curve{}$ is called a $\quotl$maximal geodesic on $\Man{}$
relative to the connection $\connection{}$$\quotr$ if there can not be found a geodesic
$\function{\Curve{1}}{\interval{1}}{\Man{}}$ on $\Man{}$ such that $\interval{}\subseteq\interval{1}$
and $\reS{\Curve{1}}{\interval{}}=\Curve{}$.
\endef
%%%%%%%%%%%%%%%%%%%%%%%%%%%%%%%%%%%%%%%%%%%%%%%%%%%%%%%%%%%%%%%%%%%%%%%%%%%%%%%%%%%%%%%%%%%%%%%%%%%%%%%%%%%%%%%%%%%%%%%%%%%%%%%%
\theorem
Let $v\in\tanbun{\Man{}}$. There is a unique maximal geodesic $\function{\Curve{}}{\interval{}}{\Man{}}$
on $\Man{}$ (relative to $\connection{}$) such that $0\in\interval{}$ and $\func{\dot{\Curve{}}}{0}=v$.
\proof
Let $G_{v}$ denote the set of all geodesics $\function{\Curve{}}{\interval{}}{\Man{}}$ such that $0\in\interval{}$
and $\func{\dot{\Curve{}}}{0}=v$. According to \reflem{lemrestrictionofgeodesicsaregeodesics}
and \reflem{lemcoincidenceofgeodesics}, for any $\function{\Curve{1}}{\interval{1}}{\Man{}}$
and $\function{\Curve{2}}{\interval{2}}{\Man{}}$, $\reS{\Curve{1}}{\interval{1}\cap\interval{2}}=
\reS{\Curve{2}}{\interval{1}\cap\interval{2}}$. So, by defining
\begin{equation}
\interval{0}:=\Union{\Curve{}}{G_v}{\domain{\Curve{}}},
\end{equation}
it is permissible to define a map $\function{\Curve{0}}{\interval{0}}{\Man{}}$ in the following way.
For any $\function{\Curve{}}{\interval{}}{\Man{}}$ in $G_v$, and for any $t\in\interval{}$,
\begin{equation}
\func{\Curve{0}}{t}\eqdef\func{\Curve{}}{t}.
\end{equation}
Clearly, according to \reflem{lemcoincidenceofgeodesics} and \refdef{defmaximalgeodesic},
$\Curve{}$ is a maximal geodesic on $\Man{}$ (relative to $\connection{}$)
such that $0\in\interval{}$ and $\func{\dot{\Curve{}}}{0}=v$, and it is the only such curve.
\endthm
%%%%%%%%%%%%%%%%%%%%%%%%%%%%%%%%%%%%%%%%%%%%%%%%%%%%%%%%%%%%%%%%%%%%%%%%%%%%%%%%%%%%%%%%%%%%%%%%%%%%%%%%%%%%%%%%%%%%%%%%%%%%%%%%
\definition
We will denote by $\maxgeodesic{\Man{}}{\connection{}}{v}$ the unique maximal geodesic
$\function{\Curve{}}{\interval{}}{\Man{}}$ on $\Man{}$ relative to $\connection{}$
such that $0\in\interval{}$ and $\func{\dot{\Curve{}}}{0}=v$, which will be referred to as the
$\quotl$maximal geodesic on $\Man{}$ relative to the connection $\connection{}$ with the initial value $v$$\quotr$.
\endef
%%%%%%%%%%%%%%%%%%%%%%%%%%%%%%%%%%%%%%%%%%%%%%%%%%%%%%%%%%%%%%%%%%%%%%%%%%%%%%%%%%%%%%%%%%%%%%%%%%%%%%%%%%%%%%%%%%%%%%%%%%%%%%%%
\definition\label{defgeodesicallycompletemanifold}
Given a point $\point$ of $\Man{}$, if for every $v\in\tanspace{\point}{\Man{}}$
the domain of the maximal geodesic $\maxgeodesic{\Man{}}{\connection{}}{v}$ equals $\R$,
then the manifold $\Man{}$ with the connection $\connection{}$ on it is said to be
$\quotl$geodesically complete at $\point$$\quotr$.
The manifold $\Man{}$ with the connection $\connection{}$ on it is said to be
$\quotl$geodesically complete$\quotr$ if it is geodesically complete at every point $\point$ of $\Man{}$.
\endef
%%%%%%%%%%%%%%%%%%%%%%%%%%%%%%%%%%%%%%%%%%%%%%%%%%%%%%%%%%%%%%%%%%%%%%%%%%%%%%%%%%%%%%%%%%%%%%%%%%%%%%%%%%%%%%%%%%%%%%%%%%%%%%%%
%%%%%%%%%%%%%%%%%%%%%%%%%%%%%%%%%%%%%%%%%%%%%%%%%%%%%%%%%%%%%%%%%%%%%%%%%%%%%%%%%%%%%%%%%%%%%%%%%%%%%%%%%%%%%%%%%%%%%%%%%%%%%%%%
%%%%%%%%%%%%%%%%%%%%%%%%%%%%%%%%%%%%%%%%%%%%%%%%%%%%%%%%%%%%%%%%%%%%%%%%%%%%%%%%%%%%%%%%%%%%%%%%%%%%%%%%%%%%%%%%%%%%%%%%%%%%%%%%
%%%%%%%%%%%%%%%%%%%%%%%%%%%%%%%%%%%%%%%%%%%%%%%%%%%%%%%%%%%%%%%%%%%%%%%%%%%%%%%%%%%%%%%%%%%%%%%%%%%%%%%%%%%%%%%%%%%%%%%%%%%%%%%%
%%%%%%%%%%%%%%%%%%%%%%%%%%%%%%%%%%%%%%%%%%%%%%%%%%%%%%%%%%%%%%%%%%%%%%%%%%%%%%%%%%%%%%%%%%%%%%%%%%%%%%%%%%%%%%%%%%%%%%%%%%%%%%%%
%%%%%%%%%%%%%%%%%%%%%%%%%%%%%%%%%%%%%%%%%%%%%%%%%%%%%%%%%%%%%%%%%%%%%%%%%%%%%%%%%%%%%%%%%%%%%%%%%%%%%%%%%%%%%%%%%%%%%%%%%%%%%%%%
%%%%%%%%%%%%%%%%%%%%%%%%%%%%%%%%%%%%%%%%%%%%%%%%%%%%%%%%%%%%%%%%%%%%%%%%%%%%%%%%%%%%%%%%%%%%%%%%%%%%%%%%%%%%%%%%%%%%%%%%%%%%%%%%
%%%%%%%%%%%%%%%%%%%%%%%%%%%%%%%%%%%%%%%%%%%%%%%%%%%%%%%%%%%%%%%%%%%%%%%%%%%%%%%%%%%%%%%%%%%%%%%%%%%%%%%%%%%%%%%%%%%%%%%%%%%%%%%%
%%%%%%%%%%%%%%%%%%%%%%%%%%%%%%%%%%%%%%%%%%%%%%%%%%%%%%%%%%%%%%%%%%%%%%%%%%%%%%%%%%%%%%%%%%%%%%%%%%%%%%%%%%%%%%%%%%%%%%%%%%%%%%%%
%%%%%%%%%%%%%%%%%%%%%%%%%%%%%%%%%%%%%%%%%%%%%%%%%%%%%%%%%%%%%%%%%%%%%%%%%%%%%%%%%%%%%%%%%%%%%%%%%%%%%%%%%%%%%%%%%%%%%%%%%%%%%%%%
%%%%%%%%%%%%%%%%%%%%%%%%%%%%%%%%%%%%%%%%%%%%%%%%%%%%%%%%%%%%%%%%%%%%%%%%%%%%%%%%%%%%%%%%%%%%%%%%%%%%%%%%%%%%%%%%%%%%%%%%%%%%%%%%
\section{Vector Fields Along One-Parameter Family of Curves}
%%%%%%%%%%%%%%%%%%%%%%%%%%%%%%%%%%%%%%%%%%%%%%%%%%%%%%%%%%%%%%%%%%%%%%%%%%%%%%%%%%%%%%%%%%%%%%%%%%%%%%%%%%%%%%%%%%%%%%%%%%%%%%%%
\fixed
$\Man{}$ is fixed as a manifold with dimension $m$, and $\connection{}$ as an affine connection on $\Man{}$.
\endfixed
%%%%%%%%%%%%%%%%%%%%%%%%%%%%%%%%%%%%%%%%%%%%%%%%%%%%%%%%%%%%%%%%%%%%%%%%%%%%%%%%%%%%%%%%%%%%%%%%%%%%%%%%%%%%%%%%%%%%%%%%%%%%%%%%
\definition
Let each $\interval{1}$ and $\interval{2}$ be an open interval of $\R$. Let
$\function{\pCurve{}}{\Cprod{\interval{1}}{\interval{2}}}{\Man{}}$ be a smooth map from $\Cprod{\interval{1}}{\interval{2}}$
(endowed with its canonical differentiable structure) to $\Man{}$. Then $\pCurve{}$ is called a $\quotl$one-parameter
family of smooth curves on $\Man{}$$\quotr$.\\
Furthermore, let $\function{\valongc{}}{\Cprod{\interval{1}}{\interval{2}}}{\Tanbun{\Man{}}}$ be a smooth map from
$\Cprod{\interval{1}}{\interval{2}}$ to the tangent-bundle $\Tanbun{\Man{}}$. Then, $\valongc{}$ is called a
$\quotl$smooth vector field along the one-parameter family of smooth curves $\pCurve{}$ on $\Man{}$$\quotr$
if
\begin{equation}
\Foreach{\opair{t}{s}}{\Cprod{\interval{1}}{\interval{2}}}
\func{\valongc{}}{\binary{t}{s}}\in\tanspace{\func{\pCurve{}}{\binary{t}{s}}}{\Man{}},
\end{equation}
or equivalently,
\begin{equation}
\cmp{\basep{\Man{}}}{\valongc{}}=\pCurve{},
\end{equation}
where $\basep{\Man{}}$ denotes the projection map of the tangent-bundle $\Tanbun{\Man{}}$.
We denote by $\valongcurves{\Man{}}{\pCurve{}}$ the set of all smooth vector fields along the one-parameter family
of smooth curves $\pCurve{}$.
\endef
%%%%%%%%%%%%%%%%%%%%%%%%%%%%%%%%%%%%%%%%%%%%%%%%%%%%%%%%%%%%%%%%%%%%%%%%%%%%%%%%%%%%%%%%%%%%%%%%%%%%%%%%%%%%%%%%%%%%%%%%%%%%%%%%
\definition\label{deflogitudinalandtransversecurves}
Let $\function{\pCurve{}}{\Cprod{\interval{1}}{\interval{2}}}{\Man{}}$ be a one-parameter family of smooth curves on $\Man{}$.
For every $t\in\interval{2}$, we define the map $\function{\pcurve{\pCurve{}}{1}{t}}{\interval{1}}{\Man{}}$ as,
\begin{equation}
\Foreach{s}{\interval{1}}
\func{\pcurve{\pCurve{}}{1}{t}}{s}\eqdef\func{\pCurve{}}{\binary{s}{t}},
\end{equation}
and for every $s\in\interval{1}$ we define the map $\function{\pcurve{\pCurve{}}{2}{s}}{\interval{2}}{\Man{}}$ as,
\begin{equation}
\Foreach{t}{\interval{2}}
\func{\pcurve{\pCurve{}}{2}{s}}{t}\eqdef\func{\pCurve{}}{\binary{s}{t}}.
\end{equation}
$\pcurve{\pCurve{}}{2}{s}$ is referred to as the $\quotl$longitudinal curve of $\pCurve{}$ with parameter $s$$\quotr$,
and $\pcurve{\pCurve{}}{1}{t}$ is referred to as the $\quotl$transverse curve of $\pCurve{}$ with parameter $t$$\quotr$.
\endef
%%%%%%%%%%%%%%%%%%%%%%%%%%%%%%%%%%%%%%%%%%%%%%%%%%%%%%%%%%%%%%%%%%%%%%%%%%%%%%%%%%%%%%%%%%%%%%%%%%%%%%%%%%%%%%%%%%%%%%%%%%%%%%%%
\proposition
Let $\function{\pCurve{}}{\Cprod{\interval{1}}{\interval{2}}}{\Man{}}$ be a one-parameter family of smooth curves on $\Man{}$.
Let $\valongc{}$ be a vector field along $\pCurve{}$.\\
For every $t\in\interval{2}$, $\pcurve{\pCurve{}}{1}{t}$ is a smooth curve on $\Man{}$.\\
For every $s\in\interval{1}$, $\pcurve{\pCurve{}}{2}{s}$ is a smooth curve on $\Man{}$.\\
For every $t\in\interval{2}$, $\reS{\valongc{}}{\Cprod{\interval{1}}{\seta{t}}}$ is a vector field along the curve
$\pcurve{\pCurve{}}{1}{t}$, that is an element of $\valongcurves{\Man{}}{\pcurve{\pCurve{}}{1}{t}}$.\\
For every $s\in\interval{1}$, $\reS{\valongc{}}{\Cprod{\seta{s}}{\interval{2}}}$ is a vector field along the curve
$\pcurve{\pCurve{}}{2}{s}$, that is an element of $\valongcurves{\Man{}}{\pcurve{\pCurve{}}{2}{s}}$.
\proof
It is trivial.
\endpro
%%%%%%%%%%%%%%%%%%%%%%%%%%%%%%%%%%%%%%%%%%%%%%%%%%%%%%%%%%%%%%%%%%%%%%%%%%%%%%%%%%%%%%%%%%%%%%%%%%%%%%%%%%%%%%%%%%%%%%%%%%%%%%%%
\remark
In the proposition above, we consider $\reS{\valongc{}}{\Cprod{\interval{1}}{\seta{t}}}$ as a map from $\interval{1}$ to $\Man{}$
sending $s\in\interval{1}$ to $\func{\valongc{}}{\binary{s}{t}}$. Similarly, we consider
$\reS{\valongc{}}{\Cprod{\seta{s}}{\interval{2}}}$ as a map from $\interval{2}$ to $\Man{}$
sending $t\in\interval{2}$ to $\func{\valongc{}}{\binary{s}{t}}$.
\endremark
%%%%%%%%%%%%%%%%%%%%%%%%%%%%%%%%%%%%%%%%%%%%%%%%%%%%%%%%%%%%%%%%%%%%%%%%%%%%%%%%%%%%%%%%%%%%%%%%%%%%%%%%%%%%%%%%%%%%%%%%%%%%%%%%
\lemma
Let $\function{\pCurve{}}{\Cprod{\interval{1}}{\interval{2}}}{\Man{}}$ be a one-parameter family of smooth curves on $\Man{}$.
For every $s\in\interval{1}$ and every $t\in\interval{2}$,
\begin{equation}
\func{\[\func{\covder{\connection{}}{\pcurve{\pCurve{}}{2}{s}}}{\dot{\pcurve{\pCurve{}}{2}{s}}}\]}{t}=
\func{\[\func{\covder{\connection{}}{\pcurve{\pCurve{}}{1}{t}}}{\dot{\pcurve{\pCurve{}}{1}{t}}}\]}{s},
\end{equation}
where $\dot{\pcurve{\pCurve{}}{2}{s}}$ and $\dot{\pcurve{\pCurve{}}{1}{t}}$ denote the velocity maps of the curves
$\pcurve{\pCurve{}}{2}{s}$ and $\pcurve{\pCurve{}}{1}{t}$, respectively.
\proof
It is obvious according to \refthm{thmlocalexpressionofcovariantderivativeofvectorfieldsalongcurves}
and \refdef{deflogitudinalandtransversecurves}.
\endlem
%%%%%%%%%%%%%%%%%%%%%%%%%%%%%%%%%%%%%%%%%%%%%%%%%%%%%%%%%%%%%%%%%%%%%%%%%%%%%%%%%%%%%%%%%%%%%%%%%%%%%%%%%%%%%%%%%%%%%%%%%%%%%%%%
\definition\label{defcovariantderivativesofvectorfieldsalongoneparameterfamiliesofcurves}
Let $\function{\pCurve{}}{\Cprod{\interval{1}}{\interval{2}}}{\Man{}}$ be a one-parameter family of smooth curves on $\Man{}$.
We define the maps $\function{\pcovder{\connection{}}{\pCurve{}}{1}}{\valongcurves{\Man{}}{\pCurve{}}}
{\valongcurves{\Man{}}{\pCurve{}}}$ and $\function{\pcovder{\connection{}}{\pCurve{}}{2}}{\valongcurves{\Man{}}{\pCurve{}}}
{\valongcurves{\Man{}}{\pCurve{}}}$ as,
\begin{align}
\Foreach{\valongc{}}{\valongcurves{\Man{}}{\pCurve{}}}
\Foreach{\opair{s}{t}}{\Cprod{\interval{1}}{\interval{2}}}
\func{\[\func{\pcovder{\connection{}}{\pCurve{}}{1}}{\valongc{}}\]}{\binary{s}{t}}\eqdef
\func{\[\func{\covder{\connection{}}{\pcurve{\pCurve{}}{1}{t}}}{\reS{\valongc{}}{\Cprod{\interval{1}}{\seta{t}}}}\]}{s},
\end{align}
and,
\begin{align}
\Foreach{\valongc{}}{\valongcurves{\Man{}}{\pCurve{}}}
\Foreach{\opair{s}{t}}{\Cprod{\interval{1}}{\interval{2}}}
\func{\[\func{\pcovder{\connection{}}{\pCurve{}}{2}}{\valongc{}}\]}{\binary{s}{t}}\eqdef
\func{\[\func{\covder{\connection{}}{\pcurve{\pCurve{}}{2}{s}}}{\reS{\valongc{}}{\Cprod{\seta{s}}{\interval{2}}}}\]}{t}.
\end{align}
$\pcovder{\connection{}}{\pCurve{}}{1}$ is called the $\quotl$longitudinal covariant derivative operator of
one-parameter family of curves $\pCurve{}$ on $\Man{}$, relative to the connection $\connection{}$$\quotr$.
$\pcovder{\connection{}}{\pCurve{}}{2}$ is called the $\quotl$transverse covariant derivative operator of
one-parameter family of curves $\pCurve{}$ on $\Man{}$, relative to the connection $\connection{}$$\quotr$.
\endef
%%%%%%%%%%%%%%%%%%%%%%%%%%%%%%%%%%%%%%%%%%%%%%%%%%%%%%%%%%%%%%%%%%%%%%%%%%%%%%%%%%%%%%%%%%%%%%%%%%%%%%%%%%%%%%%%%%%%%%%%%%%%%%%%
%%%%%%%%%%%%%%%%%%%%%%%%%%%%%%%%%%%%%%%%%%%%%%%%%%%%%%%%%%%%%%%%%%%%%%%%%%%%%%%%%%%%%%%%%%%%%%%%%%%%%%%%%%%%%%%%%%%%%%%%%%%%%%%%
%%%%%%%%%%%%%%%%%%%%%%%%%%%%%%%%%%%%%%%%%%%%%%%%%%%%%%%%%%%%%%%%%%%%%%%%%%%%%%%%%%%%%%%%%%%%%%%%%%%%%%%%%%%%%%%%%%%%%%%%%%%%%%%%
%%%%%%%%%%%%%%%%%%%%%%%%%%%%%%%%%%%%%%%%%%%%%%%%%%%%%%%%%%%%%%%%%%%%%%%%%%%%%%%%%%%%%%%%%%%%%%%%%%%%%%%%%%%%%%%%%%%%%%%%%%%%%%%%
%%%%%%%%%%%%%%%%%%%%%%%%%%%%%%%%%%%%%%%%%%%%%%%%%%%%%%%%%%%%%%%%%%%%%%%%%%%%%%%%%%%%%%%%%%%%%%%%%%%%%%%%%%%%%%%%%%%%%%%%%%%%%%%%
%%%%%%%%%%%%%%%%%%%%%%%%%%%%%%%%%%%%%%%%%%%%%%%%%%%%%%%%%%%%%%%%%%%%%%%%%%%%%%%%%%%%%%%%%%%%%%%%%%%%%%%%%%%%%%%%%%%%%%%%%%%%%%%%
%%%%%%%%%%%%%%%%%%%%%%%%%%%%%%%%%%%%%%%%%%%%%%%%%%%%%%%%%%%%%%%%%%%%%%%%%%%%%%%%%%%%%%%%%%%%%%%%%%%%%%%%%%%%%%%%%%%%%%%%%%%%%%%%
%%%%%%%%%%%%%%%%%%%%%%%%%%%%%%%%%%%%%%%%%%%%%%%%%%%%%%%%%%%%%%%%%%%%%%%%%%%%%%%%%%%%%%%%%%%%%%%%%%%%%%%%%%%%%%%%%%%%%%%%%%%%%%%%
%%%%%%%%%%%%%%%%%%%%%%%%%%%%%%%%%%%%%%%%%%%%%%%%%%%%%%%%%%%%%%%%%%%%%%%%%%%%%%%%%%%%%%%%%%%%%%%%%%%%%%%%%%%%%%%%%%%%%%%%%%%%%%%%
%%%%%%%%%%%%%%%%%%%%%%%%%%%%%%%%%%%%%%%%%%%%%%%%%%%%%%%%%%%%%%%%%%%%%%%%%%%%%%%%%%%%%%%%%%%%%%%%%%%%%%%%%%%%%%%%%%%%%%%%%%%%%%%%
%%%%%%%%%%%%%%%%%%%%%%%%%%%%%%%%%%%%%%%%%%%%%%%%%%%%%%%%%%%%%%%%%%%%%%%%%%%%%%%%%%%%%%%%%%%%%%%%%%%%%%%%%%%%%%%%%%%%%%%%%%%%%%%%
%%%%%%%%%%%%%%%%%%%%%%%%%%%%%%%%%%%%%%%%%%%%%%%%%%%%%%%%%%%%%%%%%%%%%%%%%%%%%%%%%%%%%%%%%%%%%%%%%%%%%%%%%%%%%%%%%%%%%%%%%%%%%%%%
\section{Covariant Differentiation of Tensor Fields}
%%%%%%%%%%%%%%%%%%%%%%%%%%%%%%%%%%%%%%%%%%%%%%%%%%%%%%%%%%%%%%%%%%%%%%%%%%%%%%%%%%%%%%%%%%%%%%%%%%%%%%%%%%%%%%%%%%%%%%%%%%%%%%%%
\textit{
For technical and notational reasons, we prefer to briefly review a dual approach for defining smooth tensor fields
on the tangent bundle of a manifold. Although this approach can generally be introduced
on smooth real vector bundles of finite ranks, we do it in the perticular framework of tangent-bundles of finite-dimensional
manifolds, because only this particular type of vector bundles will be of concern in the sequel.\\
Let $\Man{}$ be a real manifold with dimension $m$. Let $r$ and $s$ be non-negative integers with $r+s>0$.
An $\opair{r}{s}$ smooth tensor field $T$ on $\Tanbun{\Man{}}$ as introduced in \cite{ShahiVB}, is a smooth map from $\Man{}$
to the $\opair{r}{s}$-tensor-bundle of the tangent-bundle $\Tanbun{\Man{}}$, that is $\vbtensorbundle{r}{s}{\Tanbun{\Man{}}}$;
so, for every point $\point$ of $\Man{}$, $\func{T}{\point}\in\Tensors{r}{s}{\Tanspace{\point}{\Man{}}}$, that is
$\function{\func{T}{\point}}{\Cprod{\(\tanspace{\point}{\Man{}}\)^{r}}{\(\dualV{\tanspace{\point}{\Man{}}}\)^{s}}}{\R}$ is a
multi-linear map, considering the canonical linear structure of $\tanspace{\point}{\Man{}}$.
The set of all $\opair{r}{s}$ smooth tensor fields on $\Tanbun{\Man{}}$ is denoted by $\TF{r}{s}{\Tanbun{\Man{}}}$.\\
It is known that the set of all smooth vector fields on $\Man{}$, that is $\vectorfields{\Man{}}$,
and the set of all co-vector fields on (sections of the dual-bundle of $\Tanbun{\Man{}}$) $\Man{}$,
that is $\oneforms{\Man{}}$, each has a canonical $\smoothmaps{\Man{}}$-module (a module over the commutative ring of
real-valued smooth maps on $\Man{}$) structur. We denote by $\TTF{r}{s}{\vectorfields{\Man{}}}$ the set of all
$\smoothmaps{\Man{}}$-multilinear maps
$\function{\alpha}{\Cprod{\(\vectorfields{\Man{}}\)^r}{\(\oneforms{\Man{}}\)^s}}{\smoothmaps{\Man{}}}$.
There is natural (one-one and onto) correspondence between $\TF{r}{s}{\Tanbun{\Man{}}}$ and $\TTF{r}{s}{\vectorfields{\Man{}}}$.
We denote this correspondence as $\function{\Upsilon}{\TF{r}{s}{\Tanbun{\Man{}}}}{\TTF{r}{s}{\vectorfields{\Man{}}}}$,
and for every tensor field $T\in\TF{r}{s}{\Tanbun{\Man{}}}$ we will denote its corresponded element $\func{\Upsilon}{T}$
in $\TTF{r}{s}{\vectorfields{\Man{}}}$ simply by $\TFequiv{T}$, when the underying manifold is clear enough.
In addition, for every $\mathcal{T}\in\TTF{r}{s}{\vectorfields{\Man{}}}$, we denote its corresponded element
$\func{\finv{\Upsilon}}{\mathcal{T}}$
in $\TF{r}{s}{\Tanbun{\Man{}}}$ with  $\TTFequiv{\mathcal{T}}$.
As a straightforward consequence of the way in which the correpondence $\Upsilon$ is defined,
for every $T\in\TF{r}{s}{\Tanbun{\Man{}}}$, every $\suc{\avecf{1}}{\avecf{r}}\in\vectorfields{\Man{}}$,
every $\suc{\aoneform{1}}{\aoneform{s}}\in\oneforms{\Man{}}$, and every $\point\in\Man{}$,
\begin{equation}
\func{\[\func{\TFequiv{T}}{\binary{\suc{\avecf{1}}{\avecf{r}}}{\suc{\aoneform{1}}{\aoneform{s}}}}\]}{\point}=
\func{\[\func{T}{\point}\]}{\binary{\suc{\func{\avecf{1}}{\point}}{\func{\avecf{r}}{\point}}}
{\suc{\func{\aoneform{1}}{\point}}{\func{\aoneform{s}}{\point}}}}.
\end{equation}
In the light of this correspondence, each element of $\TTF{r}{s}{\vectorfields{\Man{}}}$
is also called an $\quotl$$\opair{r}{s}$ ($r$-covariant and $s$-contravariant) smooth tensor field on $\Man{}$$\quotr$.
For a detailed information about the correspondence $\Upsilon$, refer to \cite[Chap. 7]{JLee}.
}
%%%%%%%%%%%%%%%%%%%%%%%%%%%%%%%%%%%%%%%%%%%%%%%%%%%%%%%%%%%%%%%%%%%%%%%%%%%%%%%%%%%%%%%%%%%%%%%%%%%%%%%%%%%%%%%%%%%%%%%%%%%%%%%%
\fixed
$\Man{}$ is fixed as a manifold with dimension $m$, and $\connection{}$ as an affine connection on $\Man{}$.
\endfixed
%%%%%%%%%%%%%%%%%%%%%%%%%%%%%%%%%%%%%%%%%%%%%%%%%%%%%%%%%%%%%%%%%%%%%%%%%%%%%%%%%%%%%%%%%%%%%%%%%%%%%%%%%%%%%%%%%%%%%%%%%%%%%%%%
\proposition
Let
$\function{\Delta}{\Cprod{\vectorfields{\Man{}}}{\TTF{1}{0}{\vectorfields{\Man{}}}}}{\Func{\vectorfields{\Man{}}}{\smoothmaps{\Man{}}}}$
be the map characterized by the rule,
\begin{align}
&\Foreach{\avecf{}}{\vectorfields{\Man{}}}
\Foreach{\oneform{}}{\TTF{1}{0}{\vectorfields{\Man{}}}}\cr
&\Foreach{\avecff{}}{\vectorfields{\Man{}}}
\func{\[\func{\Delta}{\binary{\avecf{}}{\oneform{}}}\]}{\avecff{}}\eqdef
\lieder{\avecf{}}{\[\func{\oneform{}}{\avecff{}}\]}-
\func{\oneform{}}{\con{\avecf{}}{\avecff{}}}.
\end{align}
Then, for every $\avecf{}\in\vectorfields{\Man{}}$, and every $\oneform{}\in\TTF{1}{0}{\vectorfields{\Man{}}}$,
$\func{\Delta}{\binary{\avecf{}}{\oneform{}}}\in\TTF{1}{0}{\vectorfields{\Man{}}}$.
\proof
Fix $\avecf{}\in\vectorfields{\Man{}}$ and $\oneform{}\in\TTF{1}{0}{\vectorfields{\Man{}}}$.\\
Let $\avecff{}$ and $\avecff{1}$ be arbitrary elements of $\vectorfields{\Man{}}$, and let
$f$ be an arbitrary element of $\smoothmaps{\Man{}}$.
\begin{itemize}
\item[\myitem{pr-1.}]
Considering the linearity of Lie-derivatives and $\oneform{}$,
\begin{align}
\func{\[\func{\Delta}{\binary{\avecf{}}{\oneform{}}}\]}{\avecff{}+\avecff{1}}&=
\lieder{\avecf{}}{\[\func{\oneform{}}{\avecff{}+\avecff{1}}\]}-\func{\oneform{}}{\con{\avecf{}}{\avecff{}+\avecff{1}}}\cr
&=\lieder{\avecf{}}{\[\func{\oneform{}}{\avecff{}}+\func{\oneform{}}{\avecff{1}}\]}-
\func{\oneform{}}{\con{\avecf{}}{\avecff{}}+\con{\avecf{}}{\avecff{1}}}\cr
&=\(\lieder{\avecf{}}{\[\func{\oneform{}}{\avecff{}}\]}-\func{\oneform{}}{\con{\avecf{}}{\avecff{}}}\)+
\(\lieder{\avecf{}}{\[\func{\oneform{}}{\avecff{1}}\]}-\func{\oneform{}}{\con{\avecf{}}{\avecff{1}}}\)\cr
&=\func{\[\func{\Delta}{\binary{\avecf{}}{\oneform{}}}\]}{\avecff{}}+
\func{\[\func{\Delta}{\binary{\avecf{}}{\oneform{}}}\]}{\avecff{1}}.
\end{align}
\item[\myitem{pr-2.}]
Considering the $\smoothmaps{\Man{}}$-linearity of $\oneform{}$, and properties of
Lie-derivatives and connections,
\begin{align}
\func{\[\func{\Delta}{\binary{\avecf{}}{\oneform{}}}\]}{f\avecff{}}&=
\lieder{\avecf{}}{\[\func{\oneform{}}{f\avecff{}}\]}-\func{\oneform{}}{\con{\avecf{}}{f\avecff{}}}\cr
&=\lieder{\avecf{}}{\(f\[\func{\oneform{}}{\avecff{}}\]\)}-
\func{\oneform{}}{\(\lieder{\avecf{}}{f}\)\avecff{}+f\con{\avecf{}}{\avecff{}}}\cr
&=f\lieder{\avecf{}}{\[\func{\oneform{}}{\avecff{}}\]}+\[\func{\oneform{}}{\avecff{}}\]\(\lieder{\avecf{}}{f}\)-
\(\lieder{\avecf{}}{f}\)\[\func{\oneform{}}{\avecff{}}\]-f\func{\oneform{}}{\con{\avecf{}}{\avecff{}}}\cr
&=f\(\lieder{\avecf{}}{\[\func{\oneform{}}{\avecff{}}\]}-\func{\oneform{}}{\con{\avecf{}}{\avecff{}}}\)\cr
&=f\(\func{\[\func{\Delta}{\binary{\avecf{}}{\oneform{}}}\]}{\avecff{}}\).
\end{align}
\end{itemize}
Therefore, $\func{\Delta}{\binary{\avecf{}}{\oneform{}}}$ is $\smoothmaps{\Man{}}$-linear,
and hence an element of $\TTF{1}{0}{\vectorfields{\Man{}}}$.
\endpro
%%%%%%%%%%%%%%%%%%%%%%%%%%%%%%%%%%%%%%%%%%%%%%%%%%%%%%%%%%%%%%%%%%%%%%%%%%%%%%%%%%%%%%%%%%%%%%%%%%%%%%%%%%%%%%%%%%%%%%%%%%%%%%%%
\definition
We define the map
$\function{\tfconnection{\connection{}}{1}{0}}{\Cprod{\vectorfields{\Man{}}}{\oneforms{\Man{}}}}{\oneforms{\Man{}}}$
as,
\begin{align}
&\Foreach{\avecf{}}{\vectorfields{\Man{}}}
\Foreach{\aoneform{}}{\oneforms{\Tanbun{\Man{}}}}\cr
&\Foreach{\avecff{}}{\vectorfields{\Man{}}}
\func{\[\TFequiv{\func{\tfconnection{\connection{}}{1}{0}}{\binary{\avecf{}}{\aoneform{}}}}\]}{\avecff{}}\eqdef
\lieder{\avecf{}}{\[\func{\TFequiv{\aoneform{}}}{\avecff{}}\]}-
\func{\TFequiv{\aoneform{}}}{\con{\avecf{}}{\avecff{}}},
\end{align}
and refer to it as the $\quotl$$\opair{1}{0}$-covariant-derivative on $\Man{}$ relative to the connection
$\connection{}$$\quotr$.
$\func{\tfconnection{\connection{}}{1}{0}}{\binary{\avecf{}}{\aoneform{}}}$ can be alternatively denoted by
$\tfcon{1}{0}{\avecf{}}{\aoneform{}}$.
\endef
%%%%%%%%%%%%%%%%%%%%%%%%%%%%%%%%%%%%%%%%%%%%%%%%%%%%%%%%%%%%%%%%%%%%%%%%%%%%%%%%%%%%%%%%%%%%%%%%%%%%%%%%%%%%%%%%%%%%%%%%%%%%%%%%
\proposition
Let $r$ and $s$ be non-negative integers such that $r+s>0$, and let
$\function{\Delta}{\Cprod{\vectorfields{\Man{}}}{\TTF{r}{s}{\vectorfields{\Man{}}}}}
{\Func{\Cprod{\(\vectorfields{\Man{}}\)^{r}}{\(\oneforms{\Man{}}\)^{s}}}{\smoothmaps{\Man{}}}}$
be the map characterized by the rule,
\begin{align}
&\Foreach{\avecf{}}{\vectorfields{\Man{}}}
\Foreach{\TT}{\TTF{r}{s}{\vectorfields{\Man{}}}}\cr
&\Foreach{\mtuple{\avecff{1}}{\avecff{r}}}{\vectorfields{\Man{}}}
\Foreach{\mtuple{\aoneform{1}}{\aoneform{s}}}{\oneforms{\Man{}}}\cr
&\begin{aligned}
\func{\[\func{\Delta}{\binary{\avecf{}}{\TT}}\]}{\binary{\suc{\avecff{1}}{\avecff{r}}}{\suc{\aoneform{1}}{\aoneform{s}}}}\eqdef
&~~~\lieder{\avecf{}}{\[\func{\TT}{\binary{\suc{\avecff{1}}{\avecff{r}}}{\suc{\aoneform{1}}{\aoneform{s}}}}\]}\cr
&-\sum_{i=1}^{r}\func{\TT}{\binary{\suc{\avecff{1}}{\con{\avecf{}}{\avecff{i}}},\ldots,{\avecff{r}}}{\suc{\aoneform{1}}{\aoneform{s}}}}\cr
&-\sum_{j=1}^{s}\func{\TT}{\binary{\suc{\avecff{1}}{\avecff{r}}}{\suc{\aoneform{1}}{\tfcon{1}{0}{\avecf{}}{\aoneform{j}}},\ldots,{\aoneform{s}}}}.
\end{aligned}
\end{align}
Then, for every $\avecf{}\in\vectorfields{\Man{}}$, and every
$\TT\in\TTF{r}{s}{\vectorfields{\Man{}}}$, $\func{\Delta}{\binary{\avecf{}}{\TT}}\in\TTF{r}{s}{\vectorfields{\Man{}}}$.
\proof
It is left as an exercise.
\endpro
%%%%%%%%%%%%%%%%%%%%%%%%%%%%%%%%%%%%%%%%%%%%%%%%%%%%%%%%%%%%%%%%%%%%%%%%%%%%%%%%%%%%%%%%%%%%%%%%%%%%%%%%%%%%%%%%%%%%%%%%%%%%%%%%
\definition\label{defcovariantderivativeoftensorfields1}
Let $r$ and $s$ be non-negative integers such that $r+s>0$.
We define the map
$\function{\tfconnection{\connection{}}{r}{s}}{\Cprod{\vectorfields{\Man{}}}{\TF{r}{s}{\Tanbun{\Man{}}}}}{\TF{r}{s}{\Tanbun{\Man{}}}}$
as,
\begin{align}
&\Foreach{\avecf{}}{\vectorfields{\Man{}}}
\Foreach{T}{\TF{r}{s}{\Tanbun{\Man{}}}}\cr
&\Foreach{\mtuple{\avecff{1}}{\avecff{r}}}{\vectorfields{\Man{}}}
\Foreach{\mtuple{\aoneform{1}}{\aoneform{s}}}{\oneforms{\Man{}}}\cr
&\begin{aligned}
\func{\[\TFequiv{\func{\tfconnection{\connection{}}{r}{s}}{\binary{\avecf{}}{T}}}\]}{\binary{\suc{\avecff{1}}{\avecff{r}}}{\suc{\aoneform{1}}{\aoneform{s}}}}\eqdef
&~~~\lieder{\avecf{}}{\[\func{\TFequiv{T}}{\binary{\suc{\avecff{1}}{\avecff{r}}}{\suc{\aoneform{1}}{\aoneform{s}}}}\]}\cr
&-\sum_{i=1}^{r}\func{\TFequiv{T}}{\binary{\suc{\avecff{1}}{\con{\avecf{}}{\avecff{i}}},\ldots,{\avecff{r}}}{\suc{\aoneform{1}}{\aoneform{s}}}}\cr
&-\sum_{j=1}^{s}\func{\TFequiv{T}}{\binary{\suc{\avecff{1}}{\avecff{r}}}{\suc{\aoneform{1}}{\tfcon{1}{0}{\avecf{}}{\aoneform{j}}},\ldots,{\aoneform{s}}}},
\end{aligned}
\end{align}
and refer to it as the $\quotl$$\opair{r}{s}$-covariant-differential operator on tensor fields of $\Tanbun{\Man{}}$ relative to the connection
$\connection{}$$\quotr$.\\
Furthermore, by defining $\TF{0}{0}{\Tanbun{\Man{}}}:=\smoothmaps{\Man{}}$, we define in particular the map
$\function{\tfconnection{\connection{}}{0}{0}}{\Cprod{\vectorfields{\Man{}}}{\TF{0}{0}{\Tanbun{\Man{}}}}}{\TF{0}{0}{\Tanbun{\Man{}}}}$
to be the Lie-derivative operator, that is,
\begin{align}
\Foreach{\avecf{}}{\vectorfields{\Man{}}}
\Foreach{f}{\TF{0}{0}{\Tanbun{\Man{}}}}
\tfcon{0}{0}{\avecf{}}{f}\eqdef\lieder{\avecf{}}{f}.
\end{align}
$\func{\tfconnection{\connection{}}{r}{s}}{\binary{\avecf{}}{T}}$ can be alternatively denoted by
$\tfcon{r}{s}{\avecf{}}{T}$.
\endef
%%%%%%%%%%%%%%%%%%%%%%%%%%%%%%%%%%%%%%%%%%%%%%%%%%%%%%%%%%%%%%%%%%%%%%%%%%%%%%%%%%%%%%%%%%%%%%%%%%%%%%%%%%%%%%%%%%%%%%%%%%%%%%%%
\definition\label{defcovariantderivativeoftensorfields2}
Let $r$ and $s$ be non-negative integers such that $r+s>0$.
We define the map
$\function{\tfconnection{\connection{}}{r}{s}}{\Cprod{\vectorfields{\Man{}}}{\TTF{r}{s}{\vectorfields{\Man{}}}}}{\TTF{r}{s}{\vectorfields{\Man{}}}}$
as,
\begin{align}
&\Foreach{\avecf{}}{\vectorfields{\Man{}}}
\Foreach{T}{\TTF{r}{s}{\vectorfields{\Man{}}}}\cr
&\Foreach{\mtuple{\avecff{1}}{\avecff{r}}}{\vectorfields{\Man{}}}
\Foreach{\mtuple{\aoneform{1}}{\aoneform{s}}}{\oneforms{\Man{}}}\cr
&\begin{aligned}
\func{\[\func{\tfconnection{\connection{}}{r}{s}}{\binary{\avecf{}}{T}}\]}{\binary{\suc{\avecff{1}}{\avecff{r}}}{\suc{\aoneform{1}}{\aoneform{s}}}}\eqdef
&~~~\lieder{\avecf{}}{\[\func{\TFequiv{T}}{\binary{\suc{\avecff{1}}{\avecff{r}}}{\suc{\aoneform{1}}{\aoneform{s}}}}\]}\cr
&-\sum_{i=1}^{r}\func{\TFequiv{T}}{\binary{\suc{\avecff{1}}{\con{\avecf{}}{\avecff{i}}},\ldots,{\avecff{r}}}{\suc{\aoneform{1}}{\aoneform{s}}}}\cr
&-\sum_{j=1}^{s}\func{\TFequiv{T}}{\binary{\suc{\avecff{1}}{\avecff{r}}}{\suc{\aoneform{1}}{\tfcon{1}{0}{\avecf{}}{\aoneform{j}}},\ldots,{\aoneform{s}}}},
\end{aligned}
\end{align}
and refer to it as the $\quotl$$\opair{r}{s}$-covariant-differential operator on tensor fields on $\Man{}$ relative to the connection
$\connection{}$$\quotr$.\\
Furthermore, by defining $\TTF{0}{0}{\Tanbun{\Man{}}}:=\smoothmaps{\Man{}}$, we define in particular the map
$\function{\tfconnection{\connection{}}{0}{0}}{\Cprod{\vectorfields{\Man{}}}{\TTF{0}{0}{\vectorfields{\Man{}}}}}{\TTF{0}{0}{\vectorfields{\Man{}}}}$
to be the Lie-derivative operator, that is,
\begin{align}
\Foreach{\avecf{}}{\vectorfields{\Man{}}}
\Foreach{f}{\TTF{0}{0}{\vectorfields{\Man{}}}}
\tfcon{0}{0}{\avecf{}}{f}\eqdef\lieder{\avecf{}}{f}.
\end{align}
Note that the map $\tfconnection{\connection{}}{r}{s}$ defined here differs from that defined in the previous definition.
But since $\TTF{0}{0}{\vectorfields{\Man{}}}$ and $\TF{r}{s}{\Tanbun{\Man{}}}$ are naturally identified with together,
this difference is not a fundamental one, and that is why we denote both maps with the same notation.
\endef
%%%%%%%%%%%%%%%%%%%%%%%%%%%%%%%%%%%%%%%%%%%%%%%%%%%%%%%%%%%%%%%%%%%%%%%%%%%%%%%%%%%%%%%%%%%%%%%%%%%%%%%%%%%%%%%%%%%%%%%%%%%%%%%%
\remark
Consider the trivial smooth vector-bundle $\manprod{\Man{}}{\R}$ of rank $1$. There is a completely natural bijection
between the set of all smooth sections of $\manprod{\Man{}}{\R}$ and $\smoothmaps{\Man{}}$ (the set of all real-valued
smooth maps on $\Man{}$). So, $\TF{0}{0}{\Tanbun{\Man{}}}$ is naturally identified with
$\vbsections{\manprod{\Man{}}{\R}}$, and $\tfconnection{\connection{}}{0}{0}$ can equivalently be considered
as a map from $\Cprod{\vectorfields{\Man{}}}{\vbsections{\manprod{\Man{}}{\R}}}$ to
$\vbsections{\manprod{\Man{}}{\R}}$.\\
Furthermore, we define the $\opair{0}{0}$-tensor-bundle of the tangent-bundle
$\Tanbun{\Man{}}$, denoted by $\vbtensorbundle{0}{0}{\Tanbun{\Man{}}}$, to be $\manprod{\Man{}}{\R}$.
\endremark
%%%%%%%%%%%%%%%%%%%%%%%%%%%%%%%%%%%%%%%%%%%%%%%%%%%%%%%%%%%%%%%%%%%%%%%%%%%%%%%%%%%%%%%%%%%%%%%%%%%%%%%%%%%%%%%%%%%%%%%%%%%%%%%%
\remark
$\TF{0}{1}{\Tanbun{\Man{}}}$ is naturally identified with the set of all smooth vector fields on $\Man{}$,
that is $\vectorfields{\Man{}}$. This natural correspondence roots in the convention that the dual of
the dual of a vector-space is considered to be naturally identified with itself. Relying on this convention,
and considering that the value of any element of $\TF{0}{1}{\Tanbun{\Man{}}}$ at each point $\point\in\Man{}$
is an element of the double-dual of $\tanspace{\point}{\Man{}}$, the existance of such a natural
identification becomes apparent.
\endremark
%%%%%%%%%%%%%%%%%%%%%%%%%%%%%%%%%%%%%%%%%%%%%%%%%%%%%%%%%%%%%%%%%%%%%%%%%%%%%%%%%%%%%%%%%%%%%%%%%%%%%%%%%%%%%%%%%%%%%%%%%%%%%%%%
\textit{
Given a pair of non-negative integers $r$ and $s$, we denote by $\VTF{r}{s}{\Tanbun{\Man{}}}$ the vector-space
of the set of all $\opair{r}{s}$ tensor fields on $\Tanbun{\Man{}}$ endowed with its canonical linear structure.
We also define $\DVTF{\Tanbun{\Man{}}}:=\Dsum{\defSet{\VTF{r}{s}{\Tanbun{\Man{}}}}{r,s=0,1,\ldots}}$.
}
%%%%%%%%%%%%%%%%%%%%%%%%%%%%%%%%%%%%%%%%%%%%%%%%%%%%%%%%%%%%%%%%%%%%%%%%%%%%%%%%%%%%%%%%%%%%%%%%%%%%%%%%%%%%%%%%%%%%%%%%%%%%%%%%
\definition
We define
$\function{\tfconnection{\connection{}}{\infty}{\infty}}{\Cprod{\vectorfields{\Man{}}}{\DVTF{\Tanbun{\Man{}}}}}{\DVTF{\Tanbun{\Man{}}}}$
as the map such that for evey $\avecf{}\in\vectorfields{\Man{}}$,
$\tfcon{\infty}{\infty}{\avecf{}}{}$ is the unique linear operator on $\DVTF{\Tanbun{\Man{}}}$ satisfying,
\begin{equation}
\Foreach{T}{\TF{r}{s}{\Tanbun{\Man{}}}}
\tfcon{\infty}{\infty}{\avecf{}}{T}=\tfcon{r}{s}{\avecf{}}{T},
\end{equation}
for every pair of non-negative integers $r$ and $s$.
$\tfconnection{\connection{}}{\infty}{\infty}$ is referred to as the
$\quotl$total-covariant-differential operator on $\Man{}$ relative to the connection
$\connection{}$$\quotr$.\\
When everything is clear enough, $\tfconnection{\connection{}}{\infty}{\infty}$ or
$\tfconnection{\connection{}}{r}{s}$ can simply be denoted by $\connection{}$ itself.
\endef
%%%%%%%%%%%%%%%%%%%%%%%%%%%%%%%%%%%%%%%%%%%%%%%%%%%%%%%%%%%%%%%%%%%%%%%%%%%%%%%%%%%%%%%%%%%%%%%%%%%%%%%%%%%%%%%%%%%%%%%%%%%%%%%%
\textit{
We review here the notion of contraction in both the category of finite-dimensional vector-spaces and the category of
manifolds.\\
Let $\vectorspace{}$ be an $\R$-vector-space with dimension $m$. Let $\vsbase{}=\mtuple{\vsbase{1}}{\vsbase{m}}$ be an
ordered-basis of $\vectorspace{}$, and let $\dualvsbase{}=\mtuple{\dualvsbase{1}}{\dualvsbase{m}}$ denote the dual of
the ordered-basis $\vsbase{}$. Given a $\opair{1}{1}$-tensor ($1$-covariant and $1$-contravariant tensor) $\alpha$
on $\vectorspace{}$, that is an element $\alpha$ of $\Tensors{1}{1}{\vectorspace{}}$, the quantity
$\displaystyle\sum_{i=1}^{m}\func{\alpha}{\binary{\vsbase{i}}{\dualvsbase{i}}}$ is independent of the choice of
an ordered-basis for $\vectorspace{}$. Thus, the map
$\function{\vstrace{\vectorspace{}}}{\Tensors{1}{1}{\vectorspace{}}}{\R}$ characterized by the rule,
\begin{equation}
\Foreach{\alpha}{\Tensors{1}{1}{\vectorspace{}}}
\func{\vstrace{\vectorspace{}}}{\alpha}\eqdef
\sum_{i=1}^{m}\func{\alpha}{\binary{\vsbase{i}}{\dualvsbase{i}}},
\end{equation}
for any choice of the ordered-basis $\vsbase{}=\mtuple{\vsbase{1}}{\vsbase{m}}$, is well-defined.
$\vstrace{\vectorspace{}}$ is called the $\quotl$trace-operator of $\vectorspace{}$ on $\opair{1}{1}$-tensors$\quotr$.\\
Now, let $r$ and $s$ be a pair of non-negative integers, let $k\in\seta{\suc{1}{r+1}}$ and $l\in\seta{\suc{s}{s+1}}$,
and let $\beta\in\Tensors{r+1}{s+1}{\vectorspace{}}$.
Given $\suc{v_1}{v_r}\in\vectorspace{}$ and $\suc{\omega_1}{\omega_s}\in\Vdual{\vectorspace{}}$,
the mapping
\begin{equation*}
\maprule{\Delta}{\Cprod{\vectorspace{}}{\Vdual{\vectorspace{}}}\ni\opair{\bar{v}}{\bar{\omega}}}
{\func{\beta}{\binary{\suc{v_1}{v_{k-1}},~\bar{v},~\suc{v_k}{v_r}}
{\suc{\omega_1}{\omega_{l-1}},~\bar{\omega},~\suc{\omega_l}{\omega_s}}}\in\R}
\end{equation*}
is obviously an element of $\Tensors{1}{1}{\vectorspace{}}$, and therefore, the quantity
\begin{equation*}
\sum_{i=1}^{m}\func{\beta}{\binary{\suc{v_1}{v_{k-1}},~\vsbase{i},~\suc{v_k}{v_r}}
{\suc{\omega_1}{\omega_{l-1}},~\dualvsbase{i},~\suc{\omega_l}{\omega_s}}}
\end{equation*}
is independent of the choice of an ordered-basis $\vsbase{}=\mtuple{\vsbase{1}}{\vsbase{m}}$ of $\vectorspace{}$
(with $\dualvsbase{}=\mtuple{\dualvsbase{1}}{\dualvsbase{m}}$ denoting the dual of $\vsbase{}$), which equals
$\func{\vstrace{\vectorspace{}}}{\Delta}$. Therefore, the map
$\function{\contraction{r}{s}{k}{l}{\vectorspace{}}}{\Tensors{r+1}{s+1}{\vectorspace{}}}{\Tensors{r}{s}{\vectorspace{}}}$
characterized by the rule,
\begin{align}
&\Foreach{\beta}{\Tensors{r+1}{s+1}{\vectorspace{}}}
\Foreach{\mtuple{v_1}{v_r}}{{\vectorspace{}}^{r}}
\Foreach{\mtuple{\omega_1}{\omega_s}}{{\Vdual{\vectorspace{}}}^{s}}\cr
&\begin{aligned}
&~~~~\func{\[\func{\contraction{r}{s}{k}{l}{\vectorspace{}}}{\beta}\]}{\binary{\suc{v_1}{v_r}}{\suc{\omega_1}{\omega_s}}}\cr
&\eqdef\sum_{i=1}^{m}\func{\beta}
{\binary{\suc{v_1}{v_{k-1}},~\vsbase{i},~\suc{v_k}{v_r}}
{\suc{\omega_1}{\omega_{l-1}},~\dualvsbase{i},~\suc{\omega_l}{\omega_s}}},
\end{aligned}
\end{align}
for any choice of the ordered-basis $\vsbase{}=\mtuple{\vsbase{1}}{\vsbase{m}}$, is well-defined.
$\contraction{r}{s}{k}{l}{\vectorspace{}}$
is called the $\quotl$$\opair{\opair{r}{k}}{\opair{s}{l}}$-conraction (operator) of $\vectorspace{}$$\quotr$.
Note that in the special case of $r=s=0$, $\Tensors{0}{0}{\vectorspace{}}=\R$ by definition,
and $\contraction{r}{s}{k}{l}{\vectorspace{}}$ coincides with $\vstrace{\vectorspace{}}$.\\
Now, we generalize the notions of trace and contraction to the category of manifolds.
Let $r$ and $s$ be a pair of non-negative integers, let $k\in\seta{\suc{1}{r+1}}$ and $l\in\seta{\suc{s}{s+1}}$.
The mapping $\function{\contraction{r}{s}{k}{l}{\Man{}}}{\TF{r+1}{s+1}{\Tanbun{\Man{}}}}{\TF{r}{s}{\Tanbun{\Man{}}}}$
is defined as,
\begin{align}
\Foreach{T}{\TF{r+1}{s+1}{\Tanbun{\Man{}}}}
\Foreach{\point}{\Man{}}
\func{\[\func{\contraction{r}{s}{k}{l}{\Man{}}}{T}\]}{\point}\eqdef
\func{\contraction{r}{s}{k}{l}{\Tanspace{\point}{\Man{}}}}{\func{T}{\point}}.
\end{align}
We also define the map
$\function{\vstrace{\Man{}}}{\TF{1}{1}{\Tanbun{\Man{}}}}{\smoothmaps{\Man{}}}$ by,
\begin{align}
\Foreach{T}{\TF{1}{1}{\Tanbun{\Man{}}}}
\func{\[\func{\vstrace{\Man{}}}{T}\]}{\point}\eqdef
\func{\vstrace{\Tanspace{\point}{\Man{}}}}{\func{T}{\point}}.
\end{align}
There is a subtle problem that is implicitly assumed in the definition of the notion of contraction
on manifolds. The problem is that for every $T$ in $\TF{r+1}{s+1}{\Tanbun{\Man{}}}$,
$\func{\contraction{r}{s}{k}{l}{\Man{}}}{T}$ defined by the rule above, is also a smooth tensor-field.
The proof is not difficult and left as an exercise. Similar argument holds for the trace operator.
}
%%%%%%%%%%%%%%%%%%%%%%%%%%%%%%%%%%%%%%%%%%%%%%%%%%%%%%%%%%%%%%%%%%%%%%%%%%%%%%%%%%%%%%%%%%%%%%%%%%%%%%%%%%%%%%%%%%%%%%%%%%%%%%%%
\theorem
$\tfconnection{\connection{}}{\infty}{\infty}$ is the unique map from $\Cprod{\vectorfields{\Man{}}}{\DVTF{\Tanbun{\Man{}}}}$
to $\DVTF{\Tanbun{\Man{}}}$ that satisfies the following properties.
\begin{itemize}
\item[\myitem{1.}]
For every non-negative integers $r$ and $s$,
$\tfconnection{\connection{}}{r}{s}$ is a conection on the $\opair{r}{s}$-tensor-bundle of the tangent-bundle
$\Tanbun{\Man{}}$, that is,
\begin{equation}
\tfconnection{\connection{}}{r}{s}\in
\connections{\vbtensorbundle{r}{s}{\Tanbun{\Man{}}}}.
\end{equation}
\item[\myitem{2.}]
$\tfconnection{\connection{}}{0}{0}$ coincides with the Lie-derivative operator on $\Man{}$. That is,
\begin{equation}
\Foreach{\avecf{}}{\vectorfields{\Man{}}}
\Foreach{f}{\smoothmaps{\Man{}}}
\tfcon{0}{0}{\avecf{}}{f}=\lieder{\avecf{}}{f}.
\end{equation}
\item[\myitem{3.}]
Considering the natural identification of $\TF{0}{1}{\Tanbun{\Man{}}}$ and $\vectorfields{\Man{}}$,
\begin{equation}
\tfconnection{\connection{}}{0}{1}=\connection{},
\end{equation}
\item[\myitem{4.}]
For every $\avecf{}\in\vectorfields{\Man{}}$,
and every $T_1\in\TF{r_1}{s_1}{\Tanbun{\Man{}}}$ and every $T_2\in\TF{r_2}{s_2}{\Tanbun{\Man{}}}$ for some
non-negative integers $r_1$, $s_1$, $r_2$, and $s_2$,
\begin{equation}
\tfcon{\infty}{\infty}{\avecf{}}{\(T_1\tensor{}T_2\)}=
\(\tfcon{\infty}{\infty}{\avecf{}}{T_1}\)\tensor{}T_2+
T_1\tensor{}\(\tfcon{\infty}{\infty}{\avecf{}}{T_2}\).
\end{equation}
\item[\myitem{5.}]
$\tfconnection{\connection{}}{\infty}{\infty}$ commutes with contractions.
That is, for every non-negative integers $r$ and $s$,
and every $k\in\seta{\suc{1}{r+1}}$ and $l\in\seta{\suc{1}{s+1}}$,
\begin{equation}
\Foreach{\avecf{}}{\vectorfields{\Man{}}}
\Foreach{T}{\TF{r+1}{s+1}{\Tanbun{\Man{}}}}
\tfcon{r}{s}{\avecf{}}{\(\func{\contraction{r}{s}{k}{l}{\Man{}}}{T}\)}=
\func{\contraction{r}{s}{k}{l}{\Man{}}}{\tfcon{r}{s}{\avecf{}}{T}}.
\end{equation}
\end{itemize}
\proof
It is left as an exercise.
\endthm
\chapteR{Semi-Riemannian Manifolds}
\thispagestyle{fancy}
\section{Scalar Products on Vector-spaces}
%%%%%%%%%%%%%%%%%%%%%%%%%%%%%%%%%%%%%%%%%%%%%%%%%%%%%%%%%%%%%%%%%%%%%%%%%%%%%%%%%%%%%%%%%%%%%%%%%%%%%%%%%%%%%%%%%%%%%%%%%%%%%%%%
\textit{
In this section, we introduce the notion of a scalar-product on a finite-dimensional $\R$-vector-spaceand, and
review some standard relevant results without restating the proofs.
}
%%%%%%%%%%%%%%%%%%%%%%%%%%%%%%%%%%%%%%%%%%%%%%%%%%%%%%%%%%%%%%%%%%%%%%%%%%%%%%%%%%%%%%%%%%%%%%%%%%%%%%%%%%%%%%%%%%%%%%%%%%%%%%%%
\fixed
$\vectorspace{}$ is fixed as an $\R$-vector-space with dimension $m$.
Additionally, $\vectorspace{1}$ is also fixed as an $\R$-vector-space with dimension $m_1$.
\endfixed
%%%%%%%%%%%%%%%%%%%%%%%%%%%%%%%%%%%%%%%%%%%%%%%%%%%%%%%%%%%%%%%%%%%%%%%%%%%%%%%%%%%%%%%%%%%%%%%%%%%%%%%%%%%%%%%%%%%%%%%%%%%%%%%%
\definition
Let $\scalarprod{}$ be a $\opair{2}{0}$-tensor (a $2$-covariant-tensor, or a bilinear form) on the vector-space $\vectorspace{}$,
that is an element of $\Tensors{2}{0}{\vectorspace{}}$.
\begin{itemize}
\item[\myitem{1.}]
$\scalarprod{}$ is called a $\quotl$symmetric bilinear form on the vector-space $\vectorspace{}$$\quotr$ if,
\begin{equation}
\Foreach{\opair{v}{w}}{\Cprod{\vectorspace{}}{\vectorspace{}}}
\scalarprodf{}{v}{w}=\scalarprodf{}{w}{v}.
\end{equation}
\item[\myitem{2.}]
$\scalarprod{}$ is called a $\quotl$non-degenerate bilinear form on the vector-space $\vectorspace{}$$\quotr$ if,
\begin{equation}
\defset{w}{\vectorspace{}}{\(\Foreach{v}{\vectorspace{}}\scalarprodf{}{w}{v}=0\)}=\seta{\zerovec{}},
\end{equation}
where, $\zerovec{}$ denotes the additive neutral element of $\vectorspace{}$.
\item[\myitem{3.}]
$\scalarprod{}$ is called a $\quotl$scalar-product on the vector-space $\vectorspace{}$$\quotr$ if
it is both a symmetric and a 	non-degenerate bilinear form on $\vectorspace{}$.
\item[\myitem{4.}]
$\scalarprod{}$ is called an $\quotl$inner-product on the vector-space $\vectorspace{}$$\quotr$ if
it is a scalar-product on $\vectorspace{}$ and for every non-zero vector in $\vectorspace{}$,
$\scalarprodf{}{v}{v}>0$. Equivalently, it is enough for a bilinear form on $\vectorspace{}$
to be symmetric and satisfy $\scalarprodf{}{v}{v}>0$ for every non-zero vector in $\vectorspace{}$,
in order to be an inner-product on $\vectorspace{}$.
\end{itemize}
\endef
%%%%%%%%%%%%%%%%%%%%%%%%%%%%%%%%%%%%%%%%%%%%%%%%%%%%%%%%%%%%%%%%%%%%%%%%%%%%%%%%%%%%%%%%%%%%%%%%%%%%%%%%%%%%%%%%%%%%%%%%%%%%%%%%
\definition\label{defscalarproductflat}
Let $\scalarprod{}$ be an element of $\Tensors{2}{0}{\vectorspace{}}$. We define the map
$\function{\scalarprodflat{\scalarprod{}}}{\vectorspace{}}{\dualV{\vectorspace{}}}$ as,
\begin{equation}
\Foreach{v}{\vectorspace{}}
\Foreach{w}{\vectorspace{}}
\func{\[\func{\scalarprodflat{\scalarprod{}}}{v}\]}{w}\eqdef
\scalarprodf{}{v}{w},
\end{equation}
where, $\dualV{\vectorspace{}}$ denotes the dual of the vector-space $\vectorspace{}$.
$\scalarprodflat{\scalarprod{}}$ is called the $\quotl$flat operator of the scalar-product $\scalarprod{}$$\quotr$.
\endef
%%%%%%%%%%%%%%%%%%%%%%%%%%%%%%%%%%%%%%%%%%%%%%%%%%%%%%%%%%%%%%%%%%%%%%%%%%%%%%%%%%%%%%%%%%%%%%%%%%%%%%%%%%%%%%%%%%%%%%%%%%%%%%%%
\theorem\label{scalarproductflatislinear}
$\scalarprodflat{\scalarprod{}}$ is a linear map from $\vectorspace{}$ to $\dualV{\vectorspace{}}$.
\endthm
%%%%%%%%%%%%%%%%%%%%%%%%%%%%%%%%%%%%%%%%%%%%%%%%%%%%%%%%%%%%%%%%%%%%%%%%%%%%%%%%%%%%%%%%%%%%%%%%%%%%%%%%%%%%%%%%%%%%%%%%%%%%%%%%
\theorem
Let $\scalarprod{}$ be an element of $\Tensors{2}{0}{\vectorspace{}}$.
Let $\vsbase{}=\mtuple{\vsbase{1}}{\vsbase{m}}$ be an ordered-basis of the vector-space $\vectorspace{}$,
and let $\dualvsbase{}=\mtuple{\dualvsbase{1}}{\dualvsbase{m}}$ be the dual of the ordered-basis $\vsbase{}$.
\begin{equation}
\scalarprod{}=\sum_{i=1}^{m}\sum_{j=1}^{m}
\scalarprodf{}{\vsbase{i}}{\vsbase{j}}
\dualvsbase{i}\tensor{}\dualvsbase{j}.
\end{equation}
\endthm
%%%%%%%%%%%%%%%%%%%%%%%%%%%%%%%%%%%%%%%%%%%%%%%%%%%%%%%%%%%%%%%%%%%%%%%%%%%%%%%%%%%%%%%%%%%%%%%%%%%%%%%%%%%%%%%%%%%%%%%%%%%%%%%%
\definition
Let $\scalarprod{}$ be an element of $\Tensors{2}{0}{\vectorspace{}}$.
Let $\vsbase{}=\mtuple{\vsbase{1}}{\vsbase{m}}$ be an ordered-basis of the vector-space $\vectorspace{}$.
We define the $m\times m$ $\R$-matrix $\scalarprodmatrix{\scalarprod{}}{\vsbase{}}$ as,
\begin{equation}
\scalarprodmatrix{\scalarprod{}}{\vsbase{}}:=
\begin{pmatrix}
\scalarprodf{}{\vsbase{1}}{\vsbase{1}} & \cdots & \scalarprodf{}{\vsbase{1}}{\vsbase{m}}\\
\vdots & \ddots & \vdots\\
\scalarprodf{}{\vsbase{m}}{\vsbase{1}} & \cdots & \scalarprodf{}{\vsbase{m}}{\vsbase{m}}
\end{pmatrix},
\end{equation}
which will be referred to as the $\quotl$matrix representation of the bilinear form $\scalarprod{}$ on
$\vectorspace{}$ with respect to the ordered-basis $\vsbase{}$$\quotr$.
\endef
%%%%%%%%%%%%%%%%%%%%%%%%%%%%%%%%%%%%%%%%%%%%%%%%%%%%%%%%%%%%%%%%%%%%%%%%%%%%%%%%%%%%%%%%%%%%%%%%%%%%%%%%%%%%%%%%%%%%%%%%%%%%%%%%
\theorem\label{thmnondegeneratebilinearformsequivalentconditions}
Let $\scalarprod{}$ be an element of $\Tensors{2}{0}{\vectorspace{}}$,
and let $\vsbase{}=\mtuple{\vsbase{1}}{\vsbase{m}}$ be an ordered-basis of the vector-space $\vectorspace{}$.
The following statements are equivalent.
\begin{itemize}
\item[\myitem{1.}]
$\scalarprod{}$ is a non-degenerate bilinear form on $\vectorspace{}$.
\item[\myitem{2.}]
The square $\R$-matrix $\scalarprodmatrix{\scalarprod{}}{\vsbase{}}$ is non-singular,
that is $\func{\det}{\scalarprodmatrix{\scalarprod{}}{\vsbase{}}}\neq 0$.
\item[\myitem{3.}]
The map $\scalarprodflat{\scalarprod{}}$ is a linear-isomorphism from $\vectorspace{}$ to
$\dualV{\vectorspace{}}$.
\end{itemize}
\endthm
%%%%%%%%%%%%%%%%%%%%%%%%%%%%%%%%%%%%%%%%%%%%%%%%%%%%%%%%%%%%%%%%%%%%%%%%%%%%%%%%%%%%%%%%%%%%%%%%%%%%%%%%%%%%%%%%%%%%%%%%%%%%%%%%
\corollary
Let $\scalarprod{}$ be a scalar-product on $\vectorspace{}$. $\scalarprodflat{\scalarprod{}}$ is a linear-isomorphism
from $\vectorspace{}$ to $\dualV{\vectorspace{}}$.
\endcor
%%%%%%%%%%%%%%%%%%%%%%%%%%%%%%%%%%%%%%%%%%%%%%%%%%%%%%%%%%%%%%%%%%%%%%%%%%%%%%%%%%%%%%%%%%%%%%%%%%%%%%%%%%%%%%%%%%%%%%%%%%%%%%%%
\definition\label{defscalarproductsharp}
Let $\scalarprod{}$ be a scalar-product on $\vectorspace{}$. We define the map
$\function{\scalarprodsharp{\scalarprod{}}}{\dualV{\vectorspace{}}}{\vectorspace{}}$ as,
\begin{equation}
\scalarprodsharp{\scalarprod{}}:=
\finv{\(\scalarprodflat{\scalarprod{}}\)}.
\end{equation}
$\scalarprodsharp{\scalarprod{}}$ is called the $\quotl$sharp operator of the scalar-product $\scalarprod{}$$\quotr$.
\endef
%%%%%%%%%%%%%%%%%%%%%%%%%%%%%%%%%%%%%%%%%%%%%%%%%%%%%%%%%%%%%%%%%%%%%%%%%%%%%%%%%%%%%%%%%%%%%%%%%%%%%%%%%%%%%%%%%%%%%%%%%%%%%%%%
\proposition\label{proscalarproductsharpaction}
Let $\scalarprod{}$ be a scalar-product on $\vectorspace{}$.
\begin{equation}
\Foreach{T}{\dualV{\vectorspace{}}}
\Foreach{v}{\vectorspace{}}
\func{T}{v}=\scalarprodf{\func{\scalarprodsharp{\scalarprod{}}}{T}}{v}.
\end{equation}
\endpro
%%%%%%%%%%%%%%%%%%%%%%%%%%%%%%%%%%%%%%%%%%%%%%%%%%%%%%%%%%%%%%%%%%%%%%%%%%%%%%%%%%%%%%%%%%%%%%%%%%%%%%%%%%%%%%%%%%%%%%%%%%%%%%%%
\theorem
Let $\scalarprod{}$ be a scalar-product on $\vectorspace{}$.
\begin{itemize}
\item[\myitem{1.}]
There exists an ordered-basis $\vsbase{}$ of $\vectorspace{}$ such that
$\scalarprodmatrix{\scalarprod{}}{\vsbase{}}$ is a diagonal $m\times m$ $\R$-matrix.
Then, every diagonal enteries of $\scalarprodmatrix{\scalarprod{}}{\vsbase{}}$ must be non-zero
real numbers.
\item[\myitem{2.}]
There exists a pair $\opair{r}{s}$ of non-negative integers such that $r+s=n$ and
for every ordered-basis $\vsbase{}$ of $\vectorspace{}$ such that $\scalarprodmatrix{\scalarprod{}}{\vsbase{}}$
is diagonal, the number of positive and negative enteries of $\scalarprodmatrix{\scalarprod{}}{\vsbase{}}$
are $r$ and $s$, respectively.
\end{itemize}
\endthm
%%%%%%%%%%%%%%%%%%%%%%%%%%%%%%%%%%%%%%%%%%%%%%%%%%%%%%%%%%%%%%%%%%%%%%%%%%%%%%%%%%%%%%%%%%%%%%%%%%%%%%%%%%%%%%%%%%%%%%%%%%%%%%%%
\definition
Let $\scalarprod{}$ be a scalar-product on $\vectorspace{}$.
and let $\vsbase{}=\mtuple{\vsbase{1}}{\vsbase{m}}$ be an ordered-basis of the vector-space $\vectorspace{}$ such that
$\scalarprodmatrix{\scalarprod{}}{\vsbase{}}$ is diagonal. The number of negative enteries of
$\scalarprodmatrix{\scalarprod{}}{\vsbase{}}$ (which is independent of the choice of an ordered-basis which makes
the matrix representation of $\scalarprod{}$ diagonal) is referred to as the $\quotl$index of the scalar-product $\scalarprod{}$
on $\vectorspace{}$$\quotr$, and will be denoted by $\scalarprodindex{\vectorspace{}}{\scalarprod{}}$ or simply
by $\scalarprodindex{}{\scalarprod{}}$ when the underlying vector-space is clear enough.
\endef
%%%%%%%%%%%%%%%%%%%%%%%%%%%%%%%%%%%%%%%%%%%%%%%%%%%%%%%%%%%%%%%%%%%%%%%%%%%%%%%%%%%%%%%%%%%%%%%%%%%%%%%%%%%%%%%%%%%%%%%%%%%%%%%%
\definition
Let $r$ be a an element of $\seta{\suc{0}{m}}$.
The set of all scalar-products on $\vectorspace{}$ with index $r$ will be denoted by
$\scalarprods{\vectorspace{}}{r}$. The set of all scalar-products on $\vectorspace{}$
will be denoted by $\scalarprods{\vectorspace{}}{}$, that is,
\begin{equation}
\scalarprods{\vectorspace{}}{}:=\bigcup_{r=0}^{m}\scalarprods{\vectorspace{}}{r}.
\end{equation}
\endef
%%%%%%%%%%%%%%%%%%%%%%%%%%%%%%%%%%%%%%%%%%%%%%%%%%%%%%%%%%%%%%%%%%%%%%%%%%%%%%%%%%%%%%%%%%%%%%%%%%%%%%%%%%%%%%%%%%%%%%%%%%%%%%%%
\proposition
$\scalarprods{\vectorspace{}}{0}$ equals the set of all inner-products on $\vectorspace{}$.
\endpro
%%%%%%%%%%%%%%%%%%%%%%%%%%%%%%%%%%%%%%%%%%%%%%%%%%%%%%%%%%%%%%%%%%%%%%%%%%%%%%%%%%%%%%%%%%%%%%%%%%%%%%%%%%%%%%%%%%%%%%%%%%%%%%%%
\definition
Let $\scalarprod{}$ be a scalar-product on $\vectorspace{}$. The pair $\opair{\vectorspace{}}{\scalarprod{}}$
is referred to as a $\quotl$scalar-product-space$\quotr$. Particularly, when
$\scalarprodindex{\vectorspace{}}{\scalarprod{}}=0$,
$\opair{\vectorspace{}}{\scalarprod{}}$ is also called an $\quotl$inner-product-space$\quotr$.
\endef
%%%%%%%%%%%%%%%%%%%%%%%%%%%%%%%%%%%%%%%%%%%%%%%%%%%%%%%%%%%%%%%%%%%%%%%%%%%%%%%%%%%%%%%%%%%%%%%%%%%%%%%%%%%%%%%%%%%%%%%%%%%%%%%%
\definition
\begin{itemize}
\item
Let $\scalarprod{}$ be a scalar-product on $\vectorspace{}$. Let $A$ be a set of vectors in $\vectorspace{}$.
$A$ is called an $\quotl$orthonormal set of the scalar-product-space $\opair{\vectorspace{}}{\scalarprod{}}$
if for every $v$ and $w$ in $A$, $\abs{\scalarprodf{}{v}{w}}=\deltaf{v}{w}$, where $\delta$ denotes the
Kronecker delta function on $A$. Note that when $r=0$, that is when $\opair{\vectorspace{}}{\scalarprod{}}$
is an inner-product-space, an orthonormal set $A$ of $\opair{\vectorspace{}}{\scalarprod{}}$ satisfies
$\scalarprodf{}{v}{w}=\deltaf{v}{w}$ for every $v$ and $w$ in $A$.
\item
An ordered-basis $\vsbase{}=\mtuple{\vsbase{1}}{\vsbase{m}}$ of $\vectorspace{}$ is referred to as an
$\quotl$orthonormal-basis of $\vectorspace{}$$\quotr$ if $\seta{\suc{\vsbase{1}}{\vsbase{m}}}$
is an orthonormal set of the scalar-product-space $\opair{\vectorspace{}}{\scalarprod{}}$.
\end{itemize}
\endef
%%%%%%%%%%%%%%%%%%%%%%%%%%%%%%%%%%%%%%%%%%%%%%%%%%%%%%%%%%%%%%%%%%%%%%%%%%%%%%%%%%%%%%%%%%%%%%%%%%%%%%%%%%%%%%%%%%%%%%%%%%%%%%%%
\theorem
Let $\scalarprod{}$ be a scalar-product on $\vectorspace{}$ wih
$\scalarprodindex{\vectorspace{}}{\scalarprod{}}=r$. There exists an orthonormal-basis
$\vsbase{}=\mtuple{\vsbase{1}}{\vsbase{m}}$ of $\vectorspace{}$. In particular, the order of vectors of
the $\vsbase{}$ can be chosen in such a way that,
\begin{equation}
\scalarprod{}=-\(\sum_{i=1}^{r}
\dualvsbase{i}\tensor{}\dualvsbase{i}\)+
\sum_{i=r+1}^{m}\dualvsbase{i}\tensor{}\dualvsbase{i},
\end{equation}
where $\dualvsbase{}=\mtuple{\dualvsbase{1}}{\dualvsbase{m}}$ denotes the dual of the ordered-basis
$\vsbase{}$, and hence the matrix representation of $\scalarprod{}$ with respect to $\vsbase{}$
takes the form,
\begin{equation}
\scalarprodmatrix{\scalarprod{}}{\vsbase{}}=
\begin{pmatrix}
-\idmatrix{r} & \vline & \zeromatrix{r}{m-r}\\
\hline
\zeromatrix{m-r}{r} & \vline & \idmatrix{m-r}
\end{pmatrix},
\end{equation}
where for arbitrary non-negative integers $s$ and $\p{s}$, $\idmatrix{s}$ denotes the $s\times s$ identity
$\R$-matrix, and $\zeromatrix{s}{\p{s}}$ denotes the $s\times\p{s}$ zero $\R$-matrix.\\
In particular, when $r=0$, that is when $\scalarprod{}$ is an inner-product on $\vectorspace{}$,
there can be found an ordered-basis $\vsbase{}$ of $\vectorspace{}$ such that,
\begin{equation}
\scalarprod{}=
\sum_{i=1}^{m}\dualvsbase{i}\tensor{}\dualvsbase{i},
\end{equation}
and,
\begin{equation}
\scalarprodmatrix{\scalarprod{}}{\vsbase{}}=
\idmatrix{m}.
\end{equation}
\endthm
%%%%%%%%%%%%%%%%%%%%%%%%%%%%%%%%%%%%%%%%%%%%%%%%%%%%%%%%%%%%%%%%%%%%%%%%%%%%%%%%%%%%%%%%%%%%%%%%%%%%%%%%%%%%%%%%%%%%%%%%%%%%%%%%
\theorem
Let $\vsbase{}=\mtuple{\vsbase{1}}{\vsbase{m}}$ be an ordered-basis of the vector-space $\vectorspace{}$.
Let $r\in\seta{\suc{0}{m}}$. The bilinear form
$\displaystyle-\(\sum_{i=1}^{r}\dualvsbase{i}\tensor{}\dualvsbase{i}\)+\sum_{i=r+1}^{m}\dualvsbase{i}\tensor{}\dualvsbase{i}$
is an scalar-product on $\vectorspace{}$ with index $r$.
\endthm
%%%%%%%%%%%%%%%%%%%%%%%%%%%%%%%%%%%%%%%%%%%%%%%%%%%%%%%%%%%%%%%%%%%%%%%%%%%%%%%%%%%%%%%%%%%%%%%%%%%%%%%%%%%%%%%%%%%%%%%%%%%%%%%%
\definition
Let $\scalarprod{}$ be a scalar-product on $\vectorspace{}$, and let $\scalarprod{1}$ be a scalar-product on $\vectorspace{1}$.
Let $T$ be a linear map from $\vectorspace{}$ to $\vectorspace{}$. $T$ is called an $\quotl$isometry from the
scalar-product-space $\opair{\vectorspace{}}{\scalarprod{}}$ to the scalar-product-space
$\opair{\vectorspace{1}}{\scalarprod{1}}$$\quotr$ if it preserves the scalar-product operation, that is,
\begin{equation}
\Foreach{\opair{u}{v}}{\Cprod{\vectorspace{}}{\vectorspace{}}}
\func{\scalarprod{1}}{\binary{\func{T}{u}}{\func{T}{v}}}=
\func{\scalarprod{}}{\binary{u}{v}}.
\end{equation}
\endef
%%%%%%%%%%%%%%%%%%%%%%%%%%%%%%%%%%%%%%%%%%%%%%%%%%%%%%%%%%%%%%%%%%%%%%%%%%%%%%%%%%%%%%%%%%%%%%%%%%%%%%%%%%%%%%%%%%%%%%%%%%%%%%%%
\textit{
Let $\scalarprod{}$ be a scalar-product on $\vectorspace{}$.
The set of all isometries from the scalar-product-space $\opair{\vectorspace{}}{\scalarprod{}}$ to itself is denoted by
$\isometries{\binary{\vectorspace{}}{\scalarprod{}}}$. $\isometries{\binary{\vectorspace{}}{\scalarprod{}}}$ together with
the function-composition operation on it forms a group, which is famous for the $\quotl$isometry group of the
scalar-product-space $\opair{\vectorspace{}}{\scalarprod{}}$$\quotr$.
}
%%%%%%%%%%%%%%%%%%%%%%%%%%%%%%%%%%%%%%%%%%%%%%%%%%%%%%%%%%%%%%%%%%%%%%%%%%%%%%%%%%%%%%%%%%%%%%%%%%%%%%%%%%%%%%%%%%%%%%%%%%%%%%%%
\theorem
Let $\scalarprod{}$ be a scalar-product on $\vectorspace{}$, and let $\scalarprod{1}$ be a scalar-product on $\vectorspace{1}$.
Suppose that the scalar-product-spaces are isometric via the isometry $\function{T}{\vectorspace{}}{\vectorspace{1}}$.
Let $\mtuple{\vsbase{1}}{\vsbase{m}}$ be an orthonormal ordered-basis of $\opair{\vectorspace{}}{\scalarprod{}}$.
Then $\mtuple{\func{T}{\vsbase{1}}}{\func{T}{\vsbase{m}}}$ is an orthonormal ordered-basis of
$\opair{\vectorspace{1}}{\scalarprod{1}}$.\\
In other words, an isometry sends any orthonormal ordered-basis to an orthonormal ordered-basis.
\endthm
%%%%%%%%%%%%%%%%%%%%%%%%%%%%%%%%%%%%%%%%%%%%%%%%%%%%%%%%%%%%%%%%%%%%%%%%%%%%%%%%%%%%%%%%%%%%%%%%%%%%%%%%%%%%%%%%%%%%%%%%%%%%%%%%
%%%%%%%%%%%%%%%%%%%%%%%%%%%%%%%%%%%%%%%%%%%%%%%%%%%%%%%%%%%%%%%%%%%%%%%%%%%%%%%%%%%%%%%%%%%%%%%%%%%%%%%%%%%%%%%%%%%%%%%%%%%%%%%%
%%%%%%%%%%%%%%%%%%%%%%%%%%%%%%%%%%%%%%%%%%%%%%%%%%%%%%%%%%%%%%%%%%%%%%%%%%%%%%%%%%%%%%%%%%%%%%%%%%%%%%%%%%%%%%%%%%%%%%%%%%%%%%%%
%%%%%%%%%%%%%%%%%%%%%%%%%%%%%%%%%%%%%%%%%%%%%%%%%%%%%%%%%%%%%%%%%%%%%%%%%%%%%%%%%%%%%%%%%%%%%%%%%%%%%%%%%%%%%%%%%%%%%%%%%%%%%%%%
%%%%%%%%%%%%%%%%%%%%%%%%%%%%%%%%%%%%%%%%%%%%%%%%%%%%%%%%%%%%%%%%%%%%%%%%%%%%%%%%%%%%%%%%%%%%%%%%%%%%%%%%%%%%%%%%%%%%%%%%%%%%%%%%
%%%%%%%%%%%%%%%%%%%%%%%%%%%%%%%%%%%%%%%%%%%%%%%%%%%%%%%%%%%%%%%%%%%%%%%%%%%%%%%%%%%%%%%%%%%%%%%%%%%%%%%%%%%%%%%%%%%%%%%%%%%%%%%%
%%%%%%%%%%%%%%%%%%%%%%%%%%%%%%%%%%%%%%%%%%%%%%%%%%%%%%%%%%%%%%%%%%%%%%%%%%%%%%%%%%%%%%%%%%%%%%%%%%%%%%%%%%%%%%%%%%%%%%%%%%%%%%%%
%%%%%%%%%%%%%%%%%%%%%%%%%%%%%%%%%%%%%%%%%%%%%%%%%%%%%%%%%%%%%%%%%%%%%%%%%%%%%%%%%%%%%%%%%%%%%%%%%%%%%%%%%%%%%%%%%%%%%%%%%%%%%%%%
%%%%%%%%%%%%%%%%%%%%%%%%%%%%%%%%%%%%%%%%%%%%%%%%%%%%%%%%%%%%%%%%%%%%%%%%%%%%%%%%%%%%%%%%%%%%%%%%%%%%%%%%%%%%%%%%%%%%%%%%%%%%%%%%
%%%%%%%%%%%%%%%%%%%%%%%%%%%%%%%%%%%%%%%%%%%%%%%%%%%%%%%%%%%%%%%%%%%%%%%%%%%%%%%%%%%%%%%%%%%%%%%%%%%%%%%%%%%%%%%%%%%%%%%%%%%%%%%%
%%%%%%%%%%%%%%%%%%%%%%%%%%%%%%%%%%%%%%%%%%%%%%%%%%%%%%%%%%%%%%%%%%%%%%%%%%%%%%%%%%%%%%%%%%%%%%%%%%%%%%%%%%%%%%%%%%%%%%%%%%%%%%%%
\section{Metric Tensors on Manifolds}
%%%%%%%%%%%%%%%%%%%%%%%%%%%%%%%%%%%%%%%%%%%%%%%%%%%%%%%%%%%%%%%%%%%%%%%%%%%%%%%%%%%%%%%%%%%%%%%%%%%%%%%%%%%%%%%%%%%%%%%%%%%%%%%%
\fixed
$\Man{}$ is fixed as a manifold with dimension $m$. In addition, $\Man{1}$ is also fixed as a manifold with dimension $m_1$.
\endfixed
%%%%%%%%%%%%%%%%%%%%%%%%%%%%%%%%%%%%%%%%%%%%%%%%%%%%%%%%%%%%%%%%%%%%%%%%%%%%%%%%%%%%%%%%%%%%%%%%%%%%%%%%%%%%%%%%%%%%%%%%%%%%%%%%
\definition
Let $\metrictensor{}$ be a $\opair{2}{0}$ smooth tensor field of $\Tanbun{\Man{}}$ (the tangent bundle of $\Man{}$),
that is an element of $\TF{2}{0}{\Tanbun{\Man{}}}$.
\begin{itemize}
\item
Let $r$ be an integer in $\seta{\suc{0}{m}}$.
$\metrictensor{}$ is referred to as a $\quotl$semi-Riemannian metric (or a metric-tensor) on
the manifold $\Man{}$ with index $r$ (or with signature $\opair{r}{m-r}$)$\quotr$
and the pair $\opair{\Man{}}{\metrictensor{}}$ is called a
$\quotl$semi-Riemannian manifold with index $r$ (or with signature $\opair{r}{m-r}$)$\quotr$,
if for every point $\point$ of $\Man{}$, $\func{\metrictensor{}}{\point}$ is a
scalar-product (a symmetric and non-degenerate bilinear form) on the vector-space $\Tanspace{\point}{\Man{}}$ with index $r$,
that is,
\begin{equation}
\Foreach{\point}{\Man{}}
\func{\metrictensor{}}{\point}\in\scalarprods{\Tanspace{\point}{\Man{}}}{r}.
\end{equation}
\item
The set of all semi-Riemannian metrics on $\Man{}$ with index $r$ will be denoted by
$\metrictensors{\Man{}}{r}$. We will also denote by $\metrictensors{\Man{}}{}$ the
collection of all metric-tensors on $\Man{}$ with any index $r\in\seta{\suc{0}{m}}$.
For every $\metrictensor{}\in\metrictensors{\Man{}}{}$, $\metrictensor{}$ will be called a
$\quotr$semi-Riemannian metric (or a metric-tensor) on $\Man{}$$\quotr$, and $\opair{\Man{}}{\metrictensor{}}$
will be called a $\quotl$semi-Riemannian manifold$\quotr$.
\item
Particularly, for every $\metrictensor{}\in\metrictensors{\Man{}}{0}$, $\metrictensor{}$ is referred to as a
$\quotl$Riemannian metric on $\Man{}$$\quotr$, and $\opair{\Man{}}{\metrictensor{}}$ is referred to as a
$\quotl$Riemannian manifold$\quotr$.
\item
If $m\geq 2$, then for every $\metrictensor{}\in\metrictensors{\Man{}}{1}$, $\metrictensor{}$ is referred to as a
$\quotl$Lorentz metric on $\Man{}$$\quotr$, and $\opair{\Man{}}{\metrictensor{}}$ is referred to as a
$\quotl$Lorentz manifold$\quotr$.
\end{itemize}
\endef
%%%%%%%%%%%%%%%%%%%%%%%%%%%%%%%%%%%%%%%%%%%%%%%%%%%%%%%%%%%%%%%%%%%%%%%%%%%%%%%%%%%%%%%%%%%%%%%%%%%%%%%%%%%%%%%%%%%%%%%%%%%%%%%%
\definition\label{defisometry}
Let $\metrictensor{}$ be a semi-Riemannian metric on $\Man{}$ with index $r$, that is an element of $\metrictensors{\Man{}}{r}$,
and let $\metrictensor{1}$ be a semi-Riemannian metric on $\Man{1}$ with index $r_1$, that is an element of
$\metrictensors{\Man{1}}{r_1}$.
\begin{itemize}
\item
Let $\function{f}{\Man{}}{\Man{1}}$ be a smooth map from $\Man{}$ to $\Man{1}$,
that is an element of $\mapdifclass{\infty}{\Man{}}{\Man{1}}$. $f$ is reffered to as a
$\quotl$local isometry from the semi-Riemannian manifold $\opair{\Man{}}{\metrictensor{}}$ to the semi-Riemannian manifold
$\opair{\Man{1}}{\metrictensor{1}}$$\quotr$ if
\begin{equation}
\func{\VBpullback{f}{2}{0}}{\metrictensor{1}}=\metrictensor{},
\end{equation}
where, $\VBpullback{f}{2}{0}$ denotes the $\opair{2}{0}$-pullback of the smooth map $f$, or equivalently,
\begin{equation}
\Foreach{\point}{\Man{}}
\Foreach{\opair{v_1}{v_2}}{\Cprod{\tanspace{\point}{\Man{}}}{\tanspace{\point}{\Man{}}}}
\func{\[\func{\metrictensor{}}{\point}\]}{\binary{v_1}{v_2}}=
\func{\[\func{\metrictensor{1}}{\func{f}{\point}}\]}{\binary{\func{\derr{f}}{v_1}}{\func{\derr{f}}{v_2}}},
\end{equation}
where, $\function{\derr{f}}{\Tanbun{\Man{}}}{\Tanbun{\Man{1}}}$ denotes the tangent-map of
$f\in\mapdifclass{\infty}{\Man{}}{\Man{1}}$.\\
It is said that $\quotl$the semi-Riemannian manifold $\opair{\Man{}}{\metrictensor{}}$ is locally-isometric
to the semi-Riemannian manifold $\opair{\Man{1}}{\metrictensor{1}}$$\quotr$ if there exists a local-isometry
from $\opair{\Man{}}{\metrictensor{}}$ to $\opair{\Man{1}}{\metrictensor{1}}$.
\item
Suppose that the manifolds $\Man{}$ and $\Man{1}$ are diffeomorphic, and hence $m=m_1$.
A smooth map $f$ from $\Man{}$ to $\Man{1}$ is referred to as an $\quotl$isometry from $\Man{}$ to $\Man{1}$$\quotr$
if $f$ is both a diffeomorphism and a local-isometry from $\Man{}$ to $\Man{1}$.
Then, it is a trivial fact that the semi-Riemannian metrics $\metrictensor{}$ and $\metrictensor{1}$
must have the same index, that is $r=r_1$.\\
It is said that $\quotl$the semi-Riemannian manifold $\opair{\Man{}}{\metrictensor{}}$ is isometric
to the semi-Riemannian manifold $\opair{\Man{1}}{\metrictensor{1}}$$\quotr$ if there exists an isometry
from $\opair{\Man{}}{\metrictensor{}}$ to $\opair{\Man{1}}{\metrictensor{1}}$.
\item
The set of all isometries from the semi-Riemannian manifold $\opair{\Man{}}{\metrictensor{}}$ to itself is denoted by
$\isometries{\binary{\Man{}}{\metrictensor{}}}$. $\isometries{\binary{\Man{}}{\metrictensor{}}}$ together with
the function-composition operation on it forms a group, which is famous for the $\quotl$isometry group of the
semi-Riemannian manifold $\opair{\Man{}}{\metrictensor{}}$.
\item
If $\Man{1}$ is an immersed submanifold of $\Man{}$ via the immersion $\function{\iota}{\Man{1}}{\Man{}}$,
and if $\func{\VBpullback{\iota}{2}{0}}{\metrictensor{}}$ is a semi-Riemannian metric on $\Man{1}$, that is
an element of $\metrictensors{\Man{1}}{\nu}$ for some $\nu$, then $\Man{1}$ is called a $\quotl$semi-Riemannian
submanifold of the semi-Riemannian manifold $\opair{\Man{}}{\metrictensor{}}$ (via the immersion $\iota$),
with the induced metric $\func{\VBpullback{\iota}{2}{0}}{\metrictensor{}}$$\quotr$.
In the particular case of $\nu=0$,
$\Man{1}$ is called a $\quotl$Riemannian submanifold of the semi-Riemannian manifold $\opair{\Man{}}{\metrictensor{}}$$\quotr$.
Further, if $\metrictensor{1}=\func{\VBpullback{\iota}{2}{0}}{\metrictensor{}}$, then $\opair{\Man{1}}{\metrictensor{1}}$
is called an $\quotl$isometric immersion into $\opair{\Man{}}{\metrictensor{}}$$\quotr$, and if additionally $\iota$ is
an embedding, then $\opair{\Man{1}}{\metrictensor{1}}$
is called an $\quotl$isometric embedding into $\opair{\Man{}}{\metrictensor{}}$$\quotr$.
\end{itemize}
\endef
%%%%%%%%%%%%%%%%%%%%%%%%%%%%%%%%%%%%%%%%%%%%%%%%%%%%%%%%%%%%%%%%%%%%%%%%%%%%%%%%%%%%%%%%%%%%%%%%%%%%%%%%%%%%%%%%%%%%%%%%%%%%%%%%
\remark
In the particular case when $\metrictensor{}$ is a Riemannian metric on $\Man{}$,
every submanifold of $\Man{}$ is automatically a Riemannian submanifold of $\opair{\Man{}}{\metrictensor{}}$. This is
an immediate consequence of the positive-definiteness of the action of $\metrictensor{}$ at the tangent-space of any point
of $\Man{}$.
\endremark
%%%%%%%%%%%%%%%%%%%%%%%%%%%%%%%%%%%%%%%%%%%%%%%%%%%%%%%%%%%%%%%%%%%%%%%%%%%%%%%%%%%%%%%%%%%%%%%%%%%%%%%%%%%%%%%%%%%%%%%%%%%%%%%%
\proposition
Let $f$ be a smooth map from $\Man{}$ to $\Man{1}$. $f$ is a local-isometry from $\Man{}$ to $\Man{1}$ if and only if
for every $\point\in\Man{}$, there exists an open neighborhood $\U$ of $\point$ such that $\func{\image{f}}{\U}$
is an open set of $\Man{1}$,
and $\reS{f}{\U}$ is an isometry
from $\opair{\subman{\Man{}}{\U}}{\reS{\metrictensor{}}{\U}}$ to $\opair{\subman{\Man{1}}{\U}}{\reS{\metrictensor{1}}{\U}}$.
\endpro
%%%%%%%%%%%%%%%%%%%%%%%%%%%%%%%%%%%%%%%%%%%%%%%%%%%%%%%%%%%%%%%%%%%%%%%%%%%%%%%%%%%%%%%%%%%%%%%%%%%%%%%%%%%%%%%%%%%%%%%%%%%%%%%%
\textit{
Let $n$ be a positive integer, and let $\nu\in\seta{\suc{0}{n}}$. There exists a metric-tensor $\semiEucmetric{n}{\nu}$
of index $\nu$ on the manifold $\R^{n}$ (the set of all
$n$-tuples of real numbers endowed with its canonical differentiable structure), characterized by,
\begin{equation}
\Foreach{\point}{\Man{}}
\Foreach{\opair{u}{v}}{\Cprod{\tanspace{\point}{\R^n}}{\tanspace{\point}{\R^n}}}
\func{\[\func{\semiEucmetric{n}{\nu}}{\point}\]}{\binary{u}{v}}=-\sum_{i=1}^{\nu}u_{i}v_{i}+\sum_{i=\nu+1}^{n}u_{i}v_{i},
\end{equation}
where $u=\mtuple{u_1}{u_n}$ and $v=\mtuple{v_1}{v_n}$.
Equivalently,
\begin{equation}
\semiEucmetric{n}{\nu}=-\sum_{i=1}^{\nu}\Eucstandardframe{n}{i}\tensor{}\Eucstandardframe{n}{i}+
\sum_{i=\nu+1}^{n}\Eucstandardframe{n}{i}\tensor{}\Eucstandardframe{n}{i},
\end{equation}
where, $\mtuple{\Eucstandardframe{n}{1}}{\Eucstandardframe{n}{n}}$ denotes the standard global frame field of $\R^n$.
$\semiEucmetric{n}{\nu}$ is referred to as the $\quotl$standard metric of index $\nu$ on the manifold $\R^n$$\quotr$.
The semi-Riemannian manifold $\opair{\R^n}{\semiEucmetric{n}{\nu}}$ is referred to as the
$\quotl$semi-Euclidean space of dimension $n$ and index $\nu$$\quotr$, which is alternatively denoted by
$\semiEucspace{n}{\nu}$. In the particular case of $\nu=0$, $\semiEucmetric{n}{\nu}$ becomes a
Riemannian metric on $\R^n$, and $\semiEucspace{n}{0}$ is referred to as the $\quotl$Euclidean space of dimension $n$$\quotr$.
Note that $\semiEucspace{n}{n}$ carries almost the same structure as that of $\semiEucspace{n}{0}$ and hence is
not inherently a different construction.
In addition, the particular case $\semiEucspace{4}{1}$ is referred to as the $\quotl$Minkowski space$\quotr$ or
$\quotl$Minkowski space-time$\quotr$.\\
Note that, in general, $\semiEucmetric{n}{\nu}$ is not the only semi-Riemannian metric of index $\nu$ on the manifold $\R^n$;
for instance, in the theory of general relativity, the presence of matter (or more generally, a non-zero stress-energy tensor)
alters the metric-tensor of space-time, which is again of index $1$.
Speaking informally, the Minkowski space is actually the special structure of the empty space-time.
}
%%%%%%%%%%%%%%%%%%%%%%%%%%%%%%%%%%%%%%%%%%%%%%%%%%%%%%%%%%%%%%%%%%%%%%%%%%%%%%%%%%%%%%%%%%%%%%%%%%%%%%%%%%%%%%%%%%%%%%%%%%%%%%%%
\definition
Let $\metrictensor{}$ be a semi-Riemannian metric of index $\nu$ on $\Man{}$, that is an element of
$\metrictensors{\Man{}}{\nu}$.
The semi-Riemannian manifold $\opair{\Man{}}{\metrictensor{}}$ is called a $\quotl$flat semi-Riemannian manifold$\quotr$
if it is locally-isometric to semi-Euclidean space $\semiEucspace{m}{\nu}$.
\endef
%%%%%%%%%%%%%%%%%%%%%%%%%%%%%%%%%%%%%%%%%%%%%%%%%%%%%%%%%%%%%%%%%%%%%%%%%%%%%%%%%%%%%%%%%%%%%%%%%%%%%%%%%%%%%%%%%%%%%%%%%%%%%%%%
\lemma
Let $\metrictensor{}$ be a semi-Riemannian metric on $\Man{}$, that is an element of $\metrictensors{\Man{}}{}$.
\begin{itemize}
\item
Let $\avecf{}$ be a smooth vector field on $\Man{}$, that is an element of $\vectorfields{\Man{}}$.
The mapping $\maprule{F}{\Man{}\ni\point}{\func{\scalarprodflat{\[\func{\metrictensor{}}{\point}\]}}{\func{\avecf{}}{\point}}}$
is a smooth map from $\Man{}$ to $\vbDualbundle{\Tanbun{\Man{}}}$, and hence a section of the dual-bundle of the
tangent-bundle $\Tanbun{\Man{}}$, that is an element of
$\vbsections{\vbDualbundle{\Tanbun{\Man{}}}}=\TF{1}{0}{\Tanbun{\Man{}}}=\oneforms{\Man{}}$.
\item
Let $\avecff{}$ be a co-vector field on $\Man{}$, that is an element of
$\vbsections{\vbDualbundle{\Tanbun{\Man{}}}}=\TF{1}{0}{\Tanbun{\Man{}}}=\oneforms{\Man{}}$.
The mapping
$\maprule{\p{F}}{\Man{}\ni\point}{\func{\scalarprodsharp{\[\func{\metrictensor{}}{\point}\]}}{\func{\avecff{}}{\point}}}$
is a smooth map from $\Man{}$ to $\Tanbun{\Man{}}$, and hence a smooth vector field on $\Man{}$, that is an
element of $\vectorfields{\Man{}}$.
\end{itemize}
\proof
It is left as an exercise.
\endlem
%%%%%%%%%%%%%%%%%%%%%%%%%%%%%%%%%%%%%%%%%%%%%%%%%%%%%%%%%%%%%%%%%%%%%%%%%%%%%%%%%%%%%%%%%%%%%%%%%%%%%%%%%%%%%%%%%%%%%%%%%%%%%%%%
\definition\label{defflatanssharpoperatorsofmetrictensors}
Let $\metrictensor{}$ be a semi-Riemannian metric on $\Man{}$, that is an element of $\metrictensors{\Man{}}{}$.
\begin{itemize}
\item
The map $\function{\mtflat{\metrictensor{}}}{\vectorfields{\Man{}}}{\vbsections{\vbDualbundle{\Tanbun{\Man{}}}}}$ is defined as,
\begin{equation}
\Foreach{\avecf{}}{\vectorfields{\Man{}}}
\Foreach{\point}{\Man{}}
\func{\[\func{\mtflat{\metrictensor{}}}{\avecf{}}\]}{\point}\eqdef
\func{\scalarprodflat{\[\func{\metrictensor{}}{\point}\]}}{\func{\avecf{}}{\point}},
\end{equation}
where, $\vbsections{\vbDualbundle{\Tanbun{\Man{}}}}=\TF{1}{0}{\Tanbun{\Man{}}}=\oneforms{\Man{}}$
denotes the set of all sections of the dual-bundle of the
tangent-bundle $\Tanbun{\Man{}}$, or equivalently the set of all $\opair{1}{0}$ tensor fields ($1$-convariant tensor fields)
on $\Tanbun{\Man{}}$ (or differential forms of degree 1 on $\Man{}$, or co-vector fields on $\Man{}$).
So, based on \refdef{defscalarproductflat}, equivalently,
\begin{equation}
\Foreach{\avecf{}}{\vectorfields{\Man{}}}
\Foreach{\point}{\Man{}}
\Foreach{v}{\tanspace{\point}{\Man{}}}
\func{\(\func{\[\func{\mtflat{\metrictensor{}}}{\avecf{}}\]}{\point}\)}{v}\eqdef
\func{\[\func{\metrictensor{}}{\point}\]}{\binary{\func{\avecf{}}{\point}}{v}}.
\end{equation}
$\mtflat{\metrictensor{}}$ is called the $\quotl$flat operator of the semi-Riemannian metric $\metrictensor{}$$\quotr$.
\item
The map $\function{\mtsharp{\metrictensor{}}}{\vbsections{\vbDualbundle{\Tanbun{\Man{}}}}}{\vectorfields{\Man{}}}$ is defined as,
\begin{equation}
\Foreach{\avecff{}}{\vbsections{\vbDualbundle{\Tanbun{\Man{}}}}}
\Foreach{\point}{\Man{}}
\func{\[\func{\mtsharp{\metrictensor{}}}{\avecff{}}\]}{\point}\eqdef
\func{\scalarprodsharp{\[\func{\metrictensor{}}{\point}\]}}{\func{\avecff{}}{\point}}.
\end{equation}
$\mtsharp{\metrictensor{}}$ is called the $\quotl$sharp operator of the semi-Riemannian metric $\metrictensor{}$$\quotr$.
\end{itemize}
\endef
%%%%%%%%%%%%%%%%%%%%%%%%%%%%%%%%%%%%%%%%%%%%%%%%%%%%%%%%%%%%%%%%%%%%%%%%%%%%%%%%%%%%%%%%%%%%%%%%%%%%%%%%%%%%%%%%%%%%%%%%%%%%%%%%
\theorem
Let $\metrictensor{}$ be a semi-Riemannian metric on $\Man{}$, that is an element of $\metrictensors{\Man{}}{}$.
$\mtflat{\metrictensor{}}$ is a module-isomorphism (or a $\smoothmaps{\Man{}}$-linear-isomorphism) from
$\vectorfields{\Man{}}$ to $\vbsections{\vbDualbundle{\Tanbun{\Man{}}}}=\oneforms{\Man{}}$ endowed with their canonical
$\smoothmaps{\Man{}}$-module structures. Furthermore,
\begin{equation}
\finv{\(\mtflat{\metrictensor{}}\)}=
\mtsharp{\metrictensor{}}.
\end{equation}
\proof
\begin{itemize}
\item[\myitem{pr-1.}]
We first show that $\mtflat{\metrictensor{}}$ is injective and onto.
According to \refdef{defflatanssharpoperatorsofmetrictensors},
\refthm{thmnondegeneratebilinearformsequivalentconditions}, and \refdef{defscalarproductsharp},
for every $\avecf{}\in\vectorfields{\Man{}}$,
\begin{align}
\Foreach{\point}{\Man{}}
\func{\[\func{\mtsharp{\metrictensor{}}}{\func{\mtflat{\metrictensor{}}}{\avecf{}}}\]}{\point}&=
\func{\scalarprodsharp{\[\func{\metrictensor{}}{\point}\]}}{\func{\[\func{\mtflat{\metrictensor{}}}{\avecf{}}\]}{\point}}\cr
&=\func{\scalarprodsharp{\[\func{\metrictensor{}}{\point}\]}}{\func{\scalarprodflat{\[\func{\metrictensor{}}{\point}\]}}
{\func{\avecf{}}{\point}}}\cr
&=\func{\(\cmp{\scalarprodsharp{\[\func{\metrictensor{}}{\point}\]}}{\scalarprodflat{\[\func{\metrictensor{}}{\point}\]}}\)}
{\func{\avecf{}}{\point}}\cr
&=\func{\avecf{}}{\point}.
\end{align}
Therefore
\begin{equation}
\cmp{\mtsharp{\metrictensor{}}}{\mtflat{\metrictensor{}}}=
\identity{\vectorfields{\Man{}}}.
\end{equation}
In a similar manner, it can be shown that,
\begin{equation}
\cmp{\mtflat{\metrictensor{}}}{\mtsharp{\metrictensor{}}}=
\identity{\oneforms{\Man{}}}.
\end{equation}
\item[\myitem{pr-2.}]
Now, we show that $\mtflat{\metrictensor{}}$ is $\smoothmaps{\Man{}}$-linear.
Let $f\in\smoothmaps{\Man{}}$, and let $\avecf{1}$ and $\avecf{2}$ be smooth vector fields on $\Man{}$.
According to \refdef{defflatanssharpoperatorsofmetrictensors} and \refthm{scalarproductflatislinear},
\begin{align}
\Foreach{\point}{\Man{}}
\func{\[\func{\mtflat{\metrictensor{}}}{f\avecf{1}+\avecf{2}}\]}{\point}&=
\func{\scalarprodflat{\[\func{\metrictensor{}}{\point}\]}}{\func{\[f\avecf{1}+\avecf{2}\]}{\point}}\cr
&=\func{\scalarprodflat{\[\func{\metrictensor{}}{\point}\]}}{\func{f}{\point}\func{\avecf{1}}{\point}+
\func{\avecf{2}}{\point}}\cr
&=\func{f}{\point}\func{\scalarprodflat{\[\func{\metrictensor{}}{\point}\]}}{\func{\avecf{1}}{\point}}+
\func{\scalarprodflat{\[\func{\metrictensor{}}{\point}\]}}{\func{\avecf{2}}{\point}}\cr
&=\func{f}{\point}\func{\[\func{\mtflat{\metrictensor{}}}{\avecf{1}}\]}{\point}+
\func{\[\func{\mtflat{\metrictensor{}}}{\avecf{2}}\]}{\point},
\end{align}
and therefore,
\begin{equation}
\func{\mtflat{\metrictensor{}}}{f\avecf{1}+\avecf{2}}=
f\func{\mtflat{\metrictensor{}}}{\avecf{1}}+
\func{\mtflat{\metrictensor{}}}{\avecf{2}}.
\end{equation}
\end{itemize}
\endthm
%%%%%%%%%%%%%%%%%%%%%%%%%%%%%%%%%%%%%%%%%%%%%%%%%%%%%%%%%%%%%%%%%%%%%%%%%%%%%%%%%%%%%%%%%%%%%%%%%%%%%%%%%%%%%%%%%%%%%%%%%%%%%%%%
\definition\label{defmetricinducedproductofvectorfields}
Let $\metrictensor{}$ be a semi-Riemannian metric on $\Man{}$, that is an element of $\metrictensors{\Man{}}{}$.
We define the map $\function{\vfmetricprodmap{\Man{}}{\metrictensor{}}}{\Cprod{\vectorfields{\Man{}}}{\vectorfields{\Man{}}}}
{\smoothmaps{\Man{}}}$ as,
\begin{equation}
\Foreach{\opair{\avecf{1}}{\avecf{2}}}{\Cprod{\vectorfields{\Man{}}}{\vectorfields{\Man{}}}}
\Foreach{\point}{\Man{}}
\func{\[\vfmetricprod{\Man{}}{\metrictensor{}}{\avecf{1}}{\avecf{2}}\]}{\point}\eqdef
\func{\[\func{\metrictensor{}}{\point}\]}{\binary{\func{\avecf{1}}{\point}}{\func{\avecf{2}}{\point}}}.
\end{equation}
$\vfmetricprod{\Man{}}{\metrictensor{}}{\avecf{1}}{\avecf{2}}$ is called the $\quotl$metric-induced product
of smooth vector fields $\avecf{1}$ and $\avecf{2}$ on $\Man{}$ with respect to the semi-Riemannian metric
$\metrictensor{}$$\quotr$.
When there is no ambiguity about the underlying semi-Riemannian manifold,
$\vfmetricprod{\Man{}}{\metrictensor{}}{\avecf{1}}{\avecf{2}}$ can simply be denoted by
$\vfmetricproduct{\avecf{1}}{\avecf{2}}$.
\endef
%%%%%%%%%%%%%%%%%%%%%%%%%%%%%%%%%%%%%%%%%%%%%%%%%%%%%%%%%%%%%%%%%%%%%%%%%%%%%%%%%%%%%%%%%%%%%%%%%%%%%%%%%%%%%%%%%%%%%%%%%%%%%%%%
\remark
Given a metric-tensor $\metrictensor{}$ on $\Man{}$,
for any $\point\in\Man{}$, and every $\binary{v_1}{v_2}\in\tanspace{\point}{\Man{}}$, we will also alternatively denote
$\func{\[\func{\metrictensor{}}{\point}\]}{\binary{v_1}{v_2}}$ by $\vfmetricprod{\Man{}}{\metrictensor{}}{v_1}{v_2}$,
or simply by $\vfmetricproduct{v_1}{v_2}$ when everything is clear.
\endremark
%%%%%%%%%%%%%%%%%%%%%%%%%%%%%%%%%%%%%%%%%%%%%%%%%%%%%%%%%%%%%%%%%%%%%%%%%%%%%%%%%%%%%%%%%%%%%%%%%%%%%%%%%%%%%%%%%%%%%%%%%%%%%%%%
\proposition\label{prometricinducedproductofvectorfieldsisbilinear}
Let $\metrictensor{}$ be a semi-Riemannian metric on $\Man{}$, that is an element of $\metrictensors{\Man{}}{}$. The map
$\function{\vfmetricprodmap{\Man{}}{\metrictensor{}}}{\Cprod{\vectorfields{\Man{}}}{\vectorfields{\Man{}}}}$ is
$\smoothmaps{\Man{}}$-bilinear.
\proof
It is trivial.
\endpro
%%%%%%%%%%%%%%%%%%%%%%%%%%%%%%%%%%%%%%%%%%%%%%%%%%%%%%%%%%%%%%%%%%%%%%%%%%%%%%%%%%%%%%%%%%%%%%%%%%%%%%%%%%%%%%%%%%%%%%%%%%%%%%%%
\proposition\label{prometricinducedproductoneforms}
Let $\aoneform{}\in\oneforms{\Man{}}$. For every $\avecf{}\in\vectorfields{\Man{}}$,
\begin{equation}
\func{\TFequiv{\aoneform{}}}{\avecf{}}=
\vfmetricprod{\Man{}}{\metrictensor{}}{\func{\mtsharp{\metrictensor{}}}{\aoneform{}}}{\avecf{}}.
\end{equation}
It is an immediate consequence of \refpro{proscalarproductsharpaction}, \refdef{defflatanssharpoperatorsofmetrictensors},
and \refdef{defmetricinducedproductofvectorfields}.
\endpro
%%%%%%%%%%%%%%%%%%%%%%%%%%%%%%%%%%%%%%%%%%%%%%%%%%%%%%%%%%%%%%%%%%%%%%%%%%%%%%%%%%%%%%%%%%%%%%%%%%%%%%%%%%%%%%%%%%%%%%%%%%%%%%%%
\textit{
Given a real-valued smooth map on the manifold $\Man{}$, $\Rderivative{\Man{}}{f}$ denotes the derivative of $f$,
defined as the map $\function{\Rderivative{\Man{}}{f}}{\Man{}}{\vbDualbundle{\Tanbun{\Man{}}}}$ that sends each
point $\point$ of $\Man{}$ to the tangent map of $f$ at $\point$ which is an element of the dual of
$\tanspace{\point}{\Man{}}$. Since $\Rderivative{\Man{}}{f}$ is itself a smooth map from $\Man{}$ to
$\vbDualbundle{\Tanbun{\Man{}}}$, $\vbDualbundle{\Tanbun{\Man{}}}$ is a co-vector field on $\Man{}$,
that is an element of $\vbsections{\vbDualbundle{\Tanbun{\Man{}}}}=\TF{1}{0}{\Tanbun{\Man{}}}=\oneforms{\Man{}}$.\\
Furthermore, note that for every $\avecf{}\in\vectorfields{\Man{}}$, $\func{\TFequiv{\Rderivative{\Man{}}{f}}}{\avecf{}}=
\lieder{\avecf{}}{f}$, where $\lieder{\avecf{}}{f}$ denotes the Lie-derivative of $f$ relative to the vector field $\avecf{}$.
}
%%%%%%%%%%%%%%%%%%%%%%%%%%%%%%%%%%%%%%%%%%%%%%%%%%%%%%%%%%%%%%%%%%%%%%%%%%%%%%%%%%%%%%%%%%%%%%%%%%%%%%%%%%%%%%%%%%%%%%%%%%%%%%%%
\definition\label{defgradient}
Let $\metrictensor{}$ be a semi-Riemannian metric on $\Man{}$, that is an element of $\metrictensors{\Man{}}{}$.
We define the map $\function{\grad{\metrictensor{}}}{\smoothmaps{\Man{}}}{\vectorfields{\Man{}}}$ as,
\begin{equation}
\Foreach{f}{\smoothmaps{\Man{}}}
\func{\grad{\metrictensor{}}}{f}:=
\func{\mtsharp{\metrictensor{}}}{\Rderivative{\Man{}}{f}}.
\end{equation}
So, according to \refdef{defflatanssharpoperatorsofmetrictensors}, equivalently,
\begin{equation}
\Foreach{f}{\smoothmaps{\Man{}}}
\Foreach{\point}{\Man{}}
\func{\[\func{\grad{\metrictensor{}}}{f}\]}{\point}\eqdef
\func{\scalarprodsharp{\[\func{\metrictensor{}}{\point}\]}}{\func{\[\Rderivative{\Man{}}{f}\]}{\point}}.
\end{equation}
$\grad{\metrictensor{}}$ is called the $\quotl$gradiant operator of the semi-Riemannian manifold
$\opair{\Man{}}{\metrictensor{}}$$\quotr$, and for any $f\in\smoothmaps{\Man{}}$,
$\func{\grad{\metrictensor{}}}{f}$ is called the
$\quotl$gradiant of $f$ with respect to the semi-Riemmanian metric $\metrictensor{}$ on $\Man{}$$\quotr$.
When the underlying semi-Riemannian metric is clear enough, $\grad{\metrictensor{}}$ can simply be denoted by
$\grad{}$.
\endef
%%%%%%%%%%%%%%%%%%%%%%%%%%%%%%%%%%%%%%%%%%%%%%%%%%%%%%%%%%%%%%%%%%%%%%%%%%%%%%%%%%%%%%%%%%%%%%%%%%%%%%%%%%%%%%%%%%%%%%%%%%%%%%%%
\corollary
Let $\metrictensor{}$ be a semi-Riemannian metric on $\Man{}$, that is an element of $\metrictensors{\Man{}}{}$.
Let $f\in\smoothmaps{\Man{}}$. For every $\avecf{}\in\vectorfields{\Man{}}$,
\begin{equation}
\lieder{\avecf{}}{f}=
\vfmetricprod{\Man{}}{\metrictensor{}}{\func{\grad{\metrictensor{}}}{f}}{\avecf{}},
\end{equation}
and $\func{\grad{\metrictensor{}}}{f}$ is the unique element of $\vectorfields{\Man{}}$ with this property.
\endcor
%%%%%%%%%%%%%%%%%%%%%%%%%%%%%%%%%%%%%%%%%%%%%%%%%%%%%%%%%%%%%%%%%%%%%%%%%%%%%%%%%%%%%%%%%%%%%%%%%%%%%%%%%%%%%%%%%%%%%%%%%%%%%%%%
\definition
Let $\U$ be a non-empty open subset of $\Man{}$ and let $\mtuple{\avecf{1}}{\avecf{m}}$ be a system of local frame fields
on $\subman{\Man{}}{\U}$. $\mtuple{\avecf{1}}{\avecf{m}}$ is referred to as an $\quotl$orthonormal local frame on
(the open subset $\U$ of)
the semi-Riemannian manifold $\opair{\Man{}}{\metrictensor{}}$$\quotr$ if for every $\point\in\U$,
$\mtuple{\func{\avecf{1}}{\point}}{\func{\avecf{m}}{\point}}$ is an orthonormal ordered-basis of the scalar-product-space
$\opair{\Tanspace{\point}{\Man{}}}{\func{\metrictensor{}}{\point}}$, that is,
\begin{align}
\Foreach{\point}{\Man{}}
\abs{\func{\[\func{\metrictensor{}}{\point}\]}{\binary{\func{\avecf{i}}{\point}}{\func{\avecf{j}}{\point}}}}=
\deltaf{i}{j}.
\end{align}
\endef
%%%%%%%%%%%%%%%%%%%%%%%%%%%%%%%%%%%%%%%%%%%%%%%%%%%%%%%%%%%%%%%%%%%%%%%%%%%%%%%%%%%%%%%%%%%%%%%%%%%%%%%%%%%%%%%%%%%%%%%%%%%%%%%%
%%%%%%%%%%%%%%%%%%%%%%%%%%%%%%%%%%%%%%%%%%%%%%%%%%%%%%%%%%%%%%%%%%%%%%%%%%%%%%%%%%%%%%%%%%%%%%%%%%%%%%%%%%%%%%%%%%%%%%%%%%%%%%%%
%%%%%%%%%%%%%%%%%%%%%%%%%%%%%%%%%%%%%%%%%%%%%%%%%%%%%%%%%%%%%%%%%%%%%%%%%%%%%%%%%%%%%%%%%%%%%%%%%%%%%%%%%%%%%%%%%%%%%%%%%%%%%%%%
%%%%%%%%%%%%%%%%%%%%%%%%%%%%%%%%%%%%%%%%%%%%%%%%%%%%%%%%%%%%%%%%%%%%%%%%%%%%%%%%%%%%%%%%%%%%%%%%%%%%%%%%%%%%%%%%%%%%%%%%%%%%%%%%
%%%%%%%%%%%%%%%%%%%%%%%%%%%%%%%%%%%%%%%%%%%%%%%%%%%%%%%%%%%%%%%%%%%%%%%%%%%%%%%%%%%%%%%%%%%%%%%%%%%%%%%%%%%%%%%%%%%%%%%%%%%%%%%%
%%%%%%%%%%%%%%%%%%%%%%%%%%%%%%%%%%%%%%%%%%%%%%%%%%%%%%%%%%%%%%%%%%%%%%%%%%%%%%%%%%%%%%%%%%%%%%%%%%%%%%%%%%%%%%%%%%%%%%%%%%%%%%%%
%%%%%%%%%%%%%%%%%%%%%%%%%%%%%%%%%%%%%%%%%%%%%%%%%%%%%%%%%%%%%%%%%%%%%%%%%%%%%%%%%%%%%%%%%%%%%%%%%%%%%%%%%%%%%%%%%%%%%%%%%%%%%%%%
%%%%%%%%%%%%%%%%%%%%%%%%%%%%%%%%%%%%%%%%%%%%%%%%%%%%%%%%%%%%%%%%%%%%%%%%%%%%%%%%%%%%%%%%%%%%%%%%%%%%%%%%%%%%%%%%%%%%%%%%%%%%%%%%
%%%%%%%%%%%%%%%%%%%%%%%%%%%%%%%%%%%%%%%%%%%%%%%%%%%%%%%%%%%%%%%%%%%%%%%%%%%%%%%%%%%%%%%%%%%%%%%%%%%%%%%%%%%%%%%%%%%%%%%%%%%%%%%%
%%%%%%%%%%%%%%%%%%%%%%%%%%%%%%%%%%%%%%%%%%%%%%%%%%%%%%%%%%%%%%%%%%%%%%%%%%%%%%%%%%%%%%%%%%%%%%%%%%%%%%%%%%%%%%%%%%%%%%%%%%%%%%%%
%%%%%%%%%%%%%%%%%%%%%%%%%%%%%%%%%%%%%%%%%%%%%%%%%%%%%%%%%%%%%%%%%%%%%%%%%%%%%%%%%%%%%%%%%%%%%%%%%%%%%%%%%%%%%%%%%%%%%%%%%%%%%%%%
\section{Local Representations of Metric Tensors}
%%%%%%%%%%%%%%%%%%%%%%%%%%%%%%%%%%%%%%%%%%%%%%%%%%%%%%%%%%%%%%%%%%%%%%%%%%%%%%%%%%%%%%%%%%%%%%%%%%%%%%%%%%%%%%%%%%%%%%%%%%%%%%%%
\fixed
$\Man{}$ is fixed as a manifold with dimension $m$, and $\metrictensor{}$ is fixed as a
semi-Riemannian metric on $\Man{}$.
\endfixed
%%%%%%%%%%%%%%%%%%%%%%%%%%%%%%%%%%%%%%%%%%%%%%%%%%%%%%%%%%%%%%%%%%%%%%%%%%%%%%%%%%%%%%%%%%%%%%%%%%%%%%%%%%%%%%%%%%%%%%%%%%%%%%%%
\textit{
Let $\U$ be a non-empty open subset of the manifold $\Man{}$. The semi-Riemannian metric $\metrictensor{}$ on
$\Man{}$ induces a semi-Riemannian metric on the $\subman{\Man{}}{\U}$ (the open submanifold $\U$ of $\Man{}$ endowed with
its canonical differential structure inherited from that of $\Man{}$) in completely a natural way. We will denote this
canonically induced metric-tensor on $\subman{\Man{}}{\U}$ simply by $\reS{\Man{}}{\U}$.
Actually, for every point $\point\in\U$, the tangent-space of $\subman{\Man{}}{\U}$, that is
$\tanspace{\point}{\subman{\Man{}}{\U}}$, is naturally identified with $\tanspace{\point}{\Man{}}$, and
the value of $\func{\reS{\metrictensor{}}{\U}}{\point}$ at a pair of vectors in $\tanspace{\point}{\subman{\Man{}}{\U}}$
is defined to be the same as the value of $\func{\metrictensor{}}{\point}$ at the naturally corresponded pair of vectors in
$\tanspace{\point}{\Man{}}$. Then, clearly $\opair{\subman{\Man{}}{\U}}{\reS{\metrictensor{}}{\U}}$ is a
semi-Riemannian manifold, which is referred to as a $\quotl$open semi-Riemannian submanifold of
$\opair{\Man{}}{\metrictensor{}}$.\\
Furthermore, it can be easily checked that, given $\avecf{1}$ and $\avecf{2}$ in $\vectorfields{\Man{}}$,
\begin{equation}
\vfmetricprod{\subman{\Man{}}{\U}}{\reS{\metrictensor{}}{\U}}{\reS{\avecf{1}}{\U}}{\reS{\avecf{2}}{\U}}
=\reS{\[\vfmetricprod{\Man{}}{\metrictensor{}}{\avecf{1}}{\avecf{2}}\]}{\U}.
\end{equation}
Similar property holds for flat and sharp operators of the considered metric-tensors.\\
%%%%%%%%%%%%%%%%%%%%
Now, let $\mtuple{\avecf{1}}{\avecf{m}}$ be a local frame field of the vector-bundle
$\Tanbun{\Man{}}$. So each $\avecf{i}$ is a smooth vector field on the restricted manifold $\subman{\Man{}}{\U}$, that is
an element of $\vectorfields{\subman{\Man{}}{\U}}$, and for every $\point\in\U$,
$\mtuple{\func{\avecf{1}}{\point}}{\func{\avecf{m}}{\point}}$ is an ordered-basis of $\tanspace{\point}{\Man{}}$.
Let $\mtuple{\aoneform{1}}{\aoneform{m}}$ be the dual of the local frame field $\mtuple{\avecf{1}}{\avecf{m}}$.
So, each $\aoneform{i}$ is a $\opair{1}{0}$ tensor field (or a $1$-form, or a co-vector field) on the restricted
vector-bundle $\Tanbun{\subman{\Man{}}{\U}}$, that is an element of
$\vbsections{\vbDualbundle{\Tanbun{\subman{\Man{}}{\U}}}}=
\TF{1}{0}{\Tanbun{\subman{\Man{}}{\U}}}=\oneforms{\subman{\Man{}}{\U}}$, and
at every point $\point$ of $\U$, $\mtuple{\func{\aoneform{1}}{\point}}{\func{\aoneform{m}}{\point}}$
is the dual of the ordered-basis $\mtuple{\func{\avecf{1}}{\point}}{\func{\avecf{m}}{\point}}$ of
$\tanspace{\point}{\Man{}}$.\\
Given a semi-Riemannian metric $\metrictensor{}$ on $\Man{}$, clearly,
$\reS{\metrictensor{}}{\U}=\sum_{i=1}^{m}\sum_{j=1}^{m}g_{ij}\aoneform{i}\tensor{}\aoneform{j}$, where the unique
array of real-valued smooth maps $\function{g_{ij}}{\U}{\R}$ is given by,
\begin{equation}
\Foreach{\point}{\U}
\func{g_{ij}}{\point}=\func{\[\func{\metrictensor{}}{\point}\]}
{\binary{\func{\avecf{i}}{\point}}{\func{\avecf{j}}{\point}}}=
\func{\(\vfmetricprod{\Man{}}{\metrictensor{}}{\avecf{i}}{\avecf{j}}\)}{\point},
\end{equation}
for every $i$ and $j$ in $\seta{\suc{1}{m}}$. This is considered to be a local represention of the semi-Riemannian metric
$\metrictensor{}$ in $\U$ with respect to the local frame field $\mtuple{\avecf{1}}{\avecf{m}}$.
Also, note that considering that scalar-products are symmetric, for every $i$ and $j$ $g_{ij}=g_{ji}$.\\
%%%%%%%%%
We know that for every $\point\in\U$, $\func{\metrictensor{}}{\point}$ is a scalar-product on the vector-space
$\tanspace{\point}{\Man{}}$, and hence the matrix representation of $\func{\metrictensor{}}{\point}$
with respect to the ordered-basis $\vsbase{\point}:=\mtuple{\func{\avecf{1}}{\point}}{\func{\avecf{m}}{\point}}$
of $\tanspace{\point}{\Man{}}$, that is $\scalarprodmatrix{\func{\metrictensor{}}{\point}}{\vsbase{\point}}$,
is non-singular $m\times m$ $\R$-matrix. For each $i$ and $j$ in $\seta{\suc{1}{m}}$,
we define the map $\function{{\widetilde{g}}_{ij}}{\U}{\R}$ such that for every $\point\in\U$,
$\func{{\widetilde{g}}_{ij}}{\point}$ be the $ij$-th entry of the matrix
$\finv{\(\scalarprodmatrix{\func{\metrictensor{}}{\point}}{\vsbase{\point}}\)}$. It is a trivial fact that
each ${\widetilde{g}}_{ij}$ must be a smooth map.
We will simply refer to the matrix $\widetilde{g}$ of real-valued functions on $\U$ the
$\quotl$inverse of the matrix $g$ of real-valued functions on $\U$$\quotr$. It is clear that,
\begin{equation}
{\widetilde{g}}_{ij} g_{jk}=\deltaf{i}{k},
\end{equation}
for each $i$ and $k$, where $\function{\deltaf{i}{k}}{\U}{\R}$ denotes the map that assigns to
every point $\point$ of $\Man{}$ the number $1$ if $i=k$, and assigns $0$ to every point otherwise.
The matrix $\widetilde{g}$ is also symmetric, that is,
${\widetilde{g}}_{ij}={\widetilde{g}}_{ji}$ for every $i$ and $j$.\\
%%%%%%%%%
In particular, let $\opair{\U}{\phi}$ be a chart of $\Man{}$, with $\mtuple{\localframevecf{1}}{\localframevecf{m}}$
denoting its associated system of local frame fields, and $\mtuple{\localframeoneform{1}}{\localframeoneform{m}}$
denoting the dual of this local frame field. We then define for each $i$ and $j$ in $\seta{\suc{1}{m}}$,
the smooth maps $\function{\metrictensorchart{\phi}{i}{j}}{\U}{\R}$ as,
\begin{equation}
\Foreach{\point}{\U}
\func{\metrictensorchart{\phi}{i}{j}}{\point}=\func{\[\func{\metrictensor{}}{\point}\]}
{\binary{\func{\localframevecf{i}}{\point}}{\func{\localframevecf{j}}{\point}}},
\end{equation}
which is equivalent to,
\begin{equation}
\metrictensorchart{\phi}{i}{j}=
\vfmetricprod{\subman{\Man{}}{\U}}{\reS{\metrictensor{}}{\U}}{\localframevecf{i}}{\localframevecf{j}}.
\end{equation}
Thus,
\begin{equation}
\reS{\metrictensor{}}{\U}=\sum_{i=1}^{m}\sum_{j=1}^{m}\metrictensorchart{\phi}{i}{j}
\localframeoneform{i}\tensor{}\localframeoneform{j}.
\end{equation}
We will also denote by $\metrictensorchartinv{\phi}{}{}$ the inverse of the matrix
$\metrictensorchart{\phi}{}{}$ of real-valued smooth maps on $\U$.
}
%%%%%%%%%%%%%%%%%%%%%%%%%%%%%%%%%%%%%%%%%%%%%%%%%%%%%%%%%%%%%%%%%%%%%%%%%%%%%%%%%%%%%%%%%%%%%%%%%%%%%%%%%%%%%%%%%%%%%%%%%%%%%%%%
\proposition\label{prometrictensorflatlocalrepresentation}
Let $\U$ be a non-empty open set of $\Man{}$, let $\mtuple{\avecf{1}}{\avecf{m}}$ be a local frame field of the vector-bundle
$\Tanbun{\Man{}}$, and let $\mtuple{\aoneform{1}}{\aoneform{m}}$ be the dual of the local frame field
$\mtuple{\avecf{1}}{\avecf{m}}$.
Let $\reS{\metrictensor{}}{\U}=\sum_{i=1}^{m}\sum_{j=1}^{m}g_{ij}\aoneform{i}\tensor{}\aoneform{j}$
indicate the local representation of $\metrictensor{}$ with respect to $\mtuple{\avecf{1}}{\avecf{m}}$.
Let $\avecff{}\in\vectorfields{\Man{}}$ with the expansion
$\displaystyle\reS{\avecff{}}{\U}=\sum_{i=1}^{m}\veccoff{i}\avecf{i}$, where each $\veccoff{i}$ is a real-valued
smooth map on $\subman{\Man{}}{\U}$.
\begin{equation}
\reS{\[\func{\mtflat{\metrictensor{}}}{\avecff{}}\]}{\U}=
\sum_{i=1}^{m}\sum_{j=1}^{m}
\(g_{ij}\veccoff{i}\)\aoneform{j}.
\end{equation}
\proof
For every $\point\in\U$ and every
$v\in\tanspace{\point}{\Man{}}$ (note that $\tanspace{\point}{\subman{\Man{}}{\U}}$ is naturally
identified with $\tanspace{\point}{\Man{}}$), according to \refdef{defflatanssharpoperatorsofmetrictensors},
and considering that for each $i$,
$\func{\[\func{\aoneform{i}}{\point}\]}{\func{\avecff{}}{\point}}=\func{\veccoff{i}}{\point}$,
\begin{align}
\func{\(\func{\[\func{\mtflat{\metrictensor{}}}{\avecff{}}\]}{\point}\)}{v}&=
\func{\[\func{\metrictensor{}}{\point}\]}{\binary{\func{\avecff{}}{\point}}{v}}\cr
&=\func{\[\sum_{i=1}^{m}\sum_{j=1}^{m}\func{g_{ij}}{\point}\[\func{\aoneform{i}}{\point}\tensor{}
\func{\aoneform{j}}{\point}\]\]}{\binary{\func{\avecff{}}{\point}}{v}}\cr
&=\sum_{i=1}^{m}\sum_{j=1}^{m}\func{g_{ij}}{\point}\cdot\func{\[\func{\aoneform{i}}{\point}\]}{\func{\avecff{}}{\point}}\cdot
\func{\[\func{\aoneform{j}}{\point}\]}{v}\cr
&=\sum_{i=1}^{m}\sum_{j=1}^{m}\func{g_{ij}}{\point}\cdot\func{\veccoff{i}}{\point}\cdot
\func{\[\func{\aoneform{j}}{\point}\]}{v}\cr
&=\func{\(\sum_{i=1}^{m}\sum_{j=1}^{m}\func{\[g_{ij}\veccoff{}\aoneform{j}\]}{\point}\)}{v}\cr
&=\func{\[\func{\(\sum_{i=1}^{m}\sum_{j=1}^{m}g_{ij}\veccoff{}\aoneform{j}\)}{\point}\]}{v}.
\end{align}
Therefore, the desired result is obtained.
\endpro
%%%%%%%%%%%%%%%%%%%%%%%%%%%%%%%%%%%%%%%%%%%%%%%%%%%%%%%%%%%%%%%%%%%%%%%%%%%%%%%%%%%%%%%%%%%%%%%%%%%%%%%%%%%%%%%%%%%%%%%%%%%%%%%%
\proposition\label{prometrictensorsharplocalrepresentation}
Let $\U$ be a non-empty open set of $\Man{}$, let $\mtuple{\avecf{1}}{\avecf{m}}$ be a local frame field of the vector-bundle
$\Tanbun{\Man{}}$, and let $\mtuple{\aoneform{1}}{\aoneform{m}}$ be the dual of the local frame field
$\mtuple{\avecf{1}}{\avecf{m}}$.
Let $\reS{\metrictensor{}}{\U}=\sum_{i=1}^{m}\sum_{j=1}^{m}g_{ij}\aoneform{i}\tensor{}\aoneform{j}$
indicate the local representation of $\metrictensor{}$ with respect to $\mtuple{\avecf{1}}{\avecf{m}}$,
and let $\tilde{g}$ denote the inverse of the matrix $g$ of real-valued functions.
Let $\omega{}\in\oneforms{\Man{}}$ with the expansion
$\displaystyle\reS{\omega{}}{\U}=\sum_{i=1}^{m}\omega_{i}\aoneform{i}$, where each $\omega_{i}$ is a real-valued
smooth map on $\subman{\Man{}}{\U}$.
\begin{equation}
\reS{\[\func{\mtsharp{\metrictensor{}}}{\omega{}}\]}{\U}=
\sum_{i=1}^{m}\sum_{j=1}^{m}
\({\tilde{g}}_{ji}\omega_{i}\)\avecf{j}.
\end{equation}
\proof
Define $\avecff{}:=\func{\mtsharp{\metrictensor{}}}{\omega{}}$, and let
$\displaystyle\reS{\avecff{}}{\U}=\sum_{i=1}^{m}\veccoff{i}{\avecf{i}}$.
According to \refpro{prometrictensorflatlocalrepresentation},
\begin{equation}\label{prometrictensorsharplocalrepresentationeq1}
\reS{\[\func{\mtflat{\metrictensor{}}}{\avecff{}}\]}{\U}=
\sum_{i=1}^{m}\sum_{j=1}^{m}
\(g_{ij}\veccoff{i}\)\aoneform{j}.
\end{equation}
On the other hand, clearly,
\begin{align}\label{prometrictensorsharplocalrepresentationeq2}
\reS{\[\func{\mtflat{\metrictensor{}}}{\avecff{}}\]}{\U}=\reS{\omega}{\U}=
\sum_{i=1}^{m}\omega_{i}\aoneform{i}.
\end{align}
Therefor, for every $i$,
\begin{equation}\label{prometrictensorsharplocalrepresentationeq3}
\omega_{i}=\sum_{j=1}^{m}g_{ji}\veccoff{j},
\end{equation}
and therefore, for every $l\in\seta{\suc{1}{m}}$,
\begin{align}\label{prometrictensorsharplocalrepresentationeq4}
\sum_{i=1}^{m}{\tilde{g}}_{li}\omega_{i}&=
\sum_{j=1}^{m}\(\sum_{i=1}^{m}{\tilde{g}}_{li}g_{ji}\)\veccoff{j}\cr
&=\sum_{j=1}^{m}\deltaf{j}{l}\veccoff{j}\cr
&=\veccoff{l}.
\end{align}
Thus, combining \Ref{prometrictensorsharplocalrepresentationeq1} and
\Ref{prometrictensorsharplocalrepresentationeq4}, it becomes evident that,
\begin{equation}\label{prometrictensorsharplocalrepresentationeq5}
\reS{\[\func{\mtsharp{\metrictensor{}}}{\omega{}}\]}{\U}=
\sum_{i=1}^{m}\sum_{j=1}^{m}
\({\tilde{g}}_{ji}\omega_{i}\)\avecf{j}.
\end{equation}
\endpro
%%%%%%%%%%%%%%%%%%%%%%%%%%%%%%%%%%%%%%%%%%%%%%%%%%%%%%%%%%%%%%%%%%%%%%%%%%%%%%%%%%%%%%%%%%%%%%%%%%%%%%%%%%%%%%%%%%%%%%%%%%%%%%%%
\corollary
Let $f\in\smoothmaps{\Man{}}$.
Let $\U$ be a non-empty open set of $\Man{}$, let $\mtuple{\avecf{1}}{\avecf{m}}$ be a local frame field of the vector-bundle
$\Tanbun{\Man{}}$, and let $\mtuple{\aoneform{1}}{\aoneform{m}}$ be the dual of the local frame field
$\mtuple{\avecf{1}}{\avecf{m}}$.
Let $\reS{\metrictensor{}}{\U}=\sum_{i=1}^{m}\sum_{j=1}^{m}g_{ij}\aoneform{i}\tensor{}\aoneform{j}$
indicate the local representation of $\metrictensor{}$ with respect to $\mtuple{\avecf{1}}{\avecf{m}}$,
and let $\tilde{g}$ denote the inverse of the matrix $g$ of real-valued functions.
\begin{equation}
\reS{\[\func{\grad{\metrictensor{}}}{f}\]}{\U}=
\sum_{i=1}^{m}\sum_{j=1}^{m}
\({\tilde{g}}_{ji}\lieder{\avecf{i}}{f}\)\avecf{j}.
\end{equation}
\proof
According to \refdef{defgradient}, \refpro{prometrictensorsharplocalrepresentation},
and considering that $\displaystyle\reS{\[\Rderivative{\Man{}}{f}\]}{\U}=\sum_{i=1}^{m}
\(\lieder{\avecf{i}}{f}\)\aoneform{i}$, it is clear.
\endcor
%%%%%%%%%%%%%%%%%%%%%%%%%%%%%%%%%%%%%%%%%%%%%%%%%%%%%%%%%%%%%%%%%%%%%%%%%%%%%%%%%%%%%%%%%%%%%%%%%%%%%%%%%%%%%%%%%%%%%%%%%%%%%%%%
%%%%%%%%%%%%%%%%%%%%%%%%%%%%%%%%%%%%%%%%%%%%%%%%%%%%%%%%%%%%%%%%%%%%%%%%%%%%%%%%%%%%%%%%%%%%%%%%%%%%%%%%%%%%%%%%%%%%%%%%%%%%%%%%
%%%%%%%%%%%%%%%%%%%%%%%%%%%%%%%%%%%%%%%%%%%%%%%%%%%%%%%%%%%%%%%%%%%%%%%%%%%%%%%%%%%%%%%%%%%%%%%%%%%%%%%%%%%%%%%%%%%%%%%%%%%%%%%%
%%%%%%%%%%%%%%%%%%%%%%%%%%%%%%%%%%%%%%%%%%%%%%%%%%%%%%%%%%%%%%%%%%%%%%%%%%%%%%%%%%%%%%%%%%%%%%%%%%%%%%%%%%%%%%%%%%%%%%%%%%%%%%%%
%%%%%%%%%%%%%%%%%%%%%%%%%%%%%%%%%%%%%%%%%%%%%%%%%%%%%%%%%%%%%%%%%%%%%%%%%%%%%%%%%%%%%%%%%%%%%%%%%%%%%%%%%%%%%%%%%%%%%%%%%%%%%%%%
%%%%%%%%%%%%%%%%%%%%%%%%%%%%%%%%%%%%%%%%%%%%%%%%%%%%%%%%%%%%%%%%%%%%%%%%%%%%%%%%%%%%%%%%%%%%%%%%%%%%%%%%%%%%%%%%%%%%%%%%%%%%%%%%
%%%%%%%%%%%%%%%%%%%%%%%%%%%%%%%%%%%%%%%%%%%%%%%%%%%%%%%%%%%%%%%%%%%%%%%%%%%%%%%%%%%%%%%%%%%%%%%%%%%%%%%%%%%%%%%%%%%%%%%%%%%%%%%%
%%%%%%%%%%%%%%%%%%%%%%%%%%%%%%%%%%%%%%%%%%%%%%%%%%%%%%%%%%%%%%%%%%%%%%%%%%%%%%%%%%%%%%%%%%%%%%%%%%%%%%%%%%%%%%%%%%%%%%%%%%%%%%%%
%%%%%%%%%%%%%%%%%%%%%%%%%%%%%%%%%%%%%%%%%%%%%%%%%%%%%%%%%%%%%%%%%%%%%%%%%%%%%%%%%%%%%%%%%%%%%%%%%%%%%%%%%%%%%%%%%%%%%%%%%%%%%%%%
%%%%%%%%%%%%%%%%%%%%%%%%%%%%%%%%%%%%%%%%%%%%%%%%%%%%%%%%%%%%%%%%%%%%%%%%%%%%%%%%%%%%%%%%%%%%%%%%%%%%%%%%%%%%%%%%%%%%%%%%%%%%%%%%
%%%%%%%%%%%%%%%%%%%%%%%%%%%%%%%%%%%%%%%%%%%%%%%%%%%%%%%%%%%%%%%%%%%%%%%%%%%%%%%%%%%%%%%%%%%%%%%%%%%%%%%%%%%%%%%%%%%%%%%%%%%%%%%%
\section{Levi-Civita Connection}
%%%%%%%%%%%%%%%%%%%%%%%%%%%%%%%%%%%%%%%%%%%%%%%%%%%%%%%%%%%%%%%%%%%%%%%%%%%%%%%%%%%%%%%%%%%%%%%%%%%%%%%%%%%%%%%%%%%%%%%%%%%%%%%%
\fixed
$\Man{}$ is fixed as a manifold with dimension $m$, and $\metrictensor{}$ is fixed as a
semi-Riemannian metric on $\Man{}$, that is an element of $\metrictensors{\Man{}}{}$.
Additionally, $\Man{1}$ is fixed as a manifold with dimension $m_1$, and $\metrictensor{1}$ is fixed as a
semi-Riemannian metric on $\Man{1}$, that is an element of $\metrictensors{\Man{1}}{}$.
\endfixed
%%%%%%%%%%%%%%%%%%%%%%%%%%%%%%%%%%%%%%%%%%%%%%%%%%%%%%%%%%%%%%%%%%%%%%%%%%%%%%%%%%%%%%%%%%%%%%%%%%%%%%%%%%%%%%%%%%%%%%%%%%%%%%%%
\definition\label{deftorsionofaffineconnections}
Let $\connection{}$ be an affine connection on $\Man{}$, that is an element of $\connections{\Tanbun{\Man{}}}$.
The operator $\function{\Torsion{\connection{}}}{\Cprod{\vectorfields{\Man{}}}{\vectorfields{\Man{}}}}{\vectorfields{\Man{}}}$
is defined as,
\begin{equation}
\Foreach{\opair{\avecf{1}}{\avecf{2}}}{\Cprod{\vectorfields{\Man{}}}{\vectorfields{\Man{}}}}
\func{\Torsion{\connection{}}}{\binary{\avecf{1}}{\avecf{2}}}:=
\con{\avecf{1}}{\avecf{2}}-\con{\avecf{2}}{\avecf{1}}-\liebracket{\avecf{1}}{\avecf{2}}{\Man{}},
\end{equation}
where, $\liebracket{\avecf{1}}{\avecf{2}}{\Man{}}$ indicates the Lie-bracket of vector fields $\avecf{1}$
and $\avecf{2}$ on $\Man{}$.
$\Torsion{\connection{}}$ is referred to as the $\quotl$ torsion (operator) of the affine connection $\connection{}$
on the manifold $\Man{}$$\quotr$.
\endef
%%%%%%%%%%%%%%%%%%%%%%%%%%%%%%%%%%%%%%%%%%%%%%%%%%%%%%%%%%%%%%%%%%%%%%%%%%%%%%%%%%%%%%%%%%%%%%%%%%%%%%%%%%%%%%%%%%%%%%%%%%%%%%%%
\definition\label{defmetricconnections}
Let $\connection{}$ be an affine connection on the semi-Riemannian manifold $\opair{\Man{}}{\metrictensor{}}$,
that is an element of $\connections{\Tanbun{\Man{}}}$.
$\connection{}$ is called a $\quotl$metric-connection on the semi-Riemannian manifold $\opair{\Man{}}{\metrictensor{}}$$\quotr$
if for every smooth vector fields $\avecf{1}$, $\avecf{2}$, and $\avecf{3}$ on $\Man{}$,
\begin{equation}
\lieder{\avecf{1}}{\(\vfmetricprod{\Man{}}{\metrictensor{}}{\avecf{2}}{\avecf{3}}\)}=
\vfmetricprod{\Man{}}{\metrictensor{}}{\con{\avecf{1}}{\avecf{2}}}{\avecf{3}}+
\vfmetricprod{\Man{}}{\metrictensor{}}{\con{\avecf{1}}{\avecf{3}}}{\avecf{2}}.
\end{equation}
\endef
%%%%%%%%%%%%%%%%%%%%%%%%%%%%%%%%%%%%%%%%%%%%%%%%%%%%%%%%%%%%%%%%%%%%%%%%%%%%%%%%%%%%%%%%%%%%%%%%%%%%%%%%%%%%%%%%%%%%%%%%%%%%%%%%
\lemma\label{lemequivalentpropertiesofametricconnection}
Let $\connection{}$ be an affine connection on the semi-Riemannian manifold $\opair{\Man{}}{\metrictensor{}}$.
The following assertions are equivalent.
\begin{itemize}
\item[\myitem{1.~}]
$\connection{}$ is a metric-connection on $\opair{\Man{}}{\metrictensor{}}$.
\item[\myitem{2.~}]
For every smooth curve $\function{\Curve{}}{\interval{}}{\Man{}}$ on $\Man{}$, and every smooth vector fields
$\valongc{1}$ and $\valongc{2}$ along $\Curve{}$, that is any pair $\valongc{1}$
and $\valongc{2}$ in $\valongcurves{\Man{}}{\Curve{}}$,
\begin{align}
\timeder\vfmetricproduct{\valongc{1}}{\valongc{2}}=
\vfmetricproduct{\func{\covder{\connection{}}{\Curve{}}}{\valongc{1}}}{\valongc{2}}+
\vfmetricproduct{\valongc{1}}{\func{\covder{\connection{}}{\Curve{}}}{\valongc{2}}}.
\end{align}
Here, the smooth map $\function{\vfmetricproduct{\valongc{1}}{\valongc{2}}}{\interval{}}{\R}$ is defined by,
\begin{equation}
\Foreach{t}{\interval{}}
\func{\vfmetricproduct{\valongc{1}}{\valongc{2}}}{t}\eqdef
\func{\[\func{\metrictensor{}}{\func{\Curve{}}{t}}\]}{\binary{\func{\valongc{1}}{t}}{\func{\valongc{2}}{t}}}.
\end{equation}
Also, $\timeder$ denotes the ordinary differential operator of one-real-variable real-valued smooth functions.
\item[\myitem{3.~}]
For every smooth curve $\function{\Curve{}}{\interval{}}{\Man{}}$ on $\Man{}$, and every parallel vector fields
$\valongc{1}$ and $\valongc{2}$ along $\Curve{}$ relative to $\connection{}$, that is any pair $\valongc{1}$
and $\valongc{2}$ in $\parallelvalongcurves{\connection{}}{\Curve{}}$,
\begin{equation}
\timeder\vfmetricproduct{\valongc{1}}{\valongc{2}}=0,
\end{equation}
and hence $\vfmetricproduct{\valongc{1}}{\valongc{2}}$ is constant.
\item[\myitem{4.~}]
For any smooth curve $\function{\Curve{}}{\interval{}}{\Man{}}$ on $\Man{}$, and any $\binary{t_0}{t_1}\in\interval{}$,
the parallel transport operator of $\Curve{}$ (relative to $\connection{}$) from $t_0$ to $t_1$, that is
$\ptransport{\Curve{}}{t_0}{t_1}{\connection{}}$, is an isometry from the scalar-product space
$\opair{\Tanspace{\func{\Curve{}}{t_0}}{\Man{}}}{\func{\metrictensor{}}{\func{\Curve{}}{t_0}}}$
to the scalar-product space
$\opair{\Tanspace{\func{\Curve{}}{t_1}}{\Man{}}}{\func{\metrictensor{}}{\func{\Curve{}}{t_1}}}$.
\item[\myitem{5.~}]
For every smooth curve $\function{\Curve{}}{\interval{}}{\Man{}}$ on $\Man{}$, every $t\in\interval{}$,
and every orthonormal ordered-basis $\vsbase{}=\mtuple{\vsbase{1}}{\vsbase{m}}$ of
$\Tanspace{\func{\Curve{}}{t}}{\Man{}}$, there exists an $m$-tuple
$\mtuple{\valongc{1}}{\valongc{m}}$ of parallel vector fields along $\Curve{}$ (elative to $\connection{}$)
such that for every $i$, $\func{\valongc{i}}{t}=\vsbase{i}$ and for every $s\in\interval{}$
$\mtuple{\func{\valongc{1}}{s}}{\func{\valongc{m}}{s}}$ is an orthonormal ordered basis of
$\opair{\Tanspace{\func{\Curve{}}{s}}{\Man{}}}{\func{\metrictensor{}}{\point}}$.
\end{itemize}
%%%%%%%%%%%%%%%%%%%%%%%%%%%%%
%%%%%%%%%%%%%%%%%%%%%%%%%%%%%
\proof
\begin{itemize}
\item[\myitem{1}$\to$\myitem{2.~}]
Suppose that $\connection{}$ is a metric-connection on $\opair{\Man{}}{\metrictensor{}}$, and thus for every
$\triplet{\avecf{1}}{\avecf{2}}{\avecf{3}}\in\vectorfields{\Man{}}$,
\begin{equation}
\lieder{\avecf{1}}{\(\vfmetricprod{\Man{}}{\metrictensor{}}{\avecf{2}}{\avecf{3}}\)}=
\vfmetricprod{\Man{}}{\metrictensor{}}{\con{\avecf{1}}{\avecf{2}}}{\avecf{3}}+
\vfmetricprod{\Man{}}{\metrictensor{}}{\con{\avecf{1}}{\avecf{3}}}{\avecf{2}}.
\end{equation}
Let $t$ be an arbitrary element of $\interval{}$. Choose a chart $\opair{\U}{\phi}$ around the point
$\func{\Curve{}}{t}$, and let $\mtuple{\localframevecf{1}}{\localframevecf{m}}$ denote the system of
local frame fields corresponded to this chart. Let $O$ be an open sub-interval of $\interval{}$ containing
$t$ such that $\func{\image{\Curve{}}}{O}\subseteq\U$. Define the systems of smooth maps
$\mtuple{f_1}{f_m}$ and $\mtuple{h_1}{h_m}$ from $O$ to $\R$ as,
\begin{align}
\reS{\valongc{1}}{O}&=\sum_{i=1}^{m}f_i\(\cmp{\localframevecf{i}}{\Curve{}}\),\cr
\reS{\valongc{2}}{O}&=\sum_{i=1}^{m}h_i\(\cmp{\localframevecf{i}}{\Curve{}}\).
\end{align}
It can be easily checked that,
\begin{align}
\vfmetricproduct{\valongc{1}}{\valongc{2}}&=
\vfmetricproduct{\sum_{i=1}^{m}f_i\(\cmp{\localframevecf{i}}{\Curve{}}\)}{\sum_{i=1}^{m}h_i\(\cmp{\localframevecf{i}}{\Curve{}}\)}\cr
&=\sum_{i=1}^{m}\sum_{j=1}^{m}f_ih_j\vfmetricproduct{\localframevecf{i}}{\localframevecf{j}}.
\end{align}
Thus,
\begin{align}
\func{\(\timeder\vfmetricproduct{\valongc{1}}{\valongc{2}}\)}{t}&=
\sum_{i=1}^{m}\sum_{j=1}^{m}\func{\[\timeder\(f_ih_j\vfmetricproduct{\cmp{\localframevecf{i}}{\Curve{}}}{\cmp{\localframevecf{j}}{\Curve{}}}\)\]}{t}\cr
&=\sum_{i=1}^{m}\sum_{j=1}^{m}
\func{\[\({\dot{f}}_{i}h_j+f_i{\dot{h}}_{j}\)\vfmetricproduct{\cmp{\localframevecf{i}}{\Curve{}}}{\cmp{\localframevecf{j}}{\Curve{}}}\]}{t}+\func{\[f_ih_j
\(\timeder\vfmetricproduct{\cmp{\localframevecf{i}}{\Curve{}}}{\cmp{\localframevecf{j}}{\Curve{}}}\)\]}{t}.\cr
&{}
\end{align}
Choose a $\avecf{}\in\vectorfields{\subman{\Man{}}{\U}}$ such that $\func{\avecf{}}{\func{\Curve{}}{t}}=\func{\dot{\Curve{}}}{t}$.
Considering the pointwise action of Lie-derivatives, it can be easily seen that,
\begin{align}
\func{\(\timeder\vfmetricproduct{\cmp{\localframevecf{i}}{\Curve{}}}{\cmp{\localframevecf{j}}{\Curve{}}}\)}{t}&=
\func{\[\lieder{\avecf{}}{\vfmetricproduct{\localframevecf{i}}{\localframevecf{j}}}\]}{\func{\Curve{}}{t}}.
\end{align}
In addition, according to \refdef{defmetricconnections}, for every $i$ and $j$,
\begin{align}
\lieder{\avecf{}}{\vfmetricproduct{\localframevecf{i}}{\localframevecf{j}}}=
\vfmetricproduct{\con{\avecf{}}{\localframevecf{i}}}{\localframevecf{j}}+
\vfmetricproduct{\localframevecf{i}}{\con{\avecf{}}{\localframevecf{j}}},
\end{align}
and thus according to \refdef{defcovarianderivativealongcurves},
\begin{align}
\func{\[\lieder{\avecf{}}{\vfmetricproduct{\localframevecf{i}}{\localframevecf{j}}}\]}{\func{\Curve{}}{t}}&=
\vfmetricproduct{\func{\(\con{\avecf{}}{\localframevecf{i}}\)}{\func{\Curve{}}{t}}}{\func{\localframevecf{j}}{\func{\Curve{}}{t}}}+
\vfmetricproduct{\func{\localframevecf{i}}{\func{\Curve{}}{t}}}{\func{\(\con{\avecf{}}{\localframevecf{j}}\)}{\func{\Curve{}}{t}}}\cr
&=\vfmetricproduct{\tbcon{\func{\dot{\Curve{}}}{t}}{\localframevecf{i}}}{\func{\localframevecf{j}}{\func{\Curve{}}{t}}}+
\vfmetricproduct{\func{\localframevecf{i}}{\func{\Curve{}}{t}}}{\tbcon{\func{\dot{\Curve{}}}{t}}{\localframevecf{j}}}\cr
&=\vfmetricproduct{\func{\[\func{\covder{\connection{}}{\Curve{}}}{\cmp{\localframevecf{i}}{\Curve{}}}\]}{t}}{\func{\localframevecf{j}}{\func{\Curve{}}{t}}}+
\vfmetricproduct{\func{\localframevecf{i}}{\func{\Curve{}}{t}}}{\func{\[\func{\covder{\connection{}}{\Curve{}}}{\cmp{\localframevecf{j}}{\Curve{}}}\]}{t}}.\cr
&{}
\end{align}
Therefore, according to \refdef{defcovarianderivativealongcurves},
\begin{align}
&~~~~~\func{\(\timeder\vfmetricproduct{\valongc{1}}{\valongc{2}}\)}{t}\cr
&=\sum_{i=1}^{m}\sum_{j=1}^{m}
\func{\[\vfmetricproduct{\dot{f}_i\(\cmp{\localframevecf{i}}{\Curve{}}\)+\func{\covder{\connection{}}{\Curve{}}}{\cmp{\localframevecf{i}}{\Curve{}}}}
{h_j\(\cmp{\localframevecf{j}}{\Curve{}}\)}+
\vfmetricproduct{f_i\(\cmp{\localframevecf{i}}{\Curve{}}\)}
{\dot{h}_j\(\cmp{\localframevecf{j}}{\Curve{}}\)+\func{\covder{\connection{}}{\Curve{}}}{\cmp{\localframevecf{j}}{\Curve{}}}}\]}{t}\cr
&=\sum_{i=1}^{m}\sum_{j=1}^{m}
\func{\[\vfmetricproduct{\func{\covder{\connection{}}{\Curve{}}}{\cmp{f_i\localframevecf{i}}{\Curve{}}}}{\cmp{h_j\localframevecf{j}}{\Curve{}}}+
\vfmetricproduct{\cmp{f_i\localframevecf{i}}{\Curve{}}}{\func{\covder{\connection{}}{\Curve{}}}{\cmp{h_j\localframevecf{j}}{\Curve{}}}}\]}{t}\cr
&=\func{\[\vfmetricproduct{\func{\covder{\connection{}}{\Curve{}}}{\sum_{i=1}^{m}\cmp{f_i\localframevecf{i}}{\Curve{}}}}{\sum_{j=1}^{m}\cmp{h_j\localframevecf{j}}{\Curve{}}}+
\vfmetricproduct{\sum_{i=1}^{m}\cmp{f_i\localframevecf{i}}{\Curve{}}}{\func{\covder{\connection{}}{\Curve{}}}{\sum_{j=1}^{m}\cmp{h_j\localframevecf{j}}{\Curve{}}}}\]}{t}\cr
&=\func{\[\vfmetricproduct{\func{\covder{\connection{}}{\Curve{}}}{\valongc{1}}}{\valongc{2}}+
\vfmetricproduct{\valongc{1}}{\func{\covder{\connection{}}{\Curve{}}}{\valongc{2}}}\]}{t}.
\end{align}
%%%%%%%%%%%%%%%%%%%%%%%%
\item[\myitem{2}$\to$\myitem{3.~}]
It is an immediate consequence of the definition of the parallel vector fields along a curve relative to a connection,
that is \refdef{defparallelvectorfieldsalongcurves}.
%%%%%%%%%%%%%%%%%%%%%%%%
\item[\myitem{3}$\to$\myitem{4.~}]
Let $v_1$ and $v_2$ be arbitrary elements of $\tanspace{\func{\Curve{}}{t_0}}{\Man{}}$. Let $V_1$ and
$V_2$ denote the unique parallel vector fields along $\Curve{}$ (relative to $\connection{}$) with initial values
$\func{V_1}{t_0}=v_1$ and $\func{V_2}{t_0}=v_2$, respectively. Then, according to \refdef{defparalleltransport},
$\func{\[\ptransport{\Curve{}}{t_0}{t_1}{\connection{}}\]}{v_1}=\func{V_1}{t_1}$
and $\func{\[\ptransport{\Curve{}}{t_0}{t_1}{\connection{}}\]}{v_2}=\func{V_2}{t_1}$.
But according to the assumption of \myitem{3}, and considering that $V_1$ and $V_2$ are parallel vector fields along $\Curve{}$,
it is clear that $\vfmetricproduct{V_1}{V_2}$ is constant, and hence
\begin{align}
&~~~~\func{\[\func{\metrictensor{}}{\func{\Curve{}}{t_1}}\]}{\binary{\func{\[\ptransport{\Curve{}}{t_0}{t_1}{\connection{}}\]}{v_1}}
{\func{\[\ptransport{\Curve{}}{t_0}{t_1}{\connection{}}\]}{v_2}}}\cr
&=\func{\[\func{\metrictensor{}}{\func{\Curve{}}{t_1}}\]}{\binary{\func{V_1}{t_1}}{\func{V_2}{t_1}}}\cr
&=\func{\vfmetricproduct{V_1}{V_2}}{t_1}\cr
&=\func{\vfmetricproduct{V_1}{V_2}}{t_0}\cr
&=\func{\[\func{\metrictensor{}}{\func{\Curve{}}{t_0}}\]}{\binary{\func{V_1}{t_0}}{\func{V_2}{t_0}}}\cr
&=\func{\[\func{\metrictensor{}}{\func{\Curve{}}{t_0}}\]}{\binary{v_1}{v_2}}.
\end{align}
%%%%%%%%%%%%%%%%%%%%%%%%
\item[\myitem{4}$\to$\myitem{5.~}]
Fix a smooth curve $\function{\Curve{}}{\interval{}}{\Man{}}$ and some $t\in\interval{}$,
and let $\mtuple{\vsbase{1}}{\vsbase{m}}$ be an orthonormal ordered-basis of
$\Tanspace{\func{\Curve{}}{t}}{\Man{}}$. For every $i$, let $\avecf{i}$ be the unique parallel vector field along
$\Curve{}$ relative to $\connection{}$ with the initial condition $\func{\valongc{i}}{t}=\vsbase{i}$.
For every $s\in\interval{}$, according to the definition of the parallel transport along smooth curves,
$\mtuple{\func{\valongc{1}}{s}}{\func{\valongc{m}}{s}}$ equals
$\mtuple{\func{\[\ptransport{\Curve{}}{t}{s}{\connection{}}\]}{\vsbase{1}}}{\func{\[\ptransport{\Curve{}}{t}{s}{\connection{}}\]}{\vsbase{m}}}$.
So, based on the assumption of \myitem{4}, and considering that an isometry between scalar-product-spaces sends an
orthonormal ordered-basis to an orthonormal ordered-basis, $\mtuple{\func{\valongc{1}}{s}}{\func{\valongc{m}}{s}}$
must be an orthonormal ordered basis of $\opair{\Tanspace{\func{\Curve{}}{s}}{\Man{}}}{\func{\metrictensor{}}{\point}}$.
%%%%%%%%%%%%%%%%%%%%%%%%
\item[\myitem{5}$\to$\myitem{1.~}]
This can be split into parts \myitem{5}$\to$\myitem{2} and \myitem{2}$\to$\myitem{1}.
These parts can be easily verified,
based on \refdef{defcovariantderivativeinducedbyconnection} and \refthm{thmcovariantderivativeintermsofaparallelframe}.
\end{itemize}
\endlem
%%%%%%%%%%%%%%%%%%%%%%%%%%%%%%%%%%%%%%%%%%%%%%%%%%%%%%%%%%%%%%%%%%%%%%%%%%%%%%%%%%%%%%%%%%%%%%%%%%%%%%%%%%%%%%%%%%%%%%%%%%%%%%%%
\textit{
Let $\seta{f_{\avecff{}}\in\smoothmaps{\Man{}}}_{\avecff{}\in\vectorfields{\Man{}}}$ be a collection of real-valued smooth maps
on $\Man{}$ index by the set of all smooth vector fields on $\Man{}$.
First, we want to show that the collection of equations
\begin{equation}
\Foreach{\avecff{}}{\vectorfields{\Man{}}}
\vfmetricproduct{\avecf{}}{\avecff{}}=f_{\avecff{}},
\end{equation}
has at most one solution for $\avecf{}$.\\
Suppose that one $\avecf{}\in\vectorfields{\Man{}}$ exists that satisfies this condition, and let
$\aoneform{}:=\func{\mtflat{\metrictensor{}}}{\avecf{}}$. Then, since
$\avecf{}=\func{\mtsharp{\metrictensor{}}}{\aoneform{}}$,
\begin{equation}
\Foreach{\avecff{}}{\vectorfields{\Man{}}}
\vfmetricproduct{\func{\mtsharp{\metrictensor{}}}{\aoneform{}}}{\avecff{}}=f_{\avecff{}},
\end{equation}
and hence according to \refpro{prometricinducedproductoneforms},
\begin{equation}
\Foreach{\avecff{}}{\vectorfields{\Man{}}}
\func{\TFequiv{\aoneform{}}}{\avecff{}}=
f_{\avecff{}}.
\end{equation}
Obviously, $\aoneform{}$ must be the only element of $\oneforms{\Man{}}$ that satisfies the condition above,
and therefore $\avecf{}$ must be uniquely recognized.\\
Second, considering the arguments above, it is obvious that the mentioned collection of equations
has a solution if and only if the assignment
$\vectorfields{\Man{}}\ni\avecff{}\mapsto f_{\avecff{}}$ is $\smoothmaps{\Man{}}$-linear;
then, the unique solution is $\func{\mtsharp{\metrictensor{}}}{\aoneform{}}$, $\aoneform{}$ indicating the
unique element of $\TF{1}{0}{\Tanbun{\Man{}}}$ naturally corresponded to the element
$\vectorfields{\Man{}}\ni\avecff{}\mapsto\f_{\avecff{}}$ of $\TTF{1}{0}{\vectorfields{\Man{}}}$.
}
%%%%%%%%%%%%%%%%%%%%%%%%%%%%%%%%%%%%%%%%%%%%%%%%%%%%%%%%%%%%%%%%%%%%%%%%%%%%%%%%%%%%%%%%%%%%%%%%%%%%%%%%%%%%%%%%%%%%%%%%%%%%%%%%
\theorem\label{thmfundamentaltheoremofriemanniangeometry}
There exists a unique affine connection $\connection{}$ on $\Man{}$,
which is a metric-connection on $\Man{}$ with vanishing torsion
($\Torsion{\connection{}}=0$), that is, for every triple $\triple{\avecf{1}}{\avecf{2}}{\avecf{3}}$ of
smooth vector fields on $\Man{}$,
\begin{align}\label{thmfundamentaltheoremofriemanniangeometryeq}
&\lieder{\avecf{1}}{\(\vfmetricprod{\Man{}}{\metrictensor{}}{\avecf{2}}{\avecf{3}}\)}=
\vfmetricprod{\Man{}}{\metrictensor{}}{\con{\avecf{1}}{\avecf{2}}}{\avecf{3}}+
\vfmetricprod{\Man{}}{\metrictensor{}}{\con{\avecf{1}}{\avecf{3}}}{\avecf{2}},\cr
&\con{\avecf{1}}{\avecf{2}}-\con{\avecf{2}}{\avecf{1}}-\liebracket{\avecf{1}}{\avecf{2}}{\Man{}}=\zerovec{\vectorfields{\Man{}}}.
\end{align}
where $\zerovec{\vectorfields{\Man{}}}$ indicates the trivial vector field on $\Man{}$ which assigns to every point
$\point$ of $\Man{}$ the zero element of the tangent space $\Tanspace{\point}{\Man{}}$.
\proof
\begin{itemize}
\item[\myitem{pr-1.}]
First, we show that if there exists a metric-connection on $\Man{}$ without torsion, then it must be
the only such affine connection on $\Man{}$.\\
Suppose that $\connection{}$ is an affine connection on $\Man{}$ that satisfies the conditions of
\Ref{thmfundamentaltheoremofriemanniangeometryeq}. Then, using these conditions simultaneously,
it can be easily verified that for every triple $\triple{\avecf{1}}{\avecf{2}}{\avecf{3}}$
of smooth vector fields on $\Man{}$,
\begin{align}\label{thmfundamentaltheoremofriemanniangeometryeq1}
2\vfmetricproduct{\con{\avecf{1}}{\avecf{2}}}{\avecf{3}}=&~~~~
\lieder{\avecf{1}}{\vfmetricproduct{\avecf{2}}{\avecf{3}}}+
\lieder{\avecf{2}}{\vfmetricproduct{\avecf{3}}{\avecf{1}}}-
\lieder{\avecf{3}}{\vfmetricproduct{\avecf{1}}{\avecf{2}}}\cr
&+\vfmetricproduct{\avecf{3}}{\liebracket{\avecf{1}}{\avecf{2}}{}}-
\vfmetricproduct{\avecf{2}}{\liebracket{\avecf{1}}{\avecf{3}}{}}-
\vfmetricproduct{\avecf{1}}{\liebracket{\avecf{2}}{\avecf{3}}{}}.
\end{align}
Fixing $\avecf{1}$ and $\avecf{2}$, the above equation is satisfied for every $\avecf{3}\in\vectorfields{\Man{}}$,
whose righ-hand-side is an element of $\smoothmaps{\Man{}}$ dependent on $\avecf{3}$. Thus, $\con{\avecf{1}}{\avecf{2}}$
must be uniquely recognized. Since the choice of the pair of smooth vector fields $\avecf{1}$ and $\avecf{2}$
is arbitrary, the connection operator $\connection{}$ has to be uniquely recognized.
\item[\myitem{pr-2.}]
Fix arbitrary $\avecf{1}$ and $\avecf{2}$ in $\vectorfields{\Man{}}$. Since the right-hand-side of
\Ref{thmfundamentaltheoremofriemanniangeometryeq1} is $\smoothmaps{\Man{}}$-linear with respect to the variable
$\avecf{3}$, there must exist a unique $\con{\avecf{1}}{\avecf{2}}\in\vectorfields{\Man{}}$ satisfying
\Ref{thmfundamentaltheoremofriemanniangeometryeq1} for every $\avecf{3}\in\vectorfields{\Man{}}$.
Since the choice of the pair $\avecf{1}$ and $\avecf{2}$ of
smooth vector fields on $\Man{}$ was arbitrary, we conclude that there exists a unique operator
$\function{\connection{}}{\Cprod{\vectorfields{\Man{}}}{\vectorfields{\Man{}}}}{\vectorfields{\Man{}}}$
such that \Ref{thmfundamentaltheoremofriemanniangeometryeq1} holds for every triple $\avecf{1}$, $\avecf{2}$,
and $\avecf{3}$ of smooth vector fields on $\Man{}$.\\
Moreover, it can be easily checked that this unique operator $\connection{}$ satisfies the axioms of an affine connection
on a manifold, alongside the conditions of \Ref{thmfundamentaltheoremofriemanniangeometryeq}.
\end{itemize}
\endthm
%%%%%%%%%%%%%%%%%%%%%%%%%%%%%%%%%%%%%%%%%%%%%%%%%%%%%%%%%%%%%%%%%%%%%%%%%%%%%%%%%%%%%%%%%%%%%%%%%%%%%%%%%%%%%%%%%%%%%%%%%%%%%%%%
\definition\label{defLeviCivitaconnection}
The operator $\function{\LCconnection{\Man{}}{\metrictensor{}}}
{\Cprod{\vectorfields{\Man{}}}{\vectorfields{\Man{}}}}{\vectorfields{\Man{}}}$ is defined to be the unique
affine connection on $\Man{}$, which is a metric-connection on $\Man{}$ with
$\Torsion{\LCconnection{\Man{}}{\metrictensor{}}}=0$. $\LCconnection{\Man{}}{\metrictensor{}}$ is referred to as the
$\quotl$Levi-Civita connection on the semi-Riemannian manifold $\opair{\Man{}}{\metrictensor{}}$$\quotr$.
\endef
%%%%%%%%%%%%%%%%%%%%%%%%%%%%%%%%%%%%%%%%%%%%%%%%%%%%%%%%%%%%%%%%%%%%%%%%%%%%%%%%%%%%%%%%%%%%%%%%%%%%%%%%%%%%%%%%%%%%%%%%%%%%%%%%
\corollary\label{corKoszulFormula}
For every triple $\triple{\avecf{1}}{\avecf{2}}{\avecf{3}}$,
\begin{align}\label{KsozulFormula}
2\vfmetricproduct{\LCcon{\Man{}}{\metrictensor{}}{\avecf{1}}{\avecf{2}}}{\avecf{3}}=&~~~~
\lieder{\avecf{1}}{\vfmetricproduct{\avecf{2}}{\avecf{3}}}+
\lieder{\avecf{2}}{\vfmetricproduct{\avecf{3}}{\avecf{1}}}-
\lieder{\avecf{3}}{\vfmetricproduct{\avecf{1}}{\avecf{2}}}\cr
&+\vfmetricproduct{\avecf{3}}{\liebracket{\avecf{1}}{\avecf{2}}{}}-
\vfmetricproduct{\avecf{2}}{\liebracket{\avecf{1}}{\avecf{3}}{}}-
\vfmetricproduct{\avecf{1}}{\liebracket{\avecf{2}}{\avecf{3}}{}}.
\end{align}
This is called the $\quotl$Koszul property of the semi-Riemannian manifold $\opair{\Man{}}{\metrictensor{}}$$\quotr$.
\endcor
%%%%%%%%%%%%%%%%%%%%%%%%%%%%%%%%%%%%%%%%%%%%%%%%%%%%%%%%%%%%%%%%%%%%%%%%%%%%%%%%%%%%%%%%%%%%%%%%%%%%%%%%%%%%%%%%%%%%%%%%%%%%%%%%
\lemma\label{leminducedLeviCivitaconnectionnopensets}
Let $\U$ be a non-empty open subset of $\Man{}$. The restriction of the connection $\LCconnection{\Man{}}{\metrictensor{}}$
to $\U$ is the same as the Levi-Civita connection on the semi-Riemannian submanifold
$\opair{\subman{\Man{}}{\U}}{\reS{\metrictensor{}}{\U}}$. That is,
\begin{equation}
\rescon{\LCconnection{\Man{}}{\metrictensor{}}}{\U}=
\LCconnection{\subman{\Man{}}{\U}}{\reS{\metrictensor{}}{\U}}.
\end{equation}
\proof
It is trivial.
\endlem
%%%%%%%%%%%%%%%%%%%%%%%%%%%%%%%%%%%%%%%%%%%%%%%%%%%%%%%%%%%%%%%%%%%%%%%%%%%%%%%%%%%%%%%%%%%%%%%%%%%%%%%%%%%%%%%%%%%%%%%%%%%%%%%%
\theorem\label{thmChristoffelSymbolsofLeviCivitaconnection}
Let $\opair{\U}{\phi}$ be a chart of $\Man{}$, and conventionally,
let $\Christoffel{\phi}{i}{j}{k}$ denote the Christoffel symbols of the Levi-Civita connection
$\LCconnection{\Man{}}{\metrictensor{}}$ with respect to the chart $\opair{\U}{\phi}$.
Let $\mtuple{\localframevecf{1}}{\localframevecf{m}}$ denote the local frame field associated with the chart
$\opair{\U}{\phi}$. For every $i$, $j$, and $k$
in $\seta{\suc{1}{m}}$,
\begin{equation}
\Christoffel{\phi}{i}{j}{k}=\frac{1}{2}\sum_{\mu=1}^{m}
\metrictensorchartinv{\phi}{i}{\mu}\(\lieder{\localframevecf{j}}{\metrictensorchart{\phi}{\mu}{k}}+
\lieder{\localframevecf{k}}{\metrictensorchart{\phi}{j}{\mu}}-
\lieder{\localframevecf{\mu}}{\metrictensorchart{\phi}{k}{j}}\).
\end{equation}
\proof
We denote $\LCconnection{\Man{}}{\metrictensor{}}$ simply by $\connection{}$. Also, for the sake of
notational convenience,
we will denote the restriction of the connection $\connection{}$ to $\U$ simply by $\connection{}$ itself,
and the actual connection in use can be inferred from the context.
According to \refdef{defChristoffelSymbols}, for each $j$ and $k$,
\begin{equation}\label{thmChristoffelSymbolsofLeviCivitaconnectioneq1}
\con{\localframevecf{k}}{\localframevecf{j}}=
\sum_{i=1}^{m}\Christoffel{}{i}{j}{k}\localframevecf{i}.
\end{equation}
It is also known that $\mtuple{\localframevecf{1}}{\localframevecf{m}}$ is a local frame field which particularly
behaves like the Euclidean coordinate  fields, when the problem of differentiating smooth maps is under consideration.
Precisely, the second order Lie-derivatives with respect to $\localframevecf{i}$-s are symmetric, that is,
for each $i$ and $j$,
\begin{equation}
\liebracket{\localframevecf{i}}{\localframevecf{j}}{}=\zerovec{}.
\end{equation}
Therefore, according to \refcor{corKoszulFormula}, \reflem{leminducedLeviCivitaconnectionnopensets},
and considering all other relevant naturalities with regard to restrictions to open subsets,
for every $j$ and $k$, and $\mu$ in $\seta{\suc{1}{m}}$,
\begin{align}
2\vfmetricproduct{\con{\localframevecf{k}}{\localframevecf{j}}}{\localframevecf{\mu}}=
\lieder{\localframevecf{k}}{\vfmetricproduct{\localframevecf{j}}{\localframevecf{\mu}}}+
\lieder{\localframevecf{j}}{\vfmetricproduct{\localframevecf{\mu}}{\localframevecf{k}}}-
\lieder{\localframevecf{\mu}}{\vfmetricproduct{\localframevecf{k}}{\localframevecf{j}}}.
\end{align}
On the other hand, we know that for every $i$ and $j$,
\begin{equation}
\metrictensorchart{}{i}{j}=
\vfmetricproduct{\localframevecf{i}}{\localframevecf{j}}.
\end{equation}
Therefore,
\begin{equation}
2\vfmetricproduct{\con{\localframevecf{k}}{\localframevecf{j}}}{\localframevecf{\mu}}=
\lieder{\localframevecf{k}}{\metrictensorchart{}{j}{\mu}}+
\lieder{\localframevecf{j}}{\metrictensorchart{}{\mu}{k}}-
\lieder{\localframevecf{\mu}}{\metrictensorchart{}{k}{j}}.
\end{equation}
On the other hand, using \Ref{thmChristoffelSymbolsofLeviCivitaconnectioneq1},
and according to \refpro{prometricinducedproductofvectorfieldsisbilinear},
for every $j$, $k$, and $\mu$,
\begin{align}
\vfmetricproduct{\con{\localframevecf{k}}{\localframevecf{j}}}{\localframevecf{\mu}}&=
\vfmetricproduct{\sum_{i=1}^{m}\Christoffel{}{i}{j}{k}\localframevecf{i}}{\localframevecf{\mu}}\cr
&=\sum_{i=1}^{m}\Christoffel{}{i}{j}{k}\vfmetricproduct{\localframevecf{i}}{\localframevecf{\mu}}\cr
&=\sum_{i=1}^{m}\Christoffel{}{i}{j}{k}\metrictensorchart{}{i}{\mu},
\end{align}
and hence,
\begin{equation}
2\sum_{i=1}^{m}\Christoffel{}{i}{j}{k}\metrictensorchart{}{i}{\mu}=
\lieder{\localframevecf{k}}{\metrictensorchart{}{j}{\mu}}+
\lieder{\localframevecf{j}}{\metrictensorchart{}{\mu}{k}}-
\lieder{\localframevecf{\mu}}{\metrictensorchart{}{k}{j}}.
\end{equation}
Therefore, for every $j$, $k$, and $l$,
\begin{align}
\sum_{\mu=1}^{m}\metrictensorchartinv{}{l}{\mu}\(
\sum_{i=1}^{m}\Christoffel{}{i}{j}{k}\metrictensorchart{}{i}{\mu}\)=\frac{1}{2}
\sum_{\mu=1}^{m}\metrictensorchartinv{}{l}{\mu}\(
\lieder{\localframevecf{k}}{\metrictensorchart{}{j}{\mu}}+
\lieder{\localframevecf{j}}{\metrictensorchart{}{\mu}{k}}-
\lieder{\localframevecf{\mu}}{\metrictensorchart{}{k}{j}}\).
\end{align}
But,
\begin{align}
\sum_{\mu=1}^{m}\metrictensorchartinv{}{l}{\mu}\(
\sum_{i=1}^{m}\Christoffel{}{i}{j}{k}\metrictensorchart{}{i}{\mu}\)&=
\sum_{i=1}^{m}\Christoffel{}{i}{j}{k}\(\sum_{l=1}^{m}\metrictensorchartinv{}{l}{\mu}
\metrictensorchart{}{i}{\mu}\)\cr
&=\sum_{i=1}^{m}\Christoffel{}{i}{j}{k}\deltaf{l}{i}\cr
&=\Christoffel{}{l}{j}{k}.
\end{align}
Thus, for every $l$, $j$, and $k$,
\begin{equation}
\Christoffel{}{l}{j}{k}=\frac{1}{2}
\sum_{\mu=1}^{m}\metrictensorchartinv{}{l}{\mu}\(
\lieder{\localframevecf{k}}{\metrictensorchart{}{j}{\mu}}+
\lieder{\localframevecf{j}}{\metrictensorchart{}{\mu}{k}}-
\lieder{\localframevecf{\mu}}{\metrictensorchart{}{k}{j}}\).
\end{equation}
\endthm
%%%%%%%%%%%%%%%%%%%%%%%%%%%%%%%%%%%%%%%%%%%%%%%%%%%%%%%%%%%%%%%%%%%%%%%%%%%%%%%%%%%%%%%%%%%%%%%%%%%%%%%%%%%%%%%%%%%%%%%%%%%%%%%%
\definition
Let $\function{\Curve{}}{\interval{}}{\Man{}}$ be a smooth curve on $\Man{}$ for some open nterval $\interval{}$ of $\R$.
$\Curve{}$ is referred to as a $\quotl$geodesic on the semi-Riemannian manifold $\opair{\Man{}}{\metrictensor{}}$
if $\Curve{}$ is a geodesic on $\Man{}$ relative to the Levi-Civita connection $\LCconnection{\Man{}}{\metrictensor{}}$
of the semi-Riemannian manifold $\opair{\Man{}}{\metrictensor{}}$.\\
%%%
We will refer to $\maxgeodesic{\Man{}}{\LCconnection{\Man{}}{\metrictensor{}}}{v}$, that is, the unique maximal geodesic
$\function{\Curve{}}{\interval{}}{\Man{}}$ on $\Man{}$ relative to $\LCconnection{\Man{}}{\metrictensor{}}$
such that $0\in\interval{}$ and $\func{\dot{\Curve{}}}{0}=v$, as the
$\quotl$maximal geodesic on the semi-Riemannian manifold $\opair{\Man{}}{\metrictensor{}}$
with the initial value $v$$\quotr$.
\endef
%%%%%%%%%%%%%%%%%%%%%%%%%%%%%%%%%%%%%%%%%%%%%%%%%%%%%%%%%%%%%%%%%%%%%%%%%%%%%%%%%%%%%%%%%%%%%%%%%%%%%%%%%%%%%%%%%%%%%%%%%%%%%%%%
\definition
Let $\function{\Curve{}}{\interval{}}{\Man{}}$ be a smooth curve on the $\Man{}$.
We define the function $\function{\speed{\Man{}}{\metrictensor{}}{\Curve{}}}{\interval{}}{\R}$ as,
\begin{equation}
\Foreach{t}{\interval{}}
\func{\speed{\Man{}}{\metrictensor{}}{\Curve{}}}{t}\eqdef
\sqrt{\abs{\func{\[\func{\metrictensor{}}{\func{\Curve{}}{t}}\]}{\binary{\func{\dot{\Curve{}}}{t}}{\func{\dot{\Curve{}}}{t}}}}}=
\abs{\vfmetricproduct{\func{\dot{\Curve{}}}{t}}{\func{\dot{\Curve{}}}{t}}}^{\frac{1}{2}}.
\end{equation}
$\speed{\Man{}}{\metrictensor{}}{\Curve{}}$ is referred to as the $\quotl$speed function of the smooth curve
$\Curve{}$ on the semi-Riemannian manifold $\opair{\Man{}}{\metrictensor{}}$$\quotr$.\\
If $\speed{\Man{}}{\metrictensor{}}{\Curve{}}$ is a constant function, then $\Curve{}$ is called a
$\quotl$unit-speed smooth curve  on the semi-Riemannian manifold $\opair{\Man{}}{\metrictensor{}}$$\quotr$.
\endef
%%%%%%%%%%%%%%%%%%%%%%%%%%%%%%%%%%%%%%%%%%%%%%%%%%%%%%%%%%%%%%%%%%%%%%%%%%%%%%%%%%%%%%%%%%%%%%%%%%%%%%%%%%%%%%%%%%%%%%%%%%%%%%%%
\theorem
Let $\function{\Curve{}}{\interval{}}{\Man{}}$ be a geodesic on the semi-Riemannian manifold
$\opair{\Man{}}{\metrictensor{}}$. $\Curve{}$ is a unit-speed smooth curve on
$\opair{\Man{}}{\metrictensor{}}$$\quotr$, that is,
\begin{equation}
\timeder\vfmetricproduct{\dot{\Curve{}}}{\dot{\Curve{}}}=0.
\end{equation}
\proof
It is an immediate consequence of \reflem{lemequivalentpropertiesofametricconnection}, definition of the
Levi-Civita connection of a semi-Riemannian manifold, and definition of a geodesic on a manifold with connection.
\endthm
%%%%%%%%%%%%%%%%%%%%%%%%%%%%%%%%%%%%%%%%%%%%%%%%%%%%%%%%%%%%%%%%%%%%%%%%%%%%%%%%%%%%%%%%%%%%%%%%%%%%%%%%%%%%%%%%%%%%%%%%%%%%%%%%
\definition
Let $\function{\Curve{}}{\interval{}}{\Man{}}$ be a smooth curve on $\Man{}$.
Let $\binary{a}{b}\in\interval{}$ and $a<b$. We define,
\begin{align}
\length{\Man{}}{\metrictensor{}}{\Curve{}}{a}{b}&:=
\int_{a}^{b}\func{\speed{\Man{}}{\metrictensor{}}{\Curve{}}}{t}~{\mathrm{d}}t\cr
&=\int_{a}^{b}\abs{\vfmetricproduct{\func{\dot{\Curve{}}}{t}}{\func{\dot{\Curve{}}}{t}}}^{\frac{1}{2}}~{\mathrm{d}}t.
\end{align}
$\length{\Man{}}{\metrictensor{}}{\Curve{}}{a}{b}$ is referred to as the $\quotl$length of the segment
$\reS{\Curve{}}{\cpair{a}{b}}$ of the smooth curve $\Curve{}$ on the semi-Riemannian manifold
$\opair{\Man{}}{\metrictensor{}}$$\quotr$.
\endef
%%%%%%%%%%%%%%%%%%%%%%%%%%%%%%%%%%%%%%%%%%%%%%%%%%%%%%%%%%%%%%%%%%%%%%%%%%%%%%%%%%%%%%%%%%%%%%%%%%%%%%%%%%%%%%%%%%%%%%%%%%%%%%%%
\theorem
Suppose that $\opair{\Man{}}{\metrictensor{}}$ is isometric to $\opair{\Man{1}}{\metrictensor{1}}$, and let
$f$ be an isometry from $\opair{\Man{}}{\metrictensor{}}$ to $\opair{\Man{1}}{\metrictensor{1}}$. The pullback of the
Levi-Civita connection of $\opair{\Man{1}}{\metrictensor{1}}$ relative to $f$ equals the Levi-Civita connection of
$\opair{\Man{}}{\metrictensor{}}$. That is,
\begin{equation}
\func{\connectionpb{f}}{\LCconnection{\Man{1}}{\metrictensor{1}}}=
\LCconnection{\Man{}}{\metrictensor{}}.
\end{equation}
\proof
We just give the sketch of the proof, and leave the details to the reader.
Since the Levi-Civita connection of a semi-Riemannian manifold is the only connection that is a
metric-connection without torsion, it suffices to check the conditions of being a metric-connection and
torsion-lessness for the connection $\func{\connectionpb{f}}{\LCconnection{\Man{1}}{\metrictensor{1}}}$.
\endthm
%%%%%%%%%%%%%%%%%%%%%%%%%%%%%%%%%%%%%%%%%%%%%%%%%%%%%%%%%%%%%%%%%%%%%%%%%%%%%%%%%%%%%%%%%%%%%%%%%%%%%%%%%%%%%%%%%%%%%%%%%%%%%%%%
%%%%%%%%%%%%%%%%%%%%%%%%%%%%%%%%%%%%%%%%%%%%%%%%%%%%%%%%%%%%%%%%%%%%%%%%%%%%%%%%%%%%%%%%%%%%%%%%%%%%%%%%%%%%%%%%%%%%%%%%%%%%%%%%
%%%%%%%%%%%%%%%%%%%%%%%%%%%%%%%%%%%%%%%%%%%%%%%%%%%%%%%%%%%%%%%%%%%%%%%%%%%%%%%%%%%%%%%%%%%%%%%%%%%%%%%%%%%%%%%%%%%%%%%%%%%%%%%%
%%%%%%%%%%%%%%%%%%%%%%%%%%%%%%%%%%%%%%%%%%%%%%%%%%%%%%%%%%%%%%%%%%%%%%%%%%%%%%%%%%%%%%%%%%%%%%%%%%%%%%%%%%%%%%%%%%%%%%%%%%%%%%%%
%%%%%%%%%%%%%%%%%%%%%%%%%%%%%%%%%%%%%%%%%%%%%%%%%%%%%%%%%%%%%%%%%%%%%%%%%%%%%%%%%%%%%%%%%%%%%%%%%%%%%%%%%%%%%%%%%%%%%%%%%%%%%%%%
%%%%%%%%%%%%%%%%%%%%%%%%%%%%%%%%%%%%%%%%%%%%%%%%%%%%%%%%%%%%%%%%%%%%%%%%%%%%%%%%%%%%%%%%%%%%%%%%%%%%%%%%%%%%%%%%%%%%%%%%%%%%%%%%
%%%%%%%%%%%%%%%%%%%%%%%%%%%%%%%%%%%%%%%%%%%%%%%%%%%%%%%%%%%%%%%%%%%%%%%%%%%%%%%%%%%%%%%%%%%%%%%%%%%%%%%%%%%%%%%%%%%%%%%%%%%%%%%%
%%%%%%%%%%%%%%%%%%%%%%%%%%%%%%%%%%%%%%%%%%%%%%%%%%%%%%%%%%%%%%%%%%%%%%%%%%%%%%%%%%%%%%%%%%%%%%%%%%%%%%%%%%%%%%%%%%%%%%%%%%%%%%%%
%%%%%%%%%%%%%%%%%%%%%%%%%%%%%%%%%%%%%%%%%%%%%%%%%%%%%%%%%%%%%%%%%%%%%%%%%%%%%%%%%%%%%%%%%%%%%%%%%%%%%%%%%%%%%%%%%%%%%%%%%%%%%%%%
%%%%%%%%%%%%%%%%%%%%%%%%%%%%%%%%%%%%%%%%%%%%%%%%%%%%%%%%%%%%%%%%%%%%%%%%%%%%%%%%%%%%%%%%%%%%%%%%%%%%%%%%%%%%%%%%%%%%%%%%%%%%%%%%
%%%%%%%%%%%%%%%%%%%%%%%%%%%%%%%%%%%%%%%%%%%%%%%%%%%%%%%%%%%%%%%%%%%%%%%%%%%%%%%%%%%%%%%%%%%%%%%%%%%%%%%%%%%%%%%%%%%%%%%%%%%%%%%%
\section*{Exercises}
\exercise
Let $\Man{}$ be a manifold of dimension $m$. Show that there exists a Riemannian metric on $\Man{}$.\\
\textit{
{\textbf{Hint}}:
Given a point $\point$ of $\Man{}$, choose a chart $\opair{\U}{\phi}$ around $\point$.
For any $v\in\tanspace{\point}{\Man{}}$, let $\mtuple{v_1}{v_m}\in\R^m$ denote the representation
of $v$ with respect to this chart. Then, define an $\R$-valued binary operation on $\tanspace{\point}{\Man{}}$
by $\displaystyle\vfmetricproduct{u}{v}\eqdef\sum_{i=1}^{m}\func{\phi}{\point}u_iv_i$ for every
$\binary{u}{v}\in\tanspace{\point}{\Man{}}$. Verify that this binary operation on $\tanspace{\point}{\Man{}}$
is an inner-product on the vector-space $\Tanspace{\point}{\Man{}}$.\\
Define a binary operation on the tangent space of every point of $\Man{}$ with respect to any chart around $\point$
as the procedure above.
Check that any pair of such binary operations on the tangent space of a point coincide
with respect to any choice of charts around that point.
Finally, make use of the partition of unity technique to glue the constructed binary operations together,
and verify the smoothness the global tensor field thus achieved.
}
\endexercise
%%%%%%%%%%%%%%%%%%%%%%%%%%%%%%%%%%%%%%%%%%%%%%%%%%%%%%%%%%%%%%%%%%%%%%%%%%%%%%%%%%%%%%%%%%%%%%%%%%%%%%%%%%%%%%%%%%%%%%%%%%%%%%%%
\exercise
Let $\Man{}$ be a manifold and let $\metrictensor{}$ be a semi-Riemannian metric on $\Man{}$ with any
index $\nu\in\seta{\suc{0}{m}}$. For every $\point\in\Man{}$, there exists an open neighborhood $\U$ of
$\point$ in $\Man{}$ such that there exists an orthonormal local frame on the open subset $\U$ of the semi-Riemannian
manifold $\opair{\Man{}}{\metrictensor{}}$.
\endexercise
\chapteR{Curvature}
\thispagestyle{fancy}
\section{Curvature Operator of a Connection}
%%%%%%%%%%%%%%%%%%%%%%%%%%%%%%%%%%%%%%%%%%%%%%%%%%%%%%%%%%%%%%%%%%%%%%%%%%%%%%%%%%%%%%%%%%%%%%%%%%%%%%%%%%%%%%%%%%%%%%%%%%%%%%%%
\fixed
$\vbundle{}=\quintuple{\vbtotal{}}{\vbprojection{}}{\vbbase{}}{\vbfiber{}}{\vbatlas{}}$ is fixed as a real smooth vector bundle
of rank $d$,
where $\vbtotal{}=\opair{\vTot{}}{\maxatlas{\vTot{}}}$ and
$\vbbase{}=\opair{\vB{}}{\maxatlas{\vB{}}}$ are $\difclass{\infty}$ manifolds
modeled on the Banach-spaces $\R^{n_{\vTot{}}}$ and $\R^{n_{\vB{}}}$, respectively.\\
In addition, $\Man{}$ is fixed as a manifold with dimension $m$.
\endfixed
%%%%%%%%%%%%%%%%%%%%%%%%%%%%%%%%%%%%%%%%%%%%%%%%%%%%%%%%%%%%%%%%%%%%%%%%%%%%%%%%%%%%%%%%%%%%%%%%%%%%%%%%%%%%%%%%%%%%%%%%%%%%%%%%
\definition\label{defcurvatureofaconnection}
Let $\connection{}$ be a connection on the vector bundle $\vbundle{}$.
The map $\function{\curvature{\connection{}}}{\Cprod{\vectorfields{\vbbase{}}}{\vectorfields{\vbbase{}}}}
{\Func{\vbsections{\vbundle{}}}{\vbsections{\vbundle{}}}}$ is defined as,
\begin{align}
&\Foreach{\opair{\avecf{1}}{\avecf{2}}}{\Cprod{\vectorfields{\vbbase{}}}{\vectorfields{\vbbase{}}}}
\Foreach{\vbsec{}}{\vbsections{\vbundle{}}}\cr
&\func{\[\func{\curvature{\connection{}}}{\binary{\avecf{1}}{\avecf{2}}}\]}{\vbsec{}}\eqdef
\con{\avecf{1}}{\con{\avecf{2}}{\vbsec{}}}-\con{\avecf{2}}{\con{\avecf{1}}{\vbsec{}}}-
\con{\liebracket{\avecf{1}}{\avecf{2}}{}}{\vbsec{}}.
\end{align}
$\curvature{\connection{}}$ is referred to as the $\quotl$curvature (operator) on the smooth vector bundle $\vbundle{}$
relative to the connection $\connection{}$$\quotr$.
\endef
%%%%%%%%%%%%%%%%%%%%%%%%%%%%%%%%%%%%%%%%%%%%%%%%%%%%%%%%%%%%%%%%%%%%%%%%%%%%%%%%%%%%%%%%%%%%%%%%%%%%%%%%%%%%%%%%%%%%%%%%%%%%%%%%
\theorem\label{thmcurvatureofaconnectionmultilinearity1}
Let $\connection{}$ be a connection on the vector bundle $\vbundle{}$.
Let $\avecf{}$ and $\avecff{}$ be elements of $\vectorfields{\vbbase{}}$.
$\func{\curvature{\connection{}}}{\binary{\avecf{}}{\avecff{}}}\in\mendomorphism{\smoothmaps{\Man{}}}{\vbsections{\vbundle{}}}$,
$\mendomorphism{\smoothmaps{\vbbase{}}}{\vbsections{\vbundle{}}}$ denoting the set of all modul-homomorphisms
from the $\smoothmaps{\vbbase{}}$-module $\vbsections{\vbundle{}}$ to itself. That is, for every $\vbsec{}$ and $\vbsec{1}$
in $\vbsections{\vbundle{}}$, and every $f\in\smoothmaps{\vbbase{}}$,
\begin{align}
\func{\[\func{\curvature{\connection{}}}{\binary{\avecf{}}{\avecff{}}}\]}{f\vbsec{}+\vbsec{1}}=
f\func{\[\func{\curvature{\connection{}}}{\binary{\avecf{}}{\avecff{}}}\]}{\vbsec{}}+
\func{\[\func{\curvature{\connection{}}}{\binary{\avecf{}}{\avecff{}}}\]}{\vbsec{1}}.
\end{align}
\proof
Let $\vbsec{},~\vbsec{1}\in\vbsections{\vbundle{}}$, and let $f\in\smoothmaps{\vbbase{}}$.
According to \refdef{defconnection}, \refdef{defcurvatureofaconnection}, and
considering the linearity of the action of Lie-derivatives on smooth maps, it is obvious that.
\begin{align}
\func{\[\func{\curvature{\connection{}}}{\binary{\avecf{}}{\avecff{}}}\]}{\vbsec{}+\vbsec{1}}&=
\func{\[\func{\curvature{\connection{}}}{\binary{\avecf{}}{\avecff{}}}\]}{\vbsec{}}+
\func{\[\func{\curvature{\connection{}}}{\binary{\avecf{}}{\avecff{}}}\]}{\vbsec{1}}.
\end{align}
Moreover, according to \refdef{defconnection}, \refdef{defcurvatureofaconnection}, and considering the
product rule of Lie-derivatives, and considering that
$\lieder{\liebracket{\avecf{}}{\avecff{}}{}}{f}=\lieder{\avecf{}}{\lieder{\avecff{}}{f}}-
\lieder{\avecff{}}{\lieder{\avecf{}}{f}}$,
\begin{align}
&~~~~\func{\[\func{\curvature{\connection{}}}{\binary{\avecf{}}{\avecff{}}}\]}{f\vbsec{}}\cr
&=\con{\avecf{}}{\con{\avecff{}}{f\vbsec{}}}-\con{\avecff{}}{\con{\avecf{}}{f\vbsec{}}}-
\con{\liebracket{\avecf{}}{\avecff{}}{}}{f\vbsec{}}\cr
&=\con{\avecf{}}{\[\(\lieder{\avecff{}}{f}\)\vbsec{}+f\con{\avecff{}}{\vbsec{}}\]}-
\con{\avecff{}}{\[\(\lieder{\avecf{}}{f}\)\vbsec{}+f\con{\avecf{}}{\vbsec{}}\]}-
\[\(\lieder{\liebracket{\avecf{}}{\avecff{}}{}}{f}\)\vbsec{}+f\con{\liebracket{\avecf{}}{\avecff{}}{}}{\vbsec{}}\]\cr
&=~~\(\lieder{\avecf{}}{\lieder{\avecff{}}{f}}\)\vbsec{}+\(\lieder{\avecff{}}{f}\)\con{\avecf{}}{\vbsec{}}
+\(\lieder{\avecf{}}{f}\)\con{\avecff{}}{\vbsec{}}+f\con{\avecf{}}{\con{\avecff{}}{\vbsec{}}}\cr
&~~~-\(\(\lieder{\avecff{}}{\lieder{\avecf{}}{f}}\)\vbsec{}+\(\lieder{\avecf{}}{f}\)\con{\avecff{}}{\vbsec{}}
+\(\lieder{\avecff{}}{f}\)\con{\avecf{}}{\vbsec{}}+f\con{\avecff{}}{\con{\avecf{}}{\vbsec{}}}\)\cr
&~~~-\[\(\lieder{\liebracket{\avecf{}}{\avecff{}}{}}{f}\)\vbsec{}+f\con{\liebracket{\avecf{}}{\avecff{}}{}}{\vbsec{}}\]\cr
&=f\(\con{\avecf{1}}{\con{\avecf{2}}{\vbsec{}}}-\con{\avecf{2}}{\con{\avecf{1}}{\vbsec{}}}-
\con{\liebracket{\avecf{1}}{\avecf{2}}{}}{\vbsec{}}\)\cr
&=f\func{\[\func{\curvature{\connection{}}}{\binary{\avecf{}}{\avecff{}}}\]}{f\vbsec{}}.
\end{align}
\endthm
%%%%%%%%%%%%%%%%%%%%%%%%%%%%%%%%%%%%%%%%%%%%%%%%%%%%%%%%%%%%%%%%%%%%%%%%%%%%%%%%%%%%%%%%%%%%%%%%%%%%%%%%%%%%%%%%%%%%%%%%%%%%%%%%
\theorem\label{thmcurvatureofaconnectionmultilinearity2}
Let $\connection{}$ be a connection on the vector bundle $\vbundle{}$, and fix $\vbsec{}\in\vbsections{\vbundle{}}$.
The map $\Cprod{\vectorfields{\vbbase{}}}{\vectorfields{\vbbase{}}}\ni\opair{\avecf{1}}{\avecf{2}}\mapsto
\func{\[\func{\curvature{\connection{}}}{\binary{\avecf{1}}{\avecf{2}}}\]}{\vbsec{}}\in\vbsections{\vbundle{}}$
is a $\smoothmaps{\vbbase{}}$-bilinear map from $\Cprod{\vectorfields{\vbbase{}}}{\vectorfields{\vbbase{}}}$ to
$\vbsections{\vbundle{}}$, considering the canonical $\smoothmaps{\vbbase{}}$-module structures of $\vectorfields{\vbbase{}}$
and $\vbsections{\vbundle{}}$.
\proof
We show the $\smoothmaps{\vbbase{}}$-linearity just with respect to the first variable, since the proof
of $\smoothmaps{\vbbase{}}$-linearity with respect to the first and second variables are completely alike.\\
Fix $\avecf{2}\in\vectorfields{\vbbase{}}$, let $\avecf{1},~\avecff{1}\in\vectorfields{\vbbase{}}$,
and let $f\in\smoothmaps{\vbbase{}}$.
According to \refdef{defconnection}, \refdef{defcurvatureofaconnection}, and considering that
$\liebracket{f\avecf{1}}{\avecf{2}}{}f\liebracket{\avecf{1}}{\avecf{2}}{}-\(\lieder{\avecf{2}}{f}\)\avecf{1}$,
\begin{align}
&~~~~\func{\[\func{\curvature{\connection{}}}{\binary{f\avecf{1}+\avecff{1}}{\avecf{2}}}\]}{\vbsec{}}\cr
&=\con{f\avecf{1}+\avecff{1}}{\con{\avecf{2}}{\vbsec{}}}-\con{\avecf{2}}{\con{f\avecf{1}+\avecff{1}}{\vbsec{}}}-
\con{\liebracket{f\avecf{1}+\avecff{1}}{\avecf{2}}{}}{\vbsec{}}\cr
&=~~\(f\con{\avecf{1}}{\con{\avecf{2}}{\vbsec{}}}+\con{\avecff{1}}{\con{\avecf{2}}{\vbsec{}}}\)-
\[\con{\avecf{2}}{\(f\con{\avecf{1}}{\vbsec{}}+\con{\avecff{1}}{\vbsec{}}\)}\]\cr
&~~~-\(f\con{\liebracket{\avecf{1}}{\avecf{2}}{}}{\vbsec{}}+\con{\liebracket{\avecff{1}}{\avecf{2}}{}}{\vbsec{}}\)\cr
&=~~\(f\con{\avecf{1}}{\con{\avecf{2}}{\vbsec{}}}+\con{\avecff{1}}{\con{\avecf{2}}{\vbsec{}}}\)-
\[\(\lieder{\avecf{2}}{f}\)\con{\avecf{1}}{\vbsec{}}+f\con{\avecf{2}}{\con{\avecf{1}}{\vbsec{}}}+
\con{\avecf{2}}{\con{\avecff{1}}{\vbsec{}}}\]\cr
&~~~-\(f\con{\liebracket{\avecf{1}}{\avecf{2}}{}}{\vbsec{}}+\con{\liebracket{\avecff{1}}{\avecf{2}}{}}{\vbsec{}}-
\(\lieder{\avecf{2}}{f}\)\con{\avecff{1}}{\vbsec{}}\)\cr
&=f\(\con{\avecf{1}}{\con{\avecf{2}}{\vbsec{}}}-\con{\avecf{2}}{\con{\avecf{1}}{\vbsec{}}}-
\con{\liebracket{\avecf{1}}{\avecf{2}}{}}{\vbsec{}}\)-
\(\con{\avecff{1}}{\con{\avecf{2}}{\vbsec{}}}-\con{\avecf{2}}{\con{\avecff{1}}{\vbsec{}}}-
\con{\liebracket{\avecff{1}}{\avecf{2}}{}}{\vbsec{}}\)\cr
&=f\func{\[\func{\curvature{\connection{}}}{\binary{\avecf{1}}{\avecf{2}}}\]}{\vbsec{}}+
\func{\[\func{\curvature{\connection{}}}{\binary{\avecff{1}}{\avecf{2}}}\]}{\vbsec{}}.
\end{align}
\endthm
%%%%%%%%%%%%%%%%%%%%%%%%%%%%%%%%%%%%%%%%%%%%%%%%%%%%%%%%%%%%%%%%%%%%%%%%%%%%%%%%%%%%%%%%%%%%%%%%%%%%%%%%%%%%%%%%%%%%%%%%%%%%%%%%
\definition\label{defcurvatureofanaffineconnection}
Let $\connection{}$ be an affine connection on the manifold $\Man{}$.
The map
$\function{\mcurvature{\connection{}}}{\Cprod{\vectorfields{\Man{}}}{\Cprod{\vectorfields{\Man{}}}{\vectorfields{\Man{}}}}}
{\vectorfields{\Man{}}}$ is defined as,
\begin{align}
\Foreach{\triple{\avecf{1}}{\avecf{2}}{\avecf{3}}}{\vectorfields{\Man{}}^3}
\func{\mcurvature{\connection{}}}{\triplet{\avecf{1}}{\avecf{2}}{\avecf{3}}}&\eqdef
\func{\[\func{\curvature{\connection{}}}{\binary{\avecf{1}}{\avecf{2}}}\]}{\avecf{3}}\cr
&=\con{\avecf{1}}{\con{\avecf{2}}{\avecf{3}}}-\con{\avecf{2}}{\con{\avecf{1}}{\avecf{3}}}-
\con{\liebracket{\avecf{1}}{\avecf{2}}{}}{\avecf{3}}.\cr
&{}
\end{align}
$\mcurvature{\connection{}}$ is referred to as the $\quotl$curvature (operator) of the manifold $\Man{}$
relative to the affine connection $\connection{}$$\quotr$.
\endef
%%%%%%%%%%%%%%%%%%%%%%%%%%%%%%%%%%%%%%%%%%%%%%%%%%%%%%%%%%%%%%%%%%%%%%%%%%%%%%%%%%%%%%%%%%%%%%%%%%%%%%%%%%%%%%%%%%%%%%%%%%%%%%%%
\textit{
Let $\vectorspace{}$ be a finite-dimensional $\R$-vector-space with dimension $n$, and let $r$ and $s$ be non-negative integers
with $r+s>0$. By definition, an $\opair{r}{s}$ tensor on $\vectorspace{}$ is a multilinear map from
$\Cprod{{\vectorspace{}}^r}{{\Vdual{\vectorspace{}}}^s}$ to $\R$. This definition can be generalized by relaxing the limitation
on the codomain of multilinear maps on $\Cprod{{\vectorspace{}}^r}{{\Vdual{\vectorspace{}}}^s}$ to be any vector-space.
As a particular case, we are concerned here with multilinear maps from $\Cprod{{\vectorspace{}}^r}{{\Vdual{\vectorspace{}}}^s}$
to the vector-space $\vectorspace{}$ itself. We will denote by $\Tensors{r}{s}{\bbinary{\vectorspace{}}{\vectorspace{}}}$
the set of all multilinear maps from $\Cprod{{\vectorspace{}}^r}{{\Vdual{\vectorspace{}}}^s}$ to $\vectorspace{}$.
Each element of $\Tensors{r}{s}{\bbinary{\vectorspace{}}{\vectorspace{}}}$ is called a $\quotl$$\vectorspace{}$-valued
$\opair{r}{s}$-tensor on $\vectorspace{}$$\quotr$. Furthermore, similar to
$\Tensors{r}{s}{\vectorspace{}}$, $\Tensors{r}{s}{\bbinary{\vectorspace{}}{\vectorspace{}}}$ possesses a natural linear
structure, and hence can be regarded as vector-space in itself. The dimension of
$\Tensors{r}{s}{\bbinary{\vectorspace{}}{\vectorspace{}}}$ can be easily computed to be $n^{r+s+1}$.\\
Now, we aim at globalizing this generalization by carrying it on manifolds. We define,
\begin{equation}
\vbtensorbundle{r}{s}{\bbinary{\Tanbun{\Man{}}}{\Tanbun{\Man{}}}}:=
\Union{\point}{\Man{}}{\Tensors{r}{s}{\bbinary{\Tanspace{\point}{\Man{}}}{\Tanspace{\point}{\Man{}}}}}.
\end{equation}
Similar to $\vbtensorbundle{r}{s}{\Tanbun{\Man{}}}$, $\vbtensorbundle{r}{s}{\bbinary{\Tanbun{\Man{}}}{\Tanbun{\Man{}}}}$
can be equipped with a canonical differentiable structure inspired by that of
the manifold $\Man{}$, which has additionally the structure of a smooth vector bundle over the base-space $\Man{}$
with rank $m^{r+s+1}$.
Any smooth section of the vector bundle $\vbtensorbundle{r}{s}{\bbinary{\Tanbun{\Man{}}}{\Tanbun{\Man{}}}}$, that is
any element of $\vbsections{\vbtensorbundle{r}{s}{\bbinary{\Tanbun{\Man{}}}{\Tanbun{\Man{}}}}}$,
is called an $\quotl$$\opair{r}{s}$ $\Tanbun{\Man{}}$-valued smooth tensor field on $\Tanbun{\Man{}}$$\quotr$.
We will use the alternative notation $\TF{r}{s}{\bbinary{\Tanbun{\Man{}}}{\Tanbun{\Man{}}}}$ for
$\vbsections{\vbtensorbundle{r}{s}{\bbinary{\Tanbun{\Man{}}}{\Tanbun{\Man{}}}}}$.
So, any $T\in\TF{r}{s}{\bbinary{\Tanbun{\Man{}}}{\Tanbun{\Man{}}}}$ is an smooth assignment of
a multilinear map $\function{\func{T}{\point}}{\Cprod{{\(\Tanspace{\point}{\Man{}}\)}^r}{{\(\Vdual{\Tanspace{\point}{\Man{}}}\)}^s}}
{\Tanspace{\point}{\Man{}}}$ to any point $\point$ of $\Man{}$.\\
Now, we aim to look at the $\Tanbun{\Man{}}$-valued tensor fields on $\Tanbun{\Man{}}$ from another perspective.
We define $\TTF{r}{s}{\bbinary{\vectorfields{\Man{}}}{\vectorfields{\Man{}}}}$ to be the set of all
$\smoothmaps{\Man{}}$-multilinear maps from $\Cprod{{\vectorfields{\Man{}}}^r}{{\oneforms{\Man{}}}^s}$ to
$\vectorfields{\Man{}}$, considering the canonical $\smoothmaps{\Man{}}$-module structures of
$\vectorfields{\Man{}}$ and $\oneforms{\Man{}}$.
Similar to the case of ordinary tensor fields on $\Tanbun{\Man{}}$, there is a natural one-one and onto
correspondence between $\TF{r}{s}{\bbinary{\Tanbun{\Man{}}}{\Tanbun{\Man{}}}}$ and
$\TTF{r}{s}{\bbinary{\vectorfields{\Man{}}}{\vectorfields{\Man{}}}}$. This correspondence takes place via the map
$\function{\p{\Upsilon}}{\TF{r}{s}{\bbinary{\Tanbun{\Man{}}}{\Tanbun{\Man{}}}}}
{\TTF{r}{s}{\bbinary{\vectorfields{\Man{}}}{\vectorfields{\Man{}}}}}$ defined so that,
for every $T\in\TF{r}{s}{\bbinary{\Tanbun{\Man{}}}{\Tanbun{\Man{}}}}$, every $\suc{\avecf{1}}{\avecf{r}}\in\vectorfields{\Man{}}$,
every $\suc{\aoneform{1}}{\aoneform{s}}\in\oneforms{\Man{}}$, and every $\point\in\Man{}$,
\begin{equation}
\func{\[\func{\(\func{\p{\Upsilon}}{T}\)}{\binary{\suc{\avecf{1}}{\avecf{r}}}{\suc{\aoneform{1}}{\aoneform{s}}}}\]}{\point}=
\func{\[\func{T}{\point}\]}{\binary{\suc{\func{\avecf{1}}{\point}}{\func{\avecf{r}}{\point}}}
{\suc{\func{\aoneform{1}}{\point}}{\func{\aoneform{s}}{\point}}}}.
\end{equation}
For every $\Tanbun{\Man{}}$-valued tensor field $T\in\TF{r}{s}{\bbinary{\Tanbun{\Man{}}}{\Tanbun{\Man{}}}}$,
we will denote its corresponded element $\func{\p{\Upsilon}}{T}$
in $\TTF{r}{s}{\bbinary{\vectorfields{\Man{}}}{\vectorfields{\Man{}}}}$ simply by $\TFequivv{T}$,
when the underying manifold is clear enough.
In addition, for every $\mathcal{T}\in\TTF{r}{s}{\bbinary{\vectorfields{\Man{}}}{\vectorfields{\Man{}}}}$,
we will also denote its corresponded element $\func{\finv{\p{\Upsilon}}}{\mathcal{T}}$
in $\TF{r}{s}{\bbinary{\Tanbun{\Man{}}}{\Tanbun{\Man{}}}}$ with  $\TTFequivv{\mathcal{T}}$.
On the basis of this natural identification of $\TF{R}{s}{\bbinary{\Tanbun{\Man{}}}{\Tanbun{\Man{}}}}$
with $\TTF{r}{s}{\bbinary{\vectorfields{\Man{}}}{\vectorfields{\Man{}}}}$, we also refer to the elements of
$\TTF{r}{s}{\bbinary{\Tanbun{\Man{}}}{\Tanbun{\Man{}}}}$ as $\quotl$$\opair{r}{s}$ $\vectorfields{\Man{}}$-valued
smooth tensor fields on $\Man{}$$\quotr$.\\
Based on the discussions above, it is evident that given a
$\mathcal{T}\in\TTF{r}{s}{\bbinary{\vectorfields{\Man{}}}{\vectorfields{\Man{}}}}$, and a point $\point$ of $\Man{}$,
for every pair of elements $\opair{\suc{\avecf{1}}{\avecf{r}}}{\suc{\aoneform{1}}{\aoneform{s}}}$ and
$\opair{\suc{\avecff{1}}{\avecff{r}}}{\suc{\omega_1}{\omega_s}}$ in
$\Cprod{{\vectorfields{\Man{}}}^r}{{\oneforms{\Man{}}}^s}$ such that
\begin{align*}
&\Foreach{i}{\seta{\suc{1}{r}}}
\func{\avecf{i}}{\point}=\func{\avecff{i}}{\point}\cr
&\Foreach{i}{\seta{\suc{1}{s}}}
\func{\aoneform{i}}{\point}=\func{\omega_i}{\point},
\end{align*}
we must have,
\begin{align}
\func{\[\func{\mathcal{T}}{\binary{\suc{\avecf{1}}{\avecf{r}}}{\suc{\aoneform{1}}{\aoneform{s}}}}\]}{\point}=
\func{\[\func{\mathcal{T}}{\binary{\suc{\avecff{1}}{\avecff{r}}}{\suc{\omega_1}{\omega_s}}}\]}{\point}.
\end{align}
}	
%%%%%%%%%%%%%%%%%%%%%%%%%%%%%%%%%%%%%%%%%%%%%%%%%%%%%%%%%%%%%%%%%%%%%%%%%%%%%%%%%%%%%%%%%%%%%%%%%%%%%%%%%%%%%%%%%%%%%%%%%%%%%%%%
\theorem
Let $\connection{}$ be an affine connection on the manifold $\Man{}$.
$\mcurvature{\connection{}}$ is a $\smoothmaps{\Man{}}$-multilinear map from $\vectorfields{\Man{}}^3$ to
$\vectorfields{\Man{}}$, considering the canonical $\smoothmaps{\Man{}}$-module structure of $\vectorfields{\Man{}}$.
That is,
\begin{equation}
\mcurvature{\connection{}}\in
\TTF{3}{0}{\bbinary{\vectorfields{\Man{}}}{\vectorfields{\Man{}}}}.
\end{equation}
\proof
It is actually shown in the proofs of \refthm{thmcurvatureofaconnectionmultilinearity1}
and \refthm{thmcurvatureofaconnectionmultilinearity2}.
\endthm
%%%%%%%%%%%%%%%%%%%%%%%%%%%%%%%%%%%%%%%%%%%%%%%%%%%%%%%%%%%%%%%%%%%%%%%%%%%%%%%%%%%%%%%%%%%%%%%%%%%%%%%%%%%%%%%%%%%%%%%%%%%%%%%%
\corollary
Let $\connection{}$ be an affine connection on the manifold $\Man{}$, and let $\point$ be a point of $\Man{}$.
For every pair $\triple{\avecf{1}}{\avecf{2}}{\avecf{3}}$ and $\triple{\avecff{1}}{\avecff{2}}{\avecff{3}}$ of
elements of ${\vectorfields{\Man{}}}^3$, if
\begin{equation}
\Foreach{i}{\seta{\triplet{1}{2}{3}}}
\func{\avecf{i}}{\point}=\func{\avecff{i}}{\point},
\end{equation}
then
\begin{align}
\func{\[\func{\mcurvature{\connection{}}}{\triplet{\avecf{1}}{\avecf{2}}{\avecf{3}}}\]}{\point}&=
\func{\[\func{\mcurvature{\connection{}}}{\triplet{\avecff{1}}{\avecff{2}}{\avecff{3}}}\]}{\point}\cr
&=\func{\[\func{\TTFequivv{\mcurvature{\connection{}}}}{\point}\]}
{\triplet{\func{\avecf{1}}{\point}}{\func{\avecf{2}}{\point}}{\func{\avecf{3}}{\point}}}.
\end{align}
\endcor
%%%%%%%%%%%%%%%%%%%%%%%%%%%%%%%%%%%%%%%%%%%%%%%%%%%%%%%%%%%%%%%%%%%%%%%%%%%%%%%%%%%%%%%%%%%%%%%%%%%%%%%%%%%%%%%%%%%%%%%%%%%%%%%%
\definition
Let $\connection{}$ be an affine connection on the manifold $\Man{}$.
We will call $\TTFequivv{\mcurvature{\connection{}}}$ the $\quotl$curvature tensor on $\Man{}$ relative to the
affine connection $\connection{}$$\quotr$.
\endef
%%%%%%%%%%%%%%%%%%%%%%%%%%%%%%%%%%%%%%%%%%%%%%%%%%%%%%%%%%%%%%%%%%%%%%%%%%%%%%%%%%%%%%%%%%%%%%%%%%%%%%%%%%%%%%%%%%%%%%%%%%%%%%%%
\lemma\label{lemrelationofcurvatureandcovariantderivatives}
Let $\connection{}$ be an affine connection on the manifold $\Man{}$.
Let $\function{\pCurve{}}{\Cprod{\interval{1}}{\interval{2}}}{\Man{}}$ be a one-parameter family of smooth curves on $\Man{}$.
Let $\valongc{}$ be a smooth vector field along $\pCurve{}$, that is an element of $\valongcurves{\Man{}}{\pCurve{}}$.
\begin{align}
\func{\(\func{\[\cmp{\pcovder{\connection{}}{\pCurve{}}{1}}{\pcovder{\connection{}}{\pCurve{}}{2}}\]}{\valongc{}}-
\func{\[\cmp{\pcovder{\connection{}}{\pCurve{}}{2}}{\pcovder{\connection{}}{\pCurve{}}{1}}\]}{\valongc{}}\)}{\binary{s}{t}}=
\func{\[\func{\TTFequivv{\mcurvature{\connection{}}}}{\func{\pCurve{}}{\binary{s}{t}}}\]}
{\binary{\binary{\func{\dot{\pcurve{\pCurve{}}{1}{t}}}{s}}{\func{\dot{\pcurve{\pCurve{}}{2}{s}}}{t}}}{\func{\valongc{}}{\binary{s}{t}}}}.
\end{align}
\proof
It is left as an exercise.
\endlem
%%%%%%%%%%%%%%%%%%%%%%%%%%%%%%%%%%%%%%%%%%%%%%%%%%%%%%%%%%%%%%%%%%%%%%%%%%%%%%%%%%%%%%%%%%%%%%%%%%%%%%%%%%%%%%%%%%%%%%%%%%%%%%%%
\theorem
Let $\connection{}$ be an affine connection on the manifold $\Man{}$.
Let $\interval{}$ be an open interval of $\R$ containing $0$, and
let $\function{\pCurve{}}{\Cprod{\interval{}}{\interval{}}}{\Man{}}$ be a one-parameter family of smooth curves on $\Man{}$.
Define $\point:=\func{\pCurve{}}{\opair{0}{0}}$, $x:=\func{\dot{\pcurve{\pCurve{}}{1}{0}}}{0}$, and
$y:=\func{\dot{\pcurve{\pCurve{}}{2}{0}}}{0}$. For any $z\in\tanspace{\point}{\Man{}}$,
\begin{align}
&~~~~\func{\[\func{\TTFequivv{\mcurvature{\connection{}}}}{\point}\]}{\triplet{x}{y}{z}}\cr
&=\lim_{\substack{{\delta_1\to 0}\\ {\delta_2\to 0}}}
\frac{\func{\[
\cmp{\ptransport{\pcurve{\pCurve{}}{1}{0}}{\delta_1}{0}{\connection{}}}
{\cmp{\ptransport{\pcurve{\pCurve{}}{2}{\delta_1}}{\delta_2}{0}{\connection{}}}
{\cmp{\ptransport{\pcurve{\pCurve{}}{1}{\delta_2}}{0}{\delta_1}{\connection{}}}
{\ptransport{\pcurve{\pCurve{}}{2}{0}}{0}{\delta_2}{\connection{}}}}}
\]}{z}-z}
{\delta_1\delta_2}.\cr
&{}
\end{align}
\proof
Fix $z\in\tanspace{\point}{\Man{}}$. Define the map
$\function{\valongc{}}{\Cprod{\interval{}}{\interval{}}}{\Tanbun{\Man{}}}$ as,
\begin{align}
&\begin{aligned}
\Foreach{t}{\interval{}}
\func{\valongc{}}{\binary{0}{t}}\eqdef
\func{\[\pvfalongc{\pcurve{\pCurve{}}{2}{0}}{0}{z}{\connection{}}\]}{t},
\end{aligned}\cr
&\begin{aligned}
\Foreach{\opair{s}{t}}{\Cprod{\interval{}}{\interval{}}}
\func{\valongc{}}{\binary{s}{t}}\eqdef
\func{\[\pvfalongc{\pcurve{\pCurve{}}{1}{t}}{0}{\func{\valongc{}}{\binary{0}{t}}}{\connection{}}\]}{s}.
\end{aligned}
\end{align}
Clearly $\func{\valongc{}}{\binary{0}{0}}=z$,
and it can be easily verified that $\valongc{}$ is a smooth map, and hence a smooth vector field along the one-parameter
family of curves $\pCurve{}$ on $\Man{}$, that is an element of $\valongcurves{\Man{}}{\pCurve{}}$.
According to \reflem{lemrelationofcurvatureandcovariantderivatives},
\begin{align}
\func{\[\func{\TTFequivv{\mcurvature{\connection{}}}}{\point}\]}{\triplet{x}{y}{z}}=
\func{\(\func{\[\cmp{\pcovder{\connection{}}{\pCurve{}}{1}}{\pcovder{\connection{}}{\pCurve{}}{2}}\]}{\valongc{}}\)}{\binary{0}{0}}-
\func{\(\func{\[\cmp{\pcovder{\connection{}}{\pCurve{}}{2}}{\pcovder{\connection{}}{\pCurve{}}{1}}\]}{\valongc{}}\)}{\binary{0}{0}}.
\end{align}
Based on the way $\valongc{}$ is defined and according to
\refdef{defcovariantderivativesofvectorfieldsalongoneparameterfamiliesofcurves}, it is completely clear that
$\func{\pcovder{\connection{}}{\pCurve{}}{1}}{\valongc{}}=0$, and thus
$\func{\[\cmp{\pcovder{\connection{}}{\pCurve{}}{2}}{\pcovder{\connection{}}{\pCurve{}}{1}}\]}{\valongc{}}=0$, and hence,
\begin{align}
\func{\[\func{\TTFequivv{\mcurvature{\connection{}}}}{\point}\]}{\triplet{x}{y}{z}}=
\func{\(\func{\[\cmp{\pcovder{\connection{}}{\pCurve{}}{1}}{\pcovder{\connection{}}{\pCurve{}}{2}}\]}{\valongc{}}\)}{\binary{0}{0}}.
\end{align}
According to \refdef{defcovariantderivativesofvectorfieldsalongoneparameterfamiliesofcurves} and
\refthm{thmcovariantderivativeintermsofparalleltransport},
\begin{align}
\Foreach{\opair{s}{t}}{\Cprod{\interval{}}{\interval{}}}
\func{\[\func{\pcovder{\connection{}}{\pCurve{}}{2}}{\valongc{}}\]}{\binary{s}{t}}&=
\func{\[\func{\covder{\connection{}}{\pcurve{\pCurve{}}{2}{s}}}{\reS{\valongc{}}{\Cprod{\seta{s}}{\interval{2}}}}\]}{t}\cr
&=\lim_{\tau\to t}\frac{\func{\[\ptransport{\pcurve{\pCurve{}}{2}{s}}{\tau}{t}{\connection{}}\]}{\func{\reS{\valongc{}}{\Cprod{\seta{s}}{\interval{2}}}}{\tau}}-
\func{\reS{\valongc{}}{\Cprod{\seta{s}}{\interval{2}}}}{t}}{\tau-t}\cr
&=\lim_{\tau\to t}\frac{\func{\[\ptransport{\pcurve{\pCurve{}}{2}{s}}{\tau}{t}{\connection{}}\]}{\func{\valongc{}}{\binary{s}{\tau}}}-
\func{\valongc{}}{\binary{s}{t}}}{\tau-t},
\end{align}
and
\begin{align}
\func{\(\func{\[\cmp{\pcovder{\connection{}}{\pCurve{}}{1}}{\pcovder{\connection{}}{\pCurve{}}{2}}\]}{\valongc{}}\)}{\binary{0}{0}}&=
\func{\[\func{\covder{\connection{}}{\pcurve{\pCurve{}}{1}{0}}}{\reS{\func{\pcovder{\connection{}}{\pCurve{}}{2}}{\valongc{}}}{\Cprod{\interval{1}}{\seta{0}}}}\]}{0}\cr
&=\lim_{s\to 0}\frac{\func{\[\ptransport{\pcurve{\pCurve{}}{1}{0}}{s}{0}{\connection{}}\]}{\func{\reS{\func{\pcovder{\connection{}}{\pCurve{}}{2}}{\valongc{}}}{\Cprod{\interval{1}}{\seta{0}}}}{s}}-
\func{\reS{\func{\pcovder{\connection{}}{\pCurve{}}{2}}{\valongc{}}}{\Cprod{\interval{1}}{\seta{0}}}}{0}}{s}\cr
&=\lim_{s\to 0}\frac{\func{\[\ptransport{\pcurve{\pCurve{}}{1}{0}}{s}{0}{\connection{}}\]}{\func{\func{\pcovder{\connection{}}{\pCurve{}}{2}}{\valongc{}}}{\binary{s}{0}}}-
\func{\func{\pcovder{\connection{}}{\pCurve{}}{2}}{\valongc{}}}{\binary{0}{0}}}{s}.
\end{align}
By combining the last two equations,
\begin{align}
&~~~~\func{\(\func{\[\cmp{\pcovder{\connection{}}{\pCurve{}}{1}}{\pcovder{\connection{}}{\pCurve{}}{2}}\]}{\valongc{}}\)}{\binary{0}{0}}\cr
&=\lim_{\substack{{s\to 0}\\ {\tau\to 0}}}\frac{1}{s\tau}\left(~
\func{\cmp{\[\ptransport{\pcurve{\pCurve{}}{1}{0}}{s}{0}{\connection{}}\]}{\[\ptransport{\pcurve{\pCurve{}}{2}{s}}{\tau}{0}{\connection{}}\]}}{\func{\valongc{}}{\binary{s}{\tau}}}-
\func{\[\ptransport{\pcurve{\pCurve{}}{1}{0}}{s}{0}{\connection{}}\]}{\func{\valongc{}}{\binary{s}{0}}}\right.\cr
&~~~~~~~~~~~~~~~~~-\left.\func{\[\ptransport{\pcurve{\pCurve{}}{2}{0}}{\tau}{0}{\connection{}}\]}{\func{\valongc{}}{\binary{0}{\tau}}}
+\func{\valongc{}}{\binary{0}{0}}\right).
\end{align}
Now, considering that $s\mapsto\func{\valongc{}}{\binary{s}{0}}$ is a parallel vector field
along the curve $\pcurve{\pCurve{}}{1}{0}$, it is clear that for every $s\in\interval{}$,
\begin{align}
\func{\[\ptransport{\pcurve{\pCurve{}}{1}{0}}{s}{0}{\connection{}}\]}{\func{\valongc{}}{\binary{s}{0}}}=
\func{\valongc{}}{\binary{0}{0}}=z,
\end{align}
and considering that $\tau\mapsto\func{\valongc{}}{\binary{0}{0\tau}}$ is a parallel vector field along the curve
$\pcurve{\pCurve{}}{2}{0}$,
\begin{align}
\func{\[\ptransport{\pcurve{\pCurve{}}{2}{0}}{\tau}{0}{\connection{}}\]}{\func{\valongc{}}{\binary{0}{\tau}}}=
\func{\valongc{}}{\binary{0}{0}}=z.
\end{align}
Moreover, according to the way $\valongc{}$ is defined, and based on \refdef{defparalleltransport}, it is clear that,
\begin{align}
\Foreach{\opair{s}{\tau}}{\Cprod{\interval{}}{\interval{}}}
\func{\valongc{}}{\binary{s}{\tau}}=
\func{\cmp{\[\ptransport{\pcurve{\pCurve{}}{1}{\tau}}{0}{s}{\connection{}}\]}{\[\ptransport{\pcurve{\pCurve{}}{2}{0}}{0}{\tau}{\connection{}}\]}}
{\func{\valongc{}}{\binary{0}{0}}}.
\end{align}
Therefore, combining these equations, it becomes evident that,
\begin{align}
&~~~~\func{\[\func{\TTFequivv{\mcurvature{\connection{}}}}{\point}\]}{\triplet{x}{y}{z}}\cr
&=\lim_{\substack{{s\to 0}\\ {\tau\to 0}}}\frac{
\func{\(\cmp{\cmp{\[\ptransport{\pcurve{\pCurve{}}{1}{0}}{s}{0}{\connection{}}\]}{\[\ptransport{\pcurve{\pCurve{}}{2}{s}}{\tau}{0}{\connection{}}\]}}
{\cmp{\[\ptransport{\pcurve{\pCurve{}}{1}{\tau}}{0}{s}{\connection{}}\]}{\[\ptransport{\pcurve{\pCurve{}}{2}{0}}{0}{\tau}{\connection{}}\]}}\)}{z}-z
}{s\tau}.\cr
&{}
\end{align}
\endthm
%%%%%%%%%%%%%%%%%%%%%%%%%%%%%%%%%%%%%%%%%%%%%%%%%%%%%%%%%%%%%%%%%%%%%%%%%%%%%%%%%%%%%%%%%%%%%%%%%%%%%%%%%%%%%%%%%%%%%%%%%%%%%%%%
%%%%%%%%%%%%%%%%%%%%%%%%%%%%%%%%%%%%%%%%%%%%%%%%%%%%%%%%%%%%%%%%%%%%%%%%%%%%%%%%%%%%%%%%%%%%%%%%%%%%%%%%%%%%%%%%%%%%%%%%%%%%%%%%
%%%%%%%%%%%%%%%%%%%%%%%%%%%%%%%%%%%%%%%%%%%%%%%%%%%%%%%%%%%%%%%%%%%%%%%%%%%%%%%%%%%%%%%%%%%%%%%%%%%%%%%%%%%%%%%%%%%%%%%%%%%%%%%%
%%%%%%%%%%%%%%%%%%%%%%%%%%%%%%%%%%%%%%%%%%%%%%%%%%%%%%%%%%%%%%%%%%%%%%%%%%%%%%%%%%%%%%%%%%%%%%%%%%%%%%%%%%%%%%%%%%%%%%%%%%%%%%%%
%%%%%%%%%%%%%%%%%%%%%%%%%%%%%%%%%%%%%%%%%%%%%%%%%%%%%%%%%%%%%%%%%%%%%%%%%%%%%%%%%%%%%%%%%%%%%%%%%%%%%%%%%%%%%%%%%%%%%%%%%%%%%%%%
%%%%%%%%%%%%%%%%%%%%%%%%%%%%%%%%%%%%%%%%%%%%%%%%%%%%%%%%%%%%%%%%%%%%%%%%%%%%%%%%%%%%%%%%%%%%%%%%%%%%%%%%%%%%%%%%%%%%%%%%%%%%%%%%
%%%%%%%%%%%%%%%%%%%%%%%%%%%%%%%%%%%%%%%%%%%%%%%%%%%%%%%%%%%%%%%%%%%%%%%%%%%%%%%%%%%%%%%%%%%%%%%%%%%%%%%%%%%%%%%%%%%%%%%%%%%%%%%%
%%%%%%%%%%%%%%%%%%%%%%%%%%%%%%%%%%%%%%%%%%%%%%%%%%%%%%%%%%%%%%%%%%%%%%%%%%%%%%%%%%%%%%%%%%%%%%%%%%%%%%%%%%%%%%%%%%%%%%%%%%%%%%%%
%%%%%%%%%%%%%%%%%%%%%%%%%%%%%%%%%%%%%%%%%%%%%%%%%%%%%%%%%%%%%%%%%%%%%%%%%%%%%%%%%%%%%%%%%%%%%%%%%%%%%%%%%%%%%%%%%%%%%%%%%%%%%%%%
%%%%%%%%%%%%%%%%%%%%%%%%%%%%%%%%%%%%%%%%%%%%%%%%%%%%%%%%%%%%%%%%%%%%%%%%%%%%%%%%%%%%%%%%%%%%%%%%%%%%%%%%%%%%%%%%%%%%%%%%%%%%%%%%
\section{Riemann, Ricci, and Scalar Curvatures}
%%%%%%%%%%%%%%%%%%%%%%%%%%%%%%%%%%%%%%%%%%%%%%%%%%%%%%%%%%%%%%%%%%%%%%%%%%%%%%%%%%%%%%%%%%%%%%%%%%%%%%%%%%%%%%%%%%%%%%%%%%%%%%%%
\fixed
$\Man{}$ is fixed as a manifold with dimension $m$, and $\metrictensor{}$ is fixed as a
semi-Riemannian metric on $\Man{}$, that is an element of $\metrictensors{\Man{}}{}$.
Here, we will use the notation $\connection{}$ for referring
to the Levi-Civita connection of the semi-Riemannian manifold $\opair{\Man{}}{\metrictensor{}}$, that is
$\LCconnection{\Man{}}{\metrictensor{}}$.
\endfixed
%%%%%%%%%%%%%%%%%%%%%%%%%%%%%%%%%%%%%%%%%%%%%%%%%%%%%%%%%%%%%%%%%%%%%%%%%%%%%%%%%%%%%%%%%%%%%%%%%%%%%%%%%%%%%%%%%%%%%%%%%%%%%%%%
\definition\label{defRiemannCurvature}
We define the map $\function{\Riemcurvature{\Man{}}{\metrictensor{}}}
{\Cprod{\Cprod{\vectorfields{\Man{}}}{\vectorfields{\Man{}}}}{\Cprod{\vectorfields{\Man{}}}{\vectorfields{\Man{}}}}}
{\smoothmaps{\Man{}}}$ as,
\begin{align}
\Foreach{\quadruple{\avecf{1}}{\avecf{2}}{\avecf{3}}{\avecf{4}}}
{\vectorfields{\Man{}}^4}
\func{\Riemcurvature{\Man{}}{\metrictensor{}}}{\Quadruple{\avecf{1}}{\avecf{2}}{\avecf{3}}{\avecf{4}}}\eqdef
\vfmetricprod{\Man{}}{\metrictensor{}}{\func{\mcurvature{\connection{}}}{\triplet{\avecf{1}}{\avecf{2}}{\avecf{3}}}}{\avecf{4}}.
\end{align}
$\Riemcurvature{\Man{}}{\metrictensor{}}$ is referred to as the $\quotl$Riemann curvature tensor of the
semi-Riemannian manifold $\opair{\Man{}}{\metrictensor{}}$$\quotr$.
\endef
%%%%%%%%%%%%%%%%%%%%%%%%%%%%%%%%%%%%%%%%%%%%%%%%%%%%%%%%%%%%%%%%%%%%%%%%%%%%%%%%%%%%%%%%%%%%%%%%%%%%%%%%%%%%%%%%%%%%%%%%%%%%%%%%
\theorem
$\Riemcurvature{\Man{}}{\metrictensor{}}$ is a $\opair{4}{0}$ ($4$-covariant) smooth tensor field on $\Man{}$. That is,
\begin{equation}
\Riemcurvature{\Man{}}{\metrictensor{}}\in
\TTF{4}{0}{\vectorfields{\Man{}}}.
\end{equation}
\proof
Considering $\smoothmaps{\Man{}}$-bilinearity of the operation $\vfmetricprodmap{\Man{}}{\metrictensor{}}$,
and $\smoothmaps{\Man{}}$-multilinearity of $\mcurvature{\connection{}}$, it is obvious.
\endthm
%%%%%%%%%%%%%%%%%%%%%%%%%%%%%%%%%%%%%%%%%%%%%%%%%%%%%%%%%%%%%%%%%%%%%%%%%%%%%%%%%%%%%%%%%%%%%%%%%%%%%%%%%%%%%%%%%%%%%%%%%%%%%%%%
%%%%%%%%%%%%%%%%%%%%%%%%%%%%%%%%%%%%%%%%%%%%%%%%%%%%%%%%%%%%%%%%%%%%%%%%%%%%%%%%%%%%%%%%%%%%%%%%%%%%%%%%%%%%%%%%%%%%%%%%%%%%%%%%
%%%%%%%%%%%%%%%%%%%%%%%%%%%%%%%%%%%%%%%%%%%%%%%%%%%%%%%%%%%%%%%%%%%%%%%%%%%%%%%%%%%%%%%%%%%%%%%%%%%%%%%%%%%%%%%%%%%%%%%%%%%%%%%%
\textit{
As discussed previoisly, $\TTF{r}{s}{\bbinary{\vectorfields{\Man{}}}{\vectorfields{\Man{}}}}$ is naturaly identified with
$\TF{r}{s}{\bbinary{\Tanbun{\Man{}}}{\Tanbun{\Man{}}}}$ for any non-negative integers $r$ and $s$ with $r+s>0$.
So, a given element $\function{\omega}{\vectorfields{\Man{}}}{\vectorfields{\Man{}}}$
of $\TTF{1}{0}{\bbinary{\vectorfields{\Man{}}}{\vectorfields{\Man{}}}}$ is naturally corresponded
to an element $\TFequivv{\omega}$ of $\TF{1}{0}{\bbinary{\Tanbun{\Man{}}}{\Tanbun{\Man{}}}}$.
Then, for every $\point\in\Man{}$,
$\function{\func{\TFequivv{\omega}}{\point}}{\tanspace{\point}{\Man{}}}{\tanspace{\point}{\Man{}}}$ is a linear map from
$\Tanspace{\point}{\Man{}}$ to $\Tanspace{\point}{\Man{}}$. The trace of $\omega$, denoted by
$\func{\Tr{\Man{}}}{\omega}$ is defined as the element of $\smoothmaps{\Man{}}$ that assigns to each $\point\in\Man{}$
the real value $\func{\Tr{\Tanspace{\point}{\Man{}}}}{\func{\TFequivv{\omega}}{\point}}$.
Considering the canonical $\smoothmaps{\Man{}}$-module structure of
$\TTF{1}{0}{\bbinary{\vectorfields{\Man{}}}{\vectorfields{\Man{}}}}$, it can be easily verified that
$\Tr{\Man{}}$ is a $\smoothmaps{\Man{}}$-linear operator on
$\TTF{1}{0}{\bbinary{\vectorfields{\Man{}}}{\vectorfields{\Man{}}}}$.
\\
%%%%%%%%%%%
In addition, a notion of trace is defined on $\TTF{2}{0}{\vectorfields{\Man{}}}$,
that is the $2$-covariant smooth tensor fields on $\Man{}$, when the manifold $\Man{}$ is equipped with a
semi-Riemannian metric $\metrictensor{}$. First, for every $\mathcal{T}\in
\TTF{2}{0}{\vectorfields{\Man{}}}$, we define the $\opair{1}{1}$ smooth tensor field
$\function{\func{\mtsharp{\metrictensor{}}}{\mathcal{T}}}{\Cprod{\vectorfields{\Man{}}}{\oneforms{\Man{}}}}{\smoothmaps{\Man{}}}$
on $\Man{}$ as,
\begin{equation}
\Foreach{\opair{\avecf{}}{\aoneform{}}}{\Cprod{\vectorfields{\Man{}}}{\oneforms{\Man{}}}}
\func{\[\func{\mtsharp{\metrictensor{}}}{\mathcal{T}}\]}{\binary{\avecf{}}{\aoneform{}}}\eqdef
\func{\mathcal{T}}{\binary{\avecf{}}{\func{\mtsharp{\metrictensor{}}}{\aoneform{}}}}.
\end{equation}
We have already stated the definition of the notion of trace on
$\TTF{1}{1}{\vectorfields{\Man{}}}$.
The trace of any $\mathcal{T}\in\TTF{2}{0}{\vectorfields{\Man{}}}$, denoted by
$\func{\vstrace{\opair{\Man{}}{\metrictensor{}}}}{\mathcal{T}}$ is defined as,
\begin{equation}
\func{\vstrace{\opair{\Man{}}{\metrictensor{}}}}{\mathcal{T}}:=
\func{\vstrace{\Man{}}}{\func{\mtsharp{\metrictensor{}}}{\mathcal{T}}},
\end{equation}
which is an element of $\smoothmaps{\Man{}}$.
Therefore, given a $\mathcal{T}\in\TTF{2}{0}{\vectorfields{\Man{}}}$,
for every point $\point$ of $\Man{}$, and any ordered-basis $\vsbase{}=\mtuple{\vsbase{1}}{\vsbase{m}}$
of $\Tanspace{\point}{\Man{}}$,
\begin{align}
\func{\[\func{\vstrace{\opair{\Man{}}{\metrictensor{}}}}{\mathcal{T}}\]}{\point}=
\sum_{i=1}^{m}\func{\mathcal{T}}{\binary{\vsbase{i}}{\func{\mtsharp{\metrictensor{}}}{\dualvsbase{i}}}},
\end{align}
where, $\dualvsbase{}=\mtuple{\dualvsbase{1}}{\dualvsbase{m}}$ denotes the dual of the ordered-basis $\vsbase{}$.
}
%%%%%%%%%%%%%%%%%%%%%%%%%%%%%%%%%%%%%%%%%%%%%%%%%%%%%%%%%%%%%%%%%%%%%%%%%%%%%%%%%%%%%%%%%%%%%%%%%%%%%%%%%%%%%%%%%%%%%%%%%%%%%%%%
\definition\label{defRicciCurvature}
The map $\function{\Ricci{\Man{}}{\metrictensor{}}}{\Cprod{\vectorfields{\Man{}}}{\vectorfields{\Man{}}}}{\smoothmaps{\Man{}}}$
is defined as,
\begin{equation}
\Foreach{\opair{\avecf{}}{\avecff{}}}{\Cprod{\vectorfields{\Man{}}}{\vectorfields{\Man{}}}}
\func{\Ricci{\Man{}}{\metrictensor{}}}{\binary{\avecf{}}{\avecff{}}}\eqdef
\func{\Tr{\Man{}}}{\func{\Delta}{\binary{\avecf{}}{\avecff{}}}},
\end{equation}
where, the $\smoothmaps{\Man{}}$-bilinear map
$\function{\Delta}{\Cprod{\vectorfields{\Man{}}}{\vectorfields{\Man{}}}}
{\TTF{1}{0}{\bbinary{\vectorfields{\Man{}}}{\vectorfields{\Man{}}}}}$ is defined as,
\begin{equation}
\Foreach{\opair{\avecf{}}{\avecff{}}}{\Cprod{\vectorfields{\Man{}}}{\vectorfields{\Man{}}}}
\Foreach{\avecf{1}}{\vectorfields{\Man{}}}
\func{\[\func{\Delta}{\binary{\avecf{}}{\avecff{}}}\]}{\avecf{1}}\eqdef
\func{\mcurvature{\connection{}}}{\triplet{\avecf{1}}{\avecf{}}{\avecff{}}}.
\end{equation}
$\Ricci{\Man{}}{\metrictensor{}}$ is referred to as the $\quotl$Ricci curvature of the semi-Riemannian manifold
$\opair{\Man{}}{\metrictensor{}}$$\quotr$.
\endef
%%%%%%%%%%%%%%%%%%%%%%%%%%%%%%%%%%%%%%%%%%%%%%%%%%%%%%%%%%%%%%%%%%%%%%%%%%%%%%%%%%%%%%%%%%%%%%%%%%%%%%%%%%%%%%%%%%%%%%%%%%%%%%%%
\theorem
$\Ricci{\Man{}}{\metrictensor{}}$ is a bilinear map from $\Cprod{\vectorfields{\Man{}}}{\vectorfields{\Man{}}}$ to
$\smoothmaps{\Man{}}$, that is an element of $\TTF{2}{0}{\vectorfields{\Man{}}}$.
\proof
Considering the $\smoothmaps{\Man{}}$-bilinearity of the map $\Delta$ defined in \refdef{defRicciCurvature},
and the $\smoothmaps{\Man{}}$-linearity of the trace operator $\Tr{\Man{}}$, it is clear.
\endthm
%%%%%%%%%%%%%%%%%%%%%%%%%%%%%%%%%%%%%%%%%%%%%%%%%%%%%%%%%%%%%%%%%%%%%%%%%%%%%%%%%%%%%%%%%%%%%%%%%%%%%%%%%%%%%%%%%%%%%%%%%%%%%%%%
\definition
The element $\scalarcurvature{\Man{}}{\metrictensor{}}$ of $\smoothmaps{\Man{}}$ is defined as,
\begin{equation}
\scalarcurvature{\Man{}}{\metrictensor{}}:=
\func{\vstrace{\opair{\Man{}}{\metrictensor{}}}}{\Ricci{\Man{}}{\metrictensor{}}}.
\end{equation}
$\scalarcurvature{\Man{}}{\metrictensor{}}$ is referred to as the
$\quotl$scalar curvature of the semi-Riemannian manifold $\opair{\Man{}}{\metrictensor{}}$$\quotr$.
\endef
%%%%%%%%%%%%%%%%%%%%%%%%%%%%%%%%%%%%%%%%%%%%%%%%%%%%%%%%%%%%%%%%%%%%%%%%%%%%%%%%%%%%%%%%%%%%%%%%%%%%%%%%%%%%%%%%%%%%%%%%%%%%%%%%
\definition
\begin{itemize}
\item
It is said that $\quotl$the Ricci curvature of the semi-Riemannian manifold $\opair{\Man{}}{\metrictensor{}}$ is
bounded from below by the factor $\kappa$$\quotr$ if $\TTFequiv{\Ricci{\Man{}}{\metrictensor{}}}>\kappa\metrictensor{}$,
where $\kappa$ is interpreted as the element of $\smoothmaps{\Man{}}$ with constant value $\kappa\in\R$.
\item
It is said that $\quotl$the Ricci curvature of the semi-Riemannian manifold $\opair{\Man{}}{\metrictensor{}}$ is
bounded from above by the factor $\kappa$$\quotr$ if $\TTFequiv{\Ricci{\Man{}}{\metrictensor{}}}<\kappa\metrictensor{}$,
where $\kappa$ is interpreted as the element of $\smoothmaps{\Man{}}$ with constant value $\kappa\in\R$.
\item
The semi-Riemannian manifold $\opair{\Man{}}{\metrictensor{}}$ is called an
$\quotl$Einstein manifold of constant $\kappa$$\quotr$
if $\TTFequiv{\Ricci{\Man{}}{\metrictensor{}}}=\kappa\metrictensor{}$, $\kappa$ denoting the constant smooth map
on $\Man{}$ wirh real value $\kappa$.
\end{itemize}
\endef
%%%%%%%%%%%%%%%%%%%%%%%%%%%%%%%%%%%%%%%%%%%%%%%%%%%%%%%%%%%%%%%%%%%%%%%%%%%%%%%%%%%%%%%%%%%%%%%%%%%%%%%%%%%%%%%%%%%%%%%%%%%%%%%%
\definition
Let $\opair{\U}{\phi}$ be a chart of $\Man{}$, and let $\mtuple{\localframevecf{1}}{\localframevecf{m}}$ denote the system
of local frame fields corresponded to $\opair{\U}{\phi}$, and let $\mtuple{\localframeoneform{1}}{\localframeoneform{m}}$
denote the dual of this local frame field. We define the system of smooth functions $\Riemcurvaturecoef{\phi}{i}{j}{k}{l}$
from $\subman{\Man{}}{\U}$ to $\R$, where each $i$, $j$, $k$, and $l$ varies over $\seta{1,~2,~3,~4}$, as
\begin{equation}
\reS{\TTFequiv{\Riemcurvature{\Man{}}{\metrictensor{}}}}{\U}=\sum_{i=1}^{m}\sum_{j=1}^{m}\sum_{k=1}^{m}\sum_{l=1}^{m}
\Riemcurvaturecoef{\phi}{i}{j}{k}{l}\localframeoneform{i}\tensor{}\localframeoneform{j}\tensor{}
\localframeoneform{k}\tensor{}\localframeoneform{l}.
\end{equation}
\endef
%%%%%%%%%%%%%%%%%%%%%%%%%%%%%%%%%%%%%%%%%%%%%%%%%%%%%%%%%%%%%%%%%%%%%%%%%%%%%%%%%%%%%%%%%%%%%%%%%%%%%%%%%%%%%%%%%%%%%%%%%%%%%%%%
\definition
Let $\opair{\U}{\phi}$ be a chart of $\Man{}$, and let $\mtuple{\localframevecf{1}}{\localframevecf{m}}$ denote the system
of local frame fields corresponded to $\opair{\U}{\phi}$, and let $\mtuple{\localframeoneform{1}}{\localframeoneform{m}}$
denote the dual of this local frame field. We define the system of smooth functions $\Riccicoef{\phi}{i}{j}$
from $\subman{\Man{}}{\U}$ to $\R$, where each $i$, $j$, $k$, and $l$ varies over $\seta{1,~2,~3,~4}$, as
\begin{equation}
\reS{\TTFequiv{\Ricci{\Man{}}{\metrictensor{}}}}{\U}=\sum_{i=1}^{m}\sum_{j=1}^{m}
\Riccicoef{\phi}{i}{j}\localframeoneform{i}\tensor{}\localframeoneform{j}.
\end{equation}
\endef
%%%%%%%%%%%%%%%%%%%%%%%%%%%%%%%%%%%%%%%%%%%%%%%%%%%%%%%%%%%%%%%%%%%%%%%%%%%%%%%%%%%%%%%%%%%%%%%%%%%%%%%%%%%%%%%%%%%%%%%%%%%%%%%%
\theorem
Let $\opair{\U}{\phi}$ be a chart of $\Man{}$, and let $\mtuple{\localframevecf{1}}{\localframevecf{m}}$ denote the system
of local frame fields corresponded to $\opair{\U}{\phi}$, and let $\mtuple{\localframeoneform{1}}{\localframeoneform{m}}$
denote the dual of this local frame field. Also let
$\reS{\metrictensor{}}{\U}=\sum_{i=1}^{m}\sum_{j=1}^{m}\metrictensorchart{\phi}{i}{j}
\localframeoneform{i}\tensor{}\localframeoneform{j}$, and let $\metrictensorchartinv{\phi}{}{}$ denote inverse of the
smooth real-valued functions $\metrictensorchart{\phi}{}{}$ on $\U$.
Let $\Christoffel{\phi}{i}{j}{k}$-s denote the Christoffel symbols of
the Levi-Civita connection $\connection{}$ of $\opair{\Man{}}{\metrictensor{}}$ with respect to the chart $\opair{\U}{\phi}$.
\begin{itemize}
\item
\begin{align}
\Riemcurvaturecoef{\phi}{i}{j}{k}{l}=
\sum_{n=1}^{m}\metrictensorchart{\phi}{l}{n}\[\lieder{\localframevecf{i}}{\Christoffel{\phi}{n}{j}{k}}-
\lieder{\localframevecf{j}}{\Christoffel{\phi}{n}{i}{k}}+
\sum_{\alpha=1}^{m}\(\Christoffel{\phi}{\alpha}{j}{k}\Christoffel{\phi}{n}{i}{\alpha}\)-
\sum_{\alpha=1}^{m}\(\Christoffel{\phi}{\alpha}{i}{k}\Christoffel{\phi}{n}{j}{\alpha}\)\].
\end{align}
\item
\begin{equation}
\Riccicoef{\phi}{i}{j}=\sum_{\alpha=1}^{m}\sum_{\beta=1}^{m}
\metrictensorchartinv{\phi}{\alpha}{\beta}\Riemcurvaturecoef{\phi}{\alpha}{i}{j}{\beta}.
\end{equation}
\item
\begin{equation}
\reS{\scalarcurvature{\Man{}}{\metrictensor{}}}{\U}=
\sum_{i=1}^{m}\sum_{j=1}^{m}
\metrictensorchartinv{\phi}{i}{j}\Riccicoef{\phi}{i}{j}.
\end{equation}
\end{itemize}
\proof
It is left as an exercise.
\endthm
%%%%%%%%%%%%%%%%%%%%%%%%%%%%%%%%%%%%%%%%%%%%%%%%%%%%%%%%%%%%%%%%%%%%%%%%%%%%%%%%%%%%%%%%%%%%%%%%%%%%%%%%%%%%%%%%%%%%%%%%%%%%%%%%
\theorem\label{thmfirstBianchiIdentity}
For every triple $\triplet{\avecf{1}}{\avecf{2}}{\avecf{3}}$ of smooth vector fields on $\Man{}$,
\begin{equation}
\func{\mcurvature{\connection{}}}{\triplet{\avecf{1}}{\avecf{2}}{\avecf{3}}}+
\func{\mcurvature{\connection{}}}{\triplet{\avecf{2}}{\avecf{3}}{\avecf{1}}}+
\func{\mcurvature{\connection{}}}{\triplet{\avecf{3}}{\avecf{1}}{\avecf{2}}}=0.
\end{equation}
\proof
According to \refdef{defcurvatureofanaffineconnection}, and considering that the torsion of Levi-Civita connection
vanishes identically, and according to \refdef{deftorsionofaffineconnections}, and further considering that the vector-space
$\vectorfields{\Man{}}$ is a Lie-algebra when endowed with the Lie-bracket operation on it (so, Lie-bracket operation
satisfies the Jacobi identity),
\begin{align}
&~~~~\func{\mcurvature{\connection{}}}{\triplet{\avecf{1}}{\avecf{2}}{\avecf{3}}}+
\func{\mcurvature{\connection{}}}{\triplet{\avecf{2}}{\avecf{3}}{\avecf{1}}}+
\func{\mcurvature{\connection{}}}{\triplet{\avecf{3}}{\avecf{1}}{\avecf{2}}}\cr
&=~~\(\con{\avecf{1}}{\con{\avecf{2}}{\avecf{3}}}-\con{\avecf{2}}{\con{\avecf{1}}{\avecf{3}}}-\con{\liebracket{\avecf{1}}{\avecf{2}}{}}{\avecf{3}}\)+
\(\con{\avecf{2}}{\con{\avecf{3}}{\avecf{1}}}-\con{\avecf{3}}{\con{\avecf{2}}{\avecf{1}}}-\con{\liebracket{\avecf{2}}{\avecf{3}}{}}{\avecf{1}}\)\cr
&~~~+\(\con{\avecf{3}}{\con{\avecf{1}}{\avecf{2}}}-\con{\avecf{1}}{\con{\avecf{3}}{\avecf{2}}}-\con{\liebracket{\avecf{3}}{\avecf{1}}{}}{\avecf{2}}\)\cr
&=~~\con{\avecf{1}}{\(\con{\avecf{2}}{\avecf{3}}-\con{\avecf{3}}{\avecf{2}}\)}+\con{\avecf{2}}{\(\con{\avecf{3}}{\avecf{1}}-\con{\avecf{1}}{\avecf{3}}\)}+
\con{\avecf{3}}{\(\con{\avecf{1}}{\avecf{2}}-\con{\avecf{2}}{\avecf{1}}\)}\cr
&~~~-\(\con{\liebracket{\avecf{1}}{\avecf{2}}{}}{\avecf{3}}+\con{\liebracket{\avecf{2}}{\avecf{3}}{}}{\avecf{1}}+
\con{\liebracket{\avecf{3}}{\avecf{1}}{}}{\avecf{2}}\)\cr
&=~~\con{\avecf{1}}{\(\liebracket{\avecf{2}}{\avecf{3}}{}\)}+\con{\avecf{2}}{\(\liebracket{\avecf{3}}{\avecf{1}}{}\)}+
\con{\avecf{3}}{\(\liebracket{\avecf{1}}{\avecf{2}}{}\)}\cr
&~~~-\(\con{\liebracket{\avecf{1}}{\avecf{2}}{}}{\avecf{3}}+\con{\liebracket{\avecf{2}}{\avecf{3}}{}}{\avecf{1}}+
\con{\liebracket{\avecf{3}}{\avecf{1}}{}}{\avecf{2}}\)\cr
&=~~\(\con{\avecf{1}}{\(\liebracket{\avecf{2}}{\avecf{3}}{}\)}-\con{\liebracket{\avecf{2}}{\avecf{3}}{}}{\avecf{1}}\)+
\(\con{\avecf{2}}{\(\liebracket{\avecf{3}}{\avecf{1}}{}\)}-\con{\liebracket{\avecf{3}}{\avecf{1}}{}}{\avecf{2}}\)\cr
&~~~+\(\con{\avecf{3}}{\(\liebracket{\avecf{1}}{\avecf{2}}{}\)}-\con{\liebracket{\avecf{1}}{\avecf{2}}{}}{\avecf{3}}\)\cr
&=\liebracket{\avecf{1}}{\liebracket{\avecf{2}}{\avecf{3}}{}}{}+\liebracket{\avecf{2}}{\liebracket{\avecf{3}}{\avecf{1}}{}}{}+
\liebracket{\avecf{3}}{\liebracket{\avecf{1}}{\avecf{2}}{}}{}\cr
&=0.
\end{align}
\endthm
%%%%%%%%%%%%%%%%%%%%%%%%%%%%%%%%%%%%%%%%%%%%%%%%%%%%%%%%%%%%%%%%%%%%%%%%%%%%%%%%%%%%%%%%%%%%%%%%%%%%%%%%%%%%%%%%%%%%%%%%%%%%%%%%
\textit{
The equality in the previous theorem is famous for the $\quotl$first (or, algebraic) Bianchi identity
on the semi-Riemannian manifold $\opair{\Man{}}{\metrictensor{}}$$\quotr$.
}
%%%%%%%%%%%%%%%%%%%%%%%%%%%%%%%%%%%%%%%%%%%%%%%%%%%%%%%%%%%%%%%%%%%%%%%%%%%%%%%%%%%%%%%%%%%%%%%%%%%%%%%%%%%%%%%%%%%%%%%%%%%%%%%%
\theorem\label{thmRiemanncurvatureproperties0}
For every quadruple of smooth vector fields $\Quadruple{\avecf{1}}{\avecf{2}}{\avecf{3}}{\avecf{4}}$ on $\Man{}$,
the following equalities are hold.
\begin{itemize}
\item[\myitem{1.~}]
\begin{equation}
\func{\Riemcurvature{\Man{}}{\metrictensor{}}}{\Quadruple{\avecf{1}}{\avecf{2}}{\avecf{3}}{\avecf{4}}}=
-\func{\Riemcurvature{\Man{}}{\metrictensor{}}}{\Quadruple{\avecf{2}}{\avecf{1}}{\avecf{3}}{\avecf{4}}}.
\end{equation}
\item[\myitem{2.~}]
\begin{equation}
\func{\Riemcurvature{\Man{}}{\metrictensor{}}}{\Quadruple{\avecf{1}}{\avecf{2}}{\avecf{3}}{\avecf{4}}}=
-\func{\Riemcurvature{\Man{}}{\metrictensor{}}}{\Quadruple{\avecf{1}}{\avecf{2}}{\avecf{4}}{\avecf{3}}}.
\end{equation}
\item[\myitem{3.~}]
\begin{equation}
\func{\Riemcurvature{\Man{}}{\metrictensor{}}}{\Quadruple{\avecf{1}}{\avecf{2}}{\avecf{3}}{\avecf{4}}}+
\func{\Riemcurvature{\Man{}}{\metrictensor{}}}{\Quadruple{\avecf{2}}{\avecf{3}}{\avecf{1}}{\avecf{4}}}+
\func{\Riemcurvature{\Man{}}{\metrictensor{}}}{\Quadruple{\avecf{3}}{\avecf{1}}{\avecf{2}}{\avecf{4}}}=0.
\end{equation}
\item[\myitem{4.~}]
\begin{equation}
\func{\Riemcurvature{\Man{}}{\metrictensor{}}}{\Quadruple{\avecf{1}}{\avecf{2}}{\avecf{3}}{\avecf{4}}}=
-\func{\Riemcurvature{\Man{}}{\metrictensor{}}}{\Quadruple{\avecf{3}}{\avecf{4}}{\avecf{1}}{\avecf{2}}}.
\end{equation}
\end{itemize}
\proof
\begin{itemize}
\item[\myitem{pr-1.}]
The first identity is an immediate consequence of the definition of $\mcurvature{\connection{}}$, definition of the
Riemann curvature, and the bilinearity of the operation $\vfmetricprodmap{\Man{}}{\metrictensor{}}$.
\item[\myitem{pr-2.}]
Fix smooth vector fields $\avecf{1}$ and $\avecf{2}$ on $\Man{}$.\\
Since $\connection{}$ is the Levi-Civita connection of the semi-Riemannian manifold $\opair{\Man{}}{\metrictensor{}}$,
it is a metri-connection on $\opair{\Man{}}{\metrictensor{}}$, and thus according to \refdef{defmetricconnections},
for every $\avecff{}\in\vectorfields{\Man{}}$,
\begin{align}
\lieder{\avecf{1}}{\lieder{\avecf{2}}{\vfmetricproduct{\avecff{}}{\avecff{}}}}&=
\lieder{\avecf{1}}{\(2\vfmetricproduct{\con{\avecf{2}}{\avecff{}}}{\avecff{}}\)}\cr
&=2\vfmetricproduct{\con{\avecf{1}}{\con{\avecf{2}}{\avecff{}}}}{\avecff{}}
+2\vfmetricproduct{\con{\avecf{1}}{\avecff{}}}{\con{\avecf{2}}{\avecff{}}},
\end{align}
and similarly,
\begin{align}
\lieder{\avecf{2}}{\lieder{\avecf{1}}{\vfmetricproduct{\avecff{}}{\avecff{}}}}=
2\vfmetricproduct{\con{\avecf{2}}{\con{\avecf{1}}{\avecff{}}}}{\avecff{}}
+2\vfmetricproduct{\con{\avecf{2}}{\avecff{}}}{\con{\avecf{1}}{\avecff{}}},
\end{align}
and
\begin{align}
\lieder{\liebracket{\avecf{1}}{\avecf{2}}{}}{\vfmetricproduct{\avecff{}}{\avecff{}}}=
2\vfmetricproduct{\con{\liebracket{\avecf{1}}{\avecf{2}}{}}{\avecff{}}}{\avecff{}}.
\end{align}
Therefore,
\begin{align}
&~~~~~\lieder{\avecf{1}}{\lieder{\avecf{2}}{\vfmetricproduct{\avecff{}}{\avecff{}}}}-
\lieder{\avecf{2}}{\lieder{\avecf{1}}{\vfmetricproduct{\avecff{}}{\avecff{}}}}-
\lieder{\liebracket{\avecf{1}}{\avecf{2}}{}}{\vfmetricproduct{\avecff{}}{\avecff{}}}\cr
&=2\(\vfmetricproduct{\con{\avecf{1}}{\con{\avecf{2}}{\avecff{}}}}{\avecff{}}-
\vfmetricproduct{\con{\avecf{2}}{\con{\avecf{1}}{\avecff{}}}}{\avecff{}}-
\vfmetricproduct{\con{\liebracket{\avecf{1}}{\avecf{2}}{}}{\avecff{}}}{\avecff{}}\)\cr
&=2\vfmetricproduct{\con{\avecf{1}}{\con{\avecf{2}}{\avecff{}}}-\con{\avecf{2}}{\con{\avecf{1}}{\avecff{}}}-
\con{\liebracket{\avecf{1}}{\avecf{2}}{}}{\avecff{}}}{\avecff{}}\cr
&=2\vfmetricproduct{\func{\mcurvature{\connection{}}}{\triplet{\avecf{1}}{\avecf{2}}{\avecff{}}}}{\avecff{}}\cr
&=2\func{\Riemcurvature{\Man{}}{\metrictensor{}}}{\Quadruple{\avecf{1}}{\avecf{2}}{\avecff{}}{\avecff{}}}.
\end{align}
On the other hand, considering that $\lieder{\liebracket{\avecf{1}}{\avecf{2}}{}}{}=
\lieder{\avecf{1}}{\lieder{\avecf{2}}{}}-\lieder{\avecf{2}}{\lieder{\avecf{1}}{}}$,
it is clear that,
\begin{equation}
\lieder{\avecf{1}}{\lieder{\avecf{2}}{\vfmetricproduct{\avecff{}}{\avecff{}}}}-
\lieder{\avecf{2}}{\lieder{\avecf{1}}{\vfmetricproduct{\avecff{}}{\avecff{}}}}-
\lieder{\liebracket{\avecf{1}}{\avecf{2}}{}}{\vfmetricproduct{\avecff{}}{\avecff{}}}=0.
\end{equation}
Thus,
\begin{equation}\label{thmRiemanncurvatureproperties0eqesp}
\func{\Riemcurvature{\Man{}}{\metrictensor{}}}{\Quadruple{\avecf{1}}{\avecf{2}}{\avecff{}}{\avecff{}}}=0.
\end{equation}
So, in particular, by fixing another pair $\avecf{3}$ and $\avecf{4}$ of elements of $\vectorfields{\Man{}}$,
\begin{equation}
\func{\Riemcurvature{\Man{}}{\metrictensor{}}}{\Quadruple{\avecf{1}}{\avecf{2}}{\avecf{1}+\avecf{2}}{\avecf{1}+\avecf{2}}}=0.
\end{equation}
Using the $\smoothmaps{\Man{}}$-multilinearity of $\Riemcurvature{\Man{}}{\metrictensor{}}$, and considering that
\Ref{thmRiemanncurvatureproperties0eqesp} holds for every $\avecff{}\in\vectorfields{\Man{}}$,
\begin{align}
\func{\Riemcurvature{\Man{}}{\metrictensor{}}}{\Quadruple{\avecf{1}}{\avecf{2}}{\avecf{3}+\avecf{4}}{\avecf{3}+\avecf{4}}}&=~~
\func{\Riemcurvature{\Man{}}{\metrictensor{}}}{\Quadruple{\avecf{1}}{\avecf{2}}{\avecf{3}}{\avecf{3}}}\cr
&~~~+\func{\Riemcurvature{\Man{}}{\metrictensor{}}}{\Quadruple{\avecf{1}}{\avecf{2}}{\avecf{3}}{\avecf{4}}}\cr
&~~~+\func{\Riemcurvature{\Man{}}{\metrictensor{}}}{\Quadruple{\avecf{1}}{\avecf{2}}{\avecf{4}}{\avecf{3}}}\cr
&~~~+\func{\Riemcurvature{\Man{}}{\metrictensor{}}}{\Quadruple{\avecf{1}}{\avecf{2}}{\avecf{4}}{\avecf{4}}}\cr
&=~~\func{\Riemcurvature{\Man{}}{\metrictensor{}}}{\Quadruple{\avecf{1}}{\avecf{2}}{\avecf{3}}{\avecf{4}}}\cr
&~~~+\func{\Riemcurvature{\Man{}}{\metrictensor{}}}{\Quadruple{\avecf{1}}{\avecf{2}}{\avecf{4}}{\avecf{3}}}.
\end{align}
Therefore,
\begin{equation}
\func{\Riemcurvature{\Man{}}{\metrictensor{}}}{\Quadruple{\avecf{1}}{\avecf{2}}{\avecf{3}}{\avecf{4}}}
+\func{\Riemcurvature{\Man{}}{\metrictensor{}}}{\Quadruple{\avecf{1}}{\avecf{2}}{\avecf{4}}{\avecf{3}}}=0.
\end{equation}
\item[\myitem{pr-3.}]
According to \refdef{defcurvatureofanaffineconnection}, \refthm{thmfirstBianchiIdentity},
and considering the $\smoothmaps{\Man{}}$-bilinearity of the operation $\vfmetricprodmap{\Man{}}{\metrictensor{}}$,
it becomes evident that,
\begin{align}
&~~~~~\func{\Riemcurvature{\Man{}}{\metrictensor{}}}{\Quadruple{\avecf{1}}{\avecf{2}}{\avecf{3}}{\avecf{4}}}+
\func{\Riemcurvature{\Man{}}{\metrictensor{}}}{\Quadruple{\avecf{2}}{\avecf{3}}{\avecf{1}}{\avecf{4}}}+
\func{\Riemcurvature{\Man{}}{\metrictensor{}}}{\Quadruple{\avecf{3}}{\avecf{1}}{\avecf{2}}{\avecf{4}}}\cr
&=\vfmetricproduct{\func{\mcurvature{\connection{}}}{\triplet{\avecf{1}}{\avecf{2}}{\avecf{3}}}}{\avecf{4}}+
\vfmetricproduct{\func{\mcurvature{\connection{}}}{\triplet{\avecf{2}}{\avecf{3}}{\avecf{1}}}}{\avecf{4}}+
\vfmetricproduct{\func{\mcurvature{\connection{}}}{\triplet{\avecf{3}}{\avecf{1}}{\avecf{2}}}}{\avecf{4}}\cr
&=\vfmetricproduct{\func{\mcurvature{\connection{}}}{\triplet{\avecf{1}}{\avecf{2}}{\avecf{3}}}+
\func{\mcurvature{\connection{}}}{\triplet{\avecf{2}}{\avecf{3}}{\avecf{1}}}+
\func{\mcurvature{\connection{}}}{\triplet{\avecf{3}}{\avecf{1}}{\avecf{2}}}}{\avecf{4}}\cr
&=\vfmetricproduct{0}{\avecf{4}}=0.
\end{align}
\item[\myitem{pr-4.}]
It is an straightforwar consequence of the equalities \myitem{1}, \myitem{2}, and \myitem{3}.
\end{itemize}
\endthm
%%%%%%%%%%%%%%%%%%%%%%%%%%%%%%%%%%%%%%%%%%%%%%%%%%%%%%%%%%%%%%%%%%%%%%%%%%%%%%%%%%%%%%%%%%%%%%%%%%%%%%%%%%%%%%%%%%%%%%%%%%%%%%%%
\theorem
For every $\binary{\Quadruple{\avecf{1}}{\avecf{2}}{\avecf{3}}{\avecff{1}}}{\avecff{2}}\in\vectorfields{\Man{}}$,
\begin{align}
\func{\[\func{\tfconnection{\connection{}}{4}{0}}{\binary{\avecf{1}}{\Riemcurvature{\Man{}}{\metrictensor{}}}}\]}
{\Quadruple{\avecf{2}}{\avecf{3}}{\avecff{1}}{\avecff{2}}}&+
\func{\[\func{\tfconnection{\connection{}}{4}{0}}{\binary{\avecf{2}}{\Riemcurvature{\Man{}}{\metrictensor{}}}}\]}
{\Quadruple{\avecf{3}}{\avecf{1}}{\avecff{1}}{\avecff{2}}}\cr
&+\func{\[\func{\tfconnection{\connection{}}{4}{0}}{\binary{\avecf{3}}{\Riemcurvature{\Man{}}{\metrictensor{}}}}\]}
{\Quadruple{\avecf{1}}{\avecf{2}}{\avecff{1}}{\avecff{2}}}=0.\cr
&{}
\end{align}
\proof
This equality can be computed according to \refdef{defcovariantderivativeoftensorfields2}, \refdef{defRiemannCurvature},
and \refdef{defLeviCivitaconnection}.
\endthm
%%%%%%%%%%%%%%%%%%%%%%%%%%%%%%%%%%%%%%%%%%%%%%%%%%%%%%%%%%%%%%%%%%%%%%%%%%%%%%%%%%%%%%%%%%%%%%%%%%%%%%%%%%%%%%%%%%%%%%%%%%%%%%%%
\textit{
The equality in the previous theorem is famous for the $\quotl$second (or, differential) Bianchi identity
on the semi-Riemannian manifold $\opair{\Man{}}{\metrictensor{}}$$\quotr$.
}
%%%%%%%%%%%%%%%%%%%%%%%%%%%%%%%%%%%%%%%%%%%%%%%%%%%%%%%%%%%%%%%%%%%%%%%%%%%%%%%%%%%%%%%%%%%%%%%%%%%%%%%%%%%%%%%%%%%%%%%%%%%%%%%%
%%%%%%%%%%%%%%%%%%%%%%%%%%%%%%%%%%%%%%%%%%%%%%%%%%%%%%%%%%%%%%%%%%%%%%%%%%%%%%%%%%%%%%%%%%%%%%%%%%%%%%%%%%%%%%%%%%%%%%%%%%%%%%%%
%%%%%%%%%%%%%%%%%%%%%%%%%%%%%%%%%%%%%%%%%%%%%%%%%%%%%%%%%%%%%%%%%%%%%%%%%%%%%%%%%%%%%%%%%%%%%%%%%%%%%%%%%%%%%%%%%%%%%%%%%%%%%%%%
%%%%%%%%%%%%%%%%%%%%%%%%%%%%%%%%%%%%%%%%%%%%%%%%%%%%%%%%%%%%%%%%%%%%%%%%%%%%%%%%%%%%%%%%%%%%%%%%%%%%%%%%%%%%%%%%%%%%%%%%%%%%%%%%
%%%%%%%%%%%%%%%%%%%%%%%%%%%%%%%%%%%%%%%%%%%%%%%%%%%%%%%%%%%%%%%%%%%%%%%%%%%%%%%%%%%%%%%%%%%%%%%%%%%%%%%%%%%%%%%%%%%%%%%%%%%%%%%%
%%%%%%%%%%%%%%%%%%%%%%%%%%%%%%%%%%%%%%%%%%%%%%%%%%%%%%%%%%%%%%%%%%%%%%%%%%%%%%%%%%%%%%%%%%%%%%%%%%%%%%%%%%%%%%%%%%%%%%%%%%%%%%%%
%%%%%%%%%%%%%%%%%%%%%%%%%%%%%%%%%%%%%%%%%%%%%%%%%%%%%%%%%%%%%%%%%%%%%%%%%%%%%%%%%%%%%%%%%%%%%%%%%%%%%%%%%%%%%%%%%%%%%%%%%%%%%%%%
%%%%%%%%%%%%%%%%%%%%%%%%%%%%%%%%%%%%%%%%%%%%%%%%%%%%%%%%%%%%%%%%%%%%%%%%%%%%%%%%%%%%%%%%%%%%%%%%%%%%%%%%%%%%%%%%%%%%%%%%%%%%%%%%
%%%%%%%%%%%%%%%%%%%%%%%%%%%%%%%%%%%%%%%%%%%%%%%%%%%%%%%%%%%%%%%%%%%%%%%%%%%%%%%%%%%%%%%%%%%%%%%%%%%%%%%%%%%%%%%%%%%%%%%%%%%%%%%%
%%%%%%%%%%%%%%%%%%%%%%%%%%%%%%%%%%%%%%%%%%%%%%%%%%%%%%%%%%%%%%%%%%%%%%%%%%%%%%%%%%%%%%%%%%%%%%%%%%%%%%%%%%%%%%%%%%%%%%%%%%%%%%%%
\section*{Exercises}
\exercise
Let $\Man{}$ be a manifold of dimension $m$, and let $\metrictensor{}$ be a semi-Riemannian metric on $\Man{}$
with index $\nu\in\seta{\suc{0}{m}}$. Let $\point\in\Man{}$, and let
$\mtuple{\vsbase{1}}{\vsbase{m}}$ be an orthonormal ordered-basis of the scalar-product space
$\opair{\Tanspace{\point}{\Man{}}}{\func{\metrictensor{}}{\point}}$. Prove that for every
$\binary{u}{v}\in\tanspace{\point}{\Man{}}$,
\begin{align}
\func{\[\func{\TTFequiv{\Ricci{\Man{}}{\metrictensor{}}}}{\point}\]}{\binary{u}{v}}&=
\sum_{i=1}^{m}\epsilon_{i}
\func{\[\func{\TTFequiv{\Riemcurvature{\Man{}}{\metrictensor{}}}}{\point}\]}{\Quadruple{u}{\vsbase{i}}{\vsbase{i}}{v}}\cr
&=\sum_{i=1}^{m}\vfmetricproduct{
\func{\[\func{\TTFequiv{\mcurvature{\connection{}}{\metrictensor{}}}}{\point}\]}{\triplet{u}{\vsbase{i}}{\vsbase{i}}}}
{v},
\end{align}
where $\seta{\binary{+1}{-1}}\ni\epsilon_{i}:=\func{\[\func{\metrictensor{}}{\point}\]}{\binary{\vsbase{i}}{\vsbase{i}}}=
\vfmetricproduct{\vsbase{i}}{\vsbase{i}}$.
\endexercise
%%%%%%%%%%%%%%%%%%%%%%%%%%%%%%%%%%%%%%%%%%%%%%%%%%%%%%%%%%%%%%%%%%%%%%%%%%%%%%%%%%%%%%%%%%%%%%%%%%%%%%%%%%%%%%%%%%%%%%%%%%%%%%%%
\exercise
Let $\opair{\Man{}}{\metrictensor{}}$ and $\opair{\Man{1}}{\metrictensor{1}}$ be locally-isometric semi-Riemannian manifolds
with $\function{f}{\Man{}}{\Man{1}}$ a local-isometry between them. Then,
$\func{\VBpullback{f}{4}{0}}{\Riemcurvature{\Man{1}}{\metrictensor{1}}}=\Riemcurvature{\Man{}}{\metrictensor{}}$.
\endexercise
%%%%%%%%%%%%%%%%%%%%%%%%%%%%%%%%%%%%%%%%%%%%%%%%%%%%%%%%%%%%%%%%%%%%%%%%%%%%%%%%%%%%%%%%%%%%%%%%%%%%%%%%%%%%%%%%%%%%%%%%%%%%%%%%
\exercise
Let $\Man{}$ be a manifold with dimension $m$, and let $\metrictensor{}$ be a semi-Riemannian metric on $\Man{}$
with index $\nu$. Prove that the semi-Riemannian manifold $\opair{\Man{}}{\metrictensor{}}$ is flat, that is locally-isometric
to $\semiEucspace{m}{\nu}$, if and only if curvature of $\Man{}$ relative to the Levi-Civita connection of
$\opair{\Man{}}{\metrictensor{}}$ vanishes, that is $\mcurvature{\connection{}}=0$ ($\connection{}$
denotes the Levi-Civita connection of $\opair{\Man{}}{\metrictensor{}}$).
\endexercise
\newpage
\Bibliography{}
\renewcommand{\addcontentsline}[3]{}

\let\addcontentsline\oldaddcontentsline
%%%%%%%%%%%%%%%%%%%%%%%%%%%%%%%%%%%%%%%%%%%%%%%%%%%%%%%%%%%%%%%%%%%%%%%%%%%%%%%%%%%%%%%%%%%%%%%%%%%%%%%

\end{document}